\documentclass[11pt,a4paper]{amsart}
 
\usepackage[margin=2.5cm,marginpar=2cm]{geometry}
 
\usepackage[bookmarksopen,bookmarksdepth=2,hidelinks,colorlinks=false]{hyperref}
\usepackage[nameinlink]{cleveref}

\usepackage{amssymb}
\usepackage{mathrsfs}
\usepackage{mathtools}
\usepackage{amsthm}
\usepackage{graphicx}
\usepackage{braket}

\usepackage{amscd}

\usepackage[alphabetic]{amsrefs}

\usepackage[all,cmtip]{xy}
\usepackage{rotating}
\usepackage{leftidx}

\usepackage[rgb,table]{xcolor}

\usepackage{pifont}
\makeatletter
\newcommand*{\circnuma}[1]{%
  \ifnum#1<1 %
    \@ctrerr
  \else
    \ifnum#1>20 %
      \@ctrerr
    \else
      \mbox{\ding{\numexpr 171+(#1)\relax}}%
     \fi
  \fi
}
\makeatother

\usepackage[centertableaux]{ytableau}

\makeatletter
\pgfkeys{/ytableau/options,
  noframe/.default = false,
  noframe/.is choice,
  noframe/true/.code = {%
    \global\let\vrule@YT=\vrule@none@YT
    \global\let\hrule@YT=\hrule@none@YT
  },
  noframe/false/.code = {%
    \global\let\vrule@YT=\vrule@normal@YT
    \global\let\hrule@YT=\hrule@normal@YT
  },
  noframe/on/.style = {noframe/true},
  noframe/off/.style = {noframe/false},
}

\def\hrule@enon@YT{%
  \hrule width  \dimexpr \boxdim@YT + \fboxrule *2 \relax
  height 0pt
}
\def\vrule@enon@YT{%
  \vrule height \dimexpr  \boxdim@YT + \fboxrule\relax
     width \fboxrule
}

\def\enon{\omit\enon@YT}
\newcommand{\enon@YT}[2][clear]{%
  \def\thisboxcolor@YT{#1}%
  \let\hrule@YT=\hrule@enon@YT
  \let\vrule@YT=\vrule@enon@YT
  \startbox@@YT#2\endbox@YT
  \nullfont
}

\makeatother
\ytableausetup{noframe=off,boxsize=1.3em}
\let\ytb=\ytableaushort

\newcommand{\tytb}[1]{{\tiny\ytb{#1}}}

\usepackage{enumitem}
\newlist{enumC}{enumerate}{1} % Conditions in Lemma/Theorem/Prop
\setlist[enumC,1]{label=(\alph*),wide,ref=(\alph*)}
\crefname{enumCi}{condition}{conditions}
\Crefname{enumCi}{Condition}{Conditions}
\newlist{enumT}{enumerate}{3} % "Theorem"=conclusions in Lemma/Theorem/Prop
\setlist[enumT]{label=(\roman*),wide}
\setlist[enumT,1]{label=(\roman*),wide}
\setlist[enumT,2]{label=(\alph*),ref ={(\roman{enumTi}.\alph*)}}
\setlist[enumT,3]{label=(\arabic*), ref ={(\roman{enumTi}.\alph{enumTii}.\alph*)}}
\crefname{enumTi}{}{}
\Crefname{enumTi}{Item}{Items}
\crefname{enumTii}{}{}
\Crefname{enumTii}{Item}{Items}
\crefname{enumTiii}{}{}
\Crefname{enumTiii}{Item}{Items}
\newlist{enumPF}{enumerate}{3}
\setlist[enumPF]{label=(\alph*),wide}
\setlist[enumPF,1]{label=(\roman*),wide}
\setlist[enumPF,2]{label=(\alph*)}
\setlist[enumPF,3]{label=\arabic*).}
\newlist{enumS}{enumerate}{3} % Statement outside Lemma/Theorem/Prop
\setlist[enumS]{label=\roman*)}
\setlist[enumS,1]{label=\roman*)}
\setlist[enumS,2]{label=\alph*)}
\setlist[enumS,3]{label=\arabic*.}
\newlist{enumI}{enumerate}{3} % items
\setlist[enumI,1]{label=\roman*),leftmargin=*}
\setlist[enumI,2]{label=\alph*), leftmargin=*}
\setlist[enumI,3]{label=\arabic*), leftmargin=*}
\newlist{enumIL}{enumerate*}{1} % inline enum
\setlist*[enumIL]{label=\roman*)}
\newlist{enumR}{enumerate}{1} % remarks
\setlist[enumR]{label=\arabic*.,wide,labelwidth=!, labelindent=0pt}
\crefname{enumRi}{remark}{remarks}

\crefname{equation}{}{}
\Crefname{equation}{Equation}{Equations}
\Crefname{lem}{Lemma}{Lemma}
\Crefname{thm}{Theorem}{Theorem}

\newlist{des}{description}{1}
\setlist[des]{font=\sffamily\bfseries}

\blendcolors{!80!black}

\long\def\delete#1{}

\newcommand{\trivial}[2][]{\if\relax\detokenize{#1}\relax
  {%\hfill\break
      \color{orange} \vspace{0em}$[$#2$]$
      \color{black}
  }
  \else
\ifx#1h
\ifcsname showtrivial\endcsname
{%\hfill\break
    \color{orange}\vspace{0em}$[$#2$]$
    \color{black}
}
\fi
\else {\red Wrong argument!} \fi
\fi
}

\newcommand{\byhide}[2][]{\if\relax\detokenize{#1}\relax
{\color{orange} \vspace{0em} Plan to delete:  #2}
\else
\ifx#1h\relax\fi
\fi
}

\newcommand{\pr}{\mathrm{pr}}

\newcommand{\AC}{\mathrm{AC}}
\newcommand{\WF}{\mathrm{WF}}
\newcommand{\AV}{\mathrm{AV}}

\def\subset{\subseteq}



\makeatletter
\def\inn#1#2{\left\langle
      \def\ta{#1}\def\tb{#2}
      \ifx\ta\@empty{\;} \else {\ta}\fi ,
      \ifx\tb\@empty{\;} \else {\tb}\fi
      \right\rangle}
\def\binn#1#2{\left\lAngle
      \def\ta{#1}\def\tb{#2}
      \ifx\ta\@empty{\;} \else {\ta}\fi ,
      \ifx\tb\@empty{\;} \else {\tb}\fi
      \right\rAngle}
\makeatother

\makeatletter
\def\binn#1#2{\overline{\inn{#1}{#2}}}
\makeatother

 \def\ckhha{{}^a \check\fhh}

\usepackage{xparse}
\def\usecsname#1{\csname #1\endcsname}
\def\useLetter#1{#1}
\def\usedbletter#1{#1#1}

\ExplSyntaxOn

\def\mydefcirc#1#2#3{\expandafter\def\csname
  circ#3{#1}\endcsname{{}^\circ {#2{#1}}}}
\def\mydefvec#1#2#3{\expandafter\def\csname
  vec#3{#1}\endcsname{\vec{#2{#1}}}}
\def\mydefdot#1#2#3{\expandafter\def\csname
  dot#3{#1}\endcsname{\dot{#2{#1}}}}

\def\mydefacute#1#2#3{\expandafter\def\csname a#3{#1}\endcsname{\acute{#2{#1}}}}
\def\mydefbr#1#2#3{\expandafter\def\csname br#3{#1}\endcsname{\breve{#2{#1}}}}
\def\mydefbar#1#2#3{\expandafter\def\csname bar#3{#1}\endcsname{\bar{#2{#1}}}}
\def\mydefhat#1#2#3{\expandafter\def\csname hat#3{#1}\endcsname{\hat{#2{#1}}}}
\def\mydefwh#1#2#3{\expandafter\def\csname wh#3{#1}\endcsname{\widehat{#2{#1}}}}
\def\mydeft#1#2#3{\expandafter\def\csname t#3{#1}\endcsname{\tilde{#2{#1}}}}
\def\mydefu#1#2#3{\expandafter\def\csname u#3{#1}\endcsname{\underline{#2{#1}}}}
\def\mydefr#1#2#3{\expandafter\def\csname r#3{#1}\endcsname{\mathrm{#2{#1}}}}
\def\mydefb#1#2#3{\expandafter\def\csname b#3{#1}\endcsname{\mathbb{#2{#1}}}}
\def\mydefwt#1#2#3{\expandafter\def\csname wt#3{#1}\endcsname{\widetilde{#2{#1}}}}
\def\mydefbf#1#2#3{\expandafter\def\csname bf#3{#1}\endcsname{\mathbf{#2{#1}}}}
\def\mydefc#1#2#3{\expandafter\def\csname c#3{#1}\endcsname{\mathcal{#2{#1}}}}
\def\mydefsf#1#2#3{\expandafter\def\csname sf#3{#1}\endcsname{\mathsf{#2{#1}}}}
\def\mydefs#1#2#3{\expandafter\def\csname s#3{#1}\endcsname{\mathscr{#2{#1}}}}
\def\mydefcks#1#2#3{\expandafter\def\csname cks#3{#1}\endcsname{{\check{
        \csname s#2{#1}\endcsname}}}}
\def\mydefckc#1#2#3{\expandafter\def\csname ckc#3{#1}\endcsname{{\check{
      \csname c#2{#1}\endcsname}}}}
\def\mydefck#1#2#3{\expandafter\def\csname ck#3{#1}\endcsname{{\check{#2{#1}}}}}

\cs_new:Npn \mydeff #1#2#3 {\cs_new:cpn {f#3{#1}} {\mathfrak{#2{#1}}}}

\cs_new:Npn \doGreek #1
{
  \clist_map_inline:nn {alpha,beta,gamma,Gamma,delta,Delta,epsilon,varepsilon,zeta,eta,theta,vartheta,Theta,iota,kappa,lambda,Lambda,mu,nu,xi,Xi,pi,Pi,rho,sigma,varsigma,Sigma,tau,upsilon,Upsilon,phi,varphi,Phi,chi,psi,Psi,omega,Omega,tG} {#1{##1}{\usecsname}{\useLetter}}
}

\cs_new:Npn \doSymbols #1
{
  \clist_map_inline:nn {otimes,boxtimes} {#1{##1}{\usecsname}{\useLetter}}
}

\cs_new:Npn \doAtZ #1
{
  \clist_map_inline:nn {A,B,C,D,E,F,G,H,I,J,K,L,M,N,O,P,Q,R,S,T,U,V,W,X,Y,Z} {#1{##1}{\useLetter}{\useLetter}}
}

\cs_new:Npn \doatz #1
{
  \clist_map_inline:nn {a,b,c,d,e,f,g,h,i,j,k,l,m,n,o,p,q,r,s,t,u,v,w,x,y,z} {#1{##1}{\useLetter}{\usedbletter}}
}

\cs_new:Npn \doallAtZ
{
\clist_map_inline:nn {mydefsf,mydeft,mydefu,mydefwh,mydefhat,mydefr,mydefwt,mydeff,mydefb,mydefbf,mydefc,mydefs,mydefck,mydefcks,mydefckc,mydefbar,mydefvec,mydefcirc,mydefdot,mydefbr,mydefacute} {\doAtZ{\csname ##1\endcsname}}
}

\cs_new:Npn \doallatz
{
\clist_map_inline:nn {mydefsf,mydeft,mydefu,mydefwh,mydefhat,mydefr,mydefwt,mydeff,mydefb,mydefbf,mydefc,mydefs,mydefck,mydefbar,mydefvec,mydefdot,mydefbr,mydefacute} {\doatz{\csname ##1\endcsname}}
}

\cs_new:Npn \doallGreek
{
\clist_map_inline:nn {mydefck,mydefwt,mydeft,mydefwh,mydefbar,mydefu,mydefvec,mydefcirc,mydefdot,mydefbr,mydefacute} {\doGreek{\csname ##1\endcsname}}
}

\cs_new:Npn \doallSymbols
{
\clist_map_inline:nn {mydefck,mydefwt,mydeft,mydefwh,mydefbar,mydefu,mydefvec,mydefcirc,mydefdot} {\doSymbols{\csname ##1\endcsname}}
}

\cs_new:Npn \doGroups #1
{
  \clist_map_inline:nn {GL,Sp,rO,rU,fgl,fsp,foo,fuu,fkk,fuu,ufkk,uK} {#1{##1}{\usecsname}{\useLetter}}
}

\cs_new:Npn \doallGroups
{
\clist_map_inline:nn {mydeft,mydefu,mydefwh,mydefhat,mydefwt,mydefck,mydefbar} {\doGroups{\csname ##1\endcsname}}
}

\cs_new:Npn \decsyms #1
{
\clist_map_inline:nn {#1} {\expandafter\DeclareMathOperator\csname ##1\endcsname{##1}}
}

\decsyms{Mp,id,SL,Sp,SU,SO,GO,GSO,GU,GSp,PGL,Pic,Lie,Mat,Ker,Hom,Ext,Ind,reg,res,inv,Isom,Det,Tr,Norm,Sym,Span,Stab,Spec,PGSp,PSL,tr,Ad,Br,Ch,Cent,End,Aut,Dvi,Frob,Gal,GL,Gr,DO,ur,vol,ab,Nil,Supp,rank,Sign}

\def\abs#1{\left|{#1}\right|}

\NewDocumentCommand\cent{o m }{
  \IfValueTF{#1}{
    \mathop{Z}_{#1}{(#2)}}
  {\mathop{Z}{(#2)}}
}

\def\fsl{\mathfrak{sl}}

\doallAtZ
\doallatz
\doallGreek
\doallGroups
\doallSymbols
\ExplSyntaxOff

\usepackage{diagbox}
\usepackage{arydshln}

\usepackage{tikz}
\usetikzlibrary{matrix,arrows,positioning,cd,backgrounds}
\usetikzlibrary{decorations.pathmorphing,decorations.pathreplacing}

\usepackage{upgreek}

\usepackage{listings}
\newcommand{\BC}{{\mathbb {C}}}

\newcommand{\BH}{{\mathbb {H}}}

\newcommand{\BN}{{\mathbb {N}}}

\newcommand{\BR}{{\mathbb {R}}}
\newcommand{\BS}{{\mathbb {S}}}

\newcommand{\CB}{{\mathcal {B}}}
\newcommand{\CC}{{\mathcal {C}}}

\newcommand{\CK}{{\mathcal {K}}}

\newcommand{\CO}{{\mathcal {O}}}
\newcommand{\CP}{{\mathcal {P}}}

\newcommand{\CS}{{\mathcal {S}}}

\newcommand{\CU}{{\mathcal {U}}}

\newcommand{\CZ}{{\mathcal {Z}}}

\newcommand{\RS}{{\mathrm {S}}}

\DeclareMathOperator{\Ann}{Ann}

\newcommand{\cusp}{{\mathrm{cusp}}}

\def\CCL{\CCL}
\def\CCLR{\CCD}

\newcommand{\sgn}{\operatorname{sgn}}

\newcommand{\oL}{\operatorname{L}}

\newcommand{\oO}{\operatorname{O}}
\newcommand{\oS}{\operatorname{S}}

\newcommand{\oU}{\operatorname{U}}

\newcommand{\g}{\mathfrak g}

\newcommand{\h}{\mathfrak h}
\newcommand{\p}{\mathfrak p}

\renewcommand{\b}{\mathfrak b}

\newcommand{\n}{\mathfrak n}
\renewcommand{\u}{\mathfrak u}

\renewcommand{\l}{\mathfrak l}
\renewcommand{\t}{\mathfrak t}
\newcommand{\s}{\mathfrak s}

\renewcommand{\o}{\mathfrak o}

\newcommand{\gl}{\mathfrak g \mathfrak l}

\newcommand{\Z}{\mathbb{Z}}
\DeclareDocumentCommand{\C}{}{\mathbb{C}}
\newcommand{\R}{\mathbb R}

\def\eDD{\mathrm{d}^{e}}
\def\DD{\nabla}

\newcommand{\la}{\langle}
\newcommand{\ra}{\rangle}

\newcommand{\be}{\begin {equation}}
\newcommand{\ee}{\end {equation}}

\numberwithin{equation}{section}

\def\flushl#1{\ifmmode\makebox[0pt][l]{${#1}$}\else\makebox[0pt][l]{#1}\fi}
\def\flushr#1{\ifmmode\makebox[0pt][r]{${#1}$}\else\makebox[0pt][r]{#1}\fi}

\newtheorem*{thm*}{Theorem}
\newtheorem{thm}{Theorem}[section]

\newtheorem{lem}[thm]{Lemma}

\newtheorem*{lem*}{Lemma}
\newtheorem{whyp}[thm]{Hypothesis}
\newtheorem{prop}[thm]{Proposition}

\newtheorem{cor}[thm]{Corollary}

\newtheorem*{claim*}{Claim}
\newtheorem{defn}[thm]{Definition}

\newtheorem*{eg*}{Example}
\newtheorem{eg}[thm]{Example}
\newtheorem{conj}[thm]{Conjecture}

\theoremstyle{remark}
\newtheorem{remark}[thm]{Remark}
\newtheorem{remarks}[thm]{Remarks}
\newtheorem*{remark*}{Remark}

\def\half{{\tfrac{1}{2}}}

\ExplSyntaxOn
\clist_map_inline:nn {J,L,C,X,Y,H,c,e,f,h,}{
  \expandafter\def\csname #1slt\endcsname{{\mathring{#1}}}}
\ExplSyntaxOff

\NewDocumentCommand{\NilP}{t'}{
\IfBooleanTF{#1}{\Nil_{\fpp'}}{\Nil_\fpp}
}

\NewDocumentCommand{\KTW}{o g}{
  \IfValueTF{#2}{
    \left.\varsigma_{\IfValueT{#1}{#1}}\right|_{#2}}{
    \varsigma_{\IfValueT{#1}{#1}}}
}

\NewDocumentCommand{\CHI}{o g}{
  \IfValueTF{#1}{
    {\chi}_{\left[#1\right]}}{
    \IfValueTF{#2}{
      {\chi}_{\left(#2\right)}}{
      {\chi}}
  }
}
\NewDocumentCommand{\PR}{g}{
  \IfValueTF{#1}{
    \mathop{\pr}_{\left(#1\right)}}{
    \mathop{\pr}}
}
\NewDocumentCommand{\XX}{g}{
  \IfValueTF{#1}{
    {\cX}_{\left(#1\right)}}{
    {\cX}}
}
\NewDocumentCommand{\PP}{g}{
  \IfValueTF{#1}{
    {\fpp}_{\left(#1\right)}}{
    {\fpp}}
}
\NewDocumentCommand{\LL}{g}{
  \IfValueTF{#1}{
    {\bfL}_{\left(#1\right)}}{
    {\bfL}}
}
\NewDocumentCommand{\ZZ}{g}{
  \IfValueTF{#1}{
    {\cZ}_{\left(#1\right)}}{
    {\cZ}}
}

\NewDocumentCommand{\WW}{g}{
  \IfValueTF{#1}{
    {\bfW}_{\left(#1\right)}}{
    {\bfW}}
}

\NewDocumentCommand\KK{g}{
\IfValueTF{#1}{K_{(#1)}}{K}}
\NewDocumentCommand\XXo{d()}{
\IfValueTF{#1}{\cX^\circ_{(#1)}}{\cX^\circ}}

\NewDocumentCommand\ZZo{g}{
\IfValueTF{#1}{\cZ^\circ_{(#1)}}{\cZ^\circ}}

\NewDocumentCommand{\bcO}{t'}{
  \overline{\cO\IfBooleanT{#1}{'}}}

\NewDocumentCommand{\oliftc}{g}{
\IfValueTF{#1}{\boldsymbol{\vartheta} (#1)}{\boldsymbol{\vartheta}}
}
\NewDocumentCommand{\oliftr}{g}{
\IfValueTF{#1}{\vartheta_\bR(#1)}{\vartheta_\bR}
}
\NewDocumentCommand{\olift}{g}{
\IfValueTF{#1}{\vartheta(#1)}{\vartheta}
}
\NewDocumentCommand{\tlift}{g}{
\IfValueTF{#1}{\wtvartheta(#1)}{\wtvartheta}
}

\DeclareDocumentCommand{\NN}{g}{
\IfValueTF{#1}{\fN(#1)}{\fN}
}

\NewDocumentCommand\lnn{t+ t- g}{
  \IfBooleanTF{#1}{{}^l n^+\IfValueT{#3}{(#3)}}{
    \IfBooleanTF{#2}{{}^l n^-\IfValueT{#3}{(#3)}}{}
  }
}

\NewDocumentCommand\LW{g}{
\IfValueTF{#1}{L_{W_{#1}}}{L_{W}}}

\def\floor#1{{\lfloor #1 \rfloor}}

\def\cf{\emph{cf.} }

\def\Irr{\mathrm{Irr}}
\def\Irrsp{\mathrm{Irr}^{\mathrm{sp}}}

\def\Unip{\mathrm{Unip}}

\makeatletter
\pgfkeys{/ytableau/options,
  noframe/.default = false,
  noframe/.is choice,
  noframe/true/.code = {%
    \global\let\vrule@YT=\vrule@none@YT
    \global\let\hrule@YT=\hrule@none@YT
  },
  noframe/false/.code = {%
    \global\let\vrule@YT=\vrule@normal@YT
    \global\let\hrule@YT=\hrule@normal@YT
  },
  noframe/on/.style = {noframe/true},
  noframe/off/.style = {noframe/false},
}
\makeatother

\def\dBV{d_{\mathrm{BV}}}
\def\CP{\mathsf{CP}}
\def\YD{\mathsf{YD}}

\def\DD{\nabla}

\def\lamck{\lambda_\ckcO}

\def\Cint#1{\Coh_{[#1]}}
\def\PP{\mathrm{PAP}}
\def\PAP{\mathrm{PAP}}
\def\BOX#1{\mathrm{Box}(#1)}

\def\hha{{}^a\fhh}
\def\ahh{\hha}

\def\Wlamck{W_{\lamck}}

\def\LC{{}^{\scriptscriptstyle L}\sC}

\def\bsgn{\overline{\sgn}}

\def\Wg{W_{\mathrm g}}

\def\nbb{n_{\mathrm b}}
\def\ngg{n_{\mathrm g}}
\def\tU{\widetilde{\rU}}

\newcommand{\cross}{\times}

\def\AND{\quad \text{and} \quad}
\def\Coh{\mathrm{Coh}}

\def\ev#1{{\mathrm{ev}_{#1}}}

\def\cuprow{{\stackrel{r}{\sqcup}}}
\def\cupcol{{\stackrel{c}{\sqcup}}}

\def\Spr{\mathrm{Springer}}

\def\imathp{\imath_{\wp}}
\def\jmathp{\jmath_{\wp}}

\def\CPP{\mathrm{PP}}
\def\CPPs{\mathrm{PP}_{\star}}

\def\tPBP{\widetilde{\PBP}}

\def\bPsi{\overline{\Psi}}

\def\leqLR{\mathrel{\mathop{\leq}\limits_{\scriptscriptstyle LR}}}

\def\approxLR{\mathrel{\mathop{\approx}\limits^{\scalebox{0.4}{$LR$}}}}
\def\approxLR{\approx}

\def\PBP{\mathrm{PBP}}

\def\ckfgg{{\check \fgg}}

\def\Inn{\mathrm{Ad}}

\def\cuprow{{\,\stackrel{r}{\sqcup}\,}}
\def\cupcol{{\,\stackrel{c}{\sqcup}\,}}

\def\ckcOp{\ckcO'}
\def\ckcOpp{\ckcO''}

\def\ckcOb{\ckcO_{\mathrm b}}
\def\ckcOpb{\ckcO'_{\mathrm b}}
\def\cOpb{\cO'_{\mathrm b}}
\def\ckcOg{\ckcO_{\mathrm g}}

\def\nng{n_{\mathrm g}}
\def\nnb{n_{\mathrm b}}
\def\Gb{G_{\mathrm b}}
\def\Gpb{G'_{\mathrm b}}
\def\Pb{P_{\mathrm b}}
\def\Gg{G_{\mathrm g}}

\def\tPBP{\widetilde{\PBP}}
\def\PBPs{\PBP_{\star}}

\def\PBPop#1#2#3#4{\PBP_{#1}^{#2}(#3,#4)}
\newcommand{\PBPOP}[1][]{\PBPop{\star}{#1}{\ckcO}{\wp}}
\def\PBPOPp{\PBPop{\star'}{}{\ckcO'}{\wp'}}
\newcommand{\PBPOPpp}[1][]{\PBPop{\star}{#1}{\ckcO''}{\wp''}}

\def\PBPGOP{\PBPop{G}{}{\ckcO}{\wp}}

\def\yrow#1{\left[#1\right]_{\mathrm{row}}}

\def\wpu{\wp_{\uparrow}}
\def\wpm{\wp_{\downarrow}}
\def\wpd{\wp} % define the done-wp to be \wp
\def\uptauu{\uptau_{\uparrow}}
\def\uptaud{\uptau} % define the done-tau to be \uptau

 \def\tnaive{\mathrm{naive}}

\def\imathwpp{\imath_{\wp'}}
\def\jmathwpp{\jmath_{\wp'}}
\def\cPpn{\cP'_\mathrm{naive}}
\def\cQpn{\cQ'_\mathrm{naive}}

\def\uptaupn{\uptau'_{\tnaive}}

\def\eDD{\mathrm{d}^{e}}
\def\ckDD{{\check\DD}}
\def\DD{\nabla}
\def\DDn{\nabla_{\tnaive}}
\def\ckDDn{{\ckDD}_{\tnaive}}

\newcommand{\Lam}{{[\lambda]}}

\newcommand{\Grt}{\cK}

\def\CCL{{\mathcal{L}}}

\def\CCD{{\mathcal{D}}}

\usepackage{xr}
\usepackage{subfiles} % Best loaded last in the preamble

\begin{document}

\title[Special unipotent representations: counting and reduction]{Special unipotent representations of real classical groups: counting and reduction}

\author [D. Barbasch] {Dan Barbasch}
\address{Department of Mathematics\\
  310 Malott Hall, Cornell University, Ithaca, New York 14853 }
\email{barbasch@math.cornell.edu}

\author [J.-J. Ma] {Jia-Jun Ma}
\address{School of Mathematical Sciences\\
  Xiamen University\\
  Xiamen, 361005, China} 
  \address{Department of Mathematics, Xiamen University Malaysia campus, Sepang, Selangor Darul Ehsan, 43900, Malaysia} 
  \email{hoxide@xmu.edu.cn}

\author [B. Sun] {Binyong Sun}
\address{Institute for Advanced Study in Mathematics \& New Cornerstone Science Laboratory, Zhejiang University,  Hangzhou, 310058, China}
\email{sunbinyong@zju.edu.cn}

\author [C.-B. Zhu] {Chen-Bo Zhu}
\address{Department of Mathematics\\
  National University of Singapore\\
  10 Lower Kent Ridge Road, Singapore 119076} \email{matzhucb@nus.edu.sg}

\subjclass[2020]{22E46, 22E47} \keywords{Classical group, special unipotent representation, coherent continuation, Weyl group representation, primitive ideal, cell, infinitesimal character, associated variety}

% \thanks{Supported by NSFC Grant 11222101}

\begin{abstract} Let $G$ be a real reductive group in Harish-Chandra's class. We derive some consequences of theory of coherent continuation representations to the counting of irreducible representations of $G$ with a given infinitesimal character and a given bound of the complex associated variety. When $G$ is a real classical group (including the real metaplectic group), we investigate the set of special unipotent representations of $G$ attached to $\check{\mathcal O}$, in the sense of Arthur and Barbasch-Vogan. Here $\check{\mathcal O}$ is a nilpotent adjoint orbit in the Langlands dual of $G$ (or the metaplectic dual of $G$ when $G$ is a real metaplectic group). We give a precise count for the number of special unipotent representations of $G$ attached to $\check{ \mathcal O}$. We also reduce the problem of constructing special unipotent representations attached to $\check{\mathcal O}$ to the case when $\check{\mathcal O}$ is analytically even (equivalently for a real classical group, has good parity in the sense of M{\oe}glin).  The paper is the first in a series of two papers on the classification of special unipotent representations of real classical groups.
\end{abstract}

\maketitle

\tableofcontents

\section{Introduction}\label{sec:intro}

\subsection{Background and goals} In \cite{ArPro,ArUni}, Arthur introduced certain families of representations of a reductive  algebraic group $G$ over $\R$ or $\C$, in connection with his conjecture on square-integrable automorphic forms. Arthur’s representations, the special unipotent representations in the title of this paper, were made precise in the work of Barbasch-Vogan \cite{BVUni} (for groups over $\BC$; the same works for groups over $\BR$, see \cite[Chapter 27]{ABV}). They are defined in the language of primitive ideals and are attached to a nilpotent adjoint orbit $\check \CO$ in the Langlands dual of $G$.  

Apart from their clear interest for the theory of automorphic forms \cite{ArPro, ArUni,ArEnd}, special unipotent representations belong to a fundamental class of unitary representations which are associated to nilpotent coadjoint orbits in
the Kirillov philosophy (the orbit method; see \cite{Ki62,Ko70,VoBook}). These are known informally as unipotent representations, which are expected to play a central role in the classification problem of the unitary dual of a real reductive group and therefore have been a subject of great interest. See 
%{\cite{V.GL,B89}}
\cite{VoICM,VoBook,Vo89}. The aforementioned work of Barbasch-Vogan \cite{BVUni} was also motivated by this classification problem. 

In the context alluded to earlier, a long standing problem in representation theory, known as the Arthur-Barbasch-Vogan conjecture (\cite[Section 4]{ArUni}, \cite[Introduction]{ABV}), is to show that every special unipotent representation is unitary.

In a series of two papers, the authors will construct and classify special unipotent representations of real classical groups (including real metaplectic groups), and will prove the Arthur-Barbasch-Vogan conjecture for these groups as a consequence of the classification. For quasi-split classical groups, the unitarity of special unipotent representations is independently established in \cite{AAM,AM}, from the perspective of the endoscopic classification of representations \cite{ArEnd, Mok}.

The current paper is the first in the series, which has the following three goals. 
\begin{enumerate}[label=(\alph*)]
\item 
The first is to derive a  consequence of the theory of coherent continuation representations
 to the counting, in the form of an inequality, of irreducible representations of $G$ with a given infinitesimal character and a given bound of the complex associated variety. Largely this is achieved by building on existing results in  Kazhdan-Lusztig theory.   
 \item The second is to prove that the inequality just alluded to is in fact an equality under certain technical hypothesis on $G$, which holds for a real classical group and yields precise and explicit counting of special unipotent representations attached to $\check \CO$. This involves some new technical tools on $\tau$-invariants for a general $G$, as well as elaborate combinatorial constructions/arguments in the case of a real classical group $G$. 
 \item The third is to reduce
 the problem, in the case of a real classical group $G$, of constructing special unipotent representations attached to $\check \CO$ to the case when $\check \CO$ is analytically even (\Cref{def:ae}), or equivalently, has good parity in the sense of M\oe glin \cite{Mo11}. 
 This may be viewed as reduction to integral
infinitesimal character which holds for the more general setting as in \cite{ABV}. Therefore the core of the classification problem of special unipotent representations reduces to the case when $\check \CO$ has good parity. We remark that results of similar nature were proved by M\oe glin and Renard in \cite{MR.U}.

 \end{enumerate}

\subsection{Approach}
\label{subsec:appro}
We will work in the category of Casselman-Wallach representations \cite[Chapter 11]{Wa2}. The main tool for the counting of irreducible representations in general and special unipotent representations in particular is the coherent continuation representation (the idea first appeared in \cite{BV.W}). 
This is a representation of a certain integral Weyl group, and it can be compared, through a certain injective map of $\CK$ groups (due to Casian), with an analogous coherent continuation representation for the category of highest weight modules. The latter has been intensively studied in Kazhdan-Lusztig theory and is amenable to detailed analysis through the theory of primitive ideals (as in the work of Joseph and Barbasch-Vogan), as well as the theory of cells and special representations (in the sense of Lusztig). In particular one may derive precise information on what representations of the integral Weyl group may contribute to the counting, in terms of the Springer correspondence. 
(This was done in \cite[Section 5]{BVUni} for complex semisimple groups.) In effectively carrying out the mentioned steps, we build on earlier ideas of several authors including Joseph \cite{J1,J2,J.av}, Vogan \cite{Vg}, 
Barbasch-Vogan \cite{BV1,BV2,BVUni}, and Casian \cite{Cas}.
Weaving things together, we arrive at a counting inequality on irreducible representations of $G$ with a given infinitesimal character and a given bound of the complex associated variety. The inequality is valid 
for any real reductive group $G$ in Harish-Chandra's class.

For the precise counting of irreducible representations, a key technical issue is a certain relationship (expected by Vogan) between cell representations in the 
Casselman-Wallach setting (Harish-Chandra cells) and the highest weight module setting (double cells). (This was known to hold for certain $G$, notably for unitary groups \cite{BV.W} and for complex semisimple groups \cite[Section 3]{BVUni}.) We resolve this issue assuming a certain equality of $\tau$-invariants and (a weak form of) Vogan duality, both of which hold true for classical groups (including the real metaplectic group).
Together with an explicit formula of the coherent continuation representation (due to Barbasch-Vogan) and explicit branching rules of Weyl group representations (via the Littlewood-Richardson rule), this will finally yield the counting of special unipotent representations of real classical groups, explicitly described in terms of certain combinatorial constructs called painted bipartitions. It is worthwhile to note that in explicating the counting formula, we need to sum over several Weyl group representations in a certain Lusztig left cell, which we prove to have the same contribution by capitalizing on the induction mechanism provided by a certain notion of descent for painted bipartitions (defined in this article).

We now outline our approach to the third goal. The case of type $A$ groups is relatively easy and so we focus on groups of type $B$, $C$, and $D$.  The basic idea is to appeal to the theory of endoscopy and its relation to the Kazhdan-Lusztig-Vogan algorithm for nonintegral infinitesimal character (see \cite[Chapter 15]{ABV}). Specifically when $\check \CO$ is not of good parity, we introduce an $\check \CO$-relevant parabolic subgroup $P$, whose Levi component is of the form
 $\Gpb\times \Gg $, where $\Gpb$ is a general linear group, and $\Gg$ is of the same classical type as $G$. In their respective Langlands duals, we have nilpotent adjoint orbits
$\check \CO'_\mathrm b$ and $\check \CO_\mathrm g$ (determined by $\check \CO$). We will show that the normalized smooth parabolic induction from $P$ to $G$ yields a bijective map from a pair of special unipotent representations of $\Gpb$ (attached to $\check \CO'_\mathrm b$) and $\Gg$ (attached to $\check \CO_\mathrm g$) to special unipotent representations of $G$ (attached to $\check \CO$).
Firstly, there is a classical group $\Gb$ with the following properties:
\begin{itemize}
\item $\Gb\times \Gg $ is an endoscopy group of $G$ (except for the real metaplectic group, in which case there is a metaplectic analog); 
\item $\Gpb$ is naturally isomorphic to the Levi component of an $\check \CO_\mathrm b$-relevant parabolic subgroup $P_\mathrm b$ of $\Gb$.
\end{itemize}
Here $\check \CO_\mathrm b$ is a nilpotent adjoint orbit 
in the Langlands dual of $\Gb$ 
(determined by $\ckcO$).
General considerations as well as detailed information about coherent continuation representations allow us to separate bad parity and good parity, namely we may determine special unipotent representations of $G$ in terms of those of $\Gb$ and $\Gg$. The normalized smooth parabolic induction from $P_\mathrm b$ to $\Gb$ yields a bijection from special unipotent representations of $\Gpb$ (attached to $\check \CO'_\mathrm b$) to special unipotent representations of $\Gb$ (attached to $\check \CO_\mathrm b$). 
On the other hand, the Kazhdan-Lusztig-Vogan algorithm 
\cite{V3,ABV,RT2} implies that irreducibility is preserved in the endoscopy setting. 
This will finally imply the bijectivity of the induction map from $P$ to $G$.

\newcommand{\Rep}{\mathrm{Rep}}

\subsection{Organization}
In \Cref{sec:main}, we state our main results firstly on the counting of irreducible representations in our general setup, and secondly on the explicit counting of special unipotent representations of real classical groups. For the latter, the more involved cases of groups of type $B$, $C$ and $D$, the main results are in \Cref{reduction} (reduction to good parity) and \Cref{countup} (counting in the case of good parity). 
In \Cref{sec:HWM}, we review some generalities on coherent continuation representations for highest weight modules, to be precise for the category $\Rep(\g,\b)$ (see \Cref{sub:CohHWM} for the notation), as well as the full subcategory of $\Rep(\g,\b)$ defined by a prescribed support condition. Following \cite{FJMN}, we introduce a new notion of the $\tau$-invariant of a Weyl group representation, and examine its role in a certain duality notion of double cells (\Cref{subsec:tau-DC}).  
In \Cref{sec:CW}, we examine some analogous results on coherent continuation representations for Casselman-Wallach representations, and their interplay with the corresponding results in the highest weight module setting. As mentioned in \Cref{subsec:appro}, we also establish a certain relationship between Harish-Chandra cells and double cells assuming some equality of $\tau$-invariants and a weak form of the Vogan duality (\Cref{subsec:tau-HC}). In \Cref{sec:counting}, we prove the first part of our main results on the counting of irreducible representations with a given infinitesimal character and a given bound of the complex associated variety. 
In \Cref{sec:GB}, we separate good parity and bad parity for coherent continuation representations, as a preparation for the reduction step in \Cref{sec:red}. \Cref{sec:GL} and \Cref{sec:BCD} are devoted to proof of the second part of our main results, which explicitly counts special unipotent representations, for groups of type $A$, and groups of type $B$, $C$ and $D$, respectively. In \Cref{sec:red}, we carry out the reduction step to the case of good parity and prove \Cref{reduction}. In \Cref{sec:com}, we develop combinatorics of painted bipartitions which we need for the proof of \Cref{countup} in \Cref{sec:BCD}. This includes a notion of descent for painted bipartitions, helpful in explicating combinatorially the counting formula 
(mentioned earlier), and critical in matching combinatorial parameters of the current paper with special unipotent representations to be constructed in \cite{BMSZ2} (the second paper in the series) by the method of theta lifting (\cite{Howe79,Howe89}). 

\section{The main results}\label{sec:main}

\subsection{Lie algebra notation}\label{secnot}

Let $\g$ be a reductive complex Lie algebra. Its universal enveloping algebra is denoted by
$\mathcal U(\g)$, and the center of $\mathcal U(\g)$ is denoted by $\CZ(\g)$. Let
 $\hha$ denote the  universal Cartan subalgebra of $\g$    (also called the abstract Cartan subalgebra in \cite{V4}). Recall that for every Borel subalgebra $\b$ of $\g$, there is an identification $\hha=\b/[\b,\b]$. Let $W\subset \GL(\hha)$ denote the Weyl group. In general we denote by superscipt $*$ the linear dual of a vector space.  
 By the Harish-Chandra isomorphism, there is a 1-1
   correspondence between $W$-orbits of $\nu \in \hha^*$
and algebraic characters $\chi_\nu: \CZ(\g)\rightarrow \C$. 
We say that an ideal of $\CU(\g)$ has infinitesimal character $\nu$ if it contains the kernel of $\chi_\nu$. By a result of Duflo 
(\cite{Dix}, \cite[Section 3]{Bor}), 
there is a unique maximal ideal of $\CU(\g)$ that has infinitesimal  character  $\nu$. Write $I_\nu$ for this maximal ideal. (We remark that every maximal ideal of $\CU(\g)$ is primitive.) 
\delete{Note that the maps
\[
  \hha^*\rightarrow \{\textrm{algebraic character of $\CZ(\g)$}\}, \quad \nu\mapsto \chi_\nu
\]
and
\[
  \hha^*\rightarrow \{\textrm{maximal ideal  of $\CU(\g)$}\}, \quad \nu\mapsto I_\nu
\]
are both surjective, and the fibers of these two maps are precisely the $W$-orbits.}

Let $\Nil(\g)$ (resp. $\Nil(\g^*)$) be the set of nilpotent elements in $[\g,\g]$ (resp. $[\g,\g]^*$). Denote by $\mathrm{Ad}(\g)$ the inner automorphism group of $\g$, and put
\[
 \overline{\Nil}(\g):=\Inn(\g)\backslash \Nil(\g), \qquad \overline{\Nil}(\g^*):=\Inn(\g)\backslash \Nil(\g^*), \]
the set of $\Inn(\g)$-orbits in $\Nil(\g)$ and
$\Nil(\g^*)$ (which are finite). 
The Killing form on $[\g,\g]$  yields an identification $\overline{\Nil}(\g)=\overline{\Nil}(\g^*)$. 
Since $\g$ is the direct sum of its center with $[\g, \g]$, $[\g,\g]^*$ is viewed as a subspace of $\g^*$ in the obvious way.

\subsection{Counting irreducible representations}\label{subsec:countIrr}

Let $G$ be a real reductive group in Harish-Chandra's class (\cite{Kn}) whose complexified Lie algebra equals $\g$. In the rest of  this paper,  unless  mentioned otherwise, 
we use the corresponding  lowercase Gothic letter to denote the complexified Lie algebra of a Lie group. Let $\mathrm{Rep}(G)$ be the category of Casselman-Wallach representations of $G$, whose Grothendieck group (with $\C$-coefficients) is denoted by $\CK(G)$. 
Let $\Irr(G)$ be the set of isomorphism classes of irreducible Casselman-Wallach representations of $G$, which forms a basis of $\CK(G)$.

Let $\nu\in \hha^*$ and let $\sfS$ be an $\Ad(\g)$-stable Zariski closed subset of $\Nil(\g^*)$. Let $\Irr_{\nu,\sfS}(G)$ denote the subset of $\Irr(G)$ consisting of  representations that have infinitesimal character $\nu$ and whose complex associated variety is  contained in $\sfS$.
Our aim is to count the set $\Irr_{\nu,\sfS}(G)$.

To do this, we form the subset $\Lambda$ of $\hha^*$, which is the $\nu$-translate of the root lattice. Let $W(\Lambda)$ be the integral Weyl group of $\nu$, which equals the stabilizer of $\Lambda$ in $W$. We consider the space $\Coh_{\Lambda}(\CK(G))$ of $\CK(G)$-valued  coherent families on $\Lambda$, which is a finite-dimensional representation of $W(\Lambda)$. This is called the coherent continuation representation. See \Cref{subsec:coherent}. 

Attached to $\sfS$, we will define a certain subset $\Irr_\sfS(W(\Lambda))$   of $\Irr(W(\Lambda))$  by using the Springer correspondence (see \Cref{defn:IrrS}). 

Let $\widetilde \Sp_{2n}(\R)$ ($n\in \bN:=\{0,1,2,\dots\}$) denote  the metaplectic double cover of the symplectic
 group $\Sp_{2n}(\R)$ (it does not split for $n>0$), to be called a real metaplectic group. Recall that $G$ is is said to be linear if it has a faithful finite-dimensional representation. 
 
 \begin{thm}\label{Mcounteq}
 We have the inequality 
  \begin{equation*}%\label{leq2}  
  \sharp(\Irr_{\nu,\sfS}(G)) \leq  \sum_{\sigma \in \Irr_\sfS(W(\Lambda))} [1_{W_{\nu}}: \sigma]\cdot [\sigma:\Coh_{\Lambda}(\CK(G))],
\end{equation*}
 where $1_{W_{\nu}}$ denotes the trivial representation of the stabilizer ${W_{\nu}}$ of $\nu$ in $W$. The equality holds if the Coxeter group $W(\Lambda)$ has no simple factor of type $F_4$, $E_6$, $E_7$, or $E_8$, and $G$ is linear or isomorphic to a real metaplectic group. 
  \end{thm}

Here and henceforth, $[\ : \ ]$ indicates the multiplicity of the first
(irreducible) representation in the second one, and $\sharp$ indicates the cardinality of a finite set.

For any (two-sided) ideal $I$ of $\cU(\fgg)$, let $\Irr_I(G)$ denote the subset of $\Irr(G)$ consisting of representations that are annihilated by $I$. Note that for the maximal ideal $I_{\nu}$, \[\Irr_{I_\nu}(G)=
\Irr_{\nu,\overline{\CO_\nu}}(G),\] where $\overline{\CO_\nu}$ is the associated variety of $I_\nu$. In this case, work of Lusztig \cite{Lu} 
and Barbasch-Vogan \cite{BVUni} imply that \[ [1_{W_{\nu}}: \sigma]\leq 1\]
for all $\sigma \in \Irr_{\overline{\CO_\nu}}(W(\Lambda))$, and 
the set 
\[
\{\sigma\in 
\Irr_{\overline{\CO_\nu}}(W(\Lambda))\mid [1_{W_{\nu}}: \sigma]\neq 0\}
\]
can be explicitly described by the Lusztig left cell $\LC_{\nu}$ (see \eqref{ll} for its definition).  Therefore \Cref{Mcounteq} has the following consequence.

 \begin{cor}\label{Mcountleq}
We have the inequality  
  \begin{equation*}%\label{boundc}
    \sharp(\Irr_{I_\nu}(G)) \leq \sum_{\sigma \in \LC_{\nu}}   [\sigma:\Coh_{\Lambda}(\CK(G))].
  \end{equation*}
     The equality holds if the Coxeter group $W(\Lambda)$ has no simple factor of type $F_4$, $E_6$, $E_7$, or $E_8$, and $G$ is linear or isomorphic to a real metaplectic group. 
\end{cor}

%\begin{remark} 
\Cref{Mcounteq} and \Cref{Mcountleq} will be proved in \Cref{sec:counting}. The equalities in \Cref{Mcounteq} and \Cref{Mcountleq} will be proved more generally for any $G$ 
satisfying \Cref{conjcell}.  
%\end{remark}

 \subsection{Special unipotent representations of real classical groups}
 \label{sec:defunip}

Suppose that $\star$ is one of the 10 labels, $G$ is a classical Lie group of type $\star$, and $\check G$ is a
complex Lie group, as in the following table ($n,p,q\in \bN$).
\[
  \begin{aligned}
    &\textrm{Label $\star$}&& \textrm{Classical Lie Group $G$} && \textrm{Langlands (or Metaplectic) Dual Group }\check G\\
    & A^\R&&\GL_n(\R) &&\GL_n(\C)\\
    & A^\bH&&\GL_{\frac{n}{2}}(\bH)\ \,  (n \textrm{ even}) && \GL_{n}(\C)\\
    & A &&\oU(p,q) &&\GL_{p+q}(\C)\\
    &\widetilde{A}&&\widetilde \oU(p,q) &&\GL_{p+q}(\C)/\{\pm 1_{p+q}\}\qquad(\textrm{$1_{p+q}$ is the identity matrix})\\
    & B &&\SO(p,q)\ \, (p+q\, \textrm{ odd}) && \Sp_{p+q-1}(\C)\\
    &D&&\SO(p,q)\  \, (p+q\, \textrm{ even}) &&\SO_{p+q}(\C)\\
    &C&&\Sp_{2n}(\R)&& \SO_{2n+1}(\C)\\
    &\widetilde{C}&&\widetilde \Sp_{2n}(\R) &&\Sp_{2n}(\C)\\
    &D^*&& \SO^*(2n) &&\SO_{2n}(\C)\\
    &C^* && \Sp(\frac{p}{2}, \frac{q}{2})  \ \, (p,q\, \textrm{ even}) &&\SO_{p+q+1}(\C)  
  \end{aligned}
\]

 In this table  $\tU(p,q)$ is
 the (linear) double cover of $\rU(p,q)$ defined by a square root of the
 determinant character. 
In the last column, $\check G$ is the
   the Langlands dual of  the complexification of $G$, except for  $G=\widetilde \Sp_{2n}(\R)$ or $\widetilde \oU(0,0)$.  
 In the case of $G=\widetilde \Sp_{2n}(\R)$  we replace the Langlands
dual by  $\Sp_{2n}(\C)$, call it the ``\textit{metaplectic}''
dual \cite{Weis,BMSZ0}, and still write $\check G$ (or $\check G_{\mathrm{mp}}$ with the
subscript $\mathrm{mp}$ for emphasis).
 
 Write $\check \g$ for the Lie algebra of $\check G$. Let $\check \hha$ denote the universal Cartan subalgebra of $\check \g$. In all
 cases we identify  $\hha^*$ with $\check \hha$ in the standard way,  except for $\wtC$ where the identification  is via the half of the trace form on the
 Lie algebra $\g=\s\p_{2n}(\C)$.  

 Let $\check \CO\subset \check \g$ be a nilpotent orbit, namely an $\Inn(\check \g)$-orbit in $\Nil (\check \g)$. Write $\mathbf d_{\check \CO}$ for the Young diagram attached to  $\check \CO$ (\cite{CM}). It determines the nilpotent orbit $\check \CO$, unless $\check \g=\o_{4k}(\C)$ ($k\in \bN^+:=\{1,2,3, \dots\}$) and all the row lengths are even. (In the latter case, $\check \CO$ is called very even, and there are precisely two nilpotent orbits with the same Young diagram as that of $\check \CO$.)
When there is no confusion, we will not distinguish between $\check \CO$ and  $\mathbf d_{\check \CO}$.

Let $\{\check e,\check h,\check f\}$ be an $\fsl_{2}$-triple in $\check \g$ representing $\ckcO$, namely $\check e\in \check \CO$. The element $\check h/2$ is semisimple, and its 
  $\Inn(\check \g)$-orbit is uniquely determined by $\ckcO.$ Using the identification
 \be\label{sse}
   \Inn(\check \g)\backslash  \{\textrm{semisimple element in $\check \g$}\}=W\backslash \ckhha=W\backslash \hha^*,
 \ee
we pick an element $\lambda_{\check \CO}\in \hha^*$ that represents the same element in  $ W\backslash \hha^*$  as $\check h/2$. 
As in \cite[Section 5]{BVUni}, $\lambda_{\check \CO}$ determines an infinitesimal character (namely an algebraic character of $\CZ(\g)$), also called the 
 infinitesimal character associated to $\ckcO$. Write $I_{\check \CO}:=I_{\star, \check \CO}:=I_{\lambda_{\check
    \CO}}$ for the  maximal ideal of $\CU(\g)$, only dependent on
$\check \CO$. (The subscript $\star$ is one of the labels in the table.)

 Let \[\CO:=d_{\mathrm{BV}}(\check \CO)\in \overline{\Nil}(\g^*)\] denote the Barbasch-Vogan dual of $\check \CO$, namely the unique Zariski open $\Ad(\g)$-orbit in the associated variety of $I_{\check \CO}$. See \cite[Appendix]{BVUni}, and \cite[Section 1]{BMSZ0} for the case of $\star = \widetilde{C}$. (This duality notion is a recasting in terms of representation theory of a 
   duality notion defined combinatorially by Spaltenstein in \cite{Spa}.) 
 %Note that a primitive ideal of $\CU(\g)$ equals the maximal ideal $I_{\check \CO}$ if and only if it has infinitesimal character  $\lambda_{\check \CO}$ and its associated variety  equals the Zariski closure of $\CO$ in $\g^*$.   

 Following Barbasch-Vogan \cite{BVUni}, define the set of the special unipotent representations of $G$
 attached to $\ckcO$ by
\[
 %\begin{equation}\label{eq:defuni}
   \begin{split}
     \Unip_{\ckcO}(G):=&  \Unip_{\star, \ckcO}(G) \\
     :=& \begin{cases}
             %\{\pi\in \Irr(G)\mid \pi \textrm{ is genuine  and annihilated by } I_{\check \CO}\}, 
             \Set{\pi  \in \Irr_{I_\ckcO}(G) | \pi \text{ is genuine}}, 
             & \text{if } \star\in \{\widetilde A, \widetilde C\};\\
       %\{\pi\in \Irr(G)\mid \pi \textrm{ is annihilated by } I_{\check \CO}\}, 
    \Irr_{I_\ckcO}(G), 
       & otherwise.\\
     \end{cases}
   \end{split}
% \end{equation}
\]
 Here ``genuine" means that the central subgroup $\{\pm 1\}$ of $G$, which is the kernel of the covering homomorphism $\widetilde \oU(p,q)\rightarrow  \oU(p,q)$ or $\widetilde \Sp_{2n}(\R)\rightarrow \Sp_{2n}(\R)$, acts on $\pi$ through the nontrivial character. As we noted earlier, $\Irr_{I_\ckcO}(G)=
\Irr_{\lambda _{\ckcO}, \overline{\CO}}(G)$.

\begin{remark} In representation theory literature, the set $\Unip_{\ckcO}(G)$ is also called the weak Arthur/ABV packet attached to $\check \CO$, which is the union of Arthur/ABV packets ranging over unipotent Arthur parameters whose restriction to $\SL_2(\bC)$ is determined by $\check \CO$ via the Jacobson-Morozov theorem. See \cite[Corollary 27.13]{ABV}. 
\end{remark}

The main goal of this paper and the second paper in the series \cite{BMSZ2} is to count the set $\Unip_{\check \CO}(G)$ and to construct all the representations in $\Unip_{\check \CO}(G)$, both explicitly. 

\medskip

Motivated by \cite{ABV}, we introduce, for a nilpotent orbit $\check \CO\subset \check \g$, a key property to be called  analytically even. For classical groups, this notion  coincides with M\oe glin's notion of good parity \cite{Mo11}. 
The notion will be used throughout the article and will play an important role in the classification of special unipotent representations (for unitary groups and for groups of type $BCD$).

Let $\{\check e_{\mathrm{pn}},\check h_{\mathrm{pn}},\check f_{\mathrm{pn}}\}$ be an $\fsl_{2}$-triple in $\check \g$ such that $\check e_{\mathrm{pn}}$ belongs to the principal nilpotent orbit. Then the element $\exp(\pi \sqrt{-1} \check h_{\mathrm{pn}})\in \check G$ is a central element fixed by all holomorphic automorphisms of $\check G$, and is independent of the choice of the triple  $\{\check e_{\mathrm{pn}},\check h_{\mathrm{pn}},\check f_{\mathrm{pn}}\}$. Here $\exp: \check \g\rightarrow \check G$ denotes the exponential map, and $\pi$ is the usual circumference ratio.

\begin{defn}\label{def:ae}
(a) A semisimple element $a\in \check \g$ is said to be analytically integral if 
\[
  \exp(2 \pi \sqrt{-1} a )=\exp(\pi \sqrt{-1} \check h_{\mathrm{pn}})
\]

\noindent (b) A nilpotent orbit $\check \CO\subset \check \g$ is said to be analytically even if $\check h/2\in \check \g$ is analytically integral. Here and as before, $\{\check e,\check h,\check f\}$ is an $\fsl_{2}$-triple in $\check \g$ representing $\ckcO$, namely $\check e\in \check \CO$. 

\end{defn}

\begin{remarks}
(i) \Cref{def:ae} obviously generalizes to all connected complex reductive groups.

\noindent (ii) Part (b) of the definition is independent of the choice of $\{\check e,\check h,\check f\}$. 

\noindent (iii) Every  analytically even nilpotent orbit is even, and when 
 the complex group $\check G$ is an adjoint group, every even nilpotent orbit is also analytically even. 
\end{remarks}

We also define a notion of good parity and bad parity (\`a la M\oe glin \cite[page 356]{Mo11}). 
Note that if $\star\in \{A^\R, A^\bH, C, \wtC, D^*\}$, then $n$ equals the rank of $\g$. Throughout the article, we will let $n$ denote the rank of $\g$, in all cases.

\begin{defn} \label{def:gpbp}
We call an integer to have good parity (with respect to $\star$ and $n =\rank \g$) if it has the same parity as
\begin{equation}\label{parity}
  \begin{cases}
    n, &  \text{if $\star \in \{A^\R, A^\BH,  A\}$}; \\
    1+ n, &  \text{if $\star = \wtA$}; \\
   1, & \text{if } \star \in \set{C,C^{*},D,D^{*}};\\
 \text{0}, & \text{if } \star \in \set{B,\wtC}.\\
  \end{cases}
\end{equation}
Otherwise, we say the integer has bad parity.
\end{defn}

The following lemma is easily checked. 

\begin{lem}
(a) Suppose that $\star\neq \wtA$. Then a nilpotent orbit $\check \CO\subset \check \g$ is analytically even if and only if all row lengths of its  Young diagram have good parity.

\noindent (b) Suppose that $\star=\wtA$. Then a nilpotent orbit $\check \CO\subset \check \g$ is analytically even if and only if all row lengths of its  Young diagram have the same parity (or equivalently, $\check \CO$ is even).
\end{lem}

\subsection{General linear groups}\label{sec:GLRH}

All special unipotent representations of $\GL_n(\bR)$ and $\GL_{\frac{n}{2}}(\bH)$  are obtained via normalized smooth parabolic induction from quadratic characters (see \cite{V.GL}*{Page 450}).
We will review their classifications in the framework of this article.

  For a Young diagram $\imath$, write
 \[
   \mathbf r_1(\imath)\geq \mathbf r_2(\imath)\geq \mathbf r_3(\imath)\geq \cdots
 \]
 for its row lengths, and similarly, write
 \[
   \mathbf c_1(\imath)\geq \mathbf c_2(\imath)\geq \mathbf c_3(\imath)\geq \cdots
 \]
 for its column lengths. Denote by
 $\abs{\imath}:=\sum_{i=1}^\infty \mathbf r_i(\imath)$ the total size of
 $\imath$.

For any Young diagram $\imath$, we introduce the set $\mathrm{Box}(\imath)$ of
boxes of $\imath$ %as the following subset of $\bN^+\times \bN^+$:
as follows: 
\begin{equation}\label{eq:BOX}
  \mathrm{Box}(\imath):=\Set{(i,j)\in\bN^+\times \bN^+| j\leq \bfrr_i(\imath)}.
\end{equation}
A subset of $\bN^+\times \bN^+$ of the form \eqref{eq:BOX} is also said to constitute the Young diagram $\imath$.

\renewcommand{\CP}{\mathcal{P}} We also introduce five symbols $\bullet$, $s$,
$r$, $c$ and $d$, and make the following definition.
\begin{defn}
 A painting on a Young diagram $\imath$ is an assignment (we
    place a symbol in each box)
  \[
    \mathcal P: \mathrm{Box}(\imath) \rightarrow \{\bullet, s, r, c, d \}
  \]
  with the following properties:
  \begin{enumerate}
\item If we remove the boxes painted with $\{s,r,c,d\}$, $\{r,c,d\}$, $\{c,d\}$, or $\{d\}$, the remainder still constitutes a Young diagram; 
      \item every row of $\imath$ has at most one box
          painted with $s$, and has at most one box painted with $r$; 
       \item every column of $\imath$ has at most one
          box painted with $c$, and has at most one box painted with $d$. 
  \end{enumerate}
A painted Young diagram is a pair $(\imath, \CP)$ consisting of a Young diagram $\imath$ and a painting $\CP$ on $\imath$.
\end{defn}

\begin{remark} 
%(i) The first requirement of a painted Young diagram means that the symbols $\bullet$, $s$,
%$r$, $c$ and $d$ should be painted in the order as indicated so that the painted boxes with any earlier painted symbols constitute a valid Young diagram throughout.  
(i)
The first requirement for a painted Young diagram says that 
the symbols $\bullet$, $s$, $r$, $c$ and $d$ should be painted in the order as indicated so that the painted boxes with any earlier painted symbols constitute a Young diagram throughout.  
%In other words, the boxes painted with the symbols $\set{\bullet}$, $\set{\bullet, s}$, $\set{\bullet, s, r}$, $\set{\bullet, s, r, c}$, or $\set{\bullet, s, r, c, d}$ must form a valid Young diagram. 

(ii) The specific requirements of painted Young diagrams are motivated by the branching rules of Weyl group representations (via the Littlewood-Richardson rule), and will be used to count special unipotent representations in Sections \ref{sec:GL} and \ref{sec:BCD}. 
\end{remark} 

\begin{eg}
  The following represents a painted Young diagram. 
   \[ \ytb{\bullet\bullet\bullet\bullet r,\bullet r d , sr,dd} 
  \]
    Each of the followings does not represent a painted Young diagram. 
   \[ \ytb{ dc}\qquad  \qquad \ytb{ cd,c}  \qquad \qquad \ytb{\bullet ss rd,rc}
  \]
\end{eg}

\begin{defn}\label{defpbp0}
  Suppose that $\star\in \{A^\R, A^\bH\}$. A painting $\CP$ on a Young diagram $\imath$ has type $\star$ if
  \begin{enumerate}
    \item the symbols of $\CP$ are in
          \[
          \left\{
          \begin{array}{ll}
            \{\bullet, c, d\}, &\hbox{if $\star=A^\R$}; \smallskip\\
            \{\bullet\}, &\hbox{if $\star=A^\bH$},                    \end{array}
        \right.
          \]
    \item every column of $\imath$ has an even number of boxes painted with $\bullet$.  
              \end{enumerate}
  Denote by $\PAP_\star(\imath)$ the set of paintings on $\imath^{t}$ that has type $\star$, where $\imath^{t}$
  is the transpose of $\imath$.
   \end{defn}

The middle letter $\mathrm A$ in $\PAP$ refers to the common $A$ in $\{A^\R, A^\bH\}$.

It is easy to check that if $\star=A^\R$, then
  \[
    \sharp(\PAP_\star(\check \CO))=\prod_{i\in \bN^+} (1+\textrm{the
      number of rows of length $i$ in $\check \CO$}),
  \]
  and if $\star=A^\bH$, then
  \[
    \sharp(\PAP_\star(\check \CO))= \left\{
      \begin{array}{ll}
        1, &\hbox{if all row lengths of $\check \CO$ are even}; \smallskip\\
        0, &\hbox{otherwise}.  \end{array}
    \right.
  \]
Here and henceforth $\sharp$ indicates the cardinality of a finite set.

If $\star=A^\R$, for every
$\CP\in \PAP_{\star}(\ckcO)$ we attach a representation $\pi_\CP$ of $G$ as in what follows. Let $P_\CP$ be  the standard parabolic subgroup of $G$
with Levi component 
\be\label{indgl}
\GL_{\bfrr_1(\ckcO)}(\R)\times \GL_{\bfrr_2(\ckcO)}(\R)\times\dots \times \GL_{\bfrr_k(\ckcO)}(\R),
\ee
where $k$ is the number of nonempty rows of $\check \CO$.   
On each factor $\GL_{\bfrr_j(\ckcO)}(\R)$, put the trivial character or the sign character  according to whether the $j$-th column
of $\CP$ ends in $\{\bullet,c\}$ or $\{d\}$. This yields a character of $P_\CP$ and let  $\pi_\CP$ be the resulting normalized induced representation of $G$.

Similarly, if  $\star=A^\bH$, for every
$\CP\in \PAP_{\star}(\ckcO)$ we attach a representation $\pi_\CP$ of $G$ as in  what follows. Let $P_\CP$ be  the standard parabolic subgroup of $G$
with Levi component 
$$
\GL_{\bfrr_1(\ckcO)/2}(\bH)\times \GL_{\bfrr_2(\ckcO)/2}(\bH)\times\dots \times \GL_{\bfrr_k(\ckcO)/2}(\bH),
$$
where $k$ is the number of nonempty rows of $\check \CO$.   
Put the trivial character on each factor $\GL_{\bfrr_j(\ckcO)/2}(\bH)$. This yields a character of $P_\CP$ and let  $\pi_\CP$ be the resulting normalized induced representation of $G$.

Then in both cases $\pi_\CP$ is irreducible and belongs to $\Unip_{\ckcO}(G)$ (see  \cite[Theorem 3.8]{V.GL} and \cite[Example~27.5]{ABV}).
We summarize the classifications by the following theorem. A proof using the methods of this paper is in \Cref{sec:GLRH11}. 

\begin{thm}
 [{\cf \cite{V.GL}}]\label{thm:mainR00}
Suppose that  $\star\in \{A^\R, A^\bH\}$. Then the map
  \[
    \begin{array}{ccc}
      \PAP_{\star}(\ckcO) & \rightarrow & \Unip_{\ckcO}(G),\\
      \CP & \mapsto & \pi_\CP
    \end{array}
  \]
  is bijective. 
\end{thm}

\subsection{Unitary groups}

Similar to \Cref{defpbp0}, we make the following definition.

\begin{defn}\label{defpbp1}
  Suppose that $\star\in \{A, \widetilde A\}$. A  painting $\CP$ on a Young diagram $\imath$ has type $\star$ if
  \begin{enumerate}
    \item the symbols of $\CP$ are in $\{\bullet, s, r\}$, and
                    \item every row of $\imath$ has an even number of boxes painted with $\bullet$. 
  \end{enumerate}
  Denote by $\PAP_\star(\imath)$ the set of 
 paintings on $\imath^{t}$ that has type $\star$, where $\imath^{t}$
  is the transpose of $\imath$.
   \end{defn}

Now suppose that $\imath$ is a Young diagram and $\CP$ is a painting on $\imath$
that has type $A$ or $\widetilde A$. Define the signature of $\CP$ to be the pair
\begin{equation}\label{eq:signature}
    (p_\CP, q_\cP): = \left (\frac{\sharp(\cP^{-1}(\bullet))}{2}+\sharp(\cP^{-1}(r)),\,
    \frac{ \sharp(\cP^{-1}(\bullet))}{2}+\sharp(\cP^{-1}(s))\right).
\end{equation}
\trivial[h]{ The first equation is the true definition of signature. The second
  one is an easy consequence of the definition of $\AC_\cP$. }

\begin{eg}
  Suppose
  that \[ \check \CO=\ytb{\ \ \ \ \ , \ \ \ , \ , \ , \ }\quad \textrm{and}\quad \CP=\ytb{\bullet\bullet\bullet\bullet r,\bullet\bullet , sr,s,r}\in \PAP_{A}(\check \CO) .
  \]
  Then $(p_\CP, q_\cP)=(6,5)$.

\end{eg}

Given two Young diagrams $\imath$ and $\jmath$, write $\imath\cuprow \jmath$ for
the Young diagram whose multiset of nonzero row lengths equals the union of
those of $\imath$ and $\jmath$. Also write $2\imath =\imath\cuprow \imath$.
Similarly, we write $\imath\cupcol \jmath$ for
the Young diagram whose multiset of nonzero column lengths equals the union of
those of $\imath$ and $\jmath$.

For unitary groups, we have the following counting result.
\begin{thm}\label{thmu1}
  Suppose that $\star\in \{A, \widetilde A\}$. If there is a Young diagram decomposition
  \be\label{deccou}
    \ckcO=\ckcOg \cuprow 2\ckcOpb
  \ee
 such that all nonzero row lengths of   $\ckcOg$ have good parity (namely the parity of $p+q$ when $\star = A$ and the parity of $p+q+1$ when $\star =\widetilde A$) and all nonzero row lengths of $\ckcOpb$ have bad parity,   then
  \[
    \sharp(\Unip_{\ckcO}(G))= \sharp \set{\CP\in \PAP_\star(\ckcOg)|(p_\CP+\abs{\ckcOpb}, q_\CP+\abs{\ckcOpb})=(p,q)}.
  \]
 Otherwise $\sharp(\Unip_{\check \CO}(G))=0$.

\end{thm}

In particular, when $\star=\widetilde A$ and $p+q$ is odd, the set $\Unip_{\check \CO}(\widetilde \oU(p,q))$ is empty.

Suppose that $\star\in \{A, \widetilde A\}$ so that $G=\oU(p,q)$ or $\widetilde \oU(p,q)$. By Theorem \ref{thmu1},
the set $\Unip_{\ckcO}(G)$ is empty unless there is a decomposition as in \eqref{deccou} such that
\be\label{decupq}
  p,q \geq \abs{\ckcOb'}.
\ee
Assume this holds, and define two groups
$G'_{\mathrm b}:=\GL_{\abs{\check \CO'_{\mathrm b}}}(\C)$ and
\[
  G_{\mathrm g} :=
  \begin{cases}
    \rU(p-\abs{\check \CO'_{\mathrm b}},q-\abs{\check \CO'_{\mathrm b}}),  & \text{if }\star = A;\\
    \tU(p-\abs{\check \CO'_{\mathrm b}},q-\abs{\check \CO'_{\mathrm b}}),  & \text{if }\star = \wtA.
\end{cases}
\]
Then up to conjugation there is a unique real parabolic subgroup $P$ of $G$ whose Levi quotient (namely its quotient by the unipotent radical) is naturally isomorphic to   $G'_{\mathrm b}\times G_{\mathrm g}$.
  Note that $\ckcO_{\mathrm g}$ has good parity with respect to $\star$ and $\abs{\ckcO_{\mathrm g}}$.

  Let $\pi_{\check \CO'_{\mathrm b}}$ denote the
unique element in $\Unip_{\check \CO'_{\mathrm b}}(G'_{\mathrm b})$ (see Section \ref{complex} for a review on complex classical groups). Then for
every $\pi_\mathrm g\in \Unip_{\ckcO_{\mathrm g}}(\Gg)$, the normalized smooth parabolic
induced representation (from $P$ to $G$) $\pi_{\ckcOpb}  \rtimes \pi_\mathrm g$ is irreducible by
\cite{Mat96}*{Theorem~3.2.2}, 
and is an element of $\Unip_{\ckcO}(G)$ (\cf
\cite{MR.U}*{Theorem~5.3}). (This also follows from
the calculation of the wavefront cycle of the induced representation (\cite[Corollary 5.0.10]{B.Orbit}.)
%See also \cite[Chapter 27]{ABV}. 

\begin{thm}\label{thmu2}
 With the assumptions in \eqref{deccou} and \eqref{decupq}, and notation as above, the parabolic induction map
  \begin{equation}
  \label{bij000}
    \begin{array}{rcl}
      \mathrm{Unip}_{\check{\CO}_\mathrm g}(G_\mathrm g)&\longrightarrow &\mathrm{Unip}_{\ckcO}(G),\\
      \pi_{\mathrm g}& \mapsto &  \pi_{\ckcOpb} \rtimes \pi_\mathrm g  \\
    \end{array}
  \end{equation}
  is bijective.
\end{thm}

\Cref{thmu1} and \Cref{thmu2} will be proved in Section \ref{secunit}. 
%See Section \ref{secunit} for additional results on special unipotent representations of unitary groups.

\subsection{Orthogonal and symplectic groups: relevant parabolic subgroups}
\label{secrgp0}
In this subsection and the next two subsections,  we assume that
$\star\in \set{ B, C, \wtC,C^*,D, D^*}$.
Then there are Young diagram decompositions
\be\label{decdo}
 \mathbf d_\ckcO=\mathbf d_\mathrm b\cuprow \mathbf d_\mathrm g\quad\textrm{and}\quad \mathbf d_\mathrm b=2 \mathbf d_\mathrm b'
\ee
 such that $\mathbf d_\mathrm b$ has bad parity in the sense that all its nonzero row
lengths have bad parity, and $\mathbf d_\mathrm g$ has good parity in the sense that all its
nonzero row lengths have good parity.
Put
\be\label{nb000}
  n_\mathrm b:=\abs{\mathbf d_\mathrm b'}. 
\ee

Write $V$ for the standard module of $ \g$, which is either a complex symmetric bilinear space or a complex symplectic space. Denote by $\mathsf
  G_\bC(V)$ the identity connected component of the isometry group of $V$, which is a complexification of $G$.  Recall that  $\CO\subset \g^*$ is the Barbasch-Vogan dual of $\check \CO$. Note that when $\star\in \{D, D^*\}$, $\check \CO$ is very even  (namely all its row lengths are even) if and only if  $\CO$ is very even.

\begin{lem}\label{relpa00}
Up to conjugation by $\mathsf
  G_\bC(V)$, there is a unique totally isotropic subspace $X_\mathcal \CO$ of $V$ of dimension $n_\mathrm b$ that satisfies the following condition:
   \begin{equation}\label{eq:rel} 
      \CO=\mathrm{Ind}^\g_{\gl(X_\mathcal O)\times \g(V_{\mathrm g})} \mathcal O_{\mathrm b}'\times \mathcal O_{\mathrm g}\qquad(\textrm{parabolic induction of nilpotent orbit}).
\end{equation}
 Here $V_\mathrm g$ is a non-degenerate subspace of $V$ of dimension $\dim V-2 n_\mathrm b$ that is perpendicular to $X_\mathcal O$, $\g(V_{\mathrm g})$ denotes the Lie algebra of the isometry group of $V_\mathrm g$,  $\gl(X_\mathcal O)$ denotes the general linear Lie algebra,  $\gl(X_\mathcal O)\times \g(V_{\mathrm g})$ is viewed as a Levi subalgebra of $\g$ as usual, $\mathcal O_{\mathrm b}'\subset (\gl(X_\mathcal O))^*$ is the nilpotent orbit whose Young diagram is the transpose of $\mathbf d'_\mathrm b$, and  $\mathcal O_{\mathrm g}\subset (\g(V_{\mathrm g}))^*$ is the Barbasch-Vogan dual of the nilpotent orbit corresponding to the Young diagram $\mathbf d_\mathrm g$.  
 % \end{enumerate}
 \end{lem}
\begin{proof} Excluding the case when $\star\in \{D, D^*\}$, $V\neq 0$, and $\check \CO$ is very even, there is a unique totally isotropic subspace $X_\mathcal \CO$ of $V$ of dimension $n_\mathrm b$, 
up to conjugation by $\mathsf G_\bC(V)$. The condition in \eqref{eq:rel} is automatic by \cite[Proposition A2, c)]{BVUni}. When $\star\in \{D, D^*\}$, $V\neq 0$, and $\check \CO$ is very even, there are exactly two totally isotropic subspaces of $V$ of dimension $n_\mathrm b$, 
up to conjugation by $\mathsf G_\bC(V)$. Denote by $\CO'\subset \g^*$ the other nilpotent orbit having the same Young diagram as that of $\CO$.
The two totally isotropic subspaces of $V$ of dimension $n_\mathrm b$ are distinguished by the requirement in \eqref{eq:rel}, i.e., 
relating $\CO$ with $X_{\CO}$, and 
$\CO'$ with $X_{\CO'}$, 
again by \cite[Proposition A2, c)]{BVUni}. 
\end{proof}

\begin{defn}\label{def:c-relevant2}
(a) A parabolic subalgebra
 $\p$ of $\g$ is said to be $\check \CO$-relevant if 
it is the stabilizer of a totally isotropic subspace $X_\CO$ of $V$ in  \Cref{relpa00}. 

\noindent (b) A parabolic subgroup of $G$ is said to be $\check \CO$-relevant if its complexified Lie algebra is $\check \CO$-relevant.

\noindent 
(c) The orbit $\check \CO$ is said to be $G$-relevant, if there there is a parabolic subgroup of $G$ which is $\check \CO$-relevant. 
\end{defn}

We summarize by the following proposition, which is clear. 

\begin{prop}\label{relpa0}
(a) Up to conjugation by $G$, there is at most one parabolic subgroup of $G$ that is $\check \CO$-relevant.

\noindent  (b) If $\star=D^*$, $V\neq 0$, and  $\check \CO$ is very even,   then there are precisely two orbits in $\overline \Nil(\check \g)$ that has the same Young diagram as that of $\check \CO$. Between these two orbits, exactly one is $G$-relevant. Excluding this special case, $\check \CO$ is $G$-relevant if and only if
\be\label{existgl}
  \textrm{either $\star\in \{B,D, C^*\}$ and $p,q\geq n_\mathrm b,\quad $ or  $\quad \star\in \{C, \wtC, D^*\}$}.
\ee
 \end{prop}

When \eqref{existgl} holds,  we put  \be\label{gg00}
  \Gg :=
  \begin{cases}
    \SO(p-n_\mathrm b,q-n_\mathrm b), & \textrm{if $\star\in \set{B,D}$};\\
      \SO^{*}(2n-2n_\mathrm b), &\textrm{if $\star = D^{*}$};\\
    \Sp_{2n-2n_\mathrm b}(\bR), &\textrm{if $\star = C$};\\
    \wtSp_{2n-2n_\mathrm b}(\bR), &\textrm{if $\star = \wtC$};\\
      \Sp(\frac{p-n_\mathrm b}{2},\frac{q-n_\mathrm b}{2}), &\textrm{if $\star = C^{*}$}.\\
  \end{cases}
\ee
Then the Levi quotient of  every $\check \CO$-relevant parabolic subgroup  of $G$ is naturally isomorphic to   $\Gpb\times \Gg $ (or   $(\Gpb\times \Gg)/\{\pm 1\} $ when $\star=\wtC$), where
\begin{equation}\label{Gpb}
  \Gpb := \begin{cases}
    \GL_{n_\mathrm b}(\bR), & \text{if } \star \in \set{B,C,D}; \\
       \widetilde{ \GL}_{n_\mathrm b}(\bR), & \text{if } \star =\wtC; \\
    \GL_{\frac{n_\mathrm b}{2}}(\bH), & \text{if } \star \in \set{C^{*},D^{*}}. \\
  \end{cases}
\end{equation}
Here $ \widetilde{ \GL}_{n_\mathrm b}(\bR)$ is the double cover of $ \GL_{n_\mathrm b}(\bR)$ that fits the following Cartesian diagram of Lie groups:
\begin{equation}\label{wgll}
\begin{CD}
 \widetilde{ \GL}_{n_\mathrm b}(\bR)@>>>  \GL_{n_\mathrm b}(\bR)\\
  @VVV @VV g\mapsto \textrm{ sign of $\det(g)$} V\\
  \{\pm 1, \pm \sqrt{-1}\} @> x\mapsto x^2 >> \{\pm 1\}. \\
\end{CD}
\end{equation}

Unless otherwise mentioned, we will use the corresponding lower case Gothic letter to denote the complexified Lie algebra of a Lie group. For example, $\g'_\mathrm b$ is the complexified Lie algebra of $\Gpb$. Let $\check \g'_\mathrm b$ denote the Langlands dual of $\g'_\mathrm b$, and let $\check \CO'_\mathrm b\in \overline{\Nil}(\check \g'_\mathrm b)$ denote the nilpotent orbit with Young diagram $\mathbf d'_\mathrm b$. Likewise let $\check \g_\mathrm g$ denote the Langlands dual (or the metaplectic Langlands dual when $\star=\wtC$) of $\g_\mathrm g$, and let $\check \CO_\mathrm g\in \overline{\Nil}(\check \g_\mathrm g)$ denote the (unique) nilpotent orbit with Young diagram $\mathbf d_\mathrm g$.

 \subsection{Orthogonal and symplectic groups: reduction to good parity}
\label{secrgp000}
Define %$I_{\ckcOpb}:= I_{A^\R, \ckcOpb}$ and 
 \[
      \Unip_{\ckcOpb}(\widetilde{ \GL}_{n_\mathrm b}(\bR)):=
       \{\pi\in \Irr_{I_{A^\R, \ckcOpb}}(\widetilde{ \GL}_{n_\mathrm b}(\bR))\mid \pi \textrm{ is genuine}\}.
       %\{\pi\in \Irr(\widetilde{ \GL}_{n_\mathrm b}(\bR))\mid \pi \textrm{ is genuine  and annihilated by } I_{\ckcOpb}:= I_{A^\R, \ckcOpb}\}.
       \]
        Here and as before, ``genuine" means that the central subgroup $\{\pm 1\}$ acts through the nontrivial character.  Then we have a bijective map
 \[
    \Unip_{\ckcOpb}(\GL_{n_\mathrm b}(\bR))\rightarrow  \Unip_{\ckcOpb}(\widetilde{ \GL}_{n_\mathrm b}(\bR)), \quad \pi\mapsto \pi\otimes \tilde \chi_{n_\mathrm b},
 \]
 where $\tilde \chi_{n_\mathrm b}$ is the character given by the left vertical arrow of \eqref{wgll}.

In Section~\ref{sec:proof.BCDred}, we will prove the following theorem.
\begin{thm}\label{reduction}
 If  $G$ has an $\check \CO$-relevant parabolic subgroup $P$, then the normalized smooth parabolic induction from $P$ to $G$ yields
   a bijection
   \[
 % \begin{equation}\label{eq:IND}
    \begin{array}{rccc}
 &\Unip_{\ckcO'_{\mathrm b}}( G'_{\mathrm b}) \times   \Unip_{\ckcO_{\mathrm g}}( G_{\mathrm g})  &         \longrightarrow &\Unip_{\ckcO }(G), \\
                &   (\pi',\pi_\mathrm g) & \mapsto & \pi'\rtimes \pi_\mathrm g.
    \end{array}
 % \end{equation}
 \]
  Otherwise,
  \[
    \Unip_{\ckcO}(G)=\emptyset.
  \]
\end{thm}

By \Cref{reduction}, we have the more specific results on counting as follows.
\begin{enumerate}[label=(\alph*)]
  \item Assume that $\star\in \{B,D\}$ so that $G=\SO(p,q)$. Then
        \[
        \sharp(\Unip_{\check \CO}(G))=
        \begin{cases}
          \sharp(\Unip_{\check \CO_{\mathrm g}}(G_{\mathrm g}))\times \sharp(\Unip_{\check \CO'_{\mathrm b}}(\GL_{n_\mathrm b}(\R)) ), &\hbox{if $p,q\geq n_\mathrm b$}; \smallskip\\
          0, &\hbox{otherwise.}
        \end{cases}
        \]
  \item Assume that $\star=C^*$ so that $G=\Sp(\frac{p}{2},\frac{q}{2})$. Then
        \[
        \sharp(\Unip_{\check \CO}(G))=
        \begin{cases}
          \sharp(\Unip_{\check \CO_{\mathrm g}}(G_{\mathrm g} )), &\hbox{if $p,q\geq n_\mathrm b$}; \smallskip\\
          0, &\hbox{otherwise.}
        \end{cases}
        \]

  \item Assume that $\star\in \{C,\widetilde C\}$ so that $G=\Sp_{2n}(\R)$ or
        $\widetilde \Sp_{2n}(\R)$. Then
        \[
        \sharp(\Unip_{\check \CO}(G))= \sharp(\Unip_{\check \CO_{\mathrm g}}(G_{\mathrm g}))\times \sharp(\Unip_{\check \CO'_{\mathrm b}}(\GL_{n_\mathrm b}(\R)) ). \]
  \item Assume that $\star =D^*$ so that $G=\SO^*(2n)$. Then
        \[
          \sharp(\Unip_{\check \CO}(G))=
          \begin{cases}
          0,&\quad \textrm{if $\check \CO$  is not $G$-relevant;}\\
          \sharp(\Unip_{\check \CO_{\mathrm g}}(G_{\mathrm g})),&\quad \textrm{otherwise}.
          \end{cases}
        \]
\end{enumerate}

 \subsection{Orthogonal and symplectic groups: the case of good parity}\label{secorgp0}
 We now assume that $\check \CO$ has good parity, namely
 $\check \CO=\check \CO_{\mathrm g}$. By Theorem \ref{reduction}, the counting
 problem in general is reduced to this case.

\begin{defn}\label{defn:PP}
  A $\star$-pair is a pair $(i,i+1)$ of consecutive positive integers such that
  \[
    \left\{
      \begin{array}{ll}
        i\textrm{ is odd}, \quad &\textrm{if $\star\in\{C, \widetilde{C}, C^*\}$};  \\
        i \textrm{ is even}, \quad &\textrm{if $\star\in\{B, D, D^*\}$}. \\
      \end{array}
    \right.
  \]
  A $\star$-pair $(i,i+1)$ is said to be
  \begin{itemize}
    \item vacant in $\check \CO$, if
          $\mathbf r_i(\check \CO)=\mathbf r_{i+1}(\check \CO)=0$;
    \item balanced in $\check \CO$, if
          $\mathbf r_i(\check \CO)=\mathbf r_{i+1}(\check \CO)>0$;
    \item tailed in $\check \CO$, if
          $\mathbf r_i(\check \CO)-\mathbf r_{i+1}(\check \CO)$ is positive and
          odd;
    \item primitive in $\check \CO$, if
          $\mathbf r_i(\check \CO)-\mathbf r_{i+1}(\check \CO)$ is positive and
          even.
  \end{itemize}
  Denote $\CPP_\star(\check \CO)$ the set of all $\star$-pairs that are
  primitive in $\check \CO$.
\end{defn}

\begin{remark} When $\star\neq \wtC$, the power set of $\CPP_\star(\check \CO)$ gives another description of Lusztig's canonical
  quotient attached to $\ckcO$. The set $\CPP_{\star}(\ckcO)$ appears implicitly in
\cite{So}*{Section~5}.
\end{remark}

We give algorithms of how to explicitly count special unipotent representations for each type. First we attach to $\check \CO$ a pair of Young diagrams
\be\label{ijo}
  (\imath_{\check \CO}, \jmath_{\check \CO}):=(\imath_\star(\check \CO), \jmath_\star(\check \CO)),
\ee
as follows. (This is to be viewed as a Weyl group representation. See \Cref{sec:LCBCD}.)

Recall for a Young diagram $\imath$, the row and column lengths are written as 
 \[
   \mathbf r_1(\imath)\geq \mathbf r_2(\imath)\geq \mathbf r_3(\imath)\geq \cdots \quad \text{and} \quad
 \mathbf c_1(\imath)\geq \mathbf c_2(\imath)\geq \mathbf c_3(\imath)\geq \cdots . 
 \]

\medskip

\noindent {\bf The case when $\star=B$.} In this case, the nilpotent orbit $\check \CO$ is of type $C$, and has even rows
  only. Define
\[
  \mathbf c_{1}(\jmath_{\check \CO})=\frac{\mathbf r_1(\check \CO)}{2},
\]
and for all $i\geq 1$,
\[
  \left (\mathbf c_{i}(\imath_{\check \CO}), \mathbf c_{i+1}(\jmath_{\check \CO})\right )= \left (\frac{\mathbf r_{2i}(\check \CO)}{2}, \frac{\mathbf r_{2i+1}(\check \CO)}{2}\right ).
\]
(In words: 
  The largest row of size $\mathbf r_1(\check \CO)$ contributes a column of
  size $\frac{\mathbf r_1(\check \CO)}{2}$ to $\jmath_{\check \CO}$. Pair up the remaining
  rows, adding a row of size $0$ if there are an even number of
  nonzero rows. Each pair of rows $(\mathbf r_{2i}(\check \CO),\mathbf r_{2i+1}(\check \CO))$
    contributes a column of size $\frac{\mathbf r_{2i}(\check \CO)}{2}$ to
    $\imath_{\check \CO}$ and $\frac{\mathbf r_{2i+1}(\check \CO)}{2}$ to $\jmath_{\check \CO}$.)
    
\medskip

\noindent {\bf The case when $\star=\widetilde C$.} In
this case, the nilpotent
orbit $\check \CO$ is of type $C$, and has even rows only. Define, for all $i\geq 1$,
\[
  (\mathbf c_{i}(\imath_{\check \CO}), \mathbf c_{i}(\jmath_{\check \CO}))= \left (\frac{\mathbf r_{2i-1}(\check \CO)}{2}, \frac{\mathbf r_{2i}(\check \CO)}{2}\right).
\]
(In words: Pair up rows, adding a row of size 0 if there are an odd number
of nonzero rows. Each pair of rows $(\mathbf r_{2i-1}(\check \CO),\mathbf r_{2i}(\check \CO))$
    contributes a column of size $\frac{\mathbf r_{2i-1}(\check \CO)}{2}$ to
    $\imath_{\check \CO}$ and $\frac{\mathbf r_{2i}(\check \CO)}{2}$ to $\jmath_{\check \CO}$.)

\medskip

\noindent {\bf The case when $\star=\{ C,C^*\}$.} In this case, the nilpotent orbit $\check \CO$ is of type $B$, and has odd   rows only. Define, for all
$i\geq 1$,
\[
  (\mathbf c_{i}(\jmath_{\check \CO}), \mathbf c_{i}(\imath_{\check \CO}))= \left\{
    \begin{array}{ll}
      (0,  0), &\hbox{if $(2i-1, 2i)$ is vacant  in $\check \CO$};\smallskip\\
      (\frac{\mathbf r_{2i-1}(\check \CO)-1}{2},  0), & \hbox{if $(2i-1, 2i)$ is tailed in $\check \CO$};\smallskip\\
      (\frac{\mathbf r_{2i-1}(\check \CO)-1}{2},  \frac{\mathbf r_{2i}(\check \CO)+1}{2}), &\hbox{otherwise}.\\
    \end{array}
  \right.
\]
(In words: The number of nonzero rows is odd, say $2k+1$. 
  Pair them up, by adding a row of size zero. 
  %label them  $1,2\dots$ in decreasing order, 
  For each pair of rows $(\mathbf r_{2i-1}(\check \CO),\mathbf r_{2i}(\check \CO))$, add
   a column of size $\frac{\mathbf r_{2i-1}(\check \CO)-1}{2}$ to
  $\jmath_{\check \CO}$ and  $\frac{\mathbf r_{2i}(\check \CO)+1}{2}$ to $\imath_{\check \CO}$. For the last pair $(\mathbf r_{2k+1}(\check \CO),0)$, only add a column of size $\frac{\mathbf r_{2k+1}(\check \CO)-1}{2}$ to $\jmath_{\check \CO}$.)  

\medskip

\noindent {\bf The case when $\star\in \{D,D^*\}$.} In this case, the nilpotent orbit $\check \CO$ is of type $D$, and has odd  rows only. Define
\[
  \mathbf c_{1}(\imath_{\check \CO})= \left\{
    \begin{array}{ll}
      0,  &\hbox{if $\mathbf r_1(\check \CO)=0$}; \smallskip\\
      \frac{\mathbf r_1(\check \CO)+1}{2},   &\hbox{if $\mathbf r_1(\check \CO)>0$},\\
    \end{array}
  \right.
\]
and for all $i\geq 1$,
\[
  (\mathbf c_{i}(\jmath_{\check \CO}), \mathbf c_{i+1}(\imath_{\check \CO}))= \left\{
    \begin{array}{ll}
      (0,  0), &\hbox{if $(2i, 2i+1)$ is vacant in $\check \CO$};\smallskip\\
      \left  (\frac{\mathbf r_{2i}(\check \CO)-1}{2},  0\right ), & \hbox{if $(2i, 2i+1)$ is tailed in $\check \CO$};\smallskip\\
      \left  (\frac{\mathbf r_{2i}(\check \CO)-1}{2},  \frac{\mathbf r_{2i+1}(\check \CO)+1}{2}\right ), &\hbox{otherwise}.\\
    \end{array}
  \right.
\]
(In words: The number of nonzero rows is even, say $2k$. 
  The largest row of size $\mathbf r_1(\check \CO)$ contributes a column of
  size $\frac{\mathbf r_1(\check \CO)+1}{2}$ to $\imath_{\check \CO}$. Pair up the remaining
  rows, by adding a row of size $0$. For each pair of rows $(\mathbf r_{2i}(\check \CO),\mathbf r_{2i+1}(\check \CO))$, add a column of size $\frac{\mathbf r_{2i}(\check \CO)-1}{2}$ to
    $\jmath_{\check \CO}$ and $\frac{\mathbf r_{2i+1}(\check \CO)+1}{2}$ to $\imath_{\check \CO}$. For the last pair $(\mathbf r_{2k}(\check \CO),0)$, only add a column of size 
    $\frac{\mathbf r_{2k}(\check \CO)-1}{2}$ to $\jmath_{\check \CO}$.)

\begin{eg} Suppose that $\star=C$, and $\check \CO$ is the following Young
  diagram which has good parity.
  \begin{equation*}
    \tytb{\ \ \ \ \  , \ \ \  , \ \ \ , \ \ \  , \ \ \ , \  ,\  }
  \end{equation*}
  Then
  \[
    \CPP_\star(\check \CO)=\{(1,2), (5,6)\}
  \]
  and
  \[
    (\imath_{\check \CO}, \jmath_{\check \CO})= \tytb{\ \ \ ,\ \ } \times \tytb{\ \ \ , \ }.
  \]

\end{eg}

Here and henceforth, when no confusion is possible, we write
$\alpha\times \beta$ for a pair $(\alpha, \beta)$. Likewise write
$\alpha\times \beta\times \gamma$ for a triple $(\alpha, \beta, \gamma)$.

We now introduce the key combinatorial construct of this article, called a painted bipartition. We need to introduce two more labels $B^+$ and $B^-$ (over the label $B$; they play a slightly different role when counting the signature; see \eqref{ptqt}). 

\begin{defn} \label{def:pbp1}
  A painted bipartition is a triple
  $\uptau=(\imath, \CP)\times (\jmath, \cQ)\times \gamma$, where $(\imath, \CP)$
  and $ (\jmath, \mathcal Q)$ are painted Young diagrams, and
  $\gamma\in \{B^+,B^-, C,D,\widetilde {C}, C^*, D^*\}$, subject to the
  following conditions:
  \begin{enumerate}
          \delete{\item $(\imath, \jmath)\in \mathrm{BP}_\gamma$ if
          $\gamma\notin\{B^+,B^-\}$, and $(\imath, \jmath)\in \mathrm{BP}_{B}$
          if $\alpha\in\{B^+,B^-\}$;}
 \item $\CP$ and $\mathcal Q$ have the identical set of boxes painted with $\bullet$;
    \item the symbols of $\CP$ are in 
          \[
          \left\{
          \begin{array}{ll}
            \{\bullet, c\}, &\hbox{if $\gamma=B^+$ or $B^-$}; \smallskip\\
            \{\bullet,  r, c,d\}, &\hbox{if $\gamma=C$}; \smallskip\\
            \{\bullet, s, r, c,d\}, &\hbox{if $\gamma=D$}; \smallskip\\
            \{\bullet, s, c\}, &\hbox{if $\gamma =\widetilde{ C}$}; \smallskip \\
            \{\bullet\}, &\hbox{if $\gamma=C^*$}; \smallskip \\
            \{\bullet, s\}, &\hbox{if $\gamma=D^*$},\\
          \end{array}
          \right.
          \]
    \item the symbols of $\mathcal Q$ are in           \[
          \left\{
          \begin{array}{ll}
            \{\bullet, s, r, d\}, &\hbox{if $\gamma=B^+$ or $B^-$}; \smallskip\\
            \{\bullet, s\}, &\hbox{if $\gamma=C$}; \smallskip\\
            \{\bullet\}, &\hbox{if $\gamma=D$}; \smallskip\\
            \{\bullet, r, d\}, &\hbox{if $\gamma=\widetilde{ C}$}; \smallskip\\
            \{\bullet, s,r\}, &\hbox{if $\gamma=C^*$}; \smallskip \\
            \{\bullet, r\}, &\hbox{if $\gamma=D^*$}.
          \end{array}
          \right.
          \]

  \end{enumerate}
\end{defn}

\medskip
\begin{remark} In the above definition, the pairs $(\imath,\jmath)$ parametrize irreducible
representations of the classical Weyl grous $W(B)\cong W(C)$ (of type $B$ or $C$). We give some informal explanation on how painted bipartitions
arise. For each
conjugacy class of Cartan subgroups $H$ of $G$, there is a subgroup
$W(H)=\mathbf{H}_t\times W_s\times W_r\times W_c\times W_d$ in $W(\gamma)$ (roughly the integral Weyl group of the real classical group labeled by $\gamma$), a
1-dimensional character $\chi_H$ of $W(H)$, and the painted bipartitions will count the
multiplicities of a certain $W(\gamma)$-representation in $\Ind_{W(H)}^{W(\gamma)}\chi_H$. There is another subgroup $W_{\bullet}$ containing 
$\mathbf{H}_t$ and $W_\bullet,W_s,W_r,W_c,W_d$ are all
Weyl groups of classical type. They depend on the type of the group being considered, and the Cartan subgroup $H$. The induction from $\mathbf{H}_t$
to $W_\bullet$ (of a certain quadratic character) is encoded by the painting with the symbol $\bullet$. The character
$\chi_H$ is $sgn $ on $W_s, W_r$ and trivial on $W_c,W_d$. The
painting with the symbols $s,r,c,d$ reflects (use of) the Littlewood-Richardson rule (more precisely the Pieri rule). See 
\Cref{subsec:counting} for details. 
\end{remark}

For any painted bipartition $\uptau$ as in Definition \ref{def:pbp1}, we write
\[
  \imath_\uptau:=\imath,\ \cP_\uptau:=\cP,\ \jmath_\uptau:=\jmath,\ \cQ_\tau:=\cQ,\ \gamma_\uptau:=\gamma,\ \abs{\uptau}:=\abs{\imath}+\abs{\jmath},
\]
and
\[
  \star_\uptau:= \left\{
    \begin{array}{ll}
      B, &\hbox{if $\gamma=B^+$ or $B^-$}; \smallskip\\
      \gamma, & \hbox{otherwise}.           \end{array}
  \right.
\]

We further define a pair $\Sign(\uptau):=(p_{\uptau}, q_{\uptau})$ of natural numbers given by the
following recipe:
\begin{enumerate}[label=(\alph*)]
  \item If $\star_\uptau\in \{B, D, C^*\}$, then $(p_\uptau, q_\uptau)$ is given by
        counting the various symbols appearing in $(\imath, \CP)$,
        $(\jmath, \cQ)$ and $\{\alpha\}$ :
        \begin{equation}\label{ptqt}
          \left\{
            \begin{array}{l}
              p_\uptau :=( \# \bullet)+ 2 (\# r) +(\# c )+ (\# d) + (\# B^+);\smallskip\\
              q_\uptau :=( \# \bullet)+ 2 (\# s) + (\# c) + (\# d) + (\# B^-).\\
            \end{array}
          \right.
        \end{equation}
        Here
        \[
        \#\bullet:=\#(\cP^{-1}(\bullet))+\#(\cQ^{-1}(\bullet)),
        %\qquad (\textrm{$\#$        indicates the cardinality of a finite set}),
        \]
        the total number of boxes painted with $\bullet$ in $\cP$ and $\cQ$, and the other terms are defined in the obvious way. 
  \item If $\star_\uptau\in \{C, \widetilde C, D^*\}$, then
        $p_\uptau:=q_\uptau:=\abs{\uptau}$.
\end{enumerate}
\smallskip

We also define a classical group
\begin{equation*}%\label{def:Gt}
  G_\uptau:=
  \begin{cases}
    \SO(p_\uptau, q_\uptau), &\hbox{if $\star_\uptau=B$ or $D$}; \smallskip\\
    \Sp_{2\abs{\uptau}}(\R), &\hbox{if $\star_\uptau=C$}; \smallskip\\
    \widetilde{\Sp}_{2\abs{\uptau}}(\R), &\hbox{if $\star_\uptau=\widetilde{ C}$}; \smallskip \\
    \Sp(\frac{p_\uptau}{2}, \frac{q_\uptau}{2}), &\hbox{if $\star_\uptau=C^*$}; \smallskip \\
    \SO^*(2\abs{\uptau}), &\hbox{if $\star_\uptau=D^*$}.\\
  \end{cases}
\end{equation*}

Define
\begin{equation}\label{defpbp2222}
  \PBP_\star(\check \CO) :=\set{ \uptau\textrm{ is a painted
      bipartition} \mid \star_\uptau = \star, \text{ and
    } (\imath_\uptau,\jmath_\uptau) = (\imath_{\check \CO}, \jmath_{\check \CO})},
\end{equation}
and
\begin{equation} \label{defpbp3}
    \PBP_{G}(\ckcO) :=\set{\uptau\in \PBP_{\star}(\ckcO)| G_{\uptau} = G}.
\end{equation}
Here ``$G_{\uptau} = G$" amounts to saying that $(p_\uptau, q_\uptau)=(p,q)$ if $\star\in \{B, D, C^*\}$.

\begin{eg} Suppose that $\star=B$ and
  \[
    \check \CO =\tytb{\ \ \ \ \ \ , \ \ \ \ \ \ , \ \ , \ \ , \ \ }.
  \]
  Then
  \[
    \uptau:= \tytb{\bullet \bullet ,\bullet , c } \times \tytb{\bullet \bullet d ,\bullet , d }\times B^+\in \PBP_{\star}(\check \CO),
  \]
  and
  \[
    G_\uptau=\SO(10,9).
  \]
\end{eg}

We now state our final result on the explicit counting of special unipotent representations. It will be proved in \Cref{subsec:counting}. 

\begin{thm}\label{countup}
  Assume that $\star\in \{B, C,D,\widetilde {C}, C^*, D^*\}$, and $\check \CO$ has good parity. Then
 \[
   \sharp(\Unip_{\check \CO}(G)) =
    \left\{
    \begin{array}{ll}
       \sharp (\PBP_{G}(\ckcO)),  & \hbox{if $\star\in \{C^*,D^*\}$}; \smallskip\\
       2^{\sharp(\CPPs(\check \CO))} \cdot \sharp (\PBP_{G}(\ckcO)),  &\hbox{if $\star\in \{B, C,D,\widetilde {C}\}$}.
    \end{array}
  \right.
  \]
\end{thm}

In \cite{BMSZ2}, the authors will construct all representations in
$\Unip_{\check \CO}(G)$ by the method of theta lifting, 
when $\check \CO$ has good parity. See \cite[Section 3]{BMSZ2}. 

\medskip

We shall illustrate the contents of \Cref{countup} for the case $\star=C$. 
The nilpotent orbit $\check \CO$ is of type $B$ and all rows have odd sizes. The number of nonzero rows is odd, say $2k+1$. Label them as 
$$
2r_{1}+1\ge 2r_{2}+1\ge \dots\ge 2r_{2k+1}+1>0.
$$
Associate a pair of Young diagrams $(\imath_{\check \CO}, \jmath_{\check \CO})$ with columns  
$$
(r_{2}+1,r_{4}+1,\dots ,r_{2k}+1) \times (r_1,r_3,\dots ,r_{2k-1}, r_{2k+1}).
$$
We will actually associate a larger set of such pairs of Young diagrams to $\check \CO$ (see \Cref{sec:LCBCD} for details). Line up the numbers $r_1, r_2, ..., r_{2k+1}$ as follows: 
$$
(r_{1},r_{2})\dots (r_{2k-1}, r_{2k})(r_{2k+1}). 
$$
The set $\CPPs(\check \CO)$ consists of pairs $(2i-1,2i)$ such that $r_{2i-1} >r_{2i}$ (where $1\leq i\leq k$). 
For each subset $\wp \subseteq \CPPs(\check \CO)$, form a new pair of Young diagrams by
replacing the pair of numbers $(r_{2i}+1, r_{2i-1})$ in $(\imath_{\check \CO}, \jmath_{\check \CO})$ by $(r_{2i-1}+1, r_{2i})$, where $(2i-1,2i)\in \wp$. 
In this way we obtain 
$2^{\sharp(\CPPs(\check \CO))}$ number of pairs of Young diagrams associated to $\check \CO;$
they parametrize the so-called Lusztig left cell associated to $\check \CO$ (see \Cref{lem:Lcell}). 

For a pair of Young diagrams $(\imath, \jmath)$ as in the above, we associate the set of painted bipartitions  $\uptau=(\imath, \CP)\times (\jmath, \cQ)\times C$, by the painting rule specified in \Cref{def:pbp1} for $\gamma =C$. 
It turns out that
the 
set of 
painted bipartitions associated to each of the $2^{\sharp(\CPPs(\check \CO))}$ pairs of 
Young diagrams shares the same cardinality (see \Cref{MainComb}). \Cref{countup} says that the total number of painted bipartitions described above gives the number of special unipotent representations associated to $\check \CO$, for $G=\Sp_{2n}(\bR)$. 

\begin{eg} Suppose that $\star=C$, $G=\Sp_4(\R)$,  and
  \[
    \check \CO =\tytb{\ \ \ , \  , \ }.
  \]
  Then $\mathrm{PP}_\star(\check \CO)=\{(1,2)\}$, and the set $\PBP_{G}(\ckcO)$ has $4$ elements in total: 
  \[
\tytb{\bullet} \times \tytb{\bullet }\times C, \ \ \  \tytb{r} \times \tytb{s}\times C , \ \ \  \tytb{c} \times \tytb{s}\times C , \ \ \  \tytb{d} \times \tytb{s}\times C. 
  \]
 Thus there are precisely $8$ special unipotent representations of $G$ attached to $\check \CO$. 
\end{eg}

\medskip

\subsection{The case of complex classical groups}\label{complex}
Special unipotent representations of complex classical groups are all
well-understood (\cite{BVUni}, \cite{B89}) and known to be
unitarizable (\cite{B89}). We briefly review their counting and constructions in what follows. As the methods of this paper and \cite{BMSZ2} work for complex classical groups as well, we will present the results in the complex case parallel to those of this paper and \cite{BMSZ2}. For this subsection, we introduce five more labels $A^\C, B^\C,D^\C, C^\C$, and $\widetilde C^\C$, and let $\star$ be one of them. Let $G$ be a complex classical group of type $\star$, namely
\[
G=\GL_n(\C),\quad \SO_{2n+1}(\C),\quad \SO_{2n}(\C),\quad \Sp_{2n}(\C), \quad \textrm{or}\quad \Sp_{2n}(\C)\qquad (n\in \BN),
\]
respectively.
Let $\g_0$ denote the Lie algebra of $G$, which is  a complex Lie algebra.
The Langlands dual (or metaplectic Langlands dual when $\star=\widetilde C^\C$) $\check \g_0$ of $\g_0$ is respectively defined to be
\[
\g\l_n(\C),\quad \s\p_{2n}(\C), \quad \o_{2n}(\C), \quad \o_{2n+1}(\C), \quad \textrm{or}\quad  \s\p_{2n}(\C).
\]
 Let $\check \CO\in \overline{\Nil}(\check \g_0)$. As in the real case we have a maximal ideal $I_{\check \CO}:=I_{\star, \check \CO}$ of $\CU(\g_0)$.

 Write $\overline \g_0$ for the complex Lie algebra equipped with a conjugate linear isomorphism $\bar{\phantom a} :\g_0\rightarrow \overline{\g_0}$. The latter induces a  conjugate linear isomorphism $\bar{\phantom a} :\CU(\g_0)\rightarrow \CU( \overline{\g_0})$. Note that $\g_0\times \overline{\g_0}$ equals the complexified Lie algebra $\g$ of $G$.  Define the set of special unipotent representations of $G$
 attached to $\ckcO$ by
 \[
     \Unip_{\ckcO}(G):=  \Unip_{\star, \ckcO}(G)
     :=
       \{\pi\in \Irr(G)\mid \pi \textrm{ is annihilated by } I_{\check \CO}\otimes \CU(\overline{\g_0}) + \CU(\g_0)\otimes \overline{I_{\check \CO}}\, \}.
       \]

 If $\star=A^\C$ so that $G=\GL_n(\C)$, then $\Unip_{\ckcO}(G)$ is a singleton whose unique element is given by the normalized smooth parabolic induction
 $\Ind_{P}^{G} 1_P$, where $P$ is the standard parabolic subgroup whose Levi component equals
 \[
 \GL_{\mathbf r_1(\check \CO)}(\C)\times \GL_{\mathbf r_2(\check \CO)}(\C)\times \dots \times \GL_{\mathbf r_{\mathbf c_1(\check \CO)}(\check \CO)}(\C),
 \]
 and $1_P$ denotes the trivial representation of $P$.
 %See \cite{V.GL}.

Now suppose that $\star\in \{B^\C,D^\C, C^\C, \widetilde C^\C\}$. Write
\[
 \mathbf d_\ckcO=\mathbf d_\mathrm b\cuprow \mathbf d_\mathrm g\quad\textrm{and}\quad \mathbf d_\mathrm b=2 \mathbf d_\mathrm b'
\]
as in \eqref{decdo}, and put $n_\mathrm b:=\abs{\mathbf d_\mathrm b'}$ as before. %\quad\textrm{and}\quad n_\mathrm g:=n-n_\mathrm b.
Let $\p_0$ be a parabolic subalgebra of $\g_0$ that is $\check \CO$-relevant (defined in \Cref{secrgp0}). Let $P$ be the parabolic subgroup of $G$ with Lie algebra $\p_0$. Then the Levi quotient of $P$  is naturally isomorphic to $\Gpb\times \Gg$, where $\Gpb:=\GL_{n_\mathrm b}(\C)$ and
\[
  \Gg :=
  \begin{cases}
    \SO_{2n-2n_\mathrm b+1}(\C), & \textrm{if $\star=B^\C$};\\
    \SO_{2n-2n_\mathrm b}(\C), & \textrm{if $\star=D^\C$};\\
    \Sp_{2n-2n_\mathrm b}(\C), &\textrm{if $\star \in\{ C^\C, \widetilde C^\C \}$}.
      \end{cases}
\]

Define the set $\CPP_\star(\ckcOg)$ as in the real case. Then by the work of
Barbasch-Vogan (\cite[Corollary 5.29]{BVUni}, integral case) 
and Barbasch (\cite{B89}, general case; see also \cite[Theorem~6.12 and Theorem~10.1]{MR.C}),  
we have that
   \[
    \sharp(\Unip_{\check \CO}(G))=\sharp(\Unip_{\ckcOg}(\Gg))=2^{\sharp(\CPP_\star(\ckcOg))}.
  \]
   As in the real case, every representation in $\Unip_{\check \CO}(G)$ is
   obtained through irreducible parabolic induction from $P$ to $G$ via those of
   $ \Unip_{\ckcO'_{\mathrm b}}( G'_{\mathrm b})\times \Unip_{\ckcO_{\mathrm g}}( G_{\mathrm g}) $
   (\cf Theorem \ref{reduction}), and every representation in
   $\Unip_{\ckcOg}(\Gg)$ is obtained via iterated theta lifting (see
   \cite[Theorem 3.5.1]{B17}, \cite{Mo17} and \cite{BMSZ2}). The method of \cite{BMSZ2} gives an alternative proof that all representations in $\Unip_{\ckcO}( G)$ are unitarizable.

\section{Generalities on coherent families of highest weight modules}\label{sec:HWM}

In this section, we recall some generalities on coherent continuation representations for highest weight modules, as a preparation for the study of coherent continuation representations in the Casselman-Wallach setting. (The theory of primitive ideals, developed by Joseph and Barbasch-Vogan among others, plays a key role.)
Using $\tau$-invariant, we formalize a certain duality notion of double cells (\Cref{duald}), useful in understanding their counterparts in the Casselman-Wallach setting.

We retain the notation of Section \ref{secnot}.
Write \[
\Delta^+\subset \Delta\subset \hha^*\quad \textrm{and}\quad \check \Delta^+\subset \check \Delta\subset \hha
\]
for the positive root system, the root system, the positive coroot system, and the coroot system, respectively, for the reductive complex Lie algebra $\g$.
Write $Q_\g$ and $ Q^\g$ for the root lattice and the weight group of $\g$, respectively. Namely $Q_\g$ is the subgroup of $\hha^*$ spanned by $\Delta$, and
\[
Q^\g:=\{\nu\in \hha^*\mid \langle \nu, \check \alpha\rangle \in \Z\ \textrm{ for all }\check \alpha\in \check \Delta\}.
\]

\subsection{Coherent continuation representations}\label{subsec:coherent}

Throughout this article, the coefficient
ring of all Grothendieck groups will be  $\C$. When no confusion is possible, for every object $O$ in an abelian category, we use the same symbol to indicate the Grothendieck group element represented by the object $O$.

 Let $Q$ be a $W$-stable subgroup of $Q^\g$ containing $Q_\g$:
 \[ 
 Q_\g\subseteq Q\subseteq Q^\g.
 \]
Denote by $\mathrm{Rep}(\g, Q)$ the category of all finite-dimensional
representations of $\g$ whose weights 
(viewed 
as elements of $Q^\g\subset \hha^*$) are contained in $Q$. Write $\mathcal R(\g, Q)$ for the Grothendieck group of this category, which is
 a commutative $\bC$-algebra with multiplication the tensor
product of representations: $X\cdot Y:=X\otimes Y$, where $X,Y\in \mathcal R(\g, Q)$. 

Let 
 $\Lambda\subset \hha^*$ be a $Q$-coset.  
Let $W_\Lambda$ denote the stabilizer of $\Lambda$ in $W$. Specifically, if $\Lambda=\lambda +Q$ for a
   $\lambda\in \hha^*$, then 
\[
W_\Lambda:=\{w\in W \mid w\lambda-\lambda\in Q\}. 
\]

\begin{defn}
   Given  an $\mathcal R(\g, Q)$-module $\CK$, a $\CK$-valued $\Lambda$-family   is an assignment of subspaces    $\{\CK_\nu\}_{\nu\in \Lambda}$ of $\CK$ such that
\begin{itemize}
\item
  $\CK_{w \nu}=\CK_\nu\ $  for all $w\in W_\Lambda$ and $\nu\in \Lambda$;
  \item  for all representations $F$ in  $\mathrm{Rep}(\g, Q)$ and all $\nu\in \Lambda$,
  \[
   F\cdot \CK_\nu\subset \sum_{\mu 
   %\in \hha^* 
   \textrm{ is a weight of $F$} } \CK_{\nu+\mu}.
  \]
 \end{itemize}
\end{defn}

\begin{defn}[\cite{Jan}, \cite{Sch}] \label{defcoh00}
  Let $\CK$ be  an $\mathcal R(\g, Q)$-module, and $\{\CK_\nu\}_{\nu\in \Lambda}$  a $\CK$-valued  $\Lambda$-family. 
  A $\CK$-valued coherent family on $\Lambda$ based on $\{\CK_\nu\}_{\nu\in \Lambda}$ is a map
  \[
    \Psi: \Lambda\rightarrow \CK%, \qquad \lambda\mapsto \Phi_\lambda
  \]
  satisfying the following two conditions:
  \begin{itemize}
    \item for all $\nu\in \Lambda$, $\Psi(\nu)\in \CK_\nu$;
    \item for all representations $F$ in $\mathrm{Rep}(\g, Q)$
          and all $\nu\in \Lambda$,
          \[
          F \cdot (\Psi(\nu)) = \sum_{\mu} \Psi(\nu+\mu),
          \]
          where $\mu$ runs over all weights of $F$, counted with multiplicities.%  and $F$ is viewed as an element of $\mathcal R(\g)$.
  \end{itemize}
\end{defn}

When specifying a $\CK$-valued coherent family on $\Lambda$ based on $\{\CK_\nu\}_{\nu\in \Lambda}$, we will often explicitly describe $\CK$ as a Grothendieck group, while the $\mathcal R(\g, Q)$-module structure and the $\CK$-valued  $\Lambda$-family  $\{\CK_\nu\}_{\nu\in \Lambda}$ are the ones which are clear from the context. (Often  the $\CK$-valued $\Lambda$-family is specified by the infinitesimal character).
When the $\CK$-valued $\Lambda$-family  $\{\CK_\nu\}_{\nu\in \Lambda}$ is specified or clear from the context, we will just call it a $\CK$-valued coherent family on $\Lambda$. 

Given a $\CK$-valued
$\Lambda$-family $\{\CK_\nu\}_{\nu\in \Lambda}$, let $\Coh_{\Lambda}(\CK)$ denote the vector space of all $\CK$-valued coherent families on $\Lambda$ based on $\{\CK_\nu\}_{\nu\in \Lambda}$. 
It is a
representation of $W_{\Lambda}$ under the action
\be\label{actcoh}
  (w \cdot \Psi)(\nu) = \Psi(w^{-1} \nu), \qquad \textrm{for all
  }\ w\in W_\Lambda, \ \Psi\in \Coh_{\Lambda}(\CK), \   \nu\in \Lambda.
\ee
This is called a coherent continuation representation (based on $\{\CK_\nu(\g,\b)\}_{\nu\in     \Lambda}$). 

The assignment $\CK\mapsto \Coh_{\Lambda}(\CK)$ is functorial in the following sense: suppose that $\CK'$ is another  $\mathcal R(\g, Q)$-module  with a $\CK'$-valued $\Lambda$-family  $\{\CK'_\nu\}_{\nu\in \Lambda}$, and $\eta: \CK\rightarrow \CK'$ is an $\mathcal R(\g, Q)$-homomorphism such that $\eta(\CK_\nu)\subset \CK'_\nu$ for all $\nu\in \Lambda$, then
\begin{equation}\label{eq:etas}
 \eta_*: \Coh_{\Lambda}(\CK)\rightarrow \Coh_{\Lambda}(\CK'), \quad \Psi\mapsto \eta\circ \Psi
\end{equation}
is a well-defined $W_\Lambda$-equivariant linear map.

\subsection{Coherent continuation representation for highest weight modules}
\label{sub:CohHWM} 
Coherent families for highest weight modules were introduced by Jantzen \cite{Jan}. We refer the reader to \cite[Section 7]{H} for basic facts on coherent families in this setting.  

Let $\b$ be a  Borel subalgebra of $\g$, and let $\Rep(\g,\b)$ denote the category of finitely generated $\g$-modules that are unions of finite-dimensional $\b$-submodules.
For each $\nu\in \hha^*$, let  $\Rep_\nu(\g,\b)$ denote the  full subcategory of $\Rep(\g,\b)$ consisting of modules that have generalized infinitesimal character $\nu$ (by definition, a $\g$-module has generalized infinitesimal character $\nu$ if every vector in it is annihilated by $(\ker(\chi_\nu))^k$ for some $k\in \BN^+$).

\begin{defn}
  Write $\CK(\g,\b)$ for the Grothendieck group of $\Rep(\g,\b)$. It is an $\mathcal R(\g,Q)$-module under the tensor product.
 Similarly define  
 $$
\CK_{\nu}(\g,\b):=\textrm{the Grothendieck group of $\Rep_\nu(\g,\b)$},
$$
which is a subspace of  $\CK(\g,\b)$.
Denote by $\Coh_{\Lambda}( \CK(\g,\b))$ the coherent continuation
representation based on $\{\CK_\nu(\g,\b)\}_{\nu\in     \Lambda}$, as
in \eqref{actcoh}.
 \end{defn}

Write $\rho\in \hha^*$ for the half sum of positive roots.
For each $\nu\in \hha^*$, define the Verma module
\[
  \mathrm M(\nu):=\mathrm M(\g,\b,\nu):=\CU(\g)\otimes_{\CU(\b)} \C_{\nu-\rho},
\]
where  $\C_{\nu-\rho}$ is the one-dimensional $\hha$-module corresponding to the character $\nu-\rho\in \hha^*$, and every $\hha$-module is viewed an a $\b$-module via the canonical map $\b\rightarrow \hha$. Write
$\oL(\nu):=\oL(\g,\b,\nu)$ for the unique irreducible quotient of $ \mathrm M(\nu)$.

An element $\nu\in \hha^*$ is called dominant if
\be\label{dominant}
    \la \nu, \check \alpha\ra\notin -\bN^+ \qquad\textrm{for all $\check \alpha\in \check \Delta^+$},
  \ee
  and regular if
  \[
    \la \nu, \check \alpha\ra\neq 0 \qquad\textrm{for all $\check \alpha\in \check \Delta$}.
  \]

For every $w\in W$, define a map
\[
\begin{array}{rcl}
  \Psi_{w}: \Lambda&\rightarrow  &\CK(\g,\b), \\
   \nu&\mapsto& \mathrm M(w \nu).
   \end{array}
\]
Then $\Psi_w\in \Coh_{\Lambda}( \CK(\g,\b))$, and
\be\label{basis}
\textrm{ $\{\Psi_w\}_{w\in W}$ is a basis of $\Coh_{\Lambda}( \CK(\g,\b))$. }
\ee
Similarly there is a unique coherent family $\overline \Psi_w\in \Coh_{\Lambda}( \CK(\g,\b))$ such that
\be\label{psibarw}
  \overline \Psi_w(\nu)=\oL(w \nu)\quad \textrm{for all regular dominant element $\nu\in \Lambda$}.
\ee
Then $\{\overline \Psi_w\mid w\in W\}$ is also a   
basis of the coherent continuation representation $\Coh_{\Lambda}( \CK(\g,\b))$. 

Let $W$ act on  $\CK(\g,\b)$ as $\mathcal R(\g, Q)$-module automorphisms by
\[
  w\cdot (\mathrm M(\nu))=\mathrm M(w  \nu)\quad \textrm{for all }\nu\in \hha^*.
\]
By the functority \eqref{eq:etas},  this yields an action of $W$ on  
$\Coh_{\Lambda}( \CK(\g,\b))$ as  
automorphisms of $W_\Lambda$-representations. The resulting action of $W\times W_\Lambda$ on $\Coh_{\Lambda}( \CK(\g,\b))$ is explicitly given by
 \begin{equation}\label{eq:exp}
   (w_1, w_2)\cdot  \Psi_{w}=\Psi_{w_1 w w_2^{-1}}\quad\textrm{for all $w_1\in W$, $w_2\in W_\Lambda$,  $w\in W$.}
 \end{equation}

Let $\Lambda^\g\subset \hha^*$ denote the $Q^\g$-coset containing $\Lambda$, and let  $\Lambda_\g\subset \Lambda$ be a $Q_\g$-coset.
Put
\begin{equation} \label{eq:IR}
  \Delta(\Lambda):=\{\alpha\in \Delta\mid \la \check \alpha, \nu\ra\in \Z \textrm{ for some (and all) }\nu\in \Lambda\}.
\end{equation}
Here and henceforth,  $\check \alpha\in \hha$ denotes the coroot corresponding to $\alpha$.
This is a root system with the corresponding coroots
\[
  \check \Delta(\Lambda):=\{\check \alpha\in \check \Delta \mid \la \check \alpha, \nu\ra\in \Z \textrm{ for some (and all) }\nu\in \Lambda\}.
\]
Let 
\[W(\Lambda)\subset W\] denote the Weyl group of the root system $\Delta(\Lambda)$, to be called the integral Weyl group (attached to $\Lambda$).  Then
\[
  W(\Lambda_\g)=W(\Lambda)=W(\Lambda^\g)=W_{\Lambda_\g}\subset W_{\Lambda}\subset  W_{\Lambda^\g}.
\]
\trivial[h]{
This is because $\inn{\Delta(\fgg)}{Q^\fgg}\in \bZ$. 
}

To understand the coherent continuation representation $\Coh_{\Lambda}( \CK(\g,\b))$, 
it suffices to consider the case when $\Lambda=\Lambda^\fgg$ by the following lemma.

\begin{lem}\label{restco0}
The restriction map
\[
  \Coh_{\Lambda^\g}( \CK(\g,\b))\rightarrow \Coh_{\Lambda}( \CK(\g,\b))
\]
is a $W\times W_\Lambda$-equivariant  linear isomorphism,
and the restriction map
\[
 \Coh_{\Lambda}( \CK(\g,\b))\rightarrow \Coh_{\Lambda_\g}( \CK(\g,\b))
\]
is a $W\times W(\Lambda)$-equivariant  linear isomorphism.
\end{lem}
\begin{proof}
  This follows from \eqref{basis}. Essentially a coherent family on 
    $\Lambda$ extends uniquely to one on  $\Lambda^\g$ and a coherent family
    on $\Lambda_\g$ extends uniquely to one on $\Lambda$, for any    $\Lambda_\g\subset \Lambda\subset\Lambda^\g$.
\end{proof}

\subsection{Jantzen matrix}

Define the Jantzen matrix $\{a_{\Lambda}(w_1, w_2)\}_{w_1, w_2\in W}$, where $a_{\Lambda}(w_1, w_2)\in \Z$ is specified by
the following equation in $\Coh_{\Lambda}(\cK(\fgg,\fbb))$:  
\[
  \overline{\Psi}_{w_1}=\sum_{w_2\in W}  a_{\Lambda}(w_1, w_2)\Psi_{w_2},  
\]
for each  $w_1\in W$ .

\begin{lem}[Bernstein-Gelfand-Gelfand]\label{lem33} 
 Let $w_1, w_2\in W$. If $w_1 W(\Lambda)\neq w_2 W(\Lambda)$, then $a_{\Lambda}(w_1, w_2)=0$.
\end{lem}
\begin{proof}
This is a consequence of the theorem of Bernstein-Gelfand-Gelfand on the composition factors of a Verma module. 
See \cite[Corollary 5.2]{H}.
\end{proof}

\def\D{D}
Put
\[
  \Delta^+(\Lambda):=\Delta(\Lambda)\cap \Delta^+
\]
and  
\begin{equation} \label{wprime}
  \D(\Lambda):=\Set{w'\in W\mid  w' (\Delta^+(\Lambda) )\subset \Delta^+}.
\end{equation}
Then the group multiplication yields a bijective map
\[
  \D(\Lambda)\times W(\Lambda)\rightarrow W.
\]

\begin{lem}[Joseph] \label{wprime0}
For all  $w'\in \D(\Lambda)$ and $w_1, w_2\in W(\Lambda)$,  
\begin{equation}\label{eq:a1}
a_{\Lambda}(w'w_1, w'w_2)= a_{\Lambda}(w_1, w_2).
\end{equation}

\end{lem}
\begin{proof}
  This is shown by relating highest weight modules with principal series representations of complex semisimple groups.
  See \cite[Theorem~4.12 and Theorem~5.4]{J79D}. 
  \end{proof}

\begin{remark}\label{soergel}
Equation \cref{eq:a1} is also a direct consequence of a theorem of Soergel 
(see \cite[Section 2.5, Theorem 11]{Soergel}).
 The matrix   $\set{a_{\Lambda}(w_1, w_2)}_{w_1, w_2\in W(\Lambda)}$ only depends on $W(\Lambda)$ as a Coxeter group.
More precisely, suppose that $(\g', \Lambda',  W(\Lambda'))$ is a triple of the same type as 
$(\g, \Lambda, W(\Lambda))$, and $\eta: W(\Lambda)\rightarrow W(\Lambda')$ is a group
isomorphism that restricts to a bijection between the sets of  simple reflections, then
\[
a_{\Lambda}(w_1, w_2)=a_{\Lambda'}(\eta(w_1), \eta(w_2))
\]
for all $w_1, w_2\in W(\Lambda)$.
\end{remark}

Recall that a polynomial function on $\hha^*$ or $\hha$ is said to be $W(\Lambda)$-harmonic if it is  annihilated by all the $W(\Lambda)$-invariant constant coefficient differential operators without constant term.

For every $w\in W$, define a polynomial function $\tilde p_{\Lambda, w}$ on $\hha\times \hha^*$ by
\be\label{polynomial0000}
  \tilde p_{\Lambda, w}(x,\nu):= \sum_{w_1\in W} a_{\Lambda}(w, w_1)\cdot  \la x, w_1 \nu\ra^{m_{\Lambda, w}},\quad  \textrm{for all }  \ x\in \hha, \nu\in \hha^*,
\ee
where $m_{\Lambda, w}$ is the smallest non-negative integer (which always exists) that makes the right-hand side of \eqref{polynomial0000} a nonzero polynomial function.

\begin{lem}[Joseph \cite{J.hw}*{Section 5.1}] \label{leftcell00} 
Let $w\in W$. There is a $W(w\Lambda)$-harmonic polynomial function $p'_{\Lambda, w}$ on $\hha$ and a  $W(\Lambda)$-harmonic polynomial function $p_{\Lambda, w}$ on $\hha^*$ such that
\[
  \tilde  p_{\Lambda,  w}(x,\nu)= p'_{\Lambda,w}(x)\cdot  p_{\Lambda,w}(\nu), \quad\textrm{for all }  x\in \hha,  \ \nu\in \hha^*.
\]
The two harmonic polynomial functions are nonzero, homogeneous of degree $m_{\Lambda,w}$ and are uniquely determined up to scalar multiplication.
\end{lem}
\begin{proof}
  This follows from \cite[Lemma~2.3~(i), Lemma~2.5]{J2} and \cite[Theorem~2.6~(b)]{BV2}.
\end{proof}

\subsection{Left, right, and double cells}\label{seccell}

We define a basal vector space to be a complex vector space $V$ equipped with a basis $\CB\subset V$, and call elements of $\CB$ the basal elements in $V$. A subspace of a basal space $V$ is called a basal subspace if it is spanned by a set of basal elements of $V$. Every basal subspace is obviously a basal space.

As an example, if $\CK$ is the Grothendieck group of an abelian category in which all objects have finite length, then $\CK$ is a basal vector space with the irreducible objects as the basal elements.

\begin{defn}
Let $E$ be a  finite group. A basal representation of $E$ is a basal vector space carrying a representation of $E$. A basal subrepresentation of a basal representation $V$ is a subrepresentation of $V$ that is simultaneously a basal subspace.
\end{defn}

Let $V$ be a basal representation of a finite group $E$, with basal elements
$\CB\subset V$. For each subset $\CS \subset \CB$, write $\braket{\CS}$ for
the smallest basal subrepresentation of $V$ containing $\CS$. For each
$\phi\in \CB$, write $\braket{\phi}:=\la \{\phi\}\ra$ for simplicity. We define an equivalence relation $\approx$ on $\CB$  by
\[
  \phi_1 \approx \phi_2 \quad \textrm{ if and only if
  } \quad \la \phi_1 \ra =\la \phi_2 \ra \qquad (\phi_1, \phi_2\in \CB).
\]
An equivalence class of the relation $\approx$ on the set $\CB$ is called a cell in $V$.

\begin{defn}\label{defcell}
Let $\CC$ be a cell in $V$ and put $\overline{\CC}:=\braket{\CC}\cap\CB$.  Define the cell representation attached to $\CC$ by
\[
V(\CC):=\la \overline \CC \ra/ \la \overline \CC\setminus \CC\ra. 
\]    
  \end{defn}

In the notation of \Cref{defcell}, $V(\CC)$ is a representation of $E$, and the cosets $\{\phi+ \la \overline \CC\setminus \CC\ra\}_{ \phi\in \CC}$ form a basis of $V(\CC)$.

We view $\Coh_{\Lambda}( \CK(\g,\b))$ as a basal representation of $W\times W_\Lambda$ with the basal elements
 $\{\overline \Psi_{w}\mid w\in W\}$.
Write
\begin{equation}\label{eq:LRO}
\Coh^{LR}_{\Lambda}( \CK(\g,\b)):=\Coh_{\Lambda}( \CK(\g,\b)), 
\end{equation}
 to be viewed as a basal representation of $W\times W(\Lambda)$.
Likewise, write
\[
\Coh^{L}_{\Lambda}( \CK(\g,\b)):=\Coh_{\Lambda}( \CK(\g,\b)),
\]
 to be viewed as a basal representation of $W$, and write
 \[
 \Coh^{R}_{\Lambda}( \CK(\g,\b)):=\Coh_{\Lambda}( \CK(\g,\b)),
 \]
  to be viewed as a basal representation of $W(\Lambda)$ (as a subgroup of $W_\Lambda$).

\begin{defn}\
  \begin{itemize}
\item Cells in $\Coh^{L}_{\Lambda}( \CK(\g,\b))$ are called left cells. 
  \item Cells in $\Coh^{R}_{\Lambda}( \CK(\g,\b))$ are called right cells. 
  \item Cells in $\Coh^{LR}_{\Lambda}( \CK(\g,\b))$ are called double cells. 
  \end{itemize}
  For every set $\CS$ of basal elements in $\Coh_{\Lambda}( \CK(\g,\b))$,
write $\la \CS\ra_L$ for the smallest basal subrepresentation
of $\Coh^L_{\Lambda}( \CK(\g,\b))$ containing $\CS$.
Similarly write $\la \CS\ra_R$ and $\la \CS\ra_{LR}$.
\end{defn}

It is clear that every left (right) cell is contained in a (unique) double cell.

\subsection{Left cells and classification of primitive ideals}\label{sec:primitive}
Let $\Ann(M)\subset \mathcal U(\g)$ denote the 
annihilator ideal of a $\mathcal U(\g)$-module $M$. By a result of Duflo \cite[Theorem~1]{Du77}, each primitive ideal in the universal 
enveloping algebra $\cU(\fgg)$ has the form 
$\Ann \bPsi_w(\nu)$ for a $w\in W$ and a dominant $\nu\in \ahh^*$.

For each $w\in W$, let $p_{\Lambda,w}$  be the $W(\Lambda)$-harmonic polynomial function as in \Cref{leftcell00}. 

\begin{lem}\label{leftp}
Let  $w_1, w_2\in  W$.
The following conditions are equivalent.
\begin{enumerate}[label=(\roman*),wide=0pt]
\item The basal elements $\overline \Psi_{w_1}$ and $\overline \Psi_{w_2}$ lie in a common left cell. 
\item For some regular dominant $\nu\in \Lambda$, $\Ann \bPsi_{w_1}(\nu)  = \Ann \bPsi_{w_2}(\nu)$.
\item For every regular dominant $\nu\in \Lambda$, $\Ann \bPsi_{w_1}(\nu)  = \Ann \bPsi_{w_2}(\nu)$.
\item The equality $\C\cdot p_{\Lambda, w_1}=\C\cdot p_{\Lambda, w_2}$ holds. 
\end{enumerate}
\end{lem}

\begin{proof}
The equivalence of (i) and (ii) is the equivalence of $(a)$ and $(b)$
in \cite{BV2}*{Proposition~2.9}.\trivial[h]{ The theorem is conjectured by Joseph in \cite{J79W}*{Conjecture 5.7.C}
which is the main theorem of Vogan's paper on the ordering of primitive ideals. }
The equivalence of (ii) and (iii) follows from the translation principle, 
see~ \cite[Lemma~2.7]{V1}. 
The equivalence of (ii) and (iv) is \cite[Theorem~5.1 and 5.5]{J2} when 
$w_1 W(\Lambda) = w_2 W(\Lambda)$. The general case can be deduced by relating the problem to the setting of complex semisimple groups using  \cite[Theorem~4.12]{J79D} and \cite[Proposition~3.7]{Du77}. See \Cref{wprime0}. 
\end{proof}

Using \Cref{leftp}, we attach a 
polynomial function   $p_{\CCL}$  on $\hha^*$ to every left cell $\CCL$ in $\Coh_{\Lambda}( \CK(\g,\b))$,
which is uniquely determined up to scalar multiplication.

\begin{lem}\label{lem:tauinv}
For every $w\in W$ and every simple reflection $s\in W(\Lambda)$, the following statements are equivalent. 
\begin{enumerate}[label=(\roman*),wide=0pt]
\item $s\cdot p_{\Lambda,w}=-p_{\Lambda, w}$.
\item $s\cdot  \overline \Psi_w=- \overline \Psi_w$.
\item 
$ w(\alpha_s)\in \Delta^+$, where $\alpha_s\in \Delta^+(\Lambda)$ is the simple root corresponding to $s$.
\end{enumerate}
\end{lem}
\begin{proof}
The equivalence between (i) and (ii) is due to the fact that 
$p_{\Lambda,w}$ is essentially the Bernstein degree polynomial, see \cite[Theorem~4.10]{J2} and \cite{VGK}*{p85}. 
The equivalence between (ii) and (iii) is due to Jantzen, see \cite{BV2}*{Proposition~2.20}. 
\end{proof}

\begin{defn}
The  $\tau$-invariant of a basal element $\Psi:=\overline{\Psi}_w\in \Coh_{\Lambda}(  \CK(\g,\b))$ ($w\in W$) is the set of simple reflections 
satisfying the equivalent conditions in \Cref{lem:tauinv}:
\[
 \tau_\Psi:= \tau_{\Lambda, w}:=\{s\in W(\Lambda) \textrm{ is a simple reflection }\mid   s\cdot p_{\Lambda,w}=-p_{\Lambda, w}\}.
\]    
  \end{defn}

This only depends on the left cell containing $\Psi$,
and thus the $\tau$-invariant $\tau_{\CCL}$ of a left cell $\CCL$ is well defined.

For each $\nu\in \hha^*$, write $W_\nu$ for the stabilizer of $\nu$ in $W$. Note that $W_{\nu}\subseteq W(\Lambda)$, for any $\nu \in \Lambda$. The following is a reformulation of Jantzen's result on the translation to a wall, see \cite{H}*{Theorem~7.9} and 
also \cite[Corollary~7.3.22]{Vg}.
\begin{lem}[Jantzen]\label{prop:Jantzen}
Let $w\in W$ and let  $\nu\in \Lambda$ be a dominant element. Then
\[
\overline \Psi_w(\nu)=
\left\{
   \begin{array}{ll}
    \oL(w\nu),\quad & \text{if $\tau_{\Lambda, w}\cap W_{\nu} = \emptyset$;} \\
    %\textrm{if there is no element of  $\tau_{\Lambda, w}$ that fixes $\nu$};\\
       0, \quad &\textrm{otherwise}.
   \end{array}
\right.
\]             
\end{lem}

\Cref{prop:Jantzen} and the translation principle (see \cite{V1}*{Lemma~2.7}) 
imply the following classification of primitive ideals at a general infinitesimal character.  

\begin{prop}[Joseph]\label{primitivei}
Let $\nu\in \Lambda$ be a dominant element. Then the map
\[
\begin{array}{rcl}
  \set{w\in W | \tau_{\Lambda, w}\cap W_\nu = \emptyset} & \rightarrow &
 \set{\textrm{primitive ideal of $\CU(\g)$ of infinitesimal character $\nu$}}\\
 w & \mapsto & \Ann( \oL(w \nu)) 
\end{array}
\]
is surjective. Furthermore for all $w_1, w_2$ in the domain of this map,
\[
  \Ann(\oL(w_1\nu))=\Ann(\oL(w_2\nu))\quad \textrm{if and only if}\quad \C\cdot p_{\Lambda, w_1}=\C\cdot p_{\Lambda, w_2}.
\]
\end{prop}

\begin{cor}  \label{primitive000}
% By \Cref{primitivei},
Let  $\nu\in \Lambda$ be dominant.  Then for every basal element $\Psi\in \Coh_{\Lambda}(  \CK(\g,\b))$, the ideal $\Ann(\Psi(\nu))$ only depends
on the left cell $\CCL$ that $\Psi$ belongs to. 
%Denote by
%\begin{equation*}%\label{eq:ICC}
%I_{\nu,\CCL} := \Ann(\Psi(\nu)).  
%\end{equation*}
%It is a primitive ideal if $\tau_{\CCL}\cap W_\nu = \emptyset$, and $\cU(\fgg)$ otherwise.
The map 
\begin{equation}\label{eq:primitive}
\begin{split}
 & \quad \set{\textrm{left cell $\CCL $ in $\Coh_{\Lambda}( \CK(\g,\b))$} |
 \tau_{\CCL} \cap W_\nu = \emptyset }\\
 \xrightarrow{\CCL\mapsto 
 %I_{\nu, \CCL}
 \Ann(\Psi(\nu))
 } & \quad \{\textrm{primitive ideal of $\CU(\g)$ of infinitesimal character $\nu$}\}.
 \end{split}
\end{equation}
is a well-defined bijection where $\Psi$ is an arbitrary element in $\CCL$. 
 In particular when $\nu \in \Lambda$ is regular dominant, for each basal element $\Psi$ in $\Coh_\Lambda(\cK(\fgg,\fbb))$,
the left cell corresponding to the primitive ideal $\Ann (\Psi(\nu))$ is just the left cell containing  $\Psi$.     
\end{cor}

The corollary justifies the following definition.
  \begin{defn}
  For every left cell $\CCL$ in $\Coh_\Lambda(\cK(\fgg,\fbb))$ 
and every dominant element $\nu\in \Lambda$,  define 
\begin{equation}\label{eq:ICC}
I_{\nu,\CCL} := \Ann(\Psi(\nu)), \quad \text{for (any) } \Psi \in \CCL.  
\end{equation}
Conversely, for every primitive ideal $I$ of infinitesimal character $\nu\in \Lambda$, write $\CCL_{\nu, I}$ %\mjj{$\CC^L_{\Lambda, I}$?}
for the left cell that corresponds to it under the bijection \eqref{eq:primitive}.
\end{defn}

We record the following lemma for later use. 

\begin{lem}\label{primitiveii}
Let $\nu\in \Lambda$ be a dominant element. Then for all primitive ideals $I_1$ and $I_2$ of $\CU(\g)$ of infinitesimal character $\nu$,
\[
  I_1\subset I_2\quad \textrm{if and only if}\quad \la \CCL_{\nu, I_1}\ra_L\supseteq \la \CCL_{\nu, I_2}\ra_L.
\]
\end{lem}
\begin{proof}
  This follows from \Cref{primitive000}
  %\Cref{primitivei} 
  and \cite{BV2}*{Proposition~2.9}. 
\end{proof}

\subsection{Double cells and special representations}
\label{sec:DCell}

For every $\sigma\in \Irr(W)$, its fake degree is defined to be
 \begin{equation}\label{eq:fdeg}
 a(\sigma):=\min\{k\in \BN\,|\, \sigma \textrm{ occurs in the $k$-th symmetric power $\oS^k(\hha)$}\}.
 \end{equation}
This is well-defined since every  $\sigma\in \Irr(W)$ occurs in 
the symmetric algebra $\oS(\hha)$.
The representation $\sigma$ is said to be \emph{univalent} if it occurs in $\oS^{a(\sigma)}(\hha_{\mathrm s})$ with multiplicity one, where $\hha_{\mathrm s}:=\Span(\check \Delta)$ denotes the span of the coroots.

Recall Lusztig's notion of a special representation of a Weyl group 
(\cite{Lsp}).
An irreducible representation of $W$ is said to be Springer if it corresponds to
the trivial local system on a nilpotent orbit in $\g^*$ via the Springer correspondence \cite{Spr}.
Note that every special irreducible representation is Springer, and every Springer representation is univalent \cite{BM}.

We have a decomposition
\[
  \hha=(\Delta(\Lambda))^\perp\oplus \Span( \check \Delta(\Lambda))
  \]
where
\[
(\Delta(\Lambda))^\perp:=\{x\in \hha\mid \la x , \alpha \ra=0\textrm{ for all $\alpha\in \Delta(\Lambda)$}\}.
\]
 For every univalent irreducible representation $\sigma_0$ of $W(\Lambda)$, whenever it is convenient, we view it as a subrepresentation of  $\oS^{a(\sigma_0)}(\hha)$ via the inclusions
\[
\sigma_0=  \C\otimes \sigma_0\subset \oS^0((\Delta(\Lambda))^\perp)\otimes \oS^{a(\sigma_0)}(\Span( \check \Delta(\Lambda)))\subset  \oS^{a(\sigma_0)}(\hha)
\]
By the work of Macdonald, Lusztig, and Spaltenstein (\cite[Chapter 11]{Carter}), the $W$ subrepresentation of $\oS^{a(\sigma_0)}(\hha)$ generated by $\sigma_0$
is irreducible and univalent, with the same fake degree as that of $\sigma_0$.  This irreducible representation of $W$ is called the $j$-induction of $\sigma_0$, to be denoted by $j_{W(\Lambda)}^W(\sigma_0)$.

If  $\sigma_0$ is special, then the $j$-induction $j_{W(\Lambda)}^W \sigma_0$ 
is Springer. 
Write \begin{equation}\label{eq:OrbitSp}
\CO_{\sigma_0}\in \overline \Nil(\g^*)\quad \text{for the nilpotent orbit corresponding to 
$j_{W(\Lambda)}^W \sigma_0$},\end{equation}
via the Springer correspondence. Then  
\[
  \dim \CO_{\sigma_0}=2\cdot (\sharp (\Delta^+) - a(\sigma_0)).
\]
See \cite{Ho}*{Theorem~3}.

\begin{defn}\label{def:doublecell}
%\begin{enumerate}
%\item 
(a) Define a preorder  $\leqLR$ on $\Irr(W(\Lambda))$ by 
\[
\sigma_1\leqLR \sigma_2
\text{\ \  if and only if \ \ } \la \CCD_1\ra_{LR}\supseteq \la \CCD_2\ra_{LR}
\]
 for some double cells $\CCD_i$ in $\Coh_{\Lambda}(\CK( \g, \b))$ such that $\sigma_i$ occurs in  $\Coh^{LR}_{ \Lambda}(\CK( \g, \b))(\CCD_i)$ ($i=1,2$). 
% \item 

\noindent (b) Define an equivalence relation $\approxLR$ on 
$\Irr(W(\Lambda))$ by
\[
  \sigma_1\approxLR \sigma_2\quad \textrm{if and only if}\quad 
  \sigma_1\leqLR \sigma_2\textrm{ and } \sigma_2\leqLR \sigma_1. 
\]
%\item 
An equivalence class of this equivalence relation  is called a double cell in
$\Irr(W(\Lambda))$. 
%\end{enumerate}
\end{defn}

We remark that the preorder $\leqLR$ only depends on $W(\Lambda)$ as a Coxeter group (see \cite{BV2}*{Proposition~2.28}) and  
the notion of double cells in \Cref{def:doublecell}~(b)  coincides with that of Lusztig in \cite{Lcell} (also \cite{Carter}*{Section~13.2}).

A basic property of this preorder is the following result. 

\begin{prop}[{\cite{BV2}*{Proposition~2.25}}]\label{lrd}
For all $\sigma_1, \sigma_2\in \Irr{(W(\Lambda))}$,
$\sigma_1\leqLR \sigma_2$ if and only if  $\sigma_2\otimes \sgn\leqLR \sigma_1\otimes \sgn$.
\end{prop}

Here and henceforth $\sgn$ denotes the sign character of an appropriate Weyl group.

Suppose that $\CCD$ is a double cell in $\Coh_{\Lambda}( \CK(\g,\b))$. Pick a $w\in W$ such that $\overline \Psi_{w}\in \CCD$, and 
put \begin{equation}\label{def:Mc}
m_\CCD:=m_{\Lambda, w}, \quad (\text{as in } \eqref{polynomial0000}). 
\end{equation}
Note that 
$m_\CCD$ is independent of the choice of $w$ (see \cite{BV2}*{Corollary~2.15 and 2.16}).
For every $\Psi=\sum_{w\in W} a_{w} \Psi_{w}\in \la \CCD\ra_{LR}$, define a polynomial function $\tilde p_\Psi$ on $\hha\times \hha^*$ by
\begin{equation}\label{polynomial00}
  \tilde p_{\Psi}(x,\nu):= \sum_{w\in W} a_{w}\cdot  \la x, w \nu\ra^{m_{\CCD}},\quad  \textrm{for all }  \ x\in \hha, \nu\in \hha^*.
\end{equation}

Then the linear map 
\begin{equation}\label{lara}
    \braket{\CCD}_{LR}\rightarrow \oS(\hha^*)\otimes \oS(\hha), \quad \Psi\mapsto \tilde p_{\Psi}\quad (\oS\textrm{ indicates the symmetric algebra})
\end{equation}
is $W\times W(\Lambda)$-equivariant  % In fact, section~3.5 and 3.6
and descends to a linear map
\begin{equation}\label{lara2} 
 \Coh^{LR}_{\Lambda}( \CK(\g,\b))( \CCD) \rightarrow \oS(\hha^*)\otimes \oS(\hha).
\end{equation} 
See  \cite{J2}*{\S~3} and \cite[Corollary 2.15]{BV2}.

Let  $\Irr^{\mathrm{sp}}(W(\Lambda))$ denote the subset of $\Irr(W(\Lambda))$ consisting of the special 
irreducible representations. 
We recollect some results due to Joseph and Barbasch-Vogan 
in the following proposition. See \cite[Section 5]{J2} and \cite{BV2}*{Section 2}. 

\begin{prop}\label{cell00}\label{dcrep} Let $\CCD$ be a double cell in $\Coh_{\Lambda}( \CK(\g,\b))$. 
\begin{enumerate}[label=(\alph*),wide=0pt]
\item
 The image of \eqref{lara2} %equals
 is isomorphic to 
$\left(j_{W(\Lambda)}^W\sigma_\CCD\right)\otimes \sigma_\CCD$
for a unique $\sigma_\CCD$ in $\Irr^{\mathrm{sp}}(W(\Lambda))$.
Moreover, $a(\sigma_\CCD) = m_\CCD$ and for every $\sigma'\in \Irr(W(\Lambda))$ that occurs 
in $\braket{\CCD}_{LR}$, either $\sigma'=\sigma_\CCD$ or $a(\sigma')>m_\CCD$.
\item The representation
$\sigma_\CCD$ is the unique 
element of $\Irr^{\mathrm{sp}}(W(\Lambda))$ 
that occurs in the double cell representation $\Coh^{LR}_{\Lambda}( \CK(\g,\b))( \CCD)$. Moreover, the map
\[
  \{\textrm{double cell in $ \Coh_{\Lambda}( \CK(\g,\b))$}\}\rightarrow \Irr^{\mathrm{sp}}(W(\Lambda)), \quad \CCD\mapsto \sigma_\CCD
\]
is bijective.
\item 
As a representation of $W\times W(\Lambda)$,
\[
 \Coh^{LR}_{\Lambda}( \CK(\g,\b))( \CCD)\cong \bigoplus_{
 \substack{\sigma\in \Irr(W(\Lambda)),\\
 \sigma\approx \sigma_\CCD}}
\left(\Ind_{W(\Lambda)}^W\sigma\right)
\otimes \sigma.
\]
\end{enumerate}
\end{prop}

\begin{defn}\label{def:sigmaCC} For a double cell $\CCD$ in $\Coh_{\Lambda}( \CK(\g,\b))$, call $\sigma_\CCD$ as in \Cref{cell00} the special irreducible representation
attached to $\CCD$ .
\end{defn}

Recall the definition of the polynomial function $p_{\CCL}$ after \Cref{leftp}, where $\CCL$ is a left cell in  $\Coh_{\Lambda}( \CK(\g,\b))$. 

\begin{prop}[\cite{J2}*{Theorem~5.5}]\label{dcrep33}
Let $\CCD$ be a double cell in $\Coh_{\Lambda}( \CK(\g,\b))$. Then the family 
\[
\{p_{\CCL}\}_{\CCL\textrm{ is a left cell contained in }\CCD} 
\]
is linearly independent in $\mathrm S^{m_\CCD}(\hha)$ and spans a $W(\Lambda)$-subrepresentation that is isomorphic to $\sigma_\CCD$. 
\end{prop}

\subsection{Associated varieties}\label{sec:ass}

For every  ideal $I$ of $\CU(\g)$, denote by   
$\AV (I)$ its associated variety.  This is the subvariety in $\g^*$ of zeroes of the associated graded ideal in $\textrm{S}(\g)$, the symmetric algebra on $\g$. Here the grading is defined using the standard filtration of $\CU(\g)$.

\begin{defn}\label{def:sigmaC}
Define the Goldie rank representation of a primitive ideal $I$
of  infinitesimal character $\nu$,  where $\nu \in \Lambda$ is dominant, to be
\begin{equation}\label{grrj}
\sigma_{\nu, I}:=\sigma_{\CCD_{\nu,I} }\in \Irr^{\mathrm{sp}}(W(\Lambda)).
\end{equation}
Here $\CCD_{\nu,I}$ is the double cell containing the left cell $\CCL_{\nu, I}$.
\end{defn}

For a primitive ideal $I$ of infinitesimal character $\nu$, we may thus attach the nilpotent orbit 
$\CO_{\sigma_{\nu,I}}\in \overline \Nil(\g^*)$, via the Goldie rank representation 
$\sigma_{\nu,I}$ as in \eqref{eq:OrbitSp}.

The following result of Joseph determines the associated variety of a primitive ideal.  
\begin{prop}[\cite{J.av}*{Theorem~3.10}]
\label{assv}
Let $\nu\in \Lambda$ be a dominant element and let $I$ be a primitive ideal of $\CU(\g)$  of infinitesimal character $\nu$.   Then the associated variety $\AV (I)$ equals the Zariski closure
$\overline{\CO_{\sigma_{\nu,I}}}$ of $\CO_{\sigma_{\nu,I}}$ in 
$\g^*$.  
\end{prop}

\begin{defn}\label{def:av}
Fix  $\sfS\subset \mathrm{Nil}(\g^*)$, an $\Inn(\g)$-stable Zariski closed subset.
 Let $\Rep_\sfS(\g,\b)$ be the full subcategory of $\Rep(\g,\b)$ consisting of modules $M$ such that 
 $$ \AV(\Ann(M))\subset \sfS.$$
 Denote by $\CK_\sfS(\g,\b)$ its Grothendieck group. 
\end{defn}

As before, we have a coherent continuation  representation  $\Coh_{\Lambda}(\CK_\sfS(\g,\b))$ of $W_\Lambda$, which is a subrepresentation of $\Coh_{\Lambda}(\CK(\g,\b))$.

For every $w\in W$, put \[\CO_{\Lambda, w}:=\CO_{\sigma_\CCD},\]
where $\CCD$ is the double cell  in $\Coh_{\Lambda}( \CK(\g,\b))$ containing $\overline \Psi_w$.

\begin{lem}\label{basalassv}
Let $w\in W$. The coherent family $\overline \Psi_w$ belongs to
$\Coh_{\Lambda}(\CK_\sfS(\g,\b))$ if and only if  $\CO_{\Lambda, w}\subset \sfS$.
\end{lem}
\begin{proof}
This follows from \Cref{prop:Jantzen}, \Cref{primitivei}, \Cref{assv} and 
\cite{Vg}*{Proposition~7.3.11} (for the translation outside the dominant cone). 
\end{proof}

\begin{lem}\label{assv33}
Let $w_1, w_2\in W$. If $\la \overline \Psi_{w_1}\ra_L\supsetneq \la \overline \Psi_{w_2}\ra_L$, then $\overline{\CO_{\Lambda, w_1}}\supsetneq \overline{\CO_{\Lambda, w_2}}$.
\end{lem}

\begin{proof}
Pick a  regular dominant element $\nu\in \Lambda$.
Put $I_i:=\Ann(\overline{\Psi}_{w_i}(\nu))$ ($i=1,2$).
\Cref{primitiveii} implies that  $I_1\subsetneq I_2$ and hence
\[
  \AV(I_1)\supsetneq \AV(I_2).
\]
See \cite[Korollar 3.6]{BK}. On the other hand, \Cref{assv} implies that  $\AV(I_i)=\overline{\CO_{\Lambda, w_i}}$ ($i=1,2$). Hence the lemma follows.
\end{proof}

\begin{lem}\label{basals}
The space $\Coh_{\Lambda}(\CK_\sfS(\g,\b))$ is a 
$W\times W_\Lambda$-submodule
of $\Coh_{\Lambda}(\CK(\g,\b))$ spanned by the basal elements
\be\label{spanb}
 \{\overline \Psi_{w}\mid w\in W,  \CO_{\Lambda, w}\subset \sfS\}.
\ee
\end{lem}
\begin{proof}
\Cref{basalassv}  implies that  $\Coh_{\Lambda}(\CK_\sfS(\g,\b))$
is spanned by the set \eqref{spanb}.
Since the space $\Coh_{\Lambda}(\CK_\sfS(\g,\b))$ is  clearly $W_\Lambda$-stable,
\Cref{assv33} implies that it is also $W$-stable.
This proves the lemma. 
\end{proof}

\begin{defn}\label{defn:IrrS}
 Define
\[
  \Irr_\sfS^{\mathrm{sp}}(W(\Lambda)):= \Set{\sigma_0\in \Irr^{\mathrm{sp}}(W(\Lambda))| \cO_{\sigma_0}\subset \sfS}
\]
and
\[
\Irr_\sfS(W(\Lambda))
  := \Set{\sigma\in \Irr(W(\Lambda))|
  \text{there is a $\sigma_0\in \Irr_\sfS^{\mathrm{sp}}(W(\Lambda))$ such that $\sigma\approx \sigma_0$} }.
  \]
  \end{defn}
  
By \Cref{basals}, $\Coh_{\Lambda}(\CK_\sfS(\g,\b))$ is a basal subrepresentation of  $\Coh^{LR}_{\Lambda}(\CK(\g,\b))$.

\begin{prop}\label{cohbbs}
As a representation of $W\times W(\Lambda)$, % we have
\[
  \Coh_{\Lambda}( \CK_{\sfS}(\g,\b))\cong \bigoplus_{\sigma\in \Irr_\sfS(W(\Lambda))} \left(\Ind_{W(\Lambda)}^W \sigma\right)\otimes \sigma.
\]
%Here $\Ind$ indicates the induced representation.
\end{prop}

\begin{proof}
Write $\CB_0$ for the set of basal elements in $ \Coh_{\Lambda}(\CK_\sfS(\g,\b))$.
Then it is elementary to see that we may 
 choose a filtration
\[
  \CB_0\supset \CB_1\supset \cdots \supset \CB_k=\emptyset \qquad (k\in \BN)
\]
such that
\begin{itemize}
\item  $\braket{\CB_i}_{LR}$ is spanned by $\CB_i$ ($i=0, 1,2, \cdots ,k$),
\item $\CCD_i:=\CB_i\setminus \CB_{i+1}$ is a double cell in $\Coh_{\Lambda}(\CK(\g,\b))$ for all $i=0,1, \dots, k-1$, and
\item
for all $i_1, i_2\in \{0,1, \dots, k-1\}$,
$\braket{\CCD_{i_1}}_{LR} \subset \braket{ \CCD_{i_2}}_{LR}$ implies that $i_1\geq i_2$.
\end{itemize}

Then \Cref{cell00} and  \Cref{basalassv} imply that
\[
\{\sigma_{\CCD_0}, \sigma_{\CCD_1}, \dots, \sigma_{\CCD_{k-1}}\}= \Irr_\sfS^{\mathrm{sp}}(W(\Lambda)).
\]
Thus we have that
\begin{eqnarray*}
   &&\Coh_{\Lambda}(\CK_\sfS(\g,\b))
\\
   &\cong & \bigoplus_{i=0}^{k-1}  \la \CB_i\ra_{LR}/\la \CB_{i+1}\ra_{LR} \cong \bigoplus_{i=0}^{k-1}  \Coh_{\Lambda}(\CK(\g,\b))(\CCD_i)\\
    &\cong &\bigoplus_{\sigma\in \Irr_\sfS(W(\Lambda))} \left(\Ind_{W(\Lambda)}^W\sigma\right)\otimes \sigma\qquad(\textrm{by \Cref{dcrep}~(c)}).
\end{eqnarray*}
\end{proof}

Let $\nu\in \Lambda$ be a dominant element. Recall that for a primitive ideal $I$ of infinitesimal character $\nu$, we have the corresponding left cell $\CCL_{\nu, I}$ (\Cref{sec:primitive}). In the rest of the section, we focus on left cells which come from maximal ideals. 

Recall from the introductory section the maximal ideal $I_\nu$ with the infinitesimal
  character $\nu$. Write 
  \begin{equation}\label{def:onu}
  \CO_\nu\in \overline{\mathrm{Nil}}(\g^*)\end{equation}
  for the nilpotent orbit 
 whose Zariski closure $\overline{\CO_{\nu}}\subset \g^*$ equals 
 the associated variety of $I_\nu$. Note that 
 in the notation of \Cref{assv},  $\CO_{\nu}=\CO_{\sigma _{\nu, I_{\nu}}}$.

 Recall also the $J$-induction defined in \cite{Lcell}*{\S~11} (see also \cite{Lu}*{(4.1.7)}). 
\begin{prop}[Barbasch-Vogan]\label{leftcnu}
Let $\nu\in \Lambda$ be a dominant element. As $W$-representations, 
\[
  \Coh^L_{\Lambda}( \CK(\g,\b))(\CCL_{ \nu, I_\nu})\cong \Ind_{W(\Lambda)}^W \left(\left(J_{W_{\nu}}^{W(\Lambda)} \sgn \right)\otimes \sgn\right)
  \]
and as $W(\Lambda)$-representations, 
\begin{equation}\label{leftc3}
  \left(\left(J_{W_{\nu}}^{W(\Lambda)} \sgn \right)\otimes \sgn\right) \cong \bigoplus_{\sigma\in \Irr_{\overline{\CO_{\nu}}}(W(\Lambda)), \, [1_{W_\nu} :\sigma]\neq 0 }  \sigma.
\end{equation}
Moreover, 
\[
\sigma_{\nu,I_\nu} \cong \left(j_{W_{\nu}}^{W(\Lambda)} \sgn \right)\otimes \sgn ,
\]
and $\sigma_{\nu,I_\nu}$ is the unique special irreducible representation of $W(\Lambda)$ that occurs in \eqref{leftc3}.
\end{prop}
\begin{proof}
    In view of \Cref{dcrep}~(b) and \Cref{assv33}, this follows by the same line of proof as  \cite{BVUni}*{Corollary 5.30}.
\end{proof}

The isomorphism \eqref{leftc3} and explicit formulas of the $J$-induction (see \cite{Lu}*{\S4.4-4.13}) imply that
\be\label{mulone}
 [1_{W_\nu} :\sigma]\leq 1
 \ee
for all dominant  $\nu\in \Lambda$ and all $ \sigma\in \Irr_{\overline{\CO_\nu}}(W(\Lambda))$.

\begin{defn} \label{def:ll}
We call the set 
\begin{equation} \label{ll}
\LC_{\nu}:=\left\{\sigma\in \Irr(W(\Lambda))\mid \sigma \textrm{ occurs in } \left(J_{W_{\nu}}^{W(\Lambda)} \sgn \right)\otimes \sgn\right \}
\end{equation}
the Lusztig left cell attached to $\nu\in \Lambda$ (which is not necessarily  dominant). 
\end{defn}

\begin{lem}\label{leftcnu2}
Let $\nu\in \Lambda$ be a dominant element and let $I$ be a primitive ideal of $\CU(\g)$ of infinitesimal character $\nu$.    Then $I=I_\nu$ if and only if $\sigma_{\nu, I}\cong \left(j_{W_{\nu}}^{W(\Lambda)} \sgn \right)\otimes \sgn$.
\end{lem}

\begin{proof}
The ``only if$\,$" part  follows from \Cref{leftcnu}.
For the proof of the ``if'' part, note that $I \subseteq I_\nu$. 
%we suppose that 
Then $\sigma_{\nu, I}\cong \sigma_{\nu, I_\nu}$
implies that $I$ and $I_\nu$ have the same associated variety by \Cref{assv}. 
Thus $I=I_\nu$ by the maximality of $I_\nu$ and the proof is complete.
\end{proof}

\subsection{$\tau$-invariants and duals of double cells}\label{subsec:tau-DC}

Let $\Pi(\Lambda)$ denote the set of simple reflections in $W(\Lambda)$.
For each subset $S\subset \Pi(\Lambda)$, define an element 
\[
e_S :=\frac{1}{\sharp(W_{\Pi(\Lambda)\setminus S})\cdot \sharp(W_S)}\sum_{g\in W_{\Pi(\Lambda)\setminus S}}\sum_{h\in W_{S}} \sgn(h) \,  gh
\]
in the group algebra $\bC[W(\Lambda)]$, where $W_S$ and $W_{\Pi(\Lambda)\setminus S}$ are the subgroups of $W(\Lambda)$ generated by $S$ and $\Pi(\Lambda)\setminus S$ respectively. 

Following \cite{FJMN}, for every finite-dimensional representation $\sigma$ of $W(\Lambda)$,  %$\sigma\in \Irrsp{(W(\Lambda))}$,
define its $\tau$-invariant to be the set
\[
  \tau_\sigma:=\left\{S\subset \Pi(\Lambda)\mid 
   \text{the $e_S$ action on $\sigma$ is nonzero}
  \right\}. 
\]

 For every double cell $\CCD$ in $\Coh_\Lambda(\CK(\g,\b))$, we define its $\tau$-invariant 
to be the set 
\[
  \tau_{\CCD}:=\{\tau_\Psi\mid \Psi\in \CCD\}. 
\]
Recall that we have a representation $\sigma_\CCD\in \Irrsp{(W(\Lambda))}$ attached to $\CCD$ 
(\Cref{def:sigmaCC}). 

\begin{lem}\label{eqtau0}
For every  double cell $\CCD$ in $\Coh_\Lambda(\CK(\g,\b))$, $\tau_\CCD=\tau_{\sigma_{\CCD}}$. 
\end{lem}
\begin{proof}
  Let $M$ denote the cell representation $\Coh^{LR}_{\Lambda}( \CK(\g,\b))(\CCD)$ to be viewed as a
  $W(\Lambda)$ representation. By \cite{FJMN}*{Theorem~2.10}, we know that 
  $\tau_{\CCD}=\tau_M$. Since $\sigma_\CCD$ occurs in $M$, we have that
  \[
  \tau_{\sigma_{\CCD}}\subset \tau_M =  \tau_{\CCD}.
  \]
   
  Now it remains to show that $\tau_{\CCD} \subset \tau_{\sigma_{\CCD}}$.  
  Let $\Psi\in \CCD$ and put  $S:=\tau_\Psi$. 
  Write 
  \[
    e_S=Q(S')\cdot R(S),
  \]
  where $S':=\Pi(\Lambda)\setminus S$, 
\[
R(S):=
\frac{1}{ \sharp(W_S)}\sum_{h\in W_{S}} \sgn(h)\,h
\]
  and
\[
Q(S'):=
\frac{1}{ \sharp(W_{S'})}\sum_{g\in W_{S'}} g.
\]

For each $\Phi\in \CCD$, write $[\Phi]\in M$ for the class represented by $\Phi$.
Then
\[
  R(S)\cdot [\Psi]=[\Psi],
\]
and 
 \cite{FJMN}*{Proof of Proposition~2.6} 
 implies that 
 \[
  Q(S')\cdot [\Psi]=[\Psi]+\sum_{\Psi'\in \CCD, \tau_{\Psi'}\not\subset S} c_{\Psi'} [\Psi'],
 \]
 where $c_{\Psi'}\in \C$. Thus
 \be\label{es123}
  e_S\cdot [\Psi]=[\Psi]+\sum_{\Psi'\in \CCD, \tau_{\Psi'}\not\subset S} c_{\Psi'} [\Psi'].
  \ee
 
 Through \Cref{cell00}, we view $\sigma_\CCD$ as a subspace of $\oS^{m_\CCD}(\hha)$. The same proposition also implies that there is a $W(\Lambda)$-homomorphism
 \[
 \phi: M\rightarrow \sigma_{\CCD}
 \]
 such that $\phi([\Phi])$ is a nonzero scalar multiple of $p_{\CCL_{\Phi}}$ for all $\Phi\in \CCD$, where $\CCL_{\Phi}$ is the left cell containing $\Phi$. Applying $\phi$ to the quality \eqref{es123}, we obtain that 
 \[
 e_S\cdot p_{\CCL_\Psi}=p_{\CCL_\Psi}+\sum_{\Psi'\in \CCD, \tau_{\Psi'}\not\subset S} c'_{\Psi'} p_{\CCL_{\Psi'}},
 \]
 where $c'_{\Psi'}\in \C$. 
 This is nonzero by \Cref{leftp} and \Cref{dcrep33}, and by noting that  $\CCL_{\Psi} \neq \CCL_{\Psi'}$ for all $\Psi'\in \CCD$ with $\tau_{\Psi'}\not\subset S$. 
  Therefore $S\in \tau_{\sigma_{\CCD}}$ and the proof is complete. 
 \end{proof}

\begin{prop}\label{eqtau1}
Let $\sigma_1, \sigma_2\in \Irrsp{(W(\Lambda))}$. If $\tau_{\sigma_1}=\tau_{\sigma_2}$, then $\sigma_1=\sigma_2$. 
\end{prop}

\begin{proof}
\def\SignSig{\rm{SignSig}}
    Clearly it suffices to consider the case when the root system  $\Delta(\Lambda)$ is 
    irreducible. The lemma follows from \Cref{eqtau0} and the results of \cite{FJMN} (Theorem 2.12 for classical groups and the discussion in \S~6 for exceptional groups). 
\end{proof}

Let $(\check \g, \check \b, \check \Lambda)$ be a triple with the same property as the triple $(\g,\b, \Lambda)$. Let $\check W$ denote the Weyl group of $\check \g$. As before we have the basal space  $\Coh_{\check \Lambda}(\CK(\check \g,\check \b))$ carrying a representation of $\check W\times \check W(\check \Lambda)$, where  $\check W(\check \Lambda)\subseteq\check W$ is the integral Weyl group attached to $\check \Lambda$.

Suppose we are given an identification $W(\Lambda)=\check W(\check \Lambda)$ of Coxeter groups, namely a group isomorphism that sends simple reflections to simple reflections. We say that two subsets $\tau_1$ and $\tau_2$ of the power set of 
$\Pi(\Lambda)$ are dual to each other if
\[
\tau_2=\{\Pi(\Lambda)\setminus S\mid S\in \tau_1\}.
\]

\begin{prop}\label{duald} We are given an identification $W(\Lambda)=\check W(\check \Lambda)$ of Coxeter groups. 
  For every double cell $\CCD$ in $\Coh_\Lambda(\CK(\g,\b))$, there is a unique double cell $\check \CCD$ in $\Coh_{\check \Lambda}(\CK(\check \g,\check \b))$ whose $\tau$-invariant is dual to that of $\CCD$. Moreover, the map
\[
   \Set{\sigma\in \Irr(W(\Lambda)) | \sigma \approx \sigma_\CCD}
   \xrightarrow{\sigma\mapsto \sigma\otimes \sgn}  \Set{\check \sigma\in \Irr(\check W(\check \Lambda)) | \check \sigma \approx \sigma_{\check \CCD}}
\]  
%   \begin{eqnarray*}
%   &&\{\sigma\in \Irr(W(\Lambda))\mid \sigma\textrm{ occurs in $\Coh^{LR}_{ \Lambda}(\CK( \g, \b))(\CCD)$}\}\\
%   &\xrightarrow{\sigma\mapsto \sigma\otimes \sgn} &
%    \{\check \sigma\in \Irr(\check W(\check \Lambda))\mid \check \sigma\textrm{ occurs in $\Coh^{LR}_{\check \Lambda}(\CK(\check \g,\check \b))(\check \CCD)$}\}
%    \end{eqnarray*}
  is well-defined and bijective. 
\end{prop}

\begin{proof}
\def\ckCCLR{{\check{\CCD}}}
\def\ckbPsi{\check{\bPsi}}
\def\ckfbb{{\check\fbb}}
Let $\sC:=\set{w\in W(\Lambda)| \bPsi_{w}\in \CCD}$.
Then $\sC$ is a double cell of $W(\Lambda)$ in the sense of Kazhdan-Lusztig (\cite{KL}). By \cite{BV2}*{Corollary~2.23}, 
  $\CCD$ has the form 
 \[
 \set{\bPsi_{w'w}|w'\in \D(\Lambda), w\in \sC}
 \qquad
(\textrm{see \cref{wprime} for the definition of $\D(\Lambda)$}).
\]
Let $w_0$ be the element in $W(\Lambda)$ such that $w_0 (\Delta^+(\Lambda)) \cap \Delta^+(\Lambda) = \emptyset$. Similar to $\D(\Lambda)$, we have a subgroup $\D(\ckLambda)$ of $\ckW$.

By \cite{BV2}*{Corollary~2.24}, 
\[
\check \CCD:=\set{\ckbPsi_{w' w_0 w}|w'\in D(\ckLambda),  w\in \sC}
\]
is a double cell in $\Coh_{\ckLambda}(\CK(\ckfgg, \ckfbb,\ckLambda))$, where $\ckbPsi_{w' w_0 w}$ is a basal element of $\Coh_{\check \Lambda}(\CK(\check \g,\check \b))$ defined as in \eqref{psibarw}. 
By \Cref{lem:tauinv}, we have that 
\[
\tau_{\bPsi_{w'w}} = \tau_{\bPsi_{w}}= \Pi(\Lambda)\setminus \tau_{\ckbPsi_{w_0 w }} 
= \Pi(\Lambda)\setminus \tau_{\ckbPsi_{w'' w_0 w }}
\]
for each $w'\in \D(\Lambda)$, $w''\in \D(\ckLambda)$ and  $w\in W(\Lambda)$. 
This implies that $\tau_\CCD$  is dual to $\tau_{\check \CCD}$. 

The uniqueness of $\check \CCD$ follows from \Cref{eqtau1} and the bijection between double cells and special representations as in 
\Cref{cell00}~(b).  See also \cite{FJMN}*{\S~6}. 

The assertion on the operation of tensoring with $\sgn$ is in  \cite{BV2}*{Proposition~2.25}. \end{proof}

The double cell $\check\CCD$ in \Cref{duald} will be referred as the dual (in $\Coh_{\check \Lambda}(\CK(\check \g,\check \b))$) of $\CCD$.

\section{Generalities on coherent families of Casselman-Wallach representations}\label{sec:CW}

Coherent families of group representations were introduced by Schmid \cite{Sch}. See also \cite{Zu} and \cite{SpVo}. 
We refer the reader to \cite[Chapter 7]{Vg}) as a general reference for coherent families in this setting. 

Let $G$ be a real reductive group in Harish-Chandra's class (which may be
linear or non-linear).  The following assumptions hold for the real Lie groups under consideration.

A connected reductive complex Lie group $G_\C$ is fixed, together with a Lie group homomorphism $\iota: G\rightarrow G_\C$ such that its differential $\mathrm d \iota: \Lie(G)\rightarrow \Lie(G_\C)$ has the following two properties:
\begin{itemize}
  \item the kernel of $\mathrm d \iota$ is contained in the center of $ \Lie(G)$;
  \item the image of  $\mathrm d \iota$ is a real form  of  $\Lie(G_{\C})$.
  \end{itemize}
The  analytic weight lattice of $ G_\C$ is identified with a subgroup of $\hha^*$ via $\mathrm d\iota$. We write $Q_\iota\subset \hha^*$ for this subgroup.

When $G_\C=\mathrm{Ad}(\g)$ and $\iota$ is the adjoint representation, $Q_\iota$ equals the root lattice $Q_\g$. In general $Q_\iota$ is $W$-stable and $Q_\g\subseteq Q_\iota \subseteq Q^\g$.
In the rest of this section we assume that $Q=Q_\iota$.

Recall that the complex associated variety of a representation $\pi$ in $\Rep(G)$,
denoted by $\AV_\C(\pi)$, is defined to be the associated variety of $\Ann(\pi)\subset \CU(\g)$.
Let $\sfS$ be  an $\Ad(\g)$-stable Zariski closed subset   of $\Nil(\g^*)$. Let $\Rep_\sfS(G)$ denote the
category of Casselman-Wallach representations of $G$ whose complex associated
variety is contained in $\sfS$.

For every $\mu\in \hha^*$, let $\Rep_\mu(G)$ and $\Rep_{\mu,\sfS}(G)$ respectively denote the full subcategories of  $\Rep(G)$ and $\Rep_{\sfS}(G)$ consisting 
of the Casselman-Wallach representations of generalized infinitesimal character $\mu$.
Denote by $\CK(G)$, $\CK_\mu(G)$, and  $\CK_{\mu, \sfS}(G)$ the Grothendieck groups of
$\Rep(G)$, $\Rep_{\mu}(G)$, and $\Rep_{\mu, \sfS}(G)$ respectively.
The set of irreducible objects in $\Rep_{\mu}(G)$ and $\Rep_{\mu, \sfS}(G)$ will be denoted by $\Irr_{\mu}(G)$ and $\Irr_{\mu, \sfS}(G)$, respectively. 

Recall that we fixed a $Q$-coset $\Lambda=\lambda+Q\subset \hha^*$.

The algebra $\mathcal R(\g, Q)$ is identified with  the Grothendieck group of the category of holomorphic finite-dimensional representations of $G_\C$.  Thus $\CK(G)$ is naturally an $\mathcal R(\g, Q)$-module by using the tensor product and the homomorphism $\iota$. As before, we form the coherent continuation representation $\Coh_\Lambda(\CK(G))$ of $W_\Lambda$.

 The main purpose of this section is to establish good control of $\Coh_{\Lambda}(\CK(G))$, from our knowledge of double cells.  
 The main results concern embedding of coherent continuation representations (\Cref{lem0033}) and 
 relationship between the so-called Harish-Chandra cells and the double cells (\Cref{DualHC}).

\subsection{The parameter set for the coherent continuation representation}\label{extcoh}

We begin with the basic (existence and uniqueness) theorem about coherent families in $\Grt (G)$. 

For every $\nu\in \Lambda$, we have the evaluation map at $\nu$: 
\[
    \ev{\nu} \colon  \Coh_{\Lambda}(\CK(G)) \rightarrow \Grt_{\nu}(G).
\]

 \begin{thm}[{Schmid, Zuckerman}]\label{lem21}
The map  $\ev{\nu}$ is surjective for every $\nu\in \Lambda$, and it is bijective when $\nu$ is regular.
\end{thm}
\begin{proof}
The surjectivity is due to Schmid and Zuckerman, see  \cite{Vg}*{Theorem~7.2.7}.
The injectivity (for $\nu$ regular)  is due to Schmid, see \cite{Vg}*{Proposition~7.2.23}.
\end{proof}

We proceed to describe the parameter set for the coherent continuation representation.    

Suppose that $H$ is a  Cartan subgroup of $G$. Recall that by our convention, its complexified Lie algebra is denoted by $\h$. As usual, denote by $\Delta_\h\subset \h^*$ the root system of $\g$ with respect to $\h$ and $W_\h$ the corresponding Weyl group. Write $\t$ for the complexified Lie algebra of the unique maximal compact subgroup of $H$. A root $\alpha\in \Delta_\h$ is called real if $\alpha |_{\t}=0$, and imaginary if $\check \alpha\in \t$. An imaginary root $\alpha\in \Delta_\h$ is said to be compact if the root spaces $\g_\alpha$ and $\g_{-\alpha}$ are contained in  the  complexified Lie algebra of a common compact subgroup of $G$.

Note that every Casselman-Wallach representation of $H$ is finite dimensional. Since $G$ is in Harish-Chandra's class,
 for every $\Gamma\in \Irr(H)$, 
there is a unique element $\mathrm d \Gamma\in \h^*$ such that the differential of $\Gamma$ is isomorphic to a direct sum of  one-dimensional representations of $\h$ attached to $\mathrm d  \Gamma$. 

For every Borel subalgebra $\b$ of $\g$ containing $\h$, write
\be\label{xib00}
  \xi_\b: \hha\rightarrow \h
\ee
for the linear isomorphism attached to $\b$, namely the inverse of the composition of 
\[
 \h\subset \b\rightarrow \b/[\b,\b]=\hha.
\]
The  transpose inverse of the map \eqref{xib00} is still denoted by $ \xi_\b: \hha^*\rightarrow \h^*$. Write
\[
  W(\hha^*, \h^*):=\{ \xi_\b: \hha^*\rightarrow \h^*\mid \textrm{$\b$ is a Borel subalgebra of $\g$ containing $\h$}\}.
  \]

Recall that $ \Delta^+\subset \hha^*$ denotes the set of positive roots. For every element $\xi\in  W(\hha^*, \h^*)$,  put
\be\label{deltawb}
  \delta(\xi):=\frac{1}{2}\cdot \sum_{\alpha \textrm{ is an imaginary root in $\xi \Delta^+$ }} \alpha- \sum_{\beta \textrm{ is a compact imaginary root in $\xi \Delta^+$ }}\beta \in \h^*.
\ee

\begin{defn}\label{def:parameter}
Write $\sP_{\Lambda}(G)$ for the set of all triples $\upgamma=(H, \xi, \Gamma)$,
where $H$ is a Cartan subgroup of $G$, $\xi\in  W(\hha^*, \h^*)$, and
\[
\Gamma:   \Lambda\rightarrow \Irr(H), \quad \nu\mapsto \Gamma_\nu
\]
 is a map
with the following properties:
\begin{itemize}
  \item $\Gamma_{\nu+\beta}=\Gamma_\nu\otimes \xi(\beta)$ for all $\beta\in Q$ and $\nu\in \Lambda$;
  \item $\mathrm d\Gamma_\nu = \xi(\nu)+\delta(\xi)$ for all $\nu\in  \Lambda$.
  \end{itemize}
Here $\xi(\beta)$ is naturally viewed as a character of $H$ by using the homomorphism $\iota: H\rightarrow H_\C$, and $H_\C$ is the Cartan subgroup of $G_\C$ containing $\iota(H)$.
\end{defn}

\def\olgamma{\overline\gamma}
\def\olPsi{\overline\Psi}
\def\olrX{\overline\rX}

The group $G$ acts on $\sP_{\Lambda}(G)$ in the
the standard way, and we define the set of parameters for
$\Coh_{\Lambda}(\CK(G))$
to be 
\begin{equation}\label{eq:Pram}
  \cP_{\Lambda}(G) := G\backslash  \sP_{\Lambda}(G).
\end{equation}

For each $\gamma \in\cP_{\Lambda}(G)$ that is represented by $\upgamma=(H, \xi, \Gamma)$, by \cite{Vg}*{Theorem~8.2.1}, we  have two $\CK(G)$-valued coherent families  $\Psi_{\gamma}$ and $\olPsi_{\gamma}$ on $\Lambda$  such that
\begin{equation}\label{eq:psigamma}
  \Psi_{\gamma}(\nu)=\rX{(\Gamma_\nu,\xi(\nu))}
  \AND
  \olPsi_{\gamma}(\nu)=\olrX{(\Gamma_\nu,\xi(\nu))}
\end{equation}
for all  regular dominant element $\nu\in \Lambda$.
Here $\rX{(\Gamma_\nu,\xi(\nu))}$ is the standard representation defined in \cite{Vg}*{Notational Convention~6.6.3}
  and  $\olrX(\Gamma_\nu,\xi(\nu))$ is its unique irreducible subrepresentation (see \cite{Vg}*{Theorem 6.5.12}). 
  
By Langlands classification,
$\set{\olPsi_{\gamma}}_{\gamma\in \cP_{\Lambda}(G)}$ is a basis of
$\Coh_{\Lambda}(\cK(G))$ (\cite{Vg}*{Theorem~6.6.14}), and we view $\Coh_\Lambda(\cK(G))$ as a basal representation of $W_\Lambda$ with this basis. 
The family  $\set{\Psi_{\gamma}}_{\gamma\in \cP_{\Lambda}(G)}$ is also a basis of $\Coh_{\Lambda}(\cK(G))$ (\cite{Vg}*{Proposition~6.6.7}).

\begin{remark}
    For every dominant element $\nu\in \Lambda$, the set 
    \[
    \{\gamma \in\cP_{\Lambda}(G)\mid \olPsi_{\gamma}(\nu)\neq 0\}
    \]
    is identified with the set of final characters for $G$ with infinitesimal character $\nu$, by sending  the class of $\upgamma=(H, \xi, \Gamma)\in \sP_{\Lambda}(G)$ to the $G$-conjugacy class of 
    \[
    (H,\textrm{the set of imaginary roots in $\xi(\Delta^+)$}, \Gamma_\nu). 
    \] See \cite[Chapter 11]{ABV}. 
\end{remark}

The coherent continuation representation  $\Coh_{\Lambda}(\CK(G))$ is independent of $Q$ in the sense of the following lemma.

\begin{lem}%(\cf \cite[Lemma 7.2.6]{Vg})
For a $Q_\g$-coset $\Lambda_\g$ in $\Lambda$, form the coherent continuation representation $\Coh_{\Lambda_\g}(\CK(G))$ by using the adjoint representation $G\rightarrow \Ad(\g)$. Then the restriction yields a linear isomorphism
\[
  \Coh_{\Lambda}(\CK(G))\xrightarrow{\sim} \Coh_{\Lambda_\g}(\CK(G)).
\]
\end{lem}
\begin{proof}
This is a consequence of the fact that the restriction yields a bijective map
\[
   \cP_{\Lambda}(G)\rightarrow  \cP_{\Lambda_\g}(G).
\]
\end{proof}

The cross action of $W_{\Lambda}$ 
on the set $
  \sP_{\Lambda}(G)$ is defined by (\cite{V4}*{Definition~4.2}):
\[
  w\cross(H, \xi, \Gamma)=\left(H, \xi w^{-1}, (\nu\mapsto \Gamma_\nu\otimes (\xi w^{-1}\nu+\delta(\xi w^{-1})-\xi \nu-\delta(\xi)) )\right).
\]
This commutes with the action of $G$ and thus descends to an action on $\cP_{\Lambda}(G)$: 
\begin{equation}\label{eq:cross}
 W_{\Lambda}\times   \cP_{\Lambda}(G)\rightarrow \cP_{\Lambda}(G), \quad (w, \gamma)\mapsto w\cross \gamma.
 \end{equation}

Denote by $W_{\h,\t}$  the stabilizer of the space $\t$ in the Weyl group  $W_\h$, whose detailed structure is analyzed in \cite{V4}*{Section~3}.
Write 
$
\Delta_{\h,\mathrm{im}}$ for the set of  imaginary roots in $\Delta_\h$, which is a root system. The corresponding Weyl group is denoted by $W(\Delta_{\h,\mathrm{im}})$, which is identified with a normal subgroup of  $W_{\h,\t}$.  
 Then there is a unique quadratic character \[
 \mathrm{sgn}_\mathrm{im}: W_{\h,\t}\rightarrow \C^\times
 \]
 such that 
\begin{itemize}
    \item its restriction to $W(\Delta_{\h,\mathrm{im}})$ equals the sign character, and 
    \item its restriction to $W_{\h,\t, \Delta^+_{\h,\mathrm{im}}}$ is trivial for some (and hence all) positive system $\Delta^+_{\h,\mathrm{im}}$ of $\Delta_{\h,\mathrm{im}}$, where $W_{\h,\t, \Delta^+_{\h,\mathrm{im}}}$ denotes the stabilizer of the set $\Delta^+_{\h,\mathrm{im}}$ in $W_{\h,\t}$.
\end{itemize}
Note that with the notation as above, we have that 
\[
W_{\h,\t}= W_{\h,\t, \Delta^+_{\h,\mathrm{im}}}\ltimes W(\Delta_{\h,\mathrm{im}}).
\]

Since $G$ is in Harish-Chandra's class, the real Weyl group
\[
W_H:=\textrm{(the normalizer of $H$ in $G$)}/H
\]
 is identified with a subgroup of $W_{\h,\t}$.  
Choose a representative $\upgamma=(H, \xi, \Gamma)$
      for an element $\gamma\in \cP_{\Lambda}(G)$. Write $W_{\gamma}$ for the stabilizer of $\gamma$ in $W_{\Lambda}$ under the cross action.
Then for any $w\in W_\gamma$, $\xi w
\xi^{-1}\in W_H\subset W_{\h,\t}$, and we have a quadratic character
\begin{equation}\label{eq:imaginary}
  \mathrm{sgn}_{\gamma}: W_{\gamma}\rightarrow \C^\times, \qquad w\mapsto \mathrm{sgn}_\mathrm{im}(\xi w \xi^{-1}).
\end{equation}
This quadratic character is independent of  the representative $\upgamma$.

The coherent continuation representation may be computed by using the basis of standard modules $\set{\Psi_{\gamma}}_{\gamma\in \cP_{\Lambda}(G)}$ (see \cite[Section 14]{V4}).
The following result is due to Barbasch-Vogan, in a suitably modified
form from \cite{BV.W}*{Proposition~2.4}. As its proof follows the same
line as that of \cite{BV.W}*{Proposition~2.4}, we just state the precise result.

\begin{thm}[{\cf \cite{BV.W}*{Proposition~2.4}}]
  \label{thm:cohHC}
As a representation of $W_{\Lambda}$,
  \[
    \Coh_{\Lambda}(\CK(G)) \cong \bigoplus_{\gamma}
    \Ind_{W_{\gamma}}^{W_{\Lambda}}  \mathrm{sgn}_{\gamma},
  \]
  where $\gamma$ runs over a representative set of the $W_{\Lambda}$-orbits
  in $\cP_{\Lambda}(G)$ under the cross action.
\end{thm}

\subsection{Some properties of the coherent continuation representation $\Coh_{\Lambda}(\CK(G))$}

We present some properties of the coherent continuation representation as a basal representation.   

\begin{lem}\label{annlc}
Let $\Psi$ be a basal element in $\Coh_\Lambda(\CK(G))$. Then there is a unique left cell $\CCL_\Psi$ in  $\Coh_\Lambda(\CK(\g,\b))$ such that $ \Ann(\Psi(\nu))=I_{\nu, \CCL_\Psi}$ for every dominant element $\nu\in \Lambda$.
\end{lem}
\begin{proof}
Since a coherent family restricted to the dominant cone can be constructed via the translation functor,   
 the proposition follows from \cite{V1}*{Lemma~2.7 and Lemma~2.8}.  
\end{proof}

For every  basal element $\Psi\in \Coh_\Lambda(\CK(G))$, define its $\tau$-invariant to be 
\begin{equation}\label{def:tauH} 
\tau_\Psi:=\tau_{\CCL_\Psi},\end{equation}
where $\CCL_\Psi$ is as in \Cref{annlc}. 
A basic fact about $\tau_\Psi$ is the following: For each simple reflection $s$ in $W(\Lambda)$, $s\in \tau_\Psi$ if and only if 
$s\cdot \Psi=-\Psi$. See \cite{Vtau}*{Theorem~2.4} and \cite{Vg}*{Corollary~7.3.9}.

\begin{lem} [\cite{Vg}*{Corollary~7.3.23}]
\label{lemirr}
Let $\nu\in \Lambda$ be a dominant element. 
Then for every  basal element $\Psi\in \Coh_\Lambda(\cK(G))$,  $ \Psi(\nu)$ is
\[
\left\{
\begin{array}{ll}
irreducible,\quad & \text{if $\tau_{\Psi}\cap W_\nu = \emptyset$;}\\
0, \quad &\textrm{otherwise}.
\end{array}
\right.
\]
\end{lem}

Let $\CCD$ be a double cell in $\Coh_\Lambda(\CK(\fgg,\b))$. Denote by 
\begin{equation}
    \label{eq:subCoh}
\Coh_{\Lambda, \CCD}(\CK(G))
\end{equation}
 the basal subspace of $\Coh_\Lambda(\CK(G))$ spanned by all the basal elements $\Psi$ such that $\CCL_\Psi\subseteq \la \CCD \ra_{LR}$.
Recall that $\la \CCD \ra_{LR}$ is the smallest basal subrepresentation of $\Coh^{LR}_\Lambda(\CK(\fgg,\b))$ containing $\CCD$.

The following lemma should be known in the expert community and it follows from the properties of tensoring with finite dimensional representations and the containment of annihilators. We include a proof for the sake of completeness. 
\begin{lem}\label{annl2} For every double cell $\CCD$ in $\Coh_\Lambda(\CK(\fgg,\b))$, the 
space $\Coh_{\Lambda, \CCD}(\CK(G))$ is a basal $W(\Lambda)$-subrepresentation of $\Coh_{\Lambda}(\CK(G))$. 
\end{lem}
\begin{proof}
By passing to a suitable covering group, we  assume without loss of generality that the derived group of $G_\C$ is simply connected. 

Let $\Psi$ be a basal element in $\Coh_{\Lambda, \CCD}(\CK(G))$, and let $\alpha$ be a simple root in $\Delta(\Lambda)$. We need to show that $s_\alpha\cdot \Psi\in \Coh_{\Lambda, \CCD}(\CK(G))$, where $s_\alpha\in W(\Lambda)$ is the simple reflection associated to $\alpha$. 

Our assumption of simply connectedness implies that there is a regular dominant element $\mu\in \Lambda$ such that $\la \mu, \check \alpha \ra \in 2\Z$. Put 
\[
\nu:=  \mu - \frac{ \inn{\mu}{\ckalpha}}{2} \, \alpha\in \Lambda,
\]
which is a dominant element such that $W_\nu=\{1, s_\alpha\}$. 

If $\Psi(\nu)=0$, then $s_\alpha\cdot \Psi=-\Psi\in \Coh_{\Lambda, \CCD}(\CK(G))$.
Now we assume that $\Psi(\nu)\neq 0$. Then $\Psi(\nu)$ is irreducible.

Let $F$ denote the irreducible representation of $G$ with extremal weight  $\half \inn{\mu}{\ckalpha}\alpha$ that occurs in the symmetric power $\oS(\g)$.
Then there are two nonzero $G$-homomorphisms 
\[
\Psi(\mu)\xrightarrow{\iota_1} \Psi(\nu)\otimes F\xrightarrow{\iota_2}  \Psi(\mu)
\]
such that $\iota_2\circ \iota_1=0$.
Moreover, 
$\Psi(\mu)$  occurs 
in the composition series of $\Psi(\nu)\otimes F$ with multiplicity two, and 
\be\label{sap2}
 ( s_\alpha\cdot \Psi)(\mu)=\mathrm{Pr}_\mu(\Psi(\nu)\otimes F)-\Psi(\mu)\in \CK(G),
\ee
where $\mathrm{Pr}_\mu$ denotes the projection map
\be\label{prinf}
 \mathrm{Pr}_\mu: \CK(G)=\bigoplus_{\lambda\in W\backslash \hha^*} \CK_\lambda(G)\rightarrow \CK_\mu(G).
\ee
See \cite{Vg}*{Theorem~7.3.16 and Corollary ~7.3.18}.

By \eqref{sap2}, we know that $s_\alpha\cdot \Psi$ is a positive integral linear  combination of basal elements of  $\Coh_{\Lambda}(\CK(G))$. 
Let $\Psi'$ be a basal element that occurs in the  positive integral combination. Then 
\[
  \Ann(\Psi'(\mu)) \supseteq \Ann((\Psi(\nu)\otimes F)). 
\]

Pick an element $w\in W$ such that  $\Ann(\Psi(\mu)) =\Ann(\bPsi_w(\mu))$, where $\bPsi_w$ is defined as in \eqref{psibarw}. 
The  translation principle \cite{V1}*{Lemma~2.7} implies that
$\Ann(\Psi(\nu)) = \Ann(\overline \Psi_w(\nu))$. In particular, $\bPsi_w(\nu)$ is nonzero and hence irreducible.  

As in the group case,  $\bPsi_w(\mu)$ occurs as a subrepresentation of $\bPsi_w(\nu)\otimes F$ (and also as a quotient representation), and occurs in
the composition series of $\bPsi_w(\nu)\otimes F$ with multiplicity two. Moreover, 
\be\label{sap3}
 ( s_\alpha\cdot \bPsi_w)(\mu)=\mathrm{Pr}_\mu(\bPsi_w(\nu)\otimes F)-\bPsi_w(\mu)\in \CK(\g,\b),
\ee
where $\mathrm{Pr}_\mu$ denotes  the projection map analogous to \eqref{prinf}.

We have an ideal  $I\subset \CU(\g)$ that is a finite product of primitive ideals of $\CU(\g)$ with infinitesimal character distinct from $\mu$,
and a sequence $w_1, w_2, \cdots, w_{k-1}, w_k=w$ in $W$ ($k\in \BN^+$) such that 
\[
\Ann(\bPsi_w(\nu)\otimes F)\supseteq \Ann (\bPsi_{w}(\mu))\cdot \Ann (\bPsi_{w_1}(\mu))\cdot \Ann (\bPsi_{w_2}(\mu))\cdot \dots \cdot  \Ann( \bPsi_{w_k}(\mu)) \cdot I,
\]
and
\begin{equation}\label{eq:sPsiw}
  (s_\alpha\cdot \bPsi_w)= \bPsi_{w_1}+\bPsi_{w_2} \cdots+ \bPsi_{w_k}\in \CK(\g,\b).
\end{equation}

Then 
\[
  \Ann(\Psi'(\mu))\supseteq \Ann (\bPsi_{w}(\mu))\cdot \Ann (\bPsi_{w_1}(\mu))\cdot \Ann (\bPsi_{w_2}(\mu))\cdot \dots \cdot  \Ann (\bPsi_{w_k}(\mu))\cdot I.
\]
Since the primitive ideal $\Ann(\Psi'(\mu))$ is prime, 
%and $w=w_i$ for some $i\in\{1,2,\dots, k\}$, 
the above implies that $\Ann(\Psi'(\mu))\supseteq \Ann(\bPsi_{w_i}(\mu))$
for some $i\in \{1,2, \dots,k\}$. Hence $\CCL_{\Psi'} \subset \la \bPsi_{w_i}\ra_L$ by \Cref{primitiveii}.
On the other hand, $\bPsi_{w_i}$ is in $\braket{\bPsi_w}_{R}$ by \eqref{eq:sPsiw}. 
We conclude that $\Psi'\in \Coh_{\Lambda, \CCD}(\CK(G))$. The proof is complete. 
\end{proof}

\subsection{The coherent continuation representation $\Coh_\Lambda(\CK(\fgg,H, \b'))$}

  Let $H$ be a Cartan subgroup of $G$. Since
  a Cartan subgroup can be disconnected, we need a mild
  generalization of the category of $\b$-modules considered earlier. Let $\b'$ be a Borel subalgebra
  of $\g$ containing the complexified Lie algebra $\h$ of $H$. (In the sequel, we will need to use several Cartan subgroups. We shall reserve the notation $\b$ for a fixed Borel subalgebra of $\g$ as in \Cref{sub:CohHWM}.) 
  Recall that a $(\g, H)$-module is defined to be a $\g$-module $V$ together with a locally-finite representation of $H$ on it such that
     \begin{itemize}
     \item
        $h\cdot (X\cdot (h^{-1}\cdot u))=(\Ad_h(X))\cdot u$, for all $h\in H, X\in \CU(\g), u\in V$ ($\Ad$ stands for the Adjoint representation);
        \item the differential of the representation of $H$ and the restriction of the representation of $\g$ yields the same representation of $\h$ on $V$.
     \end{itemize}

Let $\Rep(\g,H,\b')$ denote the category of finitely generated $(\g, H)$-modules that  are unions of finite-dimensional $\b'$-submodules. 
As before, we form the Grothendieck group  $\CK(\g,H, \b')$, and the coherent continuation representation $\Coh_{\Lambda}(\CK(\g,H, \b'))$ of $W(\Lambda)$.

Write $\sP_{\Lambda}(\g, H, \b')$ for the set of all pairs $(w, \Gamma)$ where $w\in W$, and 
\[
\Gamma:   \Lambda\rightarrow \Irr(H), \quad \nu\mapsto \Gamma_\nu
\]
is a map
 such that 
\begin{itemize}
  \item $\Gamma_{\nu+\beta}=\Gamma_\nu\otimes \xi_{\b'}(\beta)$ for all $\beta\in Q$ and $\nu\in \Lambda$;
  \item $\mathrm d  \Gamma_\nu= \xi_{\b'}(w\nu-\rho)$ for all $\nu\in  \Lambda$.
  \end{itemize}
Here  $\xi_{\b'}$ is defined in \cref{xib00},  $\xi_{\b'}(\beta)$ is viewed as a character of $H$ by  pulling-back through
the homomorphism $\iota: H\rightarrow H_\C$,  
and $H_\C$ is the Cartan subgroup of $G_\C$ containing $\iota(H)$. 
%\xi_{\b'}

\def\olgamma{\overline\gamma}
\def\olPsi{\overline\Psi}
\def\olrX{\overline\rX}

For each $\sigma\in \Irr(H)$, put
\[
  \mathrm M(\sigma):=\CU(\g)\otimes_{\CU(\b')} \sigma,
\]
which is a module in $\Rep(\g,H,\b')$, where the $\CU(\g)$-action is given by the left multiplication, and the $H$-action is given by
\[
 h\cdot (X\otimes u):=\Ad_h(X)\otimes h \cdot u, \qquad h\in H, \, X\in \CU(\g),\, u\in \sigma.
\]
As before, the $(\g,H)$-module $\mathrm M(\sigma)$ has a unique irreducible quotient, to be denoted by $\mathrm L(\sigma)$.

For each $(w, \Gamma)\in\sP_{\Lambda}(\g,H,\b')$,  we  have two $\CK(\g,H,\b')$-valued coherent families  $\Psi_{w,\Gamma}$ and $\olPsi_{w, \Gamma}$ on $\Lambda$  such that
\begin{equation}\label{eq:psiGamma}
  \Psi_{w,\Gamma}(\nu)=\mathrm M(\Gamma_\nu) \quad\textrm{for all $\nu\in \Lambda$},
  \end{equation}
  and
  \begin{equation}
  \olPsi_{w,\Gamma}(\nu)=\mathrm L(\Gamma_\nu) \quad\textrm{for all regular dominant element $\nu\in \Lambda$}.
\end{equation}
As before, both 
$\{ \Psi_{w,\Gamma}\}_{(w,\Gamma)\in\sP_{\Lambda}(\g,H,\b')}$ and $\{ \olPsi_{w,\Gamma}\}_{(w,\Gamma)\in\sP_{\Lambda}(\g,H,\b')}$ are bases of $\Coh_{\Lambda}(\CK(\g,H, \b'))$. We view $\Coh_{\Lambda}(\CK(\g,H, \b'))$ as a basal representation of $W(\Lambda)$ by using the second basis. 

The analog of \Cref{annlc}  for   $\Coh_\Lambda(\CK(\g, H, \b'))$ still holds (with the same proof), and we define similarly the left cell $\CCL_\Psi$ in  $\Coh_\Lambda(\CK(\g,\b))$ for every basal element $\Psi (=\olPsi_{w,\Gamma})
\in \Coh_\Lambda(\CK(\g, H, \b'))$.

Let $\CCD$ be a double cell in $\Coh_\Lambda(\CK(\fgg,\b))$ as before. Define
\[\Coh_{\Lambda, \CCD}(\CK(\g,H,\b'))\] to be the basal subspace of $\Coh_\Lambda(\CK(\g,H,\b'))$ spanned by all the basal elements $\Psi$ such that $\CCL_\Psi\subseteq \la \CCD\ra_{LR}$.
By the same argument as that of \Cref{annl2}, the space $ \Coh_{\Lambda, \CCD}(\CK(\g,H,\b'))$ is a basal $W(\Lambda)$-subrepresentation of $\Coh_{\Lambda}(\CK(\g,H,\b'))$.

 \begin{lem}\label{lem0022}
The representation $\Coh_{\Lambda, \CCD}(\CK(\g,H, \b'))$ of $W(\Lambda)$ is isomorphic to a subrepresentation of $(\la \CCD\ra_{LR})^k$, for some $k\in \BN$.
     \end{lem}
\begin{proof}
Let $Q$ act on the set $\Irr(H)$ by
\[
   \beta\cdot \sigma=\sigma\otimes \xi_{\b'}(\beta), \qquad\beta\in Q, \ \sigma\in \Irr(H), 
\]
where $\xi_{\b'}(\beta)$ is as in \eqref{xib00}. 
For each $Q$-orbit $\Omega \subset \Irr(H)$, let $\Rep_{\Omega}(\g,H,\b')$ be the full subcategory of $\Rep(\g,H,\b')$
whose objects are modules $V$ such that every irreducible subquotient of $V|_H$ (viewed as a representation of $H$) belongs to $\Omega$.
Write $\CK_{\Omega}(\g,H,\b')$ for the  Grothendieck group of the category $\Rep_{\Omega}(\g,H,\b')$. Then $\Coh_{\Lambda}(\CK_{\Omega}(\g,H,\b'))$ is a basal $W(\Lambda)$-representation. 
As $W(\Lambda)$-representations, 
\[
\Coh_{\Lambda, \CCD}(\CK(\g,H, \b'))=\bigoplus_{i=1}^k 
\Coh_{\Lambda, \CCD}(\CK_{\Omega_i}(\g,H,\b')),
%\Coh_{\Lambda, \CCD, \Omega_i}(\g,H, \b'),
\]
where $\Omega_1, \Omega_2, \cdots, \Omega _k$ ($k\in \bN$) are  certain 
$Q$-orbits in $\Irr(H)$, and 
\[
\Coh_{\Lambda, \CCD}(\CK_{\Omega_i}(\g,H,\b')):=\Coh_{\Lambda}(\CK_{\Omega_i}(\g,H,\b'))\cap 
\Coh_{\Lambda, \CCD}(\CK(\g,H, \b')).
\]

Define a linear isomorphism
\[
  T_{\b, \b'}: \CK(\g,\b)\rightarrow \CK(\g,\b')
  \]
  such that $T_{\b,\b'}(\mathrm M(\g,\b,\nu))= \rM(\g,\b',\nu)$ for all $\nu\in \hha^*$.
The map  $ T_{\b, \b'}$ is $W$-equivariant and induces an isomorphism
\[
  \Coh_{\Lambda}(\CK(\g,\b))\rightarrow \Coh_{\Lambda}(\CK(\g,\b')). 
\]
By using this isomorphism we  assume without loss of generality that $\b = \b'$. 

Denote by $\cF$ the forgetful functor from  
$\Rep(\g,H,\b')$ to $\Rep(\g,\b')$.
Since $\Psi\mapsto \cF\circ \Psi$ gives an injective $W(\Lambda)$-homomorphism 
from $\Coh_{\Lambda, \CCD}(\CK_{\Omega_i}(\g,H,\b'))$
to $\braket{\CCD}_{LR}$,  the lemma follows. 
\end{proof}

\subsection{An embedding of  coherent continuation representations}
As before $\CCD$ is a double cell in $\Coh_\Lambda \CK(\fgg,\b)$. The purpose of this subsection is to prove the following proposition.

\begin{prop}\label{lem0033}
The representation $ \Coh_{\Lambda, \CCD}(\CK(G))$  of $W(\Lambda)$ is isomorphic to a subrepresentation of $(\la \CCD\ra_{LR})^k$, for some $k\in \BN$.
\end{prop}

Let $\set{H_1, H_2, \cdots, H_r}$ ($r\in \bN^+$) be a set of representatives of the
conjugacy classes of Cartan subgroups of $G$. For each $i=1,2,\cdots, r$, 
fix a Borel subalgebra $\b_i$ of $\g$ that contains the complexified Lie algebra of $H_i$.

 \begin{lem}\label{cor:HC.embed}
  % [\cite{Cas}*{Theorem~3.1}]
 There is an injective $\mathcal R(\g, Q)$-module homomorphism
 \[
\gamma_{G}: \Grt(G)\rightarrow  \bigoplus_{i=1}^{r} \Grt(\fgg,H_{i},\b_{i})
 \]
 such that
 \be\label{gammag}
   \gamma_{G}(\CK_I(G))\subset  \bigoplus_{i=1}^{r} \Grt_{I}(\fgg,H_{i},\b_{i})
 \ee
 for every ideal $I$ of $\CU(\g)$.
 Here $\CK_I(G)$ is the Grothendieck group of $\Rep_I(G)$,  $\Rep_I(G)$ is the full subcategory of $\Rep(G)$ consisting of the representations that are annihilated by  $I$, and $\Grt_{I}(\fgg,H_{i},\b_{i})$ is the similarly defined subspace of $\Grt(\fgg,H_{i},\b_{i})$.
 \end{lem}
\begin{proof}
This follows from the work of Casian (\cite{Cas}). 
Since the proposition is not explicitly stated in \cite{Cas},
we briefly recall the argument of Casian for the convenience of the reader.

 Let $\n_i$ denote the nilpotent radical of $[\g,\g]\cap \b_i$ ($i=1,2, \cdots,r$).
 For every $q\in \Z$, let $\gamma_{\n_i}^q$ denote the $q$-th right derived functor of the following left exact functor from the category of $\g$-modules to itself:
 \[
   V\rightarrow \{u\in V\,|\, \n_i^k \cdot u=0\textrm{ for some $k\in \bN^+$}\}.
 \]

Fix a Cartan involution $\theta$ of $G$ and write $K$ for its fixed point group (which is a maximal compact subgroup of $G$).  Without loss of generality we assume that all $H_i$'s are $\theta$-stable.

For every Casselman-Wallach representation $V$ of $G$, write $V_{[K]}$ for the space of $K$-finite vectors in $V$, which is a $(\g,K)$-module of finite length. Then   $\gamma_{\n_i}^q(V_{[K]})$ is naturally a representation in $\Rep(\g, H_i, \b_i)$ (\cite[Corollary 4.9]{Cas}).

We define a linear map
 \[
\gamma_{G}: \Grt(G)\rightarrow  \bigoplus_{i=1}^{r} \Grt(\fgg,H_{i},\b_{i})
 \]
 given by
 \[
   \gamma_{G}(V)= \left\{\sum_{q\in \Z} (-1)^{q} \gamma^{q}_{\n_i}(V_{[K]})\right\}_{i=1,2, \cdots, r}
 \]
for every Casselman-Wallach representation $V$ of $G$. A form of the Osborne conjecture 
(see \cite[Theorem 3.1]{Cas}) and \cite[Corollary 4.9]{Cas}) implies that the map $\gamma_{G}$ is well-defined and injective.

Proposition~4.11 of \cite{Cas} implies that the functor $\gamma_{\n_i}^q$ commutes with tensor product with finite-dimensional representations. Thus $\gamma_{G}$ is an $\mathcal R(\g, Q)$-homomorphism. Finally, \cite[Corollary 4.15]{Cas} implies that $\gamma_{G}$ satisfies the property in \eqref{gammag}.
\end{proof}

\begin{proof}[Proof of \Cref{lem0033}]
\Cref{cor:HC.embed} implies that the representation $\Coh_{\Lambda,\CCD}(\CK(G))$
of $W(\Lambda)$ is isomorphic to a subrepresentation of
$\bigoplus_{i=1}^r \Coh_{\Lambda,\CCLR}(\CK(\g,H_i,\b_i))$.
Together with  \Cref{lem0022}, this implies \Cref{lem0033}.
\end{proof}

\subsection{Harish-Chandra cells and double cells}\label{subsec:tau-HC}

We now discuss the consequences of \Cref{lem0033}.
We view $\Coh_{\Lambda}(\CK(G))$ as a basal representation of $W(\Lambda)$. A cell in $\Coh_{\Lambda}(\CK(G))$ will be called a Harish-Chandra cell. 

We fix a Harich-Chandra cell $\cC$.  

\def\CohLamLR{\Coh_\Lam^{LR}(\cK(\fgg,\fbb) }
\begin{lem} \label{lem:LRcell.HC} \label{hcass22}
There is a unique double cell $\CCLR$ in $\Coh_{\Lambda}(\CK (\g,\b))$ 
%(see \eqref{eq:LRO})
such that $\CCL_{\Psi}\subset \CCLR$ for every basal element $\Psi$ in $\CC$. 
Here $\CCL_{\Psi}$ is the left cell defined in \Cref{annlc}. 

\end{lem}
\begin{proof}
Suppose $\Psi_1$ and $\Psi_2$ are two basal elements in $\cC$. 
For $i=1,2$, let $\CCLR_{\Psi_i}$ be the double cell containing $\CCL_{\Psi_i}$. By \Cref{annl2}, we have $\braket{\CCLR_{\Psi_1}}_{LR} \subset \braket{\CCLR_{\Psi_2}}_{LR}$ since 
$\Psi_1 \in \braket{\Psi_2}$. Since the role of $\Psi_1$ and $\Psi_2$ is symmetric, we also have 
$\braket{\CCLR_{\Psi_2}}_{LR} \subset \braket{\CCLR_{\Psi_1}}_{LR}$. 
Hence $\braket{\CCD_{\Psi_1}}_{LR} = \braket{\CCD_{\Psi_2}}_{LR}$. Together with \Cref{cell00}~(a) and (b), it implies that $\CCD_{\Psi_1} = \CCD_{\Psi_2}$. This proves the existence and uniqueness of the double cell. % and proved the lemma. 
\end{proof}

\begin{defn} \label{defn:spe-HC} We call the double cell $\CCLR$ as in \Cref{lem:LRcell.HC} the double cell attached to $\CC$ and $\sigma_\CC := \sigma_{\CCLR}$ 
the special irreducible representation of $W(\Lambda)$ attached to $\CC$.
\end{defn}

The following conjecture is widely anticipated, although to our knowledge no proof has appeared in the literature. See \cite[page 1055]{V4}.

\begin{conj}\label{conjcell}
Suppose that $\sigma\in \Irr(W(\Lambda))$ occurs in 
%the cell representation 
$\Coh_{\Lambda}( \CK(G))(\CC)$. Then 
$\sigma \approx \sigma_\CC$.
%where $\sigma_{\CC}$
%is the special irreducible representation of $W(\Lambda)$ attached to $\CC$. 
%$ \in \Irrsp(W(\Lambda))$ is given by \Cref{defn:spe-HC}.
% \[
%  \{ \sigma\in \Irr(W(\Lambda))\mid \sigma\textrm{ occurs in the cell representation $\Coh_{\Lambda}( \CK(G))(\CC)$}\}
% \]
% \[
%  \Set{ \sigma\in \Irr(W(\Lambda)) | \sigma \approx \sigma_{\CC}}
% \]
% is contained in the  double cell in $\Irr(W(\Lambda))$ containing $\sigma_\CC$.
\end{conj}

The conjecture holds for complex semisimple groups (see \cite{BVUni}*{Section~3}). It also holds for unitary groups by \cite{Bo}*{Theorem~5}.  
We note McGovern's observation (\cite[Page 213]{Mc}) which amounts to the  assertion in \Cref{conjcell} for reductive linear groups in Harish-Chandra's class.
It appears to us that the argument is inadequate as presented.
In what follows we will give a proof assuming a certain equality of $\tau$-invariants (\Cref{wh:tau}) and a weak form of Vogan duality (\Cref{wh:dual}),
both of which hold true for classical groups (including the real metaplectic group).

Define the $\tau$-invariant of $\CC$ to be the set
\[
  \tau_\CC:=\set{\tau_\Psi : \Psi\in \CC}. 
\]
By \Cref{lem:LRcell.HC}, we have $\tau_\CC\subset \tau_{\CCLR}$, 
where $\CCLR$ is the  double cell attached to $\CC$.

To continue our discussion, we make the following hypothesis. 
\begin{whyp}
\label{wh:tau}
With the notation as above, 
$\tau_\CC = \tau_{\CCLR}.$ 
\end{whyp}

\begin{remark}
\begin{enumerate}[label=\arabic*., wide=0em]
In view of \Cref{eqtau0}, \Cref{wh:tau} will follow from \cite{V4}*{Corollary~14.11} 
which asserts that the special representation $\sigma_{\CCLR}$ 
occurs in the cell representation of $\cC$.  
However, \cite{V4}*{Corollary~14.11}  implicitly assumes that the Gelfand-Krillov dimension of $\Psi'$ (namely the Gelfand-Kirillov dimension of $\Psi'(\nu)$, where $\nu$ is an arbitrary regular dominant element in $\Lambda$) 
is strictly less than that of $\Psi$  
for any two basal elements $\Psi'$ and $\Psi$ in $\Coh_{\Lambda}(\CK(G))$ such that  $\braket{\Psi'} \subsetneq \braket{\Psi}$. We thank David Vogan for confirming that this has not been established in the literature.    
\end{enumerate}
\end{remark}

\begin{prop}\label{thm:tau.ABCDG}
Suppose that  the Coxeter group $W(\Lambda)$ has no simple factor of type $F_4$, $E_6$, $E_7$, or $E_8$. 
Then \Cref{wh:tau} holds.  
\end{prop}

We proceed to prove \Cref{thm:tau.ABCDG} by reducing it to 
a question on double cells, by the use of the so-called edge transport theorems (see \cite{Gar.D4}).

\def\bPsip{\overline{\Psi'}}
\def\Psip{\Psi'}
\def\cBLG{\cB_{\Lambda}(G)}
\def\cBL{\cL_{\Lambda}}
\def\CCLp{{\cC'}^L}
\def\TL{T^L}

Recall the integral root system $\Delta(\Lambda)$ from \eqref{eq:IR}. 

Let $\cBLG$ be the set of basal elements in $\Coh_{\Lambda}(\cK(G))$. 
For two simple roots $\alpha,\beta\in \Delta(\Lambda)$ spanning a root system of type $A_2$, $B_2$ (or $C_2$), or  $G_2$, define a map 
\[
T_{\alpha,\beta}\colon \cBLG \rightarrow \set{\text{subset of }\cBLG}
\]
such that, for each $\Psi\in \cBLG$,
\[
T_{\alpha,\beta}(\Psi) =\left\{ \begin{array}{l}
     \{\Psi'\in \cBLG :
 s_\alpha \in \tau_{\Psi'}, s_\beta \notin \tau_{\Psi'},
\text{ and } \\
\qquad \quad \textrm{the coefficient of $\Psi'$ in $s_\alpha\cdot \Psi$ is non-zero}\}, \qquad
\textrm{if $s_\alpha \notin \tau_{\Psi}$ and  $s_\beta \in \tau_{\Psi}$;} \\
\emptyset, \qquad \qquad \qquad \qquad  \qquad \qquad \qquad \qquad \qquad \qquad \qquad \quad \textrm{otherwise}. \end{array}
\right.
\]
Here and as usual, $s_\alpha$ denotes the simple reflection associated to $\alpha$.

Let $\cL_\Lambda$ denote the set of left cells in $\Coh_{\Lambda}(\cK(\fgg,\fbb))$. For each $w\in W$, write $\CCL_w\in \cL_\Lambda$ for the left cell containing the basal element $\overline \Psi_w$. 
The map $T_{\alpha,\beta}$ can also be defined for the category of highest weight modules (see \cite{Vtau}*{Theorem~3.2})
and descends to a map on the set of left cells (see \cite{Vtau}*{Corollary~3.6 and 3.9}).
More precisely, we have a map  

\[
\TL_{\alpha,\beta} \colon \cBL \rightarrow \set{\text{subset of }\cBL}
\]
such that, for each $w \in W(\Lambda)$, \[
\TL_{\alpha,\beta}(\CCL_w) =\left\{ \begin{array}{l}
     \{\CCL_{w'} :
w' \in W(\Lambda), s_\alpha \in \tau_{\CCL_{w'}}, s_\beta \notin \tau_{\CCL_{w'}},
\text{ and } \\
\qquad \quad \textrm{the coefficient of $\bPsi_{w'}$ in $s_\alpha\cdot \bPsi_w$ is non-zero}\}, \quad
\textrm{if $s_\alpha \notin \tau_{\CCL_w}$ and  $s_\beta \in \tau_{\CCL_w}$;} \\
\emptyset, \qquad \qquad \qquad \qquad  \qquad \qquad \qquad \qquad \qquad \qquad \qquad \quad \textrm{otherwise}. \end{array}
\right.
\]

When four simple roots $\alpha_1,\alpha_2,\alpha_3,\alpha_4$ 
in $\Delta(\Lambda)$
span a root system of type $D_4$ ,
one can define $D_4$-maps \[
T_{\alpha_i}
\colon \cBLG\rightarrow 
\set{\text{subset of }\cBLG}
\]
and 
\[
\TL_{\alpha_i}
\colon \cBL\rightarrow 
\set{\text{subset of }\cBL}
\qquad
(i=1,2,3,4) 
\]
as in \cite{GV}*{Theorem~2.15}. 

Let $(T,T^L)$ be a pair $(T_{\alpha, \beta}, T_{\alpha, \beta}^L)$ or $(T_{\alpha_i}, T_{\alpha_i}^L)$ as above. Then 
\begin{itemize}
\item
  for each basal element  $\Psi\in \cBLG$, every element  in $T(\Psi)$ is in the same Harish-Chandra cell as $\Psi$;
\item
for each left cell $\CCL\in \cBL$, every left cell in $T^L(\CCL)$  is contained in the same double cell as $\CCL$;
\item for each basal element  $\Psi\in \cBLG$,  
\begin{equation} \label{eq:TTL2}
  \set{\CCL_{\Psi'} : \Psi'\in T(\Psi)} = T^L(\CCL_{\Psi}),
   \end{equation}
\end{itemize}
where $\CCL_{\Psi}$ is the left cell associated to $\Psi$ (as in \Cref{annlc}). See \cite{Vtau}*{Theorem~3.2} and \cite{GV}*{Theorem~2.15 and  Corollaries ~2.16-2.18}.

For each $\CCL\in \cL_\Lambda$, put 
\[
[\CCL]_1 := \set{\CCL}, \quad [\CCL]_{k+1} := \bigcup_{\CCLp\in [\CCL]_k, \TL } 
\TL(\CCLp)
\quad \text{(for all $k\in \bN^+$), }
\]  
where $\TL$ runs over all possible choices of $\TL_{\alpha,\beta}$ and  $\TL_{\alpha_i}$ as above.
Let 
\[[\CCL] := \bigcup_{k=1}^\infty [\CCL]_k.\]

\begin{prop}\label{thm:LC.ABC}
Suppose that  the Coxeter group $W(\Lambda)$ has no simple factor of type $F_4$, $E_6$, $E_7$, or $E_8$. 
Then for each $\CCL\in \cL_\Lambda$, $[\CCL]$ is the set of all left cells in the double cell containing $\CCL$. 
\end{prop}
\begin{proof}
The proposition is easily reduced to the cases when the root system $\Delta(\Lambda)$ is simple. 

The proposition is a classical result of Knuth, Joseph, Vogan, and Kazhdan-Lusztig, when $\Delta(\Lambda)$ is of type $A$ (see \cite{KL}*{p.~177} and \cite{BV1}*{p.~172}), 
and a result of Devra Garfinkle Johnson when $\Delta(\Lambda)$ is of type $B$ or $C$  (see \cite{G3}*{Theorem~3.2.2}), or $D$ (\cf  \cite{Gar.D4} and \cite{GJMP}). (We thank Devra Garfinkle Johnson for confirming this result when $\Delta(\Lambda)$ is of type $D$.) Note that the root system of type $D_2$ (resp. $D_3$) is isomorphic to the root system of type $A_1\oplus A_1$ (resp. $A_3$). 

Suppose the root system $\Delta(\Lambda)$ is of type $G_2$. 
Let $\set{\alpha,\beta}$ be the set of simple roots. There are three double cells in total. Among the three double cells, two contain a single left cell in each, and the other is the union of two left cells. The $\tau$-invariant of the latter is  $\set{\set{s_\alpha},\set{s_\beta}}$. 
See \cite{Carter}*{p.~412}. Clearly it suffices to consider the double cell with two left cells, and since they are interchanged by $\TL_{\alpha,\beta}$ (or $\TL_{\beta,\alpha}$), the proposition clearly holds. % for t
\end{proof}

Now \Cref{thm:LC.ABC} implies that \Cref{wh:tau} holds via 
\eqref{eq:TTL2} as well as the two properties just above it. 
This concludes the proof of \Cref{thm:tau.ABCDG}. 

\begin{remark} \label{rem:LC.ABC}
For the root systems of type   $F_4,E_6, E_7, E_8$, one expects similar results on double cells by including some other maps defined in a similar fashion as $\TL_{\alpha,\beta}$ and  $\TL_{\alpha_i}$ (see \cite{Gar.D4}*{Section~11}). 
\end{remark}

\def\CCstar{\CC^*}
\def\Mstar{M^*}

We consider a $5$-tuple $(G,G_\C, \iota: G\rightarrow G_\C, \Lambda, \CC)$, where the conditions required of the $4$-tuple $(G,G_\C, \iota: G\rightarrow G_\C, \Lambda)$ are specified in \Cref{extcoh}, and $\CC$ is a Harish-Chandra cell in $\Coh_{\Lambda}(\CK(G))$. 

Let $(\check G,\check G_\C, \check \iota: \check G\rightarrow \check G_\C, \check \Lambda, \CCstar)$ be another $5$-tuple that has the same properties as the $5$-tuple $(G,G_\C, \iota: G\rightarrow G_\C, \Lambda, \CC)$.  
Thus $\CCstar$ is a Harish-Chandra cell in $\Coh_{\check \Lambda}(\CK(\check G))$, considered as a basal representation of $\check W(\check \Lambda)$. Here $\check W$ is the Weyl group of the complexified Lie algebra of $\check G$, and $\check W(\check \Lambda)\subseteq \check W$ is the integral Weyl group attached to $\check \Lambda$.

\begin{defn}\label{def:Vdual}
We say  that $(\check G,\check G_\C, \check \iota: \check G\rightarrow \check G_\C, \check \Lambda, \CCstar)$ is a weak Vogan dual of $(G,G_\C, \iota: G\rightarrow G_\C, \Lambda, \CC)$ if there is given an identification  $W(\Lambda)=\check W(\check \Lambda)$ of Coxeter groups such that 
    \[
    \Coh_{\check \Lambda}(\CK(\check G))(\CCstar) \cong (\Coh_{ \Lambda}(\CK( G))(\CC))\otimes \sgn
    \]
  as representations of $W(\Lambda)=\check W(\check \Lambda)$. 
 \end{defn} 
  
  \begin{lem}\label{lem:taudual}
  Suppose that $(\check G,\check G_\C, \check \iota: \check G\rightarrow \check G_\C, \check \Lambda, \CCstar)$ is a weak Vogan dual of $(G,G_\C, \iota: G\rightarrow G_\C, \Lambda, \CC)$. Then the $\tau$-invariants $\tau_{\CCstar}$ and $\tau_\CC$ are dual to each other. 
  \end{lem} 
\begin{proof}
  Write $M:=\Coh_\Lambda(\CK(G))(\CC)$, to be viewed as a representation of $W(\Lambda)$. Likewise write $\Mstar:=\Coh_{\check \Lambda}(\CK(\check G))(\CCstar)$, to be viewed as a representation of $\check W(\check \Lambda)$.
  Then by \cite{FJMN}*{Theorem~2.10}, 
  \[
  \tau_M=\tau_\CC,\quad\quad \tau_{\Mstar}=\tau_{\CCstar},
  \]
 and $\tau_M$ is dual to  $\tau_{\Mstar}$. Therefore the lemma follows.  
\end{proof}

\begin{whyp}\label{wh:dual}
The $5$-tuple $(G,G_\C, \iota: G\rightarrow G_\C, \Lambda, \CC)$ has a weak Vogan dual.
\end{whyp}

\begin{thm}\label{DualHC}
 Suppose that  the Coxeter group $W(\Lambda)$ has no simple factor of type $F_4$, $E_6$, $E_7$, or $E_8$. Assume that \Cref{wh:dual} holds,  then 
for every $\sigma\in \Irr(W(\Lambda))$ occurring in  $\Coh_{\Lambda}( \CK(G))(\CC)$,
 \[
    \sigma \approx \sigma_{\CC}.  
 \]
%  
%   then \[
%  \{ \sigma\in \Irr(W(\Lambda))\mid \sigma\textrm{ occurs in the cell representation $\Coh_{\Lambda}( \CK(G))(\CC)$}\}
% \]
% is contained in the  double cell in $\Irr(W(\Lambda))$ containing $\sigma_\CC$.
\end{thm}

\begin{proof}
Let  $(\check G,\check G_\C, \check \iota: \check G\rightarrow \check G_\C, \check \Lambda, \CCstar)$ be a weak Vogan dual of $(G,G_\C, \iota: G\rightarrow G_\C, \Lambda, \CC)$. 
Similar to $\sigma_\CC$, we have a special representation $\sigma_{\CCstar}\in \Irrsp{(\check W(\check \Lambda))}$ associated to the Harish-Chandra cell $\CCstar$.

 Suppose that  $\sigma\in \Irr(W(\Lambda))$ occurs in the cell representation $\Coh_{\Lambda}( \CK(G))(\CC)$. Then \Cref{lem0033} implies that 
 \begin{equation}\label{eq:ineq11}
 \sigma_\CC\leqLR  \sigma.
 \end{equation}
Since $\sigma\otimes \sgn$ occurs in the cell representation  $\Coh_{\check \Lambda}(\CK(\check G))(\CCstar)$, \Cref{lem0033} also implies that 
\[
\sigma_{\CCstar}\leqLR  \sigma\otimes\sgn.
\]
\Cref{lrd} then implies that
\begin{equation}\label{eq:ineq22}
\sigma\leqLR  \sigma_{\CCstar} \otimes\sgn.
 \end{equation}

\def\ckcCLR{{\check \CCD}}
\def\CCDstar{\CCD^*}
Let $\CCLR$ (resp. $\CCD^*$) be the double cell attached to $\CC$ (resp. $\CCstar$).  
%Recall that $\CCLR$ is the double cell attached to $\CC$. Likewise let $\CCD^*$ be the double cell attached to $\CCstar$.  
 \Cref{thm:tau.ABCDG} implies that
 \[
 \tau_{\CC} = \tau_{\CCLR}
  \qquad \textrm{and}\qquad 
 \tau_{\CCstar} = \tau_{\CCDstar}. 
 \]
 
 Since $\tau_{\CCstar}$ is dual to $\tau_{\CC}$ by \Cref{lem:taudual},  
  we know that $ \tau_{\CCDstar}$ is dual to $ \tau_{\CCLR}$.
 \Cref{duald} then implies that $\CCDstar = \check\CCD$ and so
 %By %definition of $\sigma_{\CC} = \sigma_{\CCD}$ and $\sigma_{\CCstar} = \sigma_{\CCDstar}$.   By 
 %\Cref{duald}, we have  
 \[
   \sigma_{\CCstar}\otimes \sgn 
 = \sigma_{\CCDstar}\otimes \sgn 
 = \sigma_{\check\CCD} \otimes \sgn   \approxLR \sigma_{\CCD} = \sigma_\CC. 
 \]
Together with \eqref{eq:ineq11} and \eqref{eq:ineq22}, we finally conclude that $\sigma \approxLR \sigma_\CC$ and the theorem is proved. 
\end{proof}

In the above proof we have used the fact that $\sigma_{\check\CCD}\otimes \sgn \approxLR \sigma_{\CCD}$.
We remark that  $\sigma_{\check\CCD}\otimes \sgn = \sigma_{\CCD}$ in most cases including the cases when $W(\Lambda)$ is of classical type (see \cite{Lu}*{Chapter~4}).  

 Since  weak Vogan duals exist for linear groups (\cite{V4}*{Corollary~14.9}) and  the real metaplectic group (\cite{RT1}, \cite{RT2}*{Theorem~5.2}), we obtain the following corollary.

\begin{cor}\label{HCLU}
 Suppose that the Coxeter group $W(\Lambda)$ has no simple factor of type $F_4$, $E_6$, $E_7$, or $E_8$, and $G$ is linear or isomorphic to a real metaplectic group.   
 Then 
for every $\sigma\in \Irr(W(\Lambda))$ occurring in  $\Coh_{\Lambda}( \CK(G))(\CC)$,
 \[
    \sigma \approx \sigma_{\CC}.  
 \]
%  Then   
%  \[
%  \{ \sigma\in \Irr(W(\Lambda))\mid \sigma\textrm{ occurs in the cell representation $\Coh_{\Lambda}( \CK(G))(\CC)$}\}
% \]
% is contained in the  double cell in $\Irr(W(\Lambda))$ containing $\sigma_\CC$.
\end{cor}

\section{Counting irreducible representations}\label{sec:counting}

The purpose of this section is to prove two counting results on irreducible representations, valid for any $G$ satisfying \Cref{conjcell}.
Special versions of the counting results are stated as \Cref{Mcounteq} and \Cref{Mcountleq} in \Cref{subsec:countIrr}, in view of \Cref{HCLU}.

For a basal element $\Psi\in \Coh_\Lambda(\CK(G))$, let \[\cO_\Psi:=\cO_{\sigma_{\Psi}}\in \overline \Nil(\g^*)\] be the nilpotent orbit corresponding to $j_{W(\Lambda)}^{W}\sigma_\Psi$ (as in \eqref{eq:OrbitSp}). Here $\sigma _{\Psi}$ is the special irreducible representation of $W(\Lambda)$ attached to the Harish-Chandra cell containing $\Psi$, as in \Cref{defn:spe-HC}. By \Cref{assv} and the translation principle~\cite{V1}*{Lemma~2.7}, 
 \[
 \AV(\Ann(\Psi(\nu))) = \overline{\CO_{\Psi}}\qquad (\textrm{the Zariski closure of $\CO_\Psi$})
 \]
 for all dominant $\nu\in \Lambda$ such that   $\Psi(\nu)\neq 0$.
Using the properties of translation outside the dominant cone as in \cite{Vg}*{Proposition~7.3.11}, we conclude that 
\[
\Psi\in\Coh_\Lambda( \cK_{\overline{\CO_{\Psi}}}(G)).
\]

Recall that $\sfS$ is a fixed $\Inn(\g)$-stable Zariski closed subset of $\mathrm{Nil}(\g^*)$.  
Recall also that $\Irr_{\nu,\sfS}(G)$ is the set of isomorphism classes of irreducible representations in $\Rep_{\nu, \sfS}(G)$, where
$\nu\in \hha^*$.  

\begin{prop}[Vogan]
\label{prop:ev00000}
\label{count1}
  For every $\nu\in \Lambda$, the evaluation map
  \[
   \mathrm{ev}_\nu:  \Coh_{\Lambda}(\CK_\sfS(G))\rightarrow \CK_{\nu, \sfS}(G), \quad \Psi\mapsto \Psi(\nu)
  \]
  descends to a linear isomorphism
  \[
     \Coh_{\Lambda}(\CK_\sfS(G))_{W_{\nu}} \xrightarrow{\sim}\CK_{\nu, \sfS}(G), 
     \]
     where the subscript group indicates the coinvariant space. Consequently, 
  \[\sharp(\Irr_{\nu,\sfS}(G)) = [1_{W_{\nu}}:\Coh_{\Lambda}(\CK_\sfS(G))].
  \]
  \end{prop}
\begin{proof}
  \def\BS{\cB_{\cS}} Without loss of generality we assume that $\nu$ is
  dominant. Put
  \[
\CB_{\sfS}:=\set{\Psi | \text{$\Psi$ is a basal element of  $\Coh_{\Lambda}(\CK(G))$ such that  $\cO_\Psi\subseteq \sfS$}},
  \]
  which is a basis of  $\Coh_{\Lambda}(\CK_\sfS(G))$ by \Cref{annlc}. 
  
It follows that 
  \begin{eqnarray*}
    \ker (\ev{\nu})& = & \Span \Set{\Psi\in \CB_{\sfS}\mid  \Psi(\nu)= 0 }\\ 
    %\quad (\textrm{by \Cref{lem21}}) \\
    &=& \Span \Set{\Psi - s\cdot\Psi \mid          \Psi\in  \CB_{\sfS},  s  \text{ is a simple reflection in $W_\nu$}} \\
    & &\quad (\textrm{by \Cref{lemirr}})\\
    & \subseteq &\Span\Set{\Psi- w\cdot \Psi | \Psi\in  \Coh_{\Lambda}(\CK_{\sfS}(G)), \, w\in W_{\nu}} \\
    &  \subseteq &  \ker (\ev{\nu}). 
  \end{eqnarray*}
Therefore
  \[
    \ker (\ev{\nu})=\Span\Set{\Psi- w \cdot\Psi \mid \Phi\in \Coh_{\Lambda}(\CK_{\sfS}(G)), \, w\in W_\nu}.
  \]
Together with \Cref{lem21}, this implies the first part. The second part is immediate since $ \sharp(\Irr_{\nu,\sfS}(G))=\dim \CK_{\nu, \sfS}(G)$.  
\end{proof}

\begin{prop}\label{count2}
For all $\sigma\in \Irr(W(\Lambda))\setminus \Irr_\sfS(W(\Lambda))$,  %$\CO_\sigma\nsubset \sfS$.
  \[
    [\sigma:\Coh_{\Lambda}(\CK_\sfS(G))]=0.
  \]
\end{prop}

\begin{proof}
This is implied by \Cref{lem0033}
%\Cref{lem0022} 
and \Cref{cohbbs}.
\end{proof}

\Cref{count1} and \Cref{count2} imply that  for every $\nu\in \Lambda$,
\begin{equation}\label{leq002}
  \sharp(\Irr_{\nu,\sfS}(G)) = \sum_{\sigma \in \Irr_\sfS(W(\Lambda))} [1_{W_{\nu}}: \sigma]\cdot [\sigma:\Coh_{\Lambda}(\CK_\sfS(G))].
  \end{equation}
  
 \begin{thm}\label{counteq}   
  For all $\nu\in \Lambda$, 
  \begin{equation*}%\label{leq2}  
  \sharp(\Irr_{\nu,\sfS}(G)) \leq  \sum_{\sigma \in \Irr_\sfS(W(\Lambda))} [1_{W_{\nu}}: \sigma]\cdot [\sigma:\Coh_{\Lambda}(\CK(G))], 
\end{equation*}
  and the equality holds if \Cref{conjcell} holds for $G$.   
    \end{thm}

\begin{proof} The first assertion is immediate from \eqref{leq002}. 
As a  representation of $W(\Lambda)$, 
\[
  \Coh_{\Lambda}(\CK(G))\\
   \cong \bigoplus_{\CC \textrm{ is a cell in $\Coh_{\Lambda}(\CK(G))$}} \Coh_{\Lambda}(\CK(G))(\CC). 
\]
\Cref{hcass22} implies that 
 \[
   \Coh_{\Lambda}(\CK_\sfS(G))\cong  \bigoplus_{\CC \textrm{ is a cell in $\Coh_{\Lambda}(\CK(G))$, $\sigma_\CC\in \Irr^{\mathrm{sp}}_\sfS(W(\Lambda))$}} \Coh_{\Lambda}(\CK(G))(\CC),
 \]
 where $\sigma_\CC$ is the special irreducible representation attached to $\CC$, as in \Cref{defn:spe-HC}. Under the assumption that \Cref{conjcell} holds for $G$,  
\begin{equation*}%\label{eqsigma}
[ \sigma: \Coh_{\Lambda}(\CK_\sfS(G))]=[\sigma:\Coh_{\Lambda}(\CK(G))]\quad \textrm{ for all  $\sigma\in \Irr_\sfS(W(\Lambda))$.}
\end{equation*}
Together with \eqref{leq002}, this implies the second assertion. 
\end{proof}

Recall that $\CO_\nu $ is the nilpotent orbit 
 whose Zariski closure $\overline{\CO_{\nu}}\subset \g^*$ equals 
 the associated variety of the maximal ideal $I_\nu$. 
 Recall also the Lusztig left cell
$\LC_{\nu}$ from \eqref{ll}.

 \begin{cor}\label{countleq}  
   For  all $\nu\in \Lambda$,
  \begin{equation*}%\label{boundc}
    \sharp(\Irr_{\nu,\overline{\CO_{\nu}}}(G)) \leq \sum_{\sigma \in \LC_{\nu}}   [\sigma:\Coh_{\Lambda}(\CK(G))],
  \end{equation*}
     and the equality holds if  \Cref{conjcell} holds for $G$.
    \end{cor}
  \begin{proof}
 In view of \eqref{mulone} and \Cref{leftcnu}, this follows from \Cref{counteq}.  \end{proof}

We prove the following proposition for later use. 

\begin{prop}\label{hcass222}
Let $\Psi$ be a basal element of $\Coh_\Lambda(\CK(G))$, and let  $\nu\in \Lambda$ be a dominant element. Then $\Psi(\nu)\in \Irr_{\nu, \overline{\CO_\nu}}(G)$ if and only if $\sigma_\Psi\cong \left( j_{W_\nu}^{W(\Lambda)} \sgn\right)\otimes \sgn $ and no element of   $ \tau_\Psi$  fixes $\nu$. 
\end{prop}

\begin{proof}
First suppose that $\sigma_\Psi\cong \left( j_{W_\nu}^{W(\Lambda)} \sgn\right)\otimes \sgn$ and there is no element of   $ \tau_\Psi$  that fixes $\nu$.
Then \Cref{lemirr} implies that $\Psi(\nu)\in \Irr(G)$ and \Cref{hcass22}  implies that
\[
\sigma_{\nu, \mathrm{Ann}(\Psi(\nu))}\cong  \left( j_{W_\nu}^{W(\Lambda)} \sgn\right)\otimes \sgn.
\]
Thus $ \mathrm{Ann}(\Psi(\nu))=I_\nu$ by \Cref{leftcnu2}, and hence $ \Psi(\nu)\in \Irr_{\nu, \overline{\CO_\nu}}(G)$.

Now suppose that $\Psi(\nu)\in \Irr_{\nu, \overline{\CO_\nu}}(G)$. \Cref{lemirr} implies that no element of $ \tau_\Psi$ fixes $\nu$.
As $ \mathrm{Ann}(\Psi(\nu))=I_\nu$, \Cref{leftcnu2} and  \Cref{hcass22} imply that $\sigma_{\Psi}\cong \left( j_{W_\nu}^{W(\Lambda)} \sgn\right)\otimes \sgn$. This completes the proof of the proposition.
\end{proof}

\section{Separating good parity and bad parity for coherent continuation  representations}\label{sec:GB}

From now on, $\star$ will be one of the 10 labels, and $G$ will be a classical group of type $\star$, as in Sections \ref{sec:defunip}-\ref{secorgp0}. We have the complex Lie group 
\be\label{gc}
  G_\C :=
  \begin{cases}
   \GL_{n}(\C), & \textrm{if $\star\in \{A^\R, A^\bH\}$};\\
     \GL_{p+q}(\C), & \textrm{if $\star\in \{A, \widetilde A\}$};\\
    \SO_{p+q}(\C), & \textrm{if $\star\in \set{B,D}$};\\
     \SO_{2n}(\C), &\textrm{if $\star = D^{*}$};\\
    \Sp_{2n}(\C), &\textrm{if $\star \in \{C, \wtC\}$};\\
    \Sp_{p+q}(\C), &\textrm{if $\star = C^{*}$},\\
  \end{cases}
\ee
%If $\star=\wtC$, then $G_\C=\check G$, otherwise $G_\C$ and $\check G$
and $\iota: G\rightarrow G_\C$ is the usual complexification homomorphism except when $\star\in \{\widetilde A, \widetilde C\}$ where it factors through the cover homomorphism
(see \Cref {sec:CW} for the general setup). 

The main purpose of this section is to separate good parity and bad parity for coherent continuation representations. 
(The broad ideas are from \cite{ABV}. We give explicit details in order to bring out the key role of standard representations. A reader who is familiar with \cite{ABV} may wish to skip this section.)
We will first describe structure data of the classical group $G$ (such as the universal Cartan subalgebra and the analytic weight lattice) in terms of the standard representation of $G$. In the Langlands dual Lie algebra $\check \g$, we define the good parity and bad parity parts, also in terms of the standard representation of $\check \g$. We then proceed to separate good parity and bad parity in all structure data necessary for the description of coherent continuation representations. This is not difficult, but requires care for the sake of clarity.

We will freely use  the notation of Sections \ref{sec:defunip}-\ref{secorgp0}. One additional convention is as follows.

If $\star\in \{\widetilde A, \widetilde C\}$, $\CK'(G)$ will denote the Grothendieck group of the category of genuine Casselman-Wallach representations of $G$. Otherwise, put $\CK'(G):=\CK(G)$. Similar notations such as $\cP'_\Lambda(G)$, $\CK'_{\nu}(G)$ and $\Irr'_{\nu}(G)$
%$ \CK'_{\sfS}(G)$, $\CK'_{\nu,\sfS}(G)$
will be used without further explanation.

\delete{As an obvious variant of \Cref{counteq} and \Cref{countleq} for $\nu=\lambda_{\check \CO}$, we have  that
 \begin{equation}\label{boundc22}
     \sharp(\Unip_{\check \CO}(G)) =\sum_{\sigma\in \LC_{\lambda_{\check \CO}}} [\sigma: \Coh_{\Lambda}(\CK_{\overline{\CO}}'(G))]  = \sum_{\sigma\in \LC_{\lambda_{\check \CO}}} [\sigma: \Coh_{\Lambda}(\CK'(G))]. 
   \end{equation}
   (The second equality holds because of \Cref{HCLU}.)
Recall that $\CO:=\dBV(\check \CO)$ and $\overline{\CO}$ denotes its Zariski closure.}

\subsection{Standard representations of  classical groups}\label{subsec:standard}

For (mostly) efficiency reasons, we would like to lump certain cases (such as $B$ and $D$) together when discussing standard representations. Partition the set of
  labels into subsets 
  \[
\{\{A^\R\}, \{A^\bH\}, \{A, \wtA\}, \{B,D\}, \{C, \wtC\}, \{C^*\},\{D^*\} \},
\]
and denote by  $[\star]$ the equivalence class containing $\star$.

Given a label $\star$, we define the notion of a $[\star]$-structure, which consists of a finite dimensional complex vector space $V$, together with some additional data. 
The group
  $\mathsf G(V)$ will be defined as the (real) subgroup of $\GL(V)$ 
  fixing the
  $[\star]-$structure. Its Zariski closure in $\GL(V)$ will be denoted by $\mathsf
  G_\bC(V).$ For cases $\widetilde{A}$ and
$\widetilde{C}$ there will be a cover, and our classical group (specified by $\star$ and $V$) will be  denoted by $\mathsf G_\star(V)$.

\medskip
\noindent {\bf The case when $\star\in \set{A^\R, A^\BH}$.} In this case, a $[\star]$-structure consists of a finite dimensional complex vector space $V$, and a conjugate linear automorphism $\mathbf j: V\rightarrow V$ such
that
\[
  \mathbf j^2= \begin{cases}
  1,  &  \text{if $\star=A^\R$};\\
  -1,    &  \text{if $\star=A^\bH$}.
    \end{cases}
\]
The group $\mathsf G(V)$ is $\GL_n(\bR)$ when $\mathbf j^2=1$, and $\GL_\frac{n}{2}(\bH)$  when $\mathbf j^2=-1$. Here $n=\dim V$. 

\medskip
\noindent {\bf The case when $\star\in \set{A, \wtA}$.} In this case, a $[\star]$-structure consists of a finite dimensional complex vector space $V$, and a non-degenerate Hermitian form
$\la\,,\,\ra: V\times V\rightarrow \C$ (which is by convention linear on the first variable and conjugate linear on the second variable). The Hermitian form has signature $(p,q)$ so $\mathsf G(V)$
is $\oU(p,q)$ with $p+q=\dim V$.

\medskip
\noindent {\bf The case when $\star\in \set{C, \wtC, C^*}$.} In this case, a $[\star]$-structure consists of a finite dimensional complex vector space $V$, a symplectic form
$\la\,,\,\ra: V\times V\rightarrow \C$ together with a conjugate linear automorphism $\mathbf j: V\rightarrow V$ such
that
\[
  \mathbf j^2= \begin{cases}
  1,  &  \text{if $\star\in\{C, \wtC\}$};\\
  -1,    &  \text{if $\star=C^*$},
    \end{cases}
\]
 and for all $u,v\in V$,
 \[
   \la \mathbf j(u), \mathbf j(v)\ra =\overline{\la u, v\ra}, \quad \text{(bar indicates the complex conjugate)}. 
 \]
When $\mathbf j^2=1,$ the group $\mathsf G(V)$ is $\Sp_{2n}(\bR)$ ($n=\frac{\dim V}{2})$. When
  $\mathbf j^2=-1,$ there are several $\mathsf G_\bC(V)$ conjugacy classes of $\mathbf j$'s and $\mathsf
  G(V)$ is $\Sp(\frac{p}{2},\frac{q}{2})$ with $p,q$ even and $p+q=\dim V$. 
  
\medskip
\noindent {\bf The case when $\star\in \set{B,D, D^*}$.} In this case, a $[\star]$-structure consists of a finite dimensional complex vector space $V$,  %$(\la\,,\,\ra, \mathbf j, \omega)$ where
a non-degenerate symmetric bilinear form $ \la\,,\,\ra: V\times V\rightarrow \C $, a conjugate linear automorphism $\mathbf j: V\rightarrow V$, and $\omega\in \wedge^{\dim V} V$ (to be called the orientation), subject to the following five  conditions: 
 \begin{itemize}
 \item
 \be\label{condition1}
  \mathbf j^2= \begin{cases}
  1,  &  \text{if $\star\in\{B,D\}$};\\
  -1,    &  \text{if $\star=D^*$}.
    \end{cases}
\ee
\item For all $u,v\in V$, 
 \be\label{condition2}
   \la \mathbf j(u), \mathbf j(v)\ra =\overline{\la u, v\ra}.
 \ee
The subgroup of $\GL(V)$ stabilizing $(\la\,,\,\ra, \mathbf j)$, to be denoted by $\mathsf G(\la\,,\,\ra, \mathbf j)$, is a real orthogonal group $\oO (p,q)$ (with $p+q=\dim V$) when $\star\in \{B,D\}$, and the quaternionic orthogonal group $\SO^*(2n)$ (with $n=\frac{\dim V}{2}$) when $\star=D^*$.
The data $\omega$ satisfies the additional conditions below. 
\item With respect to the symmetric bilinear form $\la\,,\,\ra: \wedge^{\dim V} V\times \wedge^{\dim V} V\rightarrow \C$ given by
\[
  \la u_1\wedge u_2\wedge \dots \wedge u_{\dim V}, v_1\wedge v_2\wedge \dots \wedge v_{\dim V} \ra=\det \left( [\la u_i, v_j \ra]_{i,j=1,2, \dots, \dim V}\right), 
\]
we have 
\be\label{laraw}
\la \omega, \omega\ra=1.\ee
 \item
if $\star=D^*$, then
\be\label{dstaro}
 \omega=\sqrt{-1} u_1\wedge \mathbf j(u_1)\wedge \sqrt{-1} u_2\wedge \mathbf j(u_2)\wedge \dots \wedge \sqrt{-1} u_{\frac{\dim V}{2}} \wedge \mathbf j(u_{\frac{\dim V}{2}}),
 \ee
 for some  vectors $u_1, u_2, \dots, u_{\frac{\dim V}{2}}\in V$;
 \item
If $\star\in \{B,D\}$ and $\dim V=0$, then 
\be\label{BDo}
\omega=1 \qquad (\text{as an element of $\wedge^{\dim V} V=\C$}).
\ee 
The group $\mathsf G(V)$ is a real special orthogonal group or a quaternionic orthogonal group fixing a preferred orientation. 
\end{itemize}

\begin{remarks}
    (a) In the case of $\star\in \{C,\wtC\}$, the existence of the symplectic form implies that $\dim V$ is even. In the case of $\star\in \{A^\BH, C^*,D^*\}$, the existence of $\mathbf j$ also implies that $\dim V$ is even.

    (b) We comment on the three conditions on the orientation $\omega$ in the case of $\star\in \set{B,D, D^*}$. Note that there are precisely two elements $\omega\in \wedge^{\dim V} V$ (negative of each other) satisfying $\la \omega, \omega\ra=1$, the condition in \eqref{laraw}. When $\star\in \{B,D\}$ and $\dim V\neq 0$, these two elements are conjugate under $\mathsf G(\la\,,\,\ra, \mathbf j)$. However, when $\star\in \{B,D\}$ and $\dim V=0$, or $\star=D^*$, the group $\mathsf G(\la\,,\,\ra, \mathbf j)$ is connected and acts trivially on $\wedge^{\dim V} V$, and the aforementioned  two elements are not conjugate under $\mathsf G(\la\,,\,\ra, \mathbf j)$. The two conditions \eqref{dstaro} and \eqref{BDo} thus ensure that the orientation $\omega$, as part of the $[\star]$-structure on $V$, is unique up to conjugation by $\mathsf G(\la\,,\,\ra, \mathbf j)$ in all cases.  
    
    (c) When the group $\mathsf G(V)$ is a real special orthogonal group of type $D$, or a quaternionic orthogonal group, we will need to use the orientation $\omega$ to identify the universal Cartan subalgebra with $\C^n$, where $n$ is the rank of the group. See \eqref{idabs}. 
\end{remarks}

For a $[\star]$-structure, 
by abuse of terminology we will call the underlying vector space $V$ a $[\star]$-space. We will also write 
  $\la\,,\,\ra_V:=\la\,,\,\ra$, $\mathbf j_V:=\mathbf j$, and
  $\omega_V:=\omega$ (when the additional structures are present). Given two $[\star]$-spaces $V_1$ and $V_2$, a linear isomorphism $V_1\rightarrow V_2$ is said to be an isomorphism of $[\star]$-spaces if it preserves the $[\star]$-structure. 
 The product $V_1\times V_2$ is a $[\star]$-space in the obvious way. 
 \delete{We remark that when $[\star]=\{B,D\}$ and both $\dim V_1$ and $\dim V_2$ are odd, the switching map $V_1\times V_2\rightarrow V_2\times V_1$ does not preserve the orientation.} 
% and  $V_1\times V_2$ 
% and $V_2\times V_1$ are not canonically isomorphic to each other as $[\star]$-spaces (although they are isomorphic to each other).

For a $[\star]$-space $V$, put
\[
  \mathsf G_\star(V):=\begin{cases}
   \textrm{the $\det^{\frac{1}{2}}$-double cover of $\mathsf G(V)$},  &  \text{if $\star=\wtA$};\\
 \textrm{the metaplectic double cover of $\mathsf G(V)$},  &  \text{if $\star=\wtC$};\\
  \mathsf G(V),    &  \text{otherwise.}
    \end{cases}
\]
The group $\mathsf G_\star(V)$ only depends on the isomorphism class of $V$ as a $[\star]$-space.

Now we assume that $G$ is identified with $\mathsf G_\star(V)$. We call $V$ the standard representation of $G$. 
Recall than $n$ is the rank of $\g$ in all cases. We fix a  flag
\be\label{flaga}
\{ 0\}= V_0 \subset V_1\subset \dots \subset  V_{n}
\ee
in $V$ such that
\begin{itemize}
 \item
$\dim V_i=i$ for all $i=1,2, \dots, n$;

\item
if $\star\in \{B,D,C,\wtC, C^*, D^*\}$, then $
  V_{n}
$ is totally isotropic (with respect to the bilinear form $\la\,,\,\ra_V: V\times V\rightarrow \C$);

\item
if $\star\in \{D, D^*\}$, then $V_n$ is $\omega_V$-compatible    in the sense that
\[
  \omega_V= u_1\wedge u_2\wedge \dots \wedge u_{2n}%\wedge u_n^* \wedge \dots \wedge u_{2}\wedge u_1
\]
for some elements $u_1, u_2, \dots, u_{2n}\in V_{2n}$ such that $u_1, u_2, \cdots, u_n\in V_n$ and
\[
\la u_i, u_{j}\ra=\begin{cases}
   1, &\quad \textrm{if $i+j=2n+1$};\\
   0,&  \quad \textrm{if $i+j\neq 2n+1$},\\
  \end{cases}
\]
for all $i,j=1,2, \cdots, 2n$.
\end{itemize}

\begin{remarks}\label{rm63}
    (a) When $\star\in \{D,D^*\}$ and $\dim V\neq 0$ so that $G_\C$ is a non-trivial complex even special orthogonal group, up to conjugation by $G_\C$ there are precisely two $n$-dimensional  totally isotropic subspaces of $V$. Among these two subspaces, exactly one is $\omega_V$-compatible.    
    
(b) In all cases, up to conjugation by $G_\C$ there exists a unique flag as in \eqref{flaga} satisfying  the aforementioned three conditions.  The stabilizer of the flag  \eqref{flaga} in $G_\C$ is a Borel subgroup of $G_\C$, and every Borel subgroup of $G_\C$ is uniquely of this form.  

(c) When $\star=D^*$  and $n$ is even,  the condition \eqref{dstaro} ensures that all 
$\mathbf j_V$-stable totally isotropic subspace of $V$ of dimension $n$ are $\omega_V$-compatible. 
\delete{Consequently, when $\star=D^*$  and $n$ is arbitrary, the orbit $\check \CO$ is not $G$-relevant (as defined in Section \ref{secrgp0}) if and only if the following conditions are satisfied:
\begin{itemize}
    \item $n$ is positive and even;
    \item $\check \CO$ has bad parity;
    \item an element of $\CO$ stabilizes an  $n$-dimensional totally isotropic subspace of $V$ that is not  $\omega_V$-compatible.
\end{itemize}
Recall that $\CO$ is the Barbasch-Vogan dual of $\check \CO$.}
    \end{remarks}
    
The stabilizer of the flag  \eqref{flaga} in $\g$ is a Borel subalgebra of $\g$. Using this Borel subalgebra, we get an identification
\be\label{idabs}
   \hha= \prod_{i=1}^{n} \g\l(V_i/  V_{i-1})=\C^n.
\ee
This identification is independent of the choice of the flag \eqref{flaga}.
As in (the beginning of) \Cref{sec:CW}, we assume that $Q=Q_\iota$ (the analytic weight lattice). Note that  
\be\label{weightl}
  Q=\Z^n\subset \C^n=(\C^n)^*=\hha^*,
\ee
and the positive roots are
\[
 \Delta^+= \begin{cases}
    \{e_i-e_j\mid 1\leq i<j\leq n\}, &  \text{if $\star \in \{A^\R, A^\BH,  A, \wtA\}$}; \\
   \{e_i\pm e_j\mid 1\leq i<j\leq n\}, &  \text{if $\star \in \{D, D^*\}$}; \\
   \{e_i\pm e_j\mid 1\leq i<j\leq n\} \cup \{e_i\mid 1\leq i\leq n\}, &  \text{if $\star =B$}; \\
   \{e_i\pm e_j\mid 1\leq i<j\leq n\} \cup \{2 e_i\mid 1\leq i\leq n\}, &  \text{if $\star\in \{C, \wtC, C^*\}$}. \\
  \end{cases}
\]
Here $e_1, e_2, \dots, e_n$ is the standard basis of $\C^n$. The Weyl group
\begin{equation}\label{def:Weyl}
    W= \begin{cases}
    \sfS_n, &  \text{if $\star \in \{A^\R, A^\BH,  A, \wtA\}$}; \\
    \sfW_n, &  \text{if $\star\in \{B, C, \wtC, C^*\}$}; \\
      \sfW_n' , &  \text{if $\star \in \{D, D^*\}$}, \\
  \end{cases}
\end{equation}
where $\sfS_n\subset \GL_n(\Z)$ is the group of the permutation matrices, $\sfW_n\subset \GL_n(\Z)$ is the subgroup generated by $\sfS_n$ and all the diagonal matrices with diagonal entries $\pm 1$, and $\sfW'_n\subset \GL_n(\Z)$ is the subgroup generated by $\sfS_n$ and all the diagonal matrices with  diagonal entries $\pm 1$ and determinant $1$.

\subsection{Good parity and bad parity in the Langlands dual}

Recall that an integer has  good parity (which depends on $\star$ and $n=\rank \check \g$) if it has the same parity as
\[
  \begin{cases}
    n, &  \text{if $\star \in \{A^\R, A^\BH,  A\}$}; \\
    1+ n, &  \text{if $\star = \wtA$}; \\
   1, & \text{if } \star \in \set{C,C^{*},D,D^{*}};\\
 \text{0}, & \text{if } \star \in \set{B,\wtC}.\\
  \end{cases}
\]
Otherwise it has bad parity.

Let $\check V$ denote the standard representation of $\check \g$ so that $\check \g$ is identified with a Lie subalgebra of $\g\l(\check V)$.
When $\star\in \{A^\R, A^\bH, A, \wtA\}$,  $\check \g=\g\l(\check V)$. When $\star\in \{B, \wtC\}$,
the space $\check V$ is equipped with a $\check \g$-invariant symplectic form $\la\,,\,\ra_{\check V}: \check V\times \check V\rightarrow \C$.  When $\star\in \{C, C^*, D, D^*\}$, the space $\check V$ is equipped with a $\check \g$-invariant non-degenerate symmetric bilinear form $\la\,,\,\ra_{\check V}: \check V\times \check V\rightarrow \C$ together with an element $\omega_{\check V}\in \wedge^{\dim \check V} \check V$ with $\la \omega_{\check V}, \omega_{\check V}\ra_{\check V}=1$ (see \eqref{laraw}).

As in \eqref{idabs} we also have that $\check \hha=\C^{n}$. Thus   we have an identification
  $ \check \hha=\hha^*$
that identifies $\check \g$ as the Langlands dual (or the metaplectic Langlands dual) of $\g$.

Recall the semisimple element
$\lambda^\circ_{\check \CO}\in\check \g$ that equals half of the
neutral element in an $\s\l_2$-triple attached to
$\check \CO$. Its conjugacy class  is uniquely determined by
    $\check \CO$ under $\Ad(\check \g)$. Fix a Lie algebra homomorphism
\begin{equation}\label{sl2}
  \phi : \s\l_2(\C)\rightarrow \check \g
\end{equation}
such that
\[
\phi \left(\left[
\begin{array}{cc}
0&1\\
0&0
\end{array}
\right]\right)\in \check \CO\quad \textrm{and}\quad
\phi \left(\left[
\begin{array}{cc}
1&0\\
0&-1
\end{array}
\right]\right)=2\lambda_{\check \CO}^\circ.
\]

We view $\check V$ as
an $\s\l_2(\C)$-module via
$\phi$.
Put
\[
 \check V_{\mathrm b}:=\textrm{sum of  irreducible $\s\l_2(\C)$-submodules of $\check V$ with bad parity dimension.}
 \]
 Define $ \check V_{\mathrm g}$ similarly by adding the $\s\l_2(\C)$-submodules of $\check V$ with good parity
   dimension. Then 
 \be\label{dgb}
 \check V= \check V_{\mathrm b}\times \check V_{\mathrm g}.
 \ee

 Write $\check \g_\mathrm b$ for the Lie subalgebra of $\check \g$ consisting of the elements that annihilate $\check V_{\mathrm g}$ and stabilize $\check V_{\mathrm b}$.
 Define $\check \g_\mathrm g$ similarly so that $\check \g_\mathrm b\times \check \g_\mathrm g$ is the stabilizer of the decomposition \eqref{dgb} in $\check \g$.
  We have natural isomorphisms
\[
  (\check \g_\mathrm b,\check \g_\mathrm g)\cong
  \begin{cases}
   (\g\l_{\nnb}(\bC),\g\l_{\nng}(\C)), &\quad  \text{if } \star \in \set{A^\R, A^\BH, A, \wtA}; \\
    (\s\p_{2\nnb}(\bC),\s\p_{2\nng}(\C)), &\quad  \text{if } \star \in \set{B,\wtC}; \\
    (\o_{2\nnb}(\bC),\o_{2\nng+1}(\bC)), & \quad  \text{if } \star \in \set{C,C^{*}}; \\
    (\o_{2\nnb}(\bC),\o_{2\nng}(\bC)), &\quad \text{if } \star \in \set{D,D^{*}},\\
  \end{cases}
\]
 where  $n_{\mathrm b}$ and $n_\mathrm g$ are ranks of $\check \g_\mathrm b$ and $\check \g_\mathrm g$ respectively so that $n_{\mathrm b}+n_{\mathrm g}=n$.
Note that when $\star\in \{B, C,\wtC, C^*, D, D^*\}$, $n_\mathrm b$ agrees with the one defined in \eqref{nb000}.

When $\star\in \{B, C,\wtC, C^*, D, D^*\}$,  \eqref{dgb} is an orthogonal decomposition, and the form $\la\,,\,\ra_{\check V}$ on $\check V$ restricts to bilinear forms on $\check V_{\mathrm b}$ and $\check V_{\mathrm g}$, to be respectively denoted by $\la\,,\,\ra_{\check V_\mathrm b}$ and $\la\,,\,\ra_{\check V_\mathrm g}$.
 When $\star\in \{C, C^*, D, D^*\}$, we fix an element $\omega_{\check V_\mathrm b}\in \wedge^{\dim \check V_\mathrm b} \check V_\mathrm b$
 and an element $\omega_{\check V_\mathrm g}\in \wedge^{\dim \check V_\mathrm g} \check V_\mathrm g$ such that
 \be\label{omega12}
   \la \omega_{\check V_\mathrm b}, \omega_{\check V_\mathrm b}\ra_{\check V_\mathrm b}=\la \omega_{\check V_\mathrm g}, \omega_{\check V_\mathrm g}\ra_{\check V_\mathrm g}=1\quad \textrm{and}\quad  \omega_{\check V_\mathrm b}\wedge  \omega_{\check V_\mathrm g}= \omega_{\check V}.
 \ee
Then  as in \eqref{idabs},   in all cases  we have identifications
   \[
   \check \hha_\mathrm b=\C^{n_\mathrm b},  \quad \check \hha_\mathrm g=\C^{n_\mathrm g}\quad  \textrm{and}\quad \check \hha= \check \hha_\mathrm b\times \check \hha_\mathrm g=\C^{n},
   \]
   where $ \check \hha_\mathrm b$ and $ \check \hha_\mathrm g$ are the universal Cartan subalgebras of $\check \g_\mathrm b$ and $\check \g_\mathrm g$, respectively.

The homomorphism $\phi$ in \eqref{sl2} has image in \eqref{dgb}, so induces Lie algebra homomorphisms
\[
  \phi _\mathrm g: \s\l_2(\C)\rightarrow \check \g_\mathrm g\quad \textrm{and}\quad  \phi_\mathrm b: \s\l_2(\C)\rightarrow \check \g_\mathrm b.
  \]
  By using these two homomorphisms, we obtain $\check \CO_\mathrm g\in \overline{\Nil}(\check \g_\mathrm g)$,  $\lambda_{\CO_\mathrm g}^\circ \in \check \g_\mathrm g$,  $\check \CO_\mathrm b\in \overline{\Nil}(\check \g_\mathrm b)$, and  $\lambda_{\CO_\mathrm b}^\circ \in \check \g_\mathrm b$
  as in \eqref{sl2}.
  Then \[
  \mathbf d_\ckcO = \mathbf d_{\ckcO_{\mathrm g}}\cuprow \mathbf d_{\ckcOb}.
  \]

 Put
  \begin{equation}\label{def:starg}
   \star_\mathrm g:=
  \begin{cases}
   A, &  \text{if $\star=\wtA$ and $p+q$ is odd}; \\
   \star, &  \text{otherwise},\\
  \end{cases}
\end{equation}
and
 \begin{equation}\label{def:starb}
   \star_\mathrm b:=
  \begin{cases}
   \wtA, &  \text{if $\star=A$ and $p+q$ is odd}; \\
     D, &  \text{if $\star=C$}; \\
       D^*, &  \text{if $\star=C^*$}. \\
   \star, &  \text{otherwise}.\\
  \end{cases}
\end{equation}
The orbit $\ckcO_{\mathrm g}$ has good parity (i.e. all its nonzero row lengths have good parity) with respect to $\star_\mathrm g$ and $n_\mathrm g$, while the orbit $\ckcOb$ has bad parity with respect to $\star_\mathrm b$ and $n_\mathrm b$.

\trivial[h]{Similar to \eqref{dominant}, we say that an element $\check \nu\in \check \hha$ is dominant if
\be\label{dominant2}
    \la \check \nu, \alpha\ra\notin -\bN^+ \qquad\textrm{for all positive root $\alpha\in \check \hha^*$ of $\check \g$}.
  \ee
The same terminology of course also applies to elements of  $ \check \hha_\mathrm b$ and  $ \check \hha_\mathrm g$.
Let $\lambda_{\check \CO_\mathrm b}\in \check \hha_\mathrm b$ be the unique  dominant element  that represents the $\Ad(\check \g_\mathrm b)$-conjugacy class of $\lambda_{\check \CO_\mathrm b}^\circ$, and likewise let $\lambda_{\check \CO_\mathrm g}\in \check \hha_\mathrm g$ be the unique  dominant element  that represents the $\Ad(\check \g_\mathrm g)$-conjugacy class of $\lambda_{\check \CO_\mathrm g}^\circ$. Put $\lambda_{\check \CO}:=(\lambda_{\check \CO_\mathrm b}, \lambda_{\check \CO_\mathrm g})\in \check \hha$, which is a dominant element that represents the $\Ad(\check \g)$-conjugacy class of $\lambda_{\check \CO}^\circ$.

Define two subsets
\[
\Lambda_\mathrm b:=\lambda_{\check \CO_\mathrm b}+\Z^{n_\mathrm b}\subset \C^{n_\mathrm b}= \check \hha_\mathrm b
\]
 and
 \[
 \Lambda_\mathrm g:=\lambda_{\check \CO_\mathrm g}+\Z^{n_\mathrm g}\subset \C^{n_\mathrm g}=\check \hha_\mathrm g.
 \]
}

Put 
\[
  \Lambda_\mathrm b:= \begin{cases}
    \Z^{n_\mathrm b}, &  \text{if $\star \in \{A^\R,  A\}$ and $n$ is even, or $\star \in \{A^\BH,  B, \wtC\}$}; \\
    (\frac{1}{2},\frac{1}{2},\dots,  \frac{1}{2})+ \Z^{n_\mathrm b}, &  \text{otherwise,}
      \end{cases}
\]
which is the $\Z^{n_\mathrm b}$-coset in $\C^{n_\mathrm b}= \check \hha_\mathrm b$ containing some (equiv. all) elements in $\check \hha_\mathrm b$ that represent the  $\Ad(\check \g_\mathrm b)$-conjugacy class of $\lambda_{\check \CO_\mathrm b}^\circ$. 
Likewise put
\[
  \Lambda_\mathrm g:= \begin{cases}
   (\frac{1}{2},\frac{1}{2},\dots,  \frac{1}{2})+ \Z^{n_\mathrm g}, &  \text{if $\star \in \{A^\R,  A\}$ and $n$ is even, or $\star \in \{A^\BH,  B, \wtC\}$}; \\
     \Z^{n_\mathrm g}, &  \text{otherwise,}
      \end{cases}
\]
which is the $\Z^{n_\mathrm g}$-coset in $\C^{n_\mathrm g}= \check \hha_\mathrm g$ containing some (equiv. all) elements in $\check \hha_\mathrm g$ that represent the $\Ad(\check \g_\mathrm g)$-conjugacy class of $\lambda_{\check \CO_\mathrm g}^\circ$. 
Finally put 
\begin{equation}\label{eq:Lambda}
  \Lambda:= \Lambda_\mathrm b\times  \Lambda_\mathrm g\subset \hha_\mathrm b^*\times  \hha_\mathrm g^*= \C^{n_\mathrm b}\times  \C^{n_\mathrm g}= \C^{n}= \check \hha=\hha^*.
\end{equation}

\subsection{Separating the universal Cartan subalgebra}\label{sec:SepCar}
When $\star\in \set{ B, C, \wtC,C^*,D, D^*}$, we have defined in Section \ref{secrgp0} the notion of an $\check \CO$-relevant parabolic subgroup of $G$ (\Cref{def:c-relevant2}). We extend this notion to the other cases: If $\star\in \{A^\R, A^\bH\}$, then a parabolic subgroup of $G$ is said to be  $\check \CO$-relevant if it is the stabilizer group of a $\mathbf j_V$-stable subspace of $V$ of dimension $n_\mathrm b$.  If $\star\in \{A,\wtA\}$, then a parabolic subgroup of $G$ is said to be  $\check \CO$-relevant if it is the stabilizer group of a totally isotropic subspace of  $V$ of dimension $\frac{n_\mathrm b}{2}$ (in particular, the existence of such a parabolic subgroup implies that $n_\mathrm b$ is even).

Recall the notions in \Cref{extcoh}, particularly
  Definition \ref{def:parameter} of the data $\upgamma=(H, \xi,
  \Gamma)$ parametrizing irreducible Casselman-Wallach representations. Recall also, \Cref{sec:defunip}, $\mathbf d_{\check \CO}$ denotes the Young diagram attached to  $\check \CO$.
  
Through a direct verification, we record the following lemma without proof.  
\begin{lem}\label{leme}
For every $\upgamma=(H, \xi, \Gamma)\in \sP'_\Lambda(G)$, $H$ is contained in a parabolic subgroup of $G$ that is ${\check \CO}'$-relevant for some ${\check \CO}'\in \overline{\Nil}(\check \g)$ with $\mathbf d_{\check \CO'}=\mathbf d_{\check \CO}$.
\end{lem}

\begin{remark} The ${\check \CO}'$ in the above lemma is nothing but ${\check \CO}$ (which we fix), except in the case when $\star\in \{D, D^*\}$, $V\neq 0$, and $\check \CO$ is very even. See \Cref{secrgp0} for more discussions on the notion of relevant parabolic subgroups. 
\end{remark}

In the rest of this paper we assume that $G$ has  a parabolic subgroup that is ${\check \CO}'$-relevant for some ${\check \CO}'\in \overline{\Nil}(\check \g)$ with $\mathbf d_{\check \CO'}=\mathbf d_{\check \CO}$.
Otherwise
\Cref{leme} and the surjectivity of the evaluation map in 
\Cref{lem21} imply that $\Unip_\ckcO(G)$ is empty.
The assumption is equivalent to saying that
\be\label{nonemp0}
\textrm{ $n_\mathrm b$ is even if $\star=\wtA$},
\ee
 and
\be\label{nonemp00}
  p,q\geq  \begin{cases}
   \frac{n_\mathrm b}{2},  &  \text{if $\star\in\{A, \wtA\}$};\\
  n_\mathrm b,    &  \text{if $\star\in\{B,C^*, D\}$}.
    \end{cases}
  \ee
See \Cref{relpa0}. In particular, $\star_\mathrm g=\star$ in all cases.
% In the rest of this paper we assume that  these two conditions are  satisfied. Consequently, $\star_\mathrm g=\star$.

Define two classical groups
\begin{equation}\label{eq:GbGg}
  (G_\mathrm b, G_{\mathrm g}) =
  \begin{cases}
   ( \GL_{n_{\mathrm b}}(\R), \GL_{n_{\mathrm g}}(\R) ), &\quad  \text{if } \star = A^\R; \\
  (\GL_{\frac{n_{\mathrm b}}{2}}(\BH), \GL_{\frac{n_{\mathrm g}}{2}}(\BH)  ), &\quad  \text{if } \star = A^\bH; \\
   (\oU(\frac{n_{\mathrm b}}{2}, \frac{n_\mathrm b}{2} ), \oU(p_{\mathrm g}, q_\mathrm g)), &\quad  \text{if $ \star = A$ and $p+q$ is even}; \\
    ( \widetilde \oU(\frac{n_{\mathrm b}}{2}, \frac{n_\mathrm b}{2} ),\oU(p_{\mathrm g}, q_\mathrm g)), &\quad  \text{if $ \star = A$ and $p+q$ is odd}; \\
     (\widetilde \oU(\frac{n_{\mathrm b}}{2}, \frac{n_\mathrm b}{2} ), \widetilde \oU(p_{\mathrm g}, q_\mathrm g)), &\quad  \text{if $ \star = \wtA$}; \\
    (\SO(n_{\mathrm b},n_{\mathrm b}+1), \SO(p_{\mathrm g},q_{\mathrm g})), &\quad  \text{if } \star = B; \\
    (\SO(n_{\mathrm b}, n_{\mathrm b}), \Sp_{2n_{\mathrm g}}(\bR)), &\quad  \text{if } \star = C; \\
    (\SO^*(2n_\mathrm b), \Sp(\frac{p_{\mathrm g}}{2},\frac{q_{\mathrm g}}{2})), &\quad  \text{if }  \star = C^{*}; \\
    (\widetilde{\Sp}_{2n_{\mathrm b}}(\bR), \widetilde{\Sp}_{2n_{\mathrm g}}(\bR)), &\quad  \text{if } \star = \wtC; \\
     (\SO(n_{\mathrm b},n_{\mathrm b}), \SO(p_{\mathrm g},q_{\mathrm g}) ), &\quad  \text{if }  \star = D;\\
   (\SO^{*}(2 n_{\mathrm b}), \SO^{*}(2n_{\mathrm g})) , &\quad  \text{if } \star = D^{*}, \\
    \end{cases}
  \end{equation}
  where
\[
  (p_\mathrm g, q_{\mathrm g}) =
  \begin{cases}
  (p-\frac{n_\mathrm b}{2}, q-\frac{n_\mathrm b}{2} ),  &\quad  \text{if } \star \in \{A, \wtA\};\\
   (p-n_\mathrm b, q-n_\mathrm b),  &\quad  \text{if } \star \in \{B, C^*, D\}.\\
  \end{cases}
\]
Then $G_\mathrm b$ has type $\star_\mathrm b$ and $G_\mathrm g$ has type $\star$. 
See \eqref{def:starg} and \eqref{def:starb}. 
We remark that  $G_\mathrm b\times G_\mathrm g$ is an endoscopic group of $G$ (\cite[Chapter 21]{ABV}), when  $\star\ne \wtC$.

\begin{defn}\label{def:VbVg}
Define $V_\mathrm b$ to be the $[\star_\mathrm b]$-space such that $G_\mathrm b$ is identified with $\mathsf G_{\star_\mathrm b}(V_\mathrm b)$.
Likewise define $V_\mathrm g$ to be the $[\star]$-space such that $G_\mathrm g$ is identified with $\mathsf G_{\star}(V_\mathrm g)$. 
\end{defn}

As in \eqref{gc}, we have the complex groups $G_{\mathrm b,\C}$ and $G_{\mathrm g,\C}$, and the complexification homomorphisms $\iota_\mathrm b: G_\mathrm b\rightarrow G_{\mathrm b,\C}$ and  the complexification homomorphisms $\iota_\mathrm g: G_\mathrm g\rightarrow G_{\mathrm g,\C}$.

Write $\hha_\mathrm b$ and $\hha_\mathrm g$ for the universal Cartan subalgebras for $\g_\mathrm b$ and $\g_\mathrm g$, respectively. Then we have identifications
 \[
   \hha_\mathrm b^*=\C^{n_\mathrm b}=\check \hha_\mathrm b\qquad\textrm{and}\qquad  \hha_\mathrm g^*=\C^{n_\mathrm g}=\check \hha_\mathrm g.
 \]
 Using these identifications, we view $\check \g_\mathrm b$ as the Langlands dual (or the metaplectic Langlands dual when $\star=\wtC$) of  $\g_\mathrm b$, and view $\check \g_\mathrm g$ as the Langlands dual (or the metaplectic Langlands dual when $\star=\wtC$) of $\g_\mathrm g$.

 We also have identifications
   \[
    \hha_\mathrm b^*\times  \hha_\mathrm g^*= \C^{n_\mathrm b}\times  \C^{n_\mathrm g}= \C^{n}= \hha^*.
 \]
 Under this identification,  \[
   \Delta_\mathrm b^+=\Delta_\mathrm b\cap \Delta^+\qquad\textrm{and}\qquad  \Delta_\mathrm g^+=\Delta_\mathrm g\cap \Delta^+,
 \]
 where $\Delta_\mathrm b^+\subset \Delta_\mathrm b \subset \hha_\mathrm b^*$ are respectively the positive root system and the root system of $\g_\mathrm b$, and likewise $ \Delta_\mathrm g^+ \subset \Delta_\mathrm g\subset \hha_\mathrm g^*$ are respectively the positive root system and the root system of $\g_\mathrm g$.

 Similar to \eqref{weightl}, we have the analytic weight lattices
\[
  Q_{\iota_\mathrm b}=\Z^{n_\mathrm b}\subset \C^{n_\mathrm b}=(\C^{n_\mathrm b})^*=\hha_\mathrm b^*,
\]
and
\[
  Q_{\iota_\mathrm g}=\Z^{n_\mathrm g}\subset \C^{n_\mathrm g}=(\C^{n_\mathrm g})^*=\hha_\mathrm g^*.
\]
Then $\Lambda_\mathrm b$ is the unique $Q_{\iota_\mathrm b}$-coset in $\hha_\mathrm b^*$ that contains some (equiv. all) elements in $\hha_\mathrm b^*$ that represents the $\Ad(\check \g_\mathrm b)$-conjugacy class of $\lambda_{\check \CO_\mathrm b}^\circ$. Likewise, $\Lambda_\mathrm g$ is the unique $Q_{\iota_\mathrm g}$-coset in $\hha_\mathrm g^*$ that contains some (equiv. all) elements in $\hha_\mathrm g^*$ that represents the $\Ad(\check \g_\mathrm g)$-conjugacy class of $\lambda_{\check \CO_\mathrm g}^\circ$.

 Write $W_\mathrm b\subset \GL(\hha_\mathrm b)$ for the  Weyl group for $\g_\mathrm b$ so that
 \begin{equation}\label{def:Weylb}
   W_\mathrm b= \begin{cases}
    \sfS_{n_\mathrm b}, &  \text{if $\star \in \{A^\R, A^\bH, A, \wtA\} $}; \\
    \sfW_{n_\mathrm b}, &  \text{if $\star \in \{B,\wtC\} $};\\
     \sfW'_{n_\mathrm b}, &  \text{if $\star \in \{C,C^*, D, D^*\} $},\\
      \end{cases}
 \end{equation}
 and write $W_\mathrm g\subset \GL(\hha_\mathrm g)$ for the Weyl group for $\g_\mathrm g$ so that
 \begin{equation}\label{def:Weylg}
    W_\mathrm g= \begin{cases}
    \sfS_{n_\mathrm g}, &  \text{if $\star \in \{A^\R, A^\bH, A, \wtA\} $}; \\
    \sfW_{n_\mathrm g}, &  \text{if $\star \in \{B, C, \wtC, C^*\} $};\\
     \sfW'_{n_\mathrm g}, &  \text{if $\star \in \{ D, D^*\} $}.\\
      \end{cases}
 \end{equation}
Let $W'_\mathrm g\subset \GL(\hha_\mathrm g)$ denote the integral Weyl group of $\Lambda_\mathrm g$, namely the Weyl group of the root system
\[
  \{\alpha\in \Delta_\mathrm g\mid \la \check \alpha, \nu\ra\in \Z\textrm{ for some (and all) $\nu\in \Lambda_\mathrm g$}\}.
\]
 Then \[
   W'_\mathrm g= \begin{cases}
     \sfW'_{n_\mathrm g}, &  \text{if $\star =\wtC$};\\
     W_\mathrm g, &\text{otherwise}.
      \end{cases}
 \]
%We remark that $W_\mathrm b$ agrees with the integral  Weyl group of $\Lambda_\mathrm g$ in all cases, and $W_\mathrm g$ agrees with the   Weyl group of $\g_\mathrm g$ except when $\star=\wtC$.
Note that
 \be\label{wl00}
   W(\Lambda)=
    W_{\mathrm b}\times W'_{\mathrm g}\subset W_{\mathrm b}\times W_{\mathrm g}\subset W_\Lambda. \ee
% and
% \[
%   W_\Lambda= \begin{cases}
%    \sfW_n'\cap (\sfW_{n_\mathrm b}\times \sfW_{n_\mathrm g}), &  \text{if $\star \in \{D, D^*\}$}; \\
%    W_\mathrm b\times W_\mathrm g, &  \text{otherwise.}
%      \end{cases}
% \]

Let $\lambda_{\check \CO_\mathrm b}\in \Lambda_\mathrm b$ be the unique  dominant element  that represents the $\Ad(\check \g_\mathrm b)$-conjugacy class of $\lambda_{\check \CO_\mathrm b}^\circ$. Also take a dominant element $\lambda_{\check \CO_\mathrm g}\in \Lambda_\mathrm g$ that represents the $\Ad(\check \g_\mathrm g)$-conjugacy class of $\lambda_{\check \CO_\mathrm g}^\circ$, which is unique unless $\star=\wtC$. Put $\lambda_{\check \CO}:=(\lambda_{\check \CO_\mathrm b}, \lambda_{\check \CO_\mathrm g})\in \Lambda$.

 \begin{lem}\label{dmlam}
 The element $\lambda_{\check \CO}\in \hha^*$ is dominant and represents the conjugacy class of $\lambda_{\check \CO}^\circ$.
 \end{lem}
 \begin{proof}
 Note that for all $\alpha\in \Delta$,
 \[
   \la \lambda_{\check \CO}, \check \alpha\ra\in \Z \qquad \textrm{implies}\qquad \alpha\in \Delta_\mathrm b\sqcup \Delta_\mathrm g.
 \]
 This implies that $ \lambda_{\check \CO}$ is dominant.

  It is easy to see that there exists an element $\lambda'_{\check \CO_\mathrm b}\in \hha^*_\mathrm b$ that represents the conjugacy class of $\lambda_{\check \CO_\mathrm b}^\circ$ and an  element $\lambda_{\check \CO_\mathrm g}'\in \hha^*_\mathrm g$ that represents the conjugacy class of $\lambda_{\check \CO_\mathrm g}^\circ$ such that
  $\lambda'_{\check \CO}:=(\lambda'_{\check \CO_\mathrm b}, \lambda'_{\check \CO_\mathrm g})$ represents the conjugacy class of $\lambda_{\check \CO}^\circ$.
  Since $\lambda_{\check \CO_\mathrm b}$ is $W_\mathrm b$-conjugate to $\lambda'_{\check \CO_\mathrm b}$, and $\lambda_{\check \CO_\mathrm g}$ is $W_\mathrm g$-conjugate
    to $\lambda'_{\check \CO_\mathrm b}$, $\lambda_{\check \CO}$ is $W$-conjugate to $\lambda'_{\check \CO}$. Thus $\lambda_{\check \CO}\in \hha^*$ also represents the conjugacy class of $\lambda_{\check \CO}^\circ$.
   \end{proof}

 \begin{defn}\label{def:split} We say that a $[\star]$-space $V_0$ is split if
\begin{itemize}
\item $\star\in \{A^\R, A^\bH, C, \wtC\}$, or

\item $\star\in \{A, \wtA\}$ and there exists a totally isotropic subspace of $V_0$ of dimension $\frac{\dim V_0}{2}$, or
\item
 $\star\in \{B, D, C^*, D^*\}$ and  there exists a $\mathbf j_{V_0}$-stable totally isotropic subspace of $V_0$ of dimension $\frac{\dim V_0}{2}$.
 \end{itemize}
 \end{defn}

\begin{remark} In the above definition, $\mathbf j_{V_0}$ is part of the structure data of the standard representation $V_0$, which fixes a real or quarternionic structure. See \Cref{subsec:standard}. 
\end{remark}

Up to isomorphism a split $[\star]$-space is  determined by its dimension.

Recall from \Cref{def:VbVg} the spaces  $V_\mathrm b$ and $V_\mathrm g$. 
We define a modification of $V_\mathrm b$ as follows. 

\begin{defn}\label{def:splitV}
  Denote by $V_{\mathrm b}^\mathrm g$ the split $[\star]$-space with 
\[
\dim V_{\mathrm b}^\mathrm g=
 \begin{cases}
   \dim V_{\mathrm b}-1,  &  \text{if $\star=B$};\\
   \dim V_{\mathrm b},    &  \text{otherwise}.
    \end{cases}
 \]
 \end{defn}
For simplicity, write $G_\mathrm b^\mathrm g:=\mathsf G_{\star}(V_\mathrm b^\mathrm g)$. Note that
\be\label{isostar}
  V\cong V_\mathrm b^\mathrm g\times V_\mathrm g
\ee
as $[\star]$-spaces.

\subsection{Matching relevant Cartan subgroups}

By a $\star$-representation of a Lie group $E$, we mean a $[\star]$-space $V_1$ together with a Lie group homomorphism $E\rightarrow \mathsf G_\star(V_1)$. An isomorphism of $\star$-representation is defined to be a $[\star]$-space  isomorphism $f: V_1\rightarrow V_2$ of  $\star$-representations that makes the diagram
\[
\begin{tikzcd}[column sep={4cm,between origins}]
      & E
      \ar[dl," "']\ar[dr,"  "]&\\
      \mathsf G_\star(V_1)\ar[rr,"\textrm{the isomorphism induced by $f$}"]& &\mathsf G_\star(V_2).\\
    \end{tikzcd}
\]
commute. The product $V_1\times V_2$ is obviously a $\star$-representation of $E$ whenever $V_1$ and $ V_2$  are $\star$-representations of $E$.

\begin{defn} We say that a Cartan subgroup $H_\mathrm b$ of $G_\mathrm b$ is relevant if
\begin{itemize}
\item
$\star_\mathrm b\in \{A^\R, A^\bH\}$, or
\item $\star\in \{A, \wtA\}$ and $H_\mathrm b$ stabilizes a totally isotropic subspace of $V_\mathrm b$ of dimension $\frac{n_\mathrm b}{2}$, or
\item
 $\star\in \{B,C,\wtC, D, C^*, D^*\}$ and $H_\mathrm b$ stabilizes  a $\mathbf j_{V_\mathrm b}$-stable totally isotropic subspace of $V_\mathrm b$ of dimension $n_\mathrm b$.
\end{itemize}
\end{defn}

For a relevant Cartan subgroup $H_\mathrm b$ of $G_\mathrm b$ as in the above definition, $V_\mathrm b$ is naturally a $\star_\mathrm b$-representation of $H_\mathrm b$.

The following lemma can be verified directly. Recall $V^g_b$ from  \Cref{def:splitV}.

\begin{lem}\label{match1} \label{match2} \label{match3}
There exits a unique (up to isomorphism) $\star$-representation  of $H_\mathrm b$ on the $[\star]$-space $V_\mathrm b^\mathrm g$ 
with the following property.
\begin{enumerate}[label=(\alph*),wide=0pt]
    \item When $\star=B$, $V_\mathrm b\cong V_1\times V_\mathrm b^\mathrm g$ 
    as $\star$-representations of $H_\mathrm b$ 
    for a one-dimensional $\star$-representation $V_1$ of $H_\mathrm b$ with trivial action. 
    \item When $\star\in \{A, \wtA, C, C^*\}$,  
    $V_\mathrm b\cong V_\mathrm b^\mathrm g$
    as representations of $H_\mathrm b$. 
    \item When $\star\in \{A^\bR, A^\bH, \wtC,  D, D^* \}$,  
    $V_\mathrm b\cong V_\mathrm b^\mathrm g$
    as $\star$-representations of $H_\mathrm b$. 
\end{enumerate}
\end{lem}

We view  $V_\mathrm b^\mathrm g$  as a $\star$-representation of  $H_\mathrm b$ as in \Cref{match1}.

  Let $H_\mathrm g$ be a Cartan subgroup of $G_\mathrm g$ so that $V_\mathrm g$ is naturally an $H_\mathrm g$-representation.
Then $V_\mathrm b^\mathrm g\times V_\mathrm g$ is naturally a $\star$-representation of $H_\mathrm b\times H_\mathrm g$. The isomorphism \eqref{isostar} yields a Lie group homomorphism
\be\label{matchmap}
H_\mathrm b\times H_\mathrm g\rightarrow G.
\ee
Write $H$ for the image of the homomorphism \eqref{matchmap}, and let $\zeta: H_\mathrm b\times H_\mathrm g\rightarrow H$ denote the homomorphism induced by \eqref{matchmap}.
Note that $H$ is a Cartan subgroup of $G$. We call the pair $(H, \zeta)$ a matching of $H_\mathrm b\times H_\mathrm g$, which is uniquely defined up to conjugation by $G$.

When $\star=A$ and $p+q$ is odd, $H_\mathrm b\cong (\C^\times)^{\frac{n_\mathrm b}{2}}\times \{\pm 1\}$. In all cases, we put
\[
  H'_\mathrm b:=\begin{cases}
      \textrm{the identity connected component of $H_\mathrm b$},\quad & \textrm{if $\star=A$ and $p+q$ is odd};\\
      H_\mathrm b,\quad & \textrm{otherwise}. 
  \end{cases}
\]
Then the homomorphism  $\zeta: H'_\mathrm b\times H_\mathrm g\rightarrow H$ is a double cover when $\star\in \{\wtA, \wtC\}$, and is an isomorphism in all other cases. 

\subsection{Matching the parameters}
Given a parameter $\gamma_\mathrm b \in\cP'_{\Lambda_\mathrm b}(G_\mathrm b)$ that is represented by 
$(H_\mathrm b, \xi_\mathrm b, \Gamma_\mathrm b)$,
and  a parameter $\gamma_\mathrm g \in\cP'_{\Lambda_\mathrm g}(G_\mathrm g)$ 
that is represented by $(H_\mathrm g, \xi_\mathrm g, \Gamma_\mathrm g)$,
we  define  a parameter $\varphi(\gamma_\mathrm b,  \gamma_\mathrm g) \in\cP'_{\Lambda}(G)$
as in what follows.
Note that \Cref{leme} implies that $H_\mathrm b$ is  relevant.
Take a matching $(H, \zeta)$ of $H_\mathrm b\times H_\mathrm g$.
Let $\xi$ be the composition of
\[
    \hha^*=\hha_\mathrm b^*\times  \hha_\mathrm g^*\xrightarrow{\xi_\mathrm b\times \xi_\mathrm g} \h_\mathrm b^*\times \h_\mathrm g^*\xrightarrow{\textrm{the transpose inverse of the complexified differential of $\zeta$}} \h^*.
\]
Let $\Gamma: \Lambda\rightarrow \Irr'(G)$  be the map
\[
  \nu=(\nu_\mathrm b, \nu_\mathrm g)\mapsto \Gamma_{\mathrm b,\nu_\mathrm b}\otimes \Gamma_{\mathrm g,\nu_\mathrm g},
\]
where  $\Gamma_{\mathrm b,\nu_\mathrm b}\otimes \Gamma_{\mathrm g,\nu_\mathrm g}$ 
(which is originally defined as an irreducible representation of $H_\mathrm b\times H_\mathrm g$)
is viewed as an irreducible representation of $H$ 
via the descent through the homomorphism 
$\xi: H'_\mathrm b\times H_\mathrm g\rightarrow H$.

\begin{lem}
The triple $(H, \xi, \Gamma)$ defined above is an element of $\sP'_{\Lambda}(G)$.
\end{lem}
\begin{proof}
Note that
every imaginary root of $G$ with respect to $H$ is either an  imaginary root of $G_\mathrm b$ with respect to $H_\mathrm b$ or an  imaginary root of $G_\mathrm g$ with respect to $H_\mathrm g$ (here $\h$ is identified with $\h_\mathrm b\times \h_\mathrm g$ via $\xi$). This  implies that $\delta(\xi)=\delta(\xi_\mathrm b)+\delta(\xi_\mathrm g)$, and the lemma then easily follows.
\end{proof}

Now we define  $\varphi(\gamma_\mathrm b, \gamma_\mathrm g) \in\cP'_{\Lambda}(G)$ to be the $G$-orbit of the triple  $(H, \xi, \Gamma)$ defined above.
This is independent of the choices of the representatives $(H_\mathrm b, \xi_\mathrm b, \Gamma_\mathrm b)$, $(H_\mathrm g, \xi_\mathrm g, \Gamma_\mathrm g)$, and the matching $(H, \zeta)$. It is easy to see that the map
\be\label{match0}
  \varphi: \cP'_{\Lambda_\mathrm b}(G_\mathrm b)\times \cP'_{\Lambda_\mathrm g}(G_\mathrm g)\rightarrow \cP'_{\Lambda}(G)
\ee
is $W_\mathrm b\times W_\mathrm g$-equivariant under the cross actions.

 \begin{prop}\label{para}
 If either $\star=C^*$ and $n_\mathrm b>0$, or $\star=D^*$ and $n_\mathrm b, n_\mathrm g>0$  (so $W_\mathrm b\times W_\mathrm g$ is of index $2$ in $W_\Lambda$), then the map \eqref{match0} is injective and
 \[
   \cP'_{\Lambda}(G)=\varphi(\cP'_{\Lambda_\mathrm b}(G_\mathrm b)\times \cP'_{\Lambda_\mathrm g}(G_\mathrm g))\sqcup w\times\left(\varphi(\cP'_{\Lambda_\mathrm b}(G_\mathrm b)\times \cP'_{\Lambda_\mathrm g}(G_\mathrm g))\right)
 \]
 for every $w\in W_\Lambda\setminus (W_\mathrm b\times W_\mathrm g)$, where ``$\,w\times$" indicates the cross action of $w$. The map \eqref{match0} is bijective in all other cases.
 \end{prop}
 \begin{proof}
  In view of \Cref{leme}, this follows by a direct verification. 
 \end{proof}

Define a linear map
\be\label{match1234}
  \varphi: \Coh_{\Lambda_\mathrm b}(\CK'(G_\mathrm b))\otimes \Coh'_{\Lambda_\mathrm g}(\CK'(G_\mathrm g))\rightarrow \Coh_{\Lambda}(\CK'(G))
\ee
such that
\[
 \varphi(\Psi_{\gamma_\mathrm b}\otimes \Psi_{\gamma_\mathrm g})=\Psi_{\varphi(\gamma_\mathrm b, \gamma_\mathrm g)}
\]
for all $\gamma_\mathrm b \in\cP'_{\Lambda_\mathrm b}(G_\mathrm b)$ and $\gamma_\mathrm g \in\cP'_{\Lambda_\mathrm g}(G_\mathrm g)$.

\begin{prop}\label{matchth}
The linear map \eqref{match1234} is $W_\mathrm b\times W_\mathrm g$-equivariant and injective. Moreover
\[
 \varphi(\overline \Psi_{\bar \gamma_\mathrm b}\otimes \overline \Psi_{\bar \gamma_\mathrm g})=\overline \Psi_{\varphi(\bar \gamma_\mathrm b, \bar \gamma_\mathrm g)}
\]
for all $\gamma_\mathrm b \in\cP'_{\Lambda_\mathrm b}(G_\mathrm b)$ and $\gamma_\mathrm g \in\cP'_{\Lambda_\mathrm g}(G_\mathrm g)$.
\end{prop}
\begin{proof}
  This follows from the Kazhdan-Lusztig-Vogan algorithm (see \cite{V3}*{Theorems~1.6 and 1,12} and \cite[Chapter 15]{ABV}). The case of the odd orthogonal group is treated in \cite{GI}*{Section~3} and 
  the case of the real metaplectic group is treated in \cite{RT2}*{Theorem~5.2}. 
\end{proof}

By \Cref{para} and  \Cref{matchth}, we have the following corollary. 

\begin{cor}\label{cor:Ind.DCstar}
 If $\star\in \{C^*, D^*\}$, then
\[
 \Coh_{\Lambda}(\CK'(G))\cong \Ind^{W_\Lambda}_{W_\mathrm b\times W_\mathrm g} (\Coh_{\Lambda_\mathrm b}(\CK'(G_\mathrm b))\otimes \Coh_{\Lambda_\mathrm g}(\CK'(G_\mathrm g)))
 \]
 as representations of $W_\Lambda$. In all the other cases,
\[
 \Coh_{\Lambda}(\CK'(G))\cong \Coh_{\Lambda_\mathrm b}(\CK'(G_\mathrm b))\otimes \Coh_{\Lambda_\mathrm g}(\CK'(G_\mathrm g))
 \]
 as representations of $W_\mathrm b\times W_\mathrm g$.
\end{cor}

The following proposition is a reformulation of results 
in \cite{Mat04}*{Lemma~4.1.3} and \cite{GI}*{Lemma~3.3}. 

\begin{prop}\label{propKL33}
Let $\nu=(\nu_\mathrm b, \nu_\mathrm g)\in \Lambda_\mathrm b\times \Lambda_\mathrm g=\Lambda$. Then there is a unique linear map $\varphi_\nu: \CK'_{\nu_\mathrm b}(G_\mathrm b)\otimes \CK'_{\nu_\mathrm g}(G_\mathrm g)\rightarrow \CK'_{\nu}(G)$ that makes the diagram
\[
 \begin{CD}
          \Coh_{ \Lambda_\mathrm b}(\CK'(G_\mathrm b))\otimes  \Coh_{ \Lambda_\mathrm g}(\CK'(G_\mathrm g))
                  @>   \varphi  >>  \Coh_{ \Lambda}(\CK'(G))\\
            @V   \mathrm{ev}_{\nu_\mathrm b} \otimes  \mathrm{ev}_{\nu_\mathrm g}  VV         @ VV  \mathrm{ev}_{\nu} V \\
      \CK'_{\nu_\mathrm b}(G_\mathrm b)\otimes \CK'_{\nu_\mathrm g}(G_\mathrm g) @> \varphi_\nu >>  \CK'_{\nu}(G) \\
  \end{CD}
\]
commutes, where $\mathrm ev$ indicates the evaluation maps. Moreover, $\varphi_\nu$ is injective,  and
\[
\varphi_\nu(\Irr'_{\nu_\mathrm b}(G_\mathrm b)\times  \Irr'_{\nu_\mathrm g}(G_\mathrm g))\subset \Irr'_{\nu}(G).
\]
\end{prop}
\begin{proof}
The first assertion follows directly from \Cref{prop:ev00000}, since the map $\varphi$ is $W(\Lambda)=W_\mathrm b\times W'_\mathrm g$-equivariant. Since $\varphi$ is injective,  \Cref{prop:ev00000} also implies that $\varphi_\nu$ is injective.

For the proof of the second assertion, we assume without loss of generality that $\nu$ is dominant. Let $\pi_\mathrm b\in \Irr'_{\nu_\mathrm b}(G_\mathrm b)$ and $\pi_\mathrm g\in \Irr'_{\nu_\mathrm g}(G_\mathrm g)$. Pick  basal elements
\[
  \Psi_\mathrm b \in \Coh_{ \Lambda_\mathrm b}(\CK'(G_\mathrm b)) \quad \textrm{and}\quad \Psi_\mathrm g\in  \Coh_{ \Lambda_\mathrm g}(\CK'(G_\mathrm g))
   \]
 such that $\Psi_\mathrm b(\nu_\mathrm b)=\pi_\mathrm b$ and $\Psi_\mathrm g(\nu_\mathrm g)=\pi_\mathrm g$.
 Then
\[
  \varphi_\nu(\pi_\mathrm b\otimes \pi_\mathrm g)=\mathrm{ev}_{\nu}(\varphi(\Psi_\mathrm b\otimes \Psi_\mathrm g)).
\]
\Cref{matchth} implies that $\varphi(\Psi_\mathrm b\otimes \Psi_\mathrm g)$ is a basel element of  $\Coh_{\Lambda}(\CK'(G))$. Thus by \Cref{lemirr}, $\mathrm{ev}_{\nu}(\varphi(\Psi_\mathrm b\otimes \Psi_\mathrm g))$ is either irreducible or zero. It is clearly nonzero since $\varphi_\nu$ is injective. This proves the last assertion of the proposition.
\end{proof}

\section{Special unipotent representations in type $A$}\label{sec:GL}

We begin the explicit counting of special unipotent representations for a real classical group $G$. 

Recall that $\CK'(G)$ is the Grothendieck group of the category of genuine Casselman-Wallach representations of $G$ if $\star\in \{\widetilde A, \widetilde C\}$, and $\CK'(G):=\CK(G)$ otherwise.  

Recall we have the Lusztig left cell $\LC_{\lambda}$ attached to $\lambda \in \Lambda $ (\Cref
{def:ll}). To simplify the notation as well as the terminology, we will also write $\LC_{\ckcO}:= \LC_{\lambda_{\ckcO}}$, and call it the Lusztig left cell attached to $\ckcO$. To reiterate, we repeat its definition: ($\lambda_{\ckcO}\in \Lambda$)
\begin{equation} \label{llO}
\LC_{\ckcO}:=\left\{\sigma\in \Irr(W(\Lambda))\mid \sigma \textrm{ occurs in } \left(J_{W_{\lambda _{\ckcO}}}^{W(\Lambda)} \sgn \right)\otimes \sgn\right \}. 
\end{equation}

As an obvious variant of \Cref{countleq} for $\nu=\lambda_{\check \CO}$, we have that
 \begin{equation}\label{boundc22}
     \sharp(\Unip_{\check \CO}(G)) 
     =\sum_{\sigma\in \LC_{\check \CO}} [\sigma: \Coh_{\Lambda}(\CK'(G))]. 
   \end{equation}
(The equality holds because of \Cref{HCLU}.) Here $\Lambda \subset \hha^*$ is defined in \eqref{eq:Lambda}.

\subsection{Some Weyl group representations}

The group $\sfS_n$ is  identified with the permutation group of the set $\{1,2, \dots, n\}$, and
 $\sfW_n$ is identified with $\sfS_n \ltimes \set{\pm 1}^n$.

 Define a quadratic character
\be\label{defep}
  \varepsilon: \sfW_n\rightarrow \{\pm 1\}, \quad (s,(x_{1}, x_{2}, \cdots, x_{n}))\mapsto x_{1}x_{2}\cdots x_{n}.
\ee
Then $\sfW_n'$ is the kernel of this character. As always, $\sgn$ denotes the sign character (of an appropriate Weyl group).
Since $\sfS_{n}$ is a quotient of $\sfW_{n}$, we may inflate the sign character of $\sfS_{n}$ to obtain a character of $\sfW_{n}$, to be denoted by $\bsgn$: 
\begin{equation}
\label{eq:bsgn}
\bsgn: \quad \sfW_n\longrightarrow \sfS_n\overset{\sgn}{\longrightarrow} \{\pm 1\}.
\end{equation} We have
\[\varepsilon = \sgn \otimes \bsgn.\]

We also have the natural embedding 
\[\sfW_n\hookrightarrow \sfS_{2n}\]
via the homomorphism determined by
\begin{equation}\label{wn1}
\sfS_n \ni (i, i+1)\mapsto (2 i-1, 2i+1)(2i, 2i+2), \qquad (1\leq i\leq n-1),
\end{equation}
and
\begin{equation}\label{wn2}
\set{\pm 1}^n \ni (1,\cdots,1, \underbrace{-1}_{j\text{-th term}}, 1, \cdots, 1) \mapsto (2j-1, 2j), \qquad (1\leq j\leq n).
\end{equation}
Here $(i, i+1)$, $(2i-1, 2i+1)$, etc., indicate the transpositions in the permutation groups.  Note that $\varepsilon $ is also the restriction of $\sgn$ of $\sfS_{2n}$ to $\sfW_n$.

 Let $\YD_{n}$ be the set of Young diagrams of total size $n$. Suppose that $\star\in  \{A^\R, A^\BH, A, \widetilde A\}$ in the rest of this section.
We identify $\YD_{n}$ with the set $\overline{\Nil}(\g)= \GL_{n}(\bC)\backslash \Nil(\g\l_n(\C))$ of complex nilpotent orbits and
also with the set  $\Irr(\sfS_{n})$ via the Springer
correspondence (see \cite{Carter}*{11.4}).
More specifically, for $ \cO' \in \overline \Nil( \g)=\YD_{n}$, the Springer
correspondence is given by Macdonald's construction for $\sfS_{n}$ via $j$-induction:
\[
  \Spr( \cO') = j_{\prod_{j}\sfS_{\bfcc_{j}(\cO')}}^{\sfS_{n}} \sgn.
\]

Recall that  $W(\Lambda)  = \sfS_{n_{\mathrm b}}\times \sfS_{n_{\mathrm g}}$. It is easy to verify that
\begin{equation}\label{eq:WA}
  \begin{split}
      W_{\lamck} & = \prod_{i\in \bN^{+}}\sfS_{\bfcc_{i}(\ckcO_{\mathrm b})}\times \prod_{i\in \bN^{+}}\sfS_{\bfcc_{i}(\ckcO_{\mathrm g})},\\
    \LC_{\lamck} & = \set{\tau_{\lamck}}, \quad \textrm{where }  \tau_{\lamck} := (j_{W_{\lamck}}^{W(\Lambda)}\sgn )\otimes \sgn =  \Spr(\ckcOb^{t})\otimes \Spr(\ckcOg^{t}).
  \end{split}
\end{equation}
Here and as before, a superscript ``$t$" indicates the transpose of a Young diagram.
Also note that the last tensor represents the external tensor product of representations of $\sfS_{n_{\mathrm b}}$ and $\sfS_{n_{\mathrm g}}$.  

The following Propositions \ref{count000}-\ref{count002} follow from Theorem \ref{thm:cohHC} by direct computation. We omit the details.

   \begin{prop} \label{count000}
 Suppose that $\star=A^\R$.  For each $l\in \BN$, put
  \[
    \cC_l := \bigoplus_{\substack{t,c,d\in \bN \\2t+c+d=l}} \Ind_{\sfW_t\times \sfS_c\times \sfS_d}^{\sfS_{l}} \varepsilon \otimes 1\otimes 1.
  \]
  Then
  \[
    \Coh_{\Lambda_\mathrm b}(\CK'(G_\mathrm b)) \cong \cC_{n_{\mathrm b}}\quad\textrm{and}\quad   \Coh_{\Lambda_\mathrm g}(\CK'(G_\mathrm g)) \cong \cC_{n_{\mathrm g}}.
  \]
 \end{prop}

\begin{prop}\label{count001}
 Suppose that $\star=A^\BH$. For each even number $l\in \BN$, put
  \[
  \cC_{l}: =
    \Ind_{\sfW_{\frac{l}{2}}}^{\sfS_{l}}\varepsilon.
  \]
  Then \[
    \Coh_{\Lambda_\mathrm b}(\CK'(G_\mathrm b)) \cong \cC_{n_{\mathrm b}}\quad\textrm{and}\quad   \Coh_{\Lambda_\mathrm g}(\CK'(G_\mathrm g)) \cong \cC_{n_{\mathrm g}}.
   \]
\end{prop}

\begin{prop}\label{count002}
Suppose that $\star\in \{A, \wtA\}$.   Then
    \[
       \Coh_{\Lambda_\mathrm b}(\CK'(G_\mathrm b)) \cong  \Ind_{\sfW_{\frac{n_{\mathrm b}}{2}}}^{\sfS_{n_{\mathrm b}}} 1
    \]
    and
    \[
       \Coh_{\Lambda_\mathrm g}(\CK'(G_\mathrm g)) \cong  \bigoplus_{\substack{t,s,r\in \bN\\t+r=p_\mathrm g, t+s = q_\mathrm g}}
    \Ind_{\sfW_{t}\times \sfS_s\times \sfS_r}^{\sfS_{n_\mathrm g}}
 1\otimes \sgn \otimes \sgn.
    \]
    \end{prop}

\subsection{Special unipotent representations of $\GL_n(\bR)$ and $\GL_n(\bH)$}\label{sec:GLRH11}

Recall from \Cref{defpbp0} and  \Cref{defpbp1} the set $\PAP_{\star}(\ckcO)$, which is the set of paintings on $\ckcO^{t}$ that has type $\star\in \{A^\R, A^\bH, A, \wtA\}$.
Recall that if $\star=A^\R$, then
\begin{equation}\label{eq:PPA.count}
  \#(\PAP_{\star}(\ckcO)) = \prod_{r\in \bN^{+}} (\#\set{i\in \bN^{+}| \bfrr_{i}(\ckcO)=r}+1).
\end{equation}
Also note that
if $\star=A^\bH$, then
\begin{equation}\label{eq:PPA.count2}
  \#(\PAP_{\star}(\ckcO)) = \begin{cases}
   1,
    & \text{ if $\check \CO=\check \CO_\mathrm g$;}  \\
      0, & \text{ otherwise.} \\
    \end{cases}
\end{equation}

\begin{prop} \label{lem:GL.count2}
  Suppose that $\star \in\{A^\R, A^\bH\}$.  Then
   \begin{equation*}%\label{eq:A.count}
    [\tau_{\lamck}: \Cint{\Lambda}(\CK'(G))] = \# (\PAP_{\star}(\ckcOg))\times
    \# (\PAP_{\star}(\ckcOb)) = \# (\PAP_{\star}(\ckcO)).
  \end{equation*}
\end{prop}
\begin{proof}
  In view of Proposition \ref{count000} and the description of the left cell representation $\tau_{\lamck}$ in \eqref{eq:WA}, the  first equality  follows from  Pieri's rule (\cite[Corollary 9.2.4]{GW}) and
  the following branching formula (see \cite{BV.W}*{Lemma~4.1~(b)}):
\begin{equation}\label{IndFor}
  \Ind_{\sfW_{t}}^{\sfS_{2t}} \varepsilon = \bigoplus_{\substack{\sigma\in \YD_{2t}\\
      \bfcc_{i}(\sigma) \text{ is even for all }i\in \bN^+}} \sigma\qquad (t\in \bN).
\end{equation}
The last equality
  follows from \eqref{eq:PPA.count} and \eqref{eq:PPA.count2}.
\end{proof}

We are in the setting of \Cref{thm:mainR00}. Recall the map
  \[
    \begin{array}{ccc}
      \PAP_{\star}(\ckcO) & \rightarrow & \Unip_{\ckcO}(G),\\
      \CP & \mapsto & \pi_{\CP{}}, 
    \end{array}
  \]
 where $\pi_{\CP}$ is defined in \Cref{sec:GLRH}. It  is proved in  \cite[Theorem 3.8]{V.GL} that the above map is injective. Then   the map is bijective by \eqref{boundc22} and \Cref{lem:GL.count2}. This proves  \Cref{thm:mainR00}.

\subsection{Special unipotent representations of unitary groups}\label{secunit}

In this subsection, we suppose  $\star \in \set{A, \wtA}$ so that
\[
  G =
  \begin{cases}
    \rU(p,q),  & \text{if }\star = A;\\
    \tU(p,q),  & \text{if }\star = \wtA.
\end{cases}
\]

For $\cP \in \PAP_{\star}(\ckcO)$, we have defined its signature $(p_\CP, q_\cP)$ in \eqref{eq:signature}. Recall that $\CO:=\dBV(\check \CO)=\ckcO^t\subset \Nil(\g^*)$.  Let $\overline \Nil_{G}(\cO)$ denote the set of $G$-orbits in $(\sqrt{-1}\g_0^*)\cap \CO$, where $\g_0$ denotes the Lie algebra of $G$ which equals $\u(p,q)$.

We first consider the case when  $\ckcO = \ckcO_{\mathrm g}$.
In this good parity setting, we will state a counting result on
$\Unip_{\ckcO}(G)$. The elements in $\Unip_{\ckcO}(G)$ can be constructed by
cohomological induction explicitly and they are irreducible and unitary. See \cite[Theorem 4.2]{BV.W}, and also 
\cite{Mat96,Tr.U}, \cite{Tr.H}*{Section~2} and \cite{MR.U}*{Section~4}.
We refer the reader to \cite[Section 12]{BMSZ2} for the construction of all elements of $\Unip_{\ckcO}(G)$ by the method of theta lifting.

\begin{thm}[\cf {\cite[Theorem 4.2]{BV.W} and \cite{Tr.H}*{Theorem~2.1}}]\label{thmunit0}
  Suppose that $\ckcO = \ckcO_{\mathrm g}$. Then
  \begin{equation} \label{sharpuni0} \sharp(\Unip_{\ckcO}(G))= \sharp \set{\CP\in \PAP_{\star}(\ckcO)|(p_\CP, q_\CP)=(p,q)}=\sharp(\overline \Nil_{G}(\cO)).
  \end{equation}
  Moreover, for every $\pi\in \Unip_{\ckcO}(G)$, its wavefront set $\WF(\pi)$ is
  the closure in $\sqrt{-1}\g_0^*$ of a unique orbit
  $\sO_\pi\in \overline \Nil_{G}(\cO)$, and the map
  \begin{equation}\label{sharpuni1}
    \begin{array}{rcl}
      \Unip_{\ckcO}(G) &\longrightarrow& \overline \Nil_{G}(\cO), \\%\Nil_{\mathrm g}(\cO),\\
      \pi & \mapsto & \sO_\pi
    \end{array}
  \end{equation}
  is bijective.
\end{thm}

\begin{proof} 
\delete{The Harish-Chandra cell representations in $\Coh_{\Lambda}(\CK'(G))$  are all irreducible and explicitly
  described in terms of the Springer correspondence by \cite{Bo}*{Lemma~4}. }
  In view of the counting formula in \eqref{boundc22}, the first equality  in \eqref{sharpuni0}
  follows from Proposition \ref{count002}, \eqref{eq:WA},  \eqref{IndFor}, and Pieri's rule (\cite[Corollary 9.2.4]{GW}).
  The second equality in \eqref{sharpuni0} will follow directly from
  the bijectivity of \eqref{sharpuni1}.
  
The assignment of wavefront set yields a bijection
 \[
   \{\textrm{cell in the basal representation $\Coh_\Lambda(\CK_{\overline{\CO}}'(G))/\Coh_\Lambda(\CK_{\overline{\CO}\setminus \CO}'(G))$}\} \rightarrow \overline \Nil_{G}(\cO),
 \]
and every cell representation in $\Coh_\Lambda(\CK_{\overline{\CO}}'(G))/\Coh_\Lambda(\CK_{\overline{\CO}\setminus \CO}'(G))$ is irreducible and isomorphic to $\tau_{\lambda_{\check \CO}}$. See \cite{BV.W}*{Theorem~4.2} and
  \cite[Theorem 5]{Bo} (we have used \cite{SV}*{Theorem~1.4} to rephrase the result in
  terms of real nilpotent orbits).

 Note that \eqref{mulone} implies   that
 \[
 [1_{W_{\lambda_{\ckcO}}} : \tau_{\lambda_{\check \CO}}]=1.
 \]
 Recall from \Cref{dmlam} that $\lambda_{\ckcO}\in \hha^*$ is dominant. As in the proof of  \Cref{counteq},  for each cell $\CC$ in $\Coh_\Lambda(\CK_{\overline{\CO}}'(G))$ that is not a cell in $\Coh_\Lambda(\CK_{\overline{\CO}\setminus \CO}'(G))$, there is a unique element $\Psi_\CC\in \CC$ such that $ \Psi_\CC(\lambda_{\ckcO})\neq 0$.
 Then
 \[
    \Unip_{\ckcO}(G) =\{ \Psi_\CC(\lambda_{\ckcO})\}_{\CC\textrm{ is a cell  in $\Coh_\Lambda(\CK_{\overline{\CO}}'(G))$ that is not a cell in $\Coh_\Lambda(\CK_{\overline{\CO}\setminus \CO}'(G))$}}.
 \]
 This implies the bijectivity assertion of the theorem.   
\end{proof}

\begin{lem}\label{thmca0000}
  If  $ \sharp(\Unip_{\ckcO}(G))\neq 0$, then
  $
    \ckcOb= 2\ckcOpb
  $
  for some Young diagram $\ckcOpb$.
\end{lem}

\begin{proof}
Suppose that  $ \sharp(\Unip_{\ckcO}(G))\neq 0$. Then
\be\label{tau5}
   [\tau_{\lambda_{\check \CO_\mathrm b}}: \Coh_{\Lambda_{\mathrm b}}(\CK'(G_\mathrm b))]\neq 0.
\ee
Recall from Proposition \ref{count002} that
  \[
       \Coh_{\Lambda_\mathrm b}(\CK'(G_\mathrm b)) \cong  \Ind_{\sfW_{t}}^{\sfS_{2t}} 1, \quad (t=\frac{n_{\mathrm b}}{2}). 
    \]
Thus \eqref{tau5}
 is the same as saying that
\[
 [ \Spr({\check \CO_\mathrm b}^t): \Ind_{\sfW_{t}}^{\sfS_{2t}} 1 ]\neq 0. 
\]
Similar to \eqref{IndFor}, we have that
\begin{equation}\label{IndFor3}
  \Ind_{\sfW_{t}}^{\sfS_{2t}} 1= \bigoplus_{\substack{\sigma\in \YD_{2t}\\
      \bfrr_{i}(\sigma) \text{ is even for all }i\in \bN^+}} \sigma\qquad (t\in \bN).
\end{equation}
This implies the lemma.
\end{proof}

Now we assume that $\ckcOb= 2\ckcOpb$
  for some Young diagram $\ckcOpb$.

The group $G$ has an $\check \CO$-relevant parabolic subgroup $P$ whose Levi quotient  is naturally isomorphic to $G'_{\mathrm b}\times G_{\mathrm g}$, where
$G'_{\mathrm b}:=\GL_{\frac{n_{\mathrm b}}{2}}(\C)$. Let $\pi_{\ckcOpb}$ denote the
unique element in $\Unip_{\ckcOpb}(\GL_{\frac{n_{\mathrm b}}{2}}(\C))$. Then for
every $\pi_\mathrm g\in \Unip_{\ckcO_{\mathrm g}}(\Gg)$, the normalized smooth parabolic
induction $\pi_{\ckcOpb}\rtimes \pi_\mathrm g$ is irreducible by
\cite{Mat96}*{Theorem~3.2.2} and is an element of $\Unip_{\ckcO}(G)$ (\cf
\cite{MR.U}*{Theorem~5.3}).

\begin{thm}\label{thmunit012}
  The equality  \begin{equation}\label{unitarred1}
    \sharp(\Unip_{\ckcO}(G)) =
    \sharp(\Unip_{\ckcO_{\mathrm g}}(\Gg))
  \end{equation}
  holds,
  and the map
  \begin{equation}
  \label{bij00}
    \begin{array}{ccc}
      \Unip_{\ckcOg}(\Gg)&\longrightarrow &\Unip_{\ckcO}(G),\\
      \pi_{\mathrm g}& \mapsto & \pi_{\ckcOpb}\rtimes \pi_{\mathrm g}\\
    \end{array}
  \end{equation}
  is bijective. 
\end{thm}

\begin{proof}
By the counting formula in \eqref{boundc22},
  \[
    \sharp(\Unip_{\check \CO}(G)) =[\tau_{\lamck}: \Coh_{\Lambda}(\CK'(G))].
  \]
  Similarly,
  \[
    \sharp(\Unip_{\check \CO_\mathrm g}(G_\mathrm g)) = [\tau_{\lambda_{\check \CO_\mathrm g}}: \Coh_{\Lambda_{\mathrm g}}(\CK'(G_\mathrm g))].
  \]
 The proof of \Cref{thmca0000} shows that
 \[
 [\tau_{\lambda_{\check \CO_\mathrm b}}:\Coh_{\Lambda_\mathrm b}(\CK'(G_\mathrm b))]=1.
 \]
  The above three equalities clearly imply \eqref{unitarred1}.

 By the  calculation of the wavefront set of the induced representations (\cite[Corollary 5.0.10]{B.Orbit}),  Theorem \ref{thmunit0} implies that all the representations $\pi_{\ckcOpb}\rtimes \pi_{\mathrm g}$, where $\pi_\mathrm g$ varies in $\Unip_{\ckcO_{\mathrm g}}(\Gg)$, have pairwise distinct wavefront sets.   Thus the map \eqref{bij00} is injective. Hence it is bijection by the counting assertion \eqref{unitarred1}.
\end{proof}

Theorems \ref{thmu1} and \ref{thmu2} then follow from Theorems \ref{thmunit0} and \ref{thmunit012}.

\begin{remark}
  When $\check \CO\neq \check \CO_{\mathrm g}$, the wavefront set $\WF(\pi)$, where $\pi\in \Unip_{\ckcO}(G)$,  may not be the closure of a single orbit in
   $\overline \Nil_{G}(\cO)$.
\end{remark}

\section{Special unipotent representations in type $BCD$: counting}\label{sec:BCD}

In this section and the next one, we assume that $\star \in \set{B,C,\wtC,C^{*},D,D^{*}}$.
We will give the detailed description of the coherent continuation representations and the Lusztig left cells, which yield the explicit counting of special unipotent representations via the Littlewood-Richardson rule.

\subsection{The coherent continuation representations}
 Let
\begin{equation}\label{eq:Ht}
\sfH_{t} := \sfW_t\ltimes \set{\pm 1}^t, \quad (t\in\bN),
\end{equation}
to be viewed as a subgroup in $\sfW_{2t}$ such
that
\begin{itemize}
  \item the first factor $\sfW_{t}$ sits in $\sfS_{2t}\subset \sfW_{2t}$ as in \eqref{wn1} and \eqref{wn2},
    \item the element $(1,\cdots,1, \underbrace{-1}_{i\text{-th
        term}}, 1, \cdots, 1)\in \set{\pm 1}^{t}$ acts on $\bC^{2t}$ by
        \[
        (x_{1},x_{2},\cdots, x_{2t} ) \mapsto (x_{1},\cdots, x_{2i-2}, -x_{2i},-x_{2i-1},x_{2i+1},\cdots, x_{2t}).
        \]
\end{itemize}
Note that $\sfH_{t}$ is also a subgroup of $\sfW'_{2t}$. Define a quadratic
character
\begin{equation}\label{eq:eta}
  \begin{array}{rcl}
    \eta : \sfH_{t}=  \sfW_{t}\ltimes \set{\pm 1}^{t}& \rightarrow & \set{\pm 1},\\
                                  (g,(a_{1},a_{2},\cdots, a_{t})) & \mapsto & a_{1}a_{2}\cdots a_{t}.
  \end{array}
\end{equation}

  Recall from \Cref{sec:SepCar} the classical groups $G_\mathrm b$ and $G_\mathrm g$, the  $\bZ^{n}$-coset
  $\Lambda= \Lambda_\mathrm b\times  \Lambda_\mathrm g$, the Weyl group $W_\mathrm b$ of $\g_\mathrm b$ and the Weyl group $W_\mathrm g$ of $\g_\mathrm g$. 
 
 The following Propositions \ref{count003} and \ref{count004}  also follow from Theorem \ref{thm:cohHC} by direct computation. 
  We omit the details. \cf  \cite[Applications]{Mc} and \cite{AC}. The following subsection gives an example in the case of type $B$ of  the structure data relevant for the computation. 
 
\begin{prop}\label{count003}
As representations of $W_\mathrm b$,
    \[
       \Coh_{\Lambda_\mathrm b}(\CK'(G_\mathrm b)) \cong
        \begin{cases}
  \bigoplus_{\substack{2t+c+d=n_\mathrm b}} \Ind_{\sfH_{t} \times \sfW_{c}\times \sfW_{d}}^{\sfW_{n_\mathrm b}}
         \eta\otimes 1\otimes 1, &  \text{if $\star \in \{B,  \wtC \}$}; \medskip\\
      \bigoplus_{\substack{2t+a=n_\mathrm b}} %\Res_{\sfW_{n}}^{\sfW'_{n}} \left(
          \Ind_{\sfH_{t} \times \sfS_a}^{\sfW_{n_\mathrm b}}\eta\otimes 1, &  \text{if $\star \in \{C, D \}$};\medskip \\
         %\Res_{\sfW_{n}}^{\sfW'_{n}} \left(
          \Ind_{\sfH_{\frac{n_\mathrm b}{2}}}^{\sfW'_{n_\mathrm b}}\eta, &  \text{if $\star \in \{C^*, D^* \}$},
     \end{cases}
       \]
       where in the case when $\star\in \{C,D\}$, the right-hand side space is viewed as a representation of $W_\mathrm b=\sfW_{n_\mathrm b}'$ by restriction.
      \end{prop}

Recall the quadratic character $\bsgn$ of $\sfW_n$ in \eqref{eq:bsgn}. 

\begin{prop}\label{count004}
As representations of $W_\mathrm g$,
    \[
       \Coh_{\Lambda_\mathrm g}(\CK'(G_\mathrm g)) \cong
        \begin{cases}
     \bigoplus_{\substack{0\leq p_\mathrm g-(2t+a+2r)\leq 1,\\0\leq q_\mathrm g - (2t+a+2s)\leq 1}} \Ind_{\sfH_{t} \times \sfS_{a}\times \sfW_s\times \sfW_r}^{\sfW_{n_\mathrm g}}
         \eta \otimes 1 \otimes \sgn \otimes \sgn, &  \text{if $\star=B$}; \medskip\\
     \bigoplus_{2t+a+c+d=n_\mathrm g}\Ind_{\sfH_{t} \times \sfS_{a} \times \sfW_c\times \sfW_d}^{W_{n_\mathrm g}}\eta \otimes
          \sgn \otimes 1 \otimes 1, &  \text{if $\star =C$};\medskip \\
          \bigoplus_{2t+a+a'=n_\mathrm g}\Ind_{\sfH_{t} \times \sfS_{a} \times \sfS_{a'}}^{\sfW_{n_\mathrm g}}\eta \otimes
          \sgn \otimes 1, &  \text{if $\star =\wtC$};\medskip \\
          \bigoplus_{\substack{(t+r,t+s)=(\frac{p_\mathrm g}{2},\frac{q_\mathrm g}{2})}} \Ind_{\sfH_{t} \times \sfW_s\times \sfW_r}^{\sfW_{n_\mathrm g}}
         \eta \otimes \sgn \otimes \sgn, &  \text{if $\star =C^*$}; \medskip\\
        \bigoplus_{\substack{2t+c+d+2r=p_\mathrm g\\2t+c+d+2s=q_\mathrm g}}
          \Ind_{\sfH_{t} \times \sfW_s\times \sfW_r\times \sfW'_{c}\times \sfW_{d} }^{\sfW_{n_\mathrm g}}\eta \otimes \bsgn \otimes \bsgn \otimes 1\otimes
          1, &  \text{if $\star =D$}; \medskip \\
          \bigoplus_{\substack{2t+a=n_\mathrm g}} \Ind_{\sfH_{t} \times \sfS_{a}}^{\sfW'_{n_\mathrm g}}
         \eta \otimes \sgn, &  \text{if $\star =D^*$}, \\
     \end{cases}
       \]
       where in the case when $\star=D$,  the right-hand side space is viewed as a representation of $W_\mathrm g=\sfW_{n_\mathrm g}'$ by restriction.
      \end{prop}

\subsection{An example: type $B$}
In this subsection, we give the relevant structure data leading to the results of \Cref{count003} and \Cref{count004} when $G = \SO(p,q)$ with $p+q=2n+1$.

We retain the notation of \Cref{sec:CW} and \Cref{sec:GB}. 
The set of the conjugate classes of Cartan subgroups is parameterized by 
\begin{equation}\label{eq:p.Cartan.B}
\set{ (t,a)\in \bN\times \bN | p - 2t - a \geq 0, q -2t-a \geq 0}.
\end{equation}
Fixing an element $(t, a)$ in \eqref{eq:p.Cartan.B}, a representative of the conjugate class is given by a Cartan subgroup  $H_{t,a}$ isomorphic to 
\[
    (\bC^\times)^t\times (\bR^\times)^a  \times \SO(2)^{n-2t-a} 
\]
with the set of real roots 
\[
\set{\pm(e_{2i-1}-e_{2i}) : 1\leq i \leq t}\cup 
\set{\pm e_{i} : 2t+1\leq i \leq 2t+a}
\]
and the set of imaginary roots
\[
\set{\pm(e_{2i-1}+e_{2i}) : 1\leq i \leq t}\cup 
\set{\pm e_{i} \pm e_j : 2t+a + 1\leq i<j\leq n}
\cup \set{e_{k} : 2t+a + 1\leq k\leq n}. 
\]

For a root $\alpha$, write $s_{\alpha}$ the corresponding reflection. Then the real Weyl group $W_{H_{t,a}}$ is given by 
\[
W_{H_{t,a}} =  \sfS_t \ltimes (W(A_1)^t \times W(A_1)^t)  \times \sfW_a\times \sfW_{s}\times \sfW_{r}, 
\]
where
 \begin{itemize}
\item $W(A_1)$ is the Weyl group of type $A_1$; 
\item $r = \floor{(p-2t-a)/2}$, and $s = \floor{(q-2t-a)/2}$  (note that $r+s=n-2t-a$);  
\item $\sfS_t$ is generated by $s_{e_{2i-1}-e_{2i+1}}s_{e_{2i}-e_{2i+2}}$ for $1\leq i <t$; \item the first $W(A_1)^t$ is generated by 
$s_{e_{2i-1}-e_{2i}}$ for $1\leq i \leq t$;
\item the second $W(A_1)^t$ is generated by 
$s_{e_{2i-1}+e_{2i}}$ for $1\leq i \leq t$;
\item $\sfW_a$, $\sfW_r$ and $\sfW_s$ are identified respectively with the subgroups of $\sfW_n$ acting on the subspaces spanned by $\set{e_i: 2t+1\leq i\leq 2t+a}$, 
$\set{e_i: 2t+a+1\leq i\leq 2t+a+r}$,
and
$\set{e_i: 2t+a+r+1\leq i\leq n}$.
\end{itemize}
Note that the imaginary roots have the same span as
\[
\set{e_1+e_2, e_3+e_4, \cdots, e_{2t-1}+e_{2t}, e_{2t+a+1},e_{2t+a+2},\cdots , e_{n}}. 
\]
The quadratic character $\sgn_{\mathrm {im}}$, when restricted to $W_{H_{t,a}}$, is 
\begin{itemize}
\item the trivial character on $\sfS_t$, the first $W(A_1)^t$ factor and  $\sfW_a$, 
\item the sign character on the second $W(A_1)^t$ factor, and $\sfW_s\times \sfW_r$. 
\end{itemize}

We identify $\sfH_t =
\sfW_t\ltimes \set{\pm 1}^t$  with 
$\sfS_t\ltimes (W(A_1)^t \times W(A_1)^t)$, viewed as a subgroup of $\sfW_{2t}$ as in \eqref{eq:Ht}. 

Let $\cP_{t,a}$ denote the set of parameters in $\cP_\Lambda(G)$ which are represented by $(H, \xi,\Gamma)\in \sP_\Lambda(G)$ with $H$ conjugate to $H_{t,a}$. (Recall 
  $\Lambda= \Lambda_\mathrm b\times  \Lambda_\mathrm g$ is the $\bZ^{n}$-coset defined in \eqref{eq:Lambda}.)
Under the cross action, the $W_\Lambda$-orbits of $\cP_{t,a}$ is parameterized by a pair of numbers $(c,d)$ in 
\[
\set{(c,d)\in \bN\times \bN | \text{$\nbb-c-d$ is a non-negative even integer and $c+d\leq a$}}.
\]
Write $t_1 = (\nbb-c-d)/2$ and $t_2 = t-t_1$. 
Set 
\[
    \nu := (\underbrace{0,\cdots, 0}_{\text{$\nbb$}},
    \underbrace{\half,\cdots, \half}_{\text{$\ngg$}})\in \Lambda_{\mathrm b}\times \Lambda_{\mathrm g} = \Lambda
\]
Then there is a unique  $W_\Lambda$-orbit $\cP_{t,a,c,d}$ in $\cP_\Lambda(G)$ represented by an element $\gamma_{t,a,c,d} =(H_{t,a},\xi, \Gamma)\in \sP_\Lambda(G)$ satisfying  the following conditions:
\begin{itemize}
\item the group $H_{t,a}$ has a decomposition 
$
H_{t,a}=H_{\mathrm b}\times H_{\mathrm g} $,
 with
 \[
H_{\mathrm b} = 
    (\bC^\times)^{t_1}\times (\bR^\times)^{c+d},\quad     H_{\mathrm g} = 
    (\bC^\times)^{t_2}\times (\bR^\times)^{a-c-d}  \times \SO(2)^{n-2t-a};
    \]
    \item 
    \[
    \xi(\hha_\mathrm b)=\h_\mathrm b\quad \textrm{  
    and } \quad \xi(\hha_\mathrm g)=\h_\mathrm g,
    \]
    where $\h_\mathrm b$ and $\h_\mathrm g$ are respectively the complexified Lie algebras of $H_\mathrm b$ and $H_\mathrm g$;
    %\item $\Gamma|_{\set{\pm 1}^{c+d}} = 1^{\otimes c} \otimes \sgn^{\otimes d}$,  
\item $\Gamma_\nu|_{\set{\pm 1}^{c+d}} = 
    \underbrace{1\otimes \cdots \otimes 1}_{c} \otimes
      \underbrace{\sgn \otimes \cdots \otimes \sgn}_{d}$,
      where $\set{\pm 1}^{c+d}$ is viewed as a subgroup of $H_{t,a}$ via the inclusions
    \[
    \set{\pm 1}^{c+d}\subset 
    (\bR^\times)^{c+d}\subset H_\mathrm b=(\bC^\times)^{t_1}\times (\bR^\times)^{c+d}\subset H_{t,a};
    \]
    \item  $\Gamma_\nu|_{\set{\pm 1}^{a-c-d}}$ is trivial, where $\set{\pm 1}^{a-c-d}$ is viewed as a subgroup of $H_{t,a}$ via the inclusions
    \[
    \set{\pm 1}^{a-c-d}\subset (\bR^\times)^{a-c-d}\subset H_\mathrm g=(\bC^\times)^{t_2}\times (\bR^\times)^{a-c-d}  \times \SO(2)^{n-2t-a}\subset H_{t,a}.
    \]
\end{itemize}

Recall that the real Weyl group $W_{H_{t,a}}$ may be identified as a subgroup of the abstract Weyl group via $\xi$. Then the cross stabilizer $W_{\gamma_{t,a,c,d}}$ of $\gamma_{t,a,c,d}$ is the following subgroup of $W_\Lambda $:  
\[
 (\sfH_{t_1}\times \sfW_c\times \sfW_d) \times (\sfH_{t_2}\times \sfS_{a-c-d} \times \sfW_s\times \sfW_r) \subset W_{\mathrm b}\times W_{\mathrm g},
\]
where
\begin{itemize}
\item $\sfH_{t_1}\times \sfH_{t_2}$ is the natural subgroup of $\sfH_t\subset \sfW_{2t}$,
\item $\sfW_c\times \sfW_d \times \sfS_{a-c-d}$ is the natural subgroup of $\sfW_c\times \sfW_d \times \sfW_{a-c-d}\subset \sfW_a$. 
\end{itemize}

With the above structure data one may deduce the two formulas on the coherent continuation representations stated in Section 8.1.

\subsection{The Lusztig  left cells}
\label{sec:LCBCD}

To ease the notation, for every sequence $a_1\geq a_2\geq \dots \geq a_k\geq 0$ ($k\geq 0$) of integers,   we let $[a_1, a_2, \cdots, a_k]_{\mathrm{col}}$ denote the Young diagram
whose $i$-th column has length $a_i$ if $1 \leq i \leq k$ and length $0$
otherwise. Likewise, we let $[a_1, a_2, \cdots, a_k]_{\mathrm{row}}$ denote the Young diagram
whose $i$-th  row has length $a_i$ if $1 \leq i \leq k$ and length $0$ otherwise. 

 As usual, we identify $\Irr(\sfW_{t})$ ($t\in \BN$) with the set of bipartitions $\tau =(\tau_{L},\tau_{R})$ of total size $t$ (\cite[Section 11.4]{Carter}). Here and henceforth, a bipartition means a pair of Young diagrams, and the total size refers to
$\abs{\tau_{L}}+\abs{\tau_{R}}$.
%Given the chosen embedding of $\sfS_{n}$ into $\sfW'_{n}$,
We also let $(\tau_{L},\tau_{R})_{I}\in \Irr(\sfW'_t)$ denote the  irreducible representation  given by
  \begin{itemize}
    \item the restriction of $(\tau_{L},\tau_{R})\in \Irr(\sfW_{t})$  if
    $\tau_{L}\neq \tau_{R}$, and
    \item
    the induced representation
    $\Ind_{\sfS_{t}}^{\sfW'_{t}} \tau_{L}$ if $\tau_{L}=\tau_{R}$.
  \end{itemize}
  Take an element $w\in \sfW_t$ such that $w \sfW_t'$ generates the group $\sfW_t/\sfW'_t$. Define   $(\tau_{L},\tau_{R})_{II}\in \Irr(\sfW'_t)$ to be the twist of $(\tau_{L},\tau_{R})_{I}$ by the conjugation by $w$, namely there is a linear isomorphism $\kappa: (\tau_{L},\tau_{R})_{I}\rightarrow (\tau_{L},\tau_{R})_{II}$ such that
  \[
    \kappa(g\cdot u)= (wgw^{-1})\cdot (\kappa(u)), \quad \textrm{for all }g\in \sfW'_t, \ u\in (\tau_{L},\tau_{R})_{I}.
  \]

  Note that
  \[
     (\tau_{L},\tau_{R})_{I}=(\tau_{R},\tau_{L})_{I}\quad \textrm{and}\quad (\tau_{L},\tau_{R})_{II}=(\tau_{R},\tau_{L})_{II}
  \]
in all cases, and  $(\tau_{L},\tau_{R})_{I}=(\tau_{R},\tau_{L})_{II}$ when     $\tau_{L}\neq \tau_{R}$.

  In the rest of this subsection, we describe the Lusztig left cell $\LC_{\lambda_{\ckcO}}$ $(\subset \Irr(W(\Lambda)))$ attached to $\lambda_{\ckcO}$ (\Cref{def:ll}). Here the integral Weyl group 
  $W(\Lambda)= W_{\mathrm b}\times W'_{\mathrm g}$, and is described in \Cref{sec:SepCar}. 

Recall that if  $\star\in \{C, C^*, D, D^*\}$, then $W_{\mathrm b}=\sfW'_{n_\mathrm b}$. In this case, we say that 
\be\label{assumi}
  \textrm{$\check \CO_\mathrm b$ has type I, if the number of negative entries of $\lambda_{\check \CO_\mathrm b}$ has the same parity as $\frac{n_\mathrm b}{2}$};
  \ee
  otherwise we say that $\check \CO_\mathrm b$ has type II.

  Recall also that if  $\star =\wtC$, then $W'_{\mathrm g}=\sfW'_{n_\mathrm g}$. In this case, we say that 
  \be\label{assumiCT}
  \textrm {$\check \CO_\mathrm g$ has type I if
 $\frac{n_\mathrm g}{2}$ is even, and 
 type II if
 $\frac{n_\mathrm g}{2}$ is odd.}
 \ee  
 
 Define two Young diagrams
 \begin{equation}\label{eq:taub1}
    \begin{split}
      \tau_{L,\mathrm b} := \begin{cases}
        \big[\half(\bfrr_{1}(\ckcO'_{\mathrm b})+1), \half(\bfrr_{2}(\ckcO'_{\mathrm b})+1), \cdots, \half(\bfrr_{c}(\ckcO'_{\mathrm b})+1)\big]_{\mathrm{col}},
               &\quad \text{if } \star \in \set{B,\wtC}; \\% \smallskip \\
         \big[\half\bfrr_{1}(\ckcO'_{\mathrm b}), \half\bfrr_{2}(\ckcO'_{\mathrm b}),\cdots, \half\bfrr_{c}(\ckcO'_{\mathrm b})\big]_{\mathrm{col}},
        &\quad  \text{if } \star \in \set{C,C^{*}, D,D^{*}},\\
      \end{cases}
    \end{split}
  \end{equation}
  and
   \begin{equation}\label{eq:taub2}
    \begin{split}
      \tau_{R,\mathrm b} := \begin{cases}
        \big(\half(\bfrr_{1}(\ckcO'_{\mathrm b})-1), \half(\bfrr_{2}(\ckcO'_{\mathrm b})-1), \cdots, \half(\bfrr_{c}(\ckcO'_{\mathrm b})-1)\big)_{\mathrm{col}},
               &\quad \text{if } \star \in \set{B,\wtC}; \\% \smallskip \\
         \big(\half\bfrr_{1}(\ckcO'_{\mathrm b}), \half\bfrr_{2}(\ckcO'_{\mathrm b}),\cdots, \half\bfrr_{c}(\ckcO'_{\mathrm b})\big)_{\mathrm{col}},
        &\quad  \text{if } \star \in \set{C,C^{*}, D,D^{*}},\\
      \end{cases}
    \end{split}
  \end{equation}
 where $c:= \bfcc_{1}(\ckcO'_{\mathrm b})$.

  Define an irreducible representation $\tau_{\mathrm b}\in \Irr(W_{\mathrm b})$ attached to $\ckcO_{\mathrm b}$ by
\begin{equation}\label{eq:taub}
 \tau_{\mathrm b}:= \left\{
     \begin{array}{ll}
       ( \tau_{L,\mathrm b}, \tau_{R,\mathrm b}), \quad
       & \text{if } \star \in \set{B,\wtC}; \medskip\\
         ( \tau_{L,\mathrm b}, \tau_{R,\mathrm b})_I, \quad & \text{if $\star \in \set{C,C^{*}, D,D^{*}}$ and $\check \CO_\mathrm b$ has type I;}\medskip\\
        ( \tau_{L,\mathrm b}, \tau_{R,\mathrm b})_{II}, \quad & \text{if $\star \in \set{C,C^{*}, D,D^{*}}$ and $\check \CO_\mathrm b$ has type II}.
\end{array}
  \right.
\end{equation}

Recall the set  $\CPPs(\ckcO_{\mathrm g})$ from   \Cref{defn:PP}.   Put
  \[
    {\mathrm A}(\ckcO) := {\mathrm A}(\ckcO_{\mathrm g}):= \textrm{the power set of $\CPPs(\ckcO_{\mathrm g})$},
    \]
    which is identified with the free $\bF_2$-vector space with free basis $\CPPs(\ckcO_{\mathrm g})$. Here $\bF_{2}:=\bZ/2\bZ$ is the field with two elements only.
Note that   $\{\emptyset, \CPPs(\ckcO_{\mathrm g})\}$ is a subgroup of ${\mathrm A}(\ckcO)$.  Define
   \begin{equation*}%\label{def:barA}
  \bar{\mathrm A}(\ckcO):= \bar{\mathrm A}(\ckcO_{\mathrm g}):=
  \begin{cases}
 {\mathrm A}(\ckcO_{\mathrm g})/\{\emptyset, \CPPs(\ckcO_{\mathrm g})\}, & \quad \text{if  } \star =\wtC;\\
 {\mathrm A}(\ckcO_{\mathrm g}),  & \quad \text{otherwise.}
  \end{cases}
  \end{equation*}

  Generalizing \eqref{ijo}, for each $\wp\in  {\mathrm A}(\ckcO)$,    we define a pair of Young diagrams \begin{equation}
      \label{def:ipjp}
(\imath_\wp, \jmath_\wp):=(\imath_\star(\check \CO, \wp), \jmath_\star(\check \CO, \wp))
\end{equation}
as in what follows.

If $\star=B$, then
 \[
   \mathbf c_{1}(\jmath_\wp)=\frac{\mathbf r_1(\check \CO_{\mathrm g})}{2},
\]
and for all $i\geq 1$,
\[
(\mathbf c_{i}(\imath_\wp), \mathbf c_{i+1}(\jmath_\wp))=
   \left\{
     \begin{array}{ll}
           (\frac{\mathbf r_{2i+1}(\check \CO_{\mathrm g})}{2},  \frac{\mathbf r_{2i}(\check \CO_{\mathrm g})}{2}), &\hbox{if $(2i, 2i+1)\in \wp$}; \smallskip\\
            (\frac{\mathbf r_{2i}(\check \CO_{\mathrm g})}{2},  \frac{\mathbf r_{2i+1}(\check \CO_{\mathrm g})}{2}), &\hbox{otherwise}.\\
            \end{array}
   \right.
\]

If $\star=\widetilde{C}$, then for all $i\geq 1$,
\[
(\mathbf c_{i}(\imath_\wp), \mathbf c_{i}(\jmath_\wp))=
   \left\{
     \begin{array}{ll}
           (\frac{\mathbf r_{2i}(\check \CO_{\mathrm g})}{2},  \frac{\mathbf r_{2i-1}(\check \CO_{\mathrm g})}{2}), &\hbox{if $(2i-1, 2i)\in \wp$}; \smallskip\\
            (\frac{\mathbf r_{2i-1}(\check \CO_{\mathrm g})}{2},  \frac{\mathbf r_{2i}(\check \CO_{\mathrm g})}{2}), &\hbox{otherwise}.\\
            \end{array}
   \right.
\]

If $\star\in\{C,C^*\}$, then for all $i\geq 1$,
\[
(\mathbf c_{i}(\jmath_\wp), \mathbf c_{i}(\imath_\wp))=
   \left\{
     \begin{array}{ll}
            (\frac{\mathbf r_{2i}(\check \CO_{\mathrm g})-1}{2},  \frac{\mathbf r_{2i-1}(\check \CO_{\mathrm g})+1}{2}), &\hbox{if $(2i-1, 2i)\in \wp$}; \smallskip\\
            (0,  0), &\hbox{if $(2i-1, 2i)$ is vacant in $\check \CO_{\mathrm g}$};\\
                    (\frac{\mathbf r_{2i-1}(\check \CO_{\mathrm g})-1}{2},  0), & \hbox{if $(2i-1, 2i)$ is tailed in $\check \CO_{\mathrm g}$};\smallskip\\
                  (\frac{\mathbf r_{2i-1}(\check \CO_{\mathrm g})-1}{2},  \frac{\mathbf r_{2i}(\check \CO_{\mathrm g})+1}{2}), &\hbox{otherwise}.\\
            \end{array}
   \right.
\]

If $\star\in\{D,D^*\}$, then
 \[
   \mathbf c_{1}(\imath_\wp)= \left\{
     \begin{array}{ll}
     0,  &\hbox{if $\mathbf r_1(\check \CO_{\mathrm g})=0$};\\
                 \frac{\mathbf r_1(\check \CO_{\mathrm g})+1}{2},   &\hbox{if $\mathbf r_1(\check \CO_{\mathrm g})>0$}, \smallskip\\
                   \end{array}
   \right.
 \]
and for all $i\geq 1$,
\[
(\mathbf c_{i}(\jmath_\wp), \mathbf c_{i+1}(\imath_\wp))=
   \left\{
     \begin{array}{ll}
            (\frac{\mathbf r_{2i+1}(\check \CO_{\mathrm g})-1}{2},  \frac{\mathbf r_{2i}(\check \CO_{\mathrm g})+1}{2}), &\hbox{if $(2i, 2i+1)\in \wp$}; \smallskip\\ 
            (0,  0), &\hbox{if $(2i, 2i+1)$ is vacant in $\check \CO_{\mathrm g}$};\\
                    (\frac{\mathbf r_{2i}(\check \CO_{\mathrm g})-1}{2},  0), & \hbox{if $(2i, 2i+1)$ is tailed in $\check \CO_{\mathrm g}$};\smallskip\\
         (\frac{\mathbf r_{2i}(\check \CO_{\mathrm g})-1}{2},  \frac{\mathbf r_{2i+1}(\check \CO_{\mathrm g})+1}{2}), &\hbox{otherwise}.\\
            \end{array}
   \right.
\]

In the notation above, we have
\[
(\imath_{\check O},  \jmath_{\check O})=(\imath_\star(\check \CO, \emptyset), \jmath_\star(\check \CO, \emptyset)). 
\]
Here $(\imath_{\check O},  \jmath_{\check O})$ is the pair of Young diagrams attached to $\check \CO$, as in \eqref{ijo}. 

We define an element $\tau_{\wp}\in \Irr(\Wg')$ by
  \begin{equation}\label{eq:tauwp}
    \tau_{\wp} :=
    \begin{cases}
      (\imathp,\jmathp),  \quad &  \text{if } \star \in \set{B,C, C^{*} }; \\
      (\imathp,\jmathp)_{I}= (\imathp,\jmathp)_{II},  \quad &  \text{if   $ \star \in \set{D,D^{*}}$};\\
       (\imathp,\jmathp)_I, \quad & \text{if $\star =\wtC$  and $\check \CO_\mathrm g$ has type I;}\medskip \\
            (\imathp,\jmathp)_{II}, \quad & \text{if $\star =\wtC$ and  $ \check \CO_\mathrm g$ has type II.}
    \end{cases}
  \end{equation}
  Note that  if $\star=\wtC$, then $\tau_{\wp} = \tau_{\wp^{c}}$, where $\wp^{c}$ is the complement of $\wp$ in $\CPPs(\ckcO_{\mathrm g})$.
Therefore in all cases,   $\tau_{\bar \wp}\in \Irr(\Wg')$ is defined for every $\bar \wp\in \bar{\mathrm A}(\check \CO)$.

Recall the Lusztig left cell $\LC_{\ckcO}$ attached to $\ckcO$, which is the set of all $\sigma\in \Irr(W(\Lambda))$  that occurs in the multiplicity free representation
   \[
    \left(J_{\Wlamck}^{W(\Lambda)} \sgn\right) \otimes \sgn.
  \]

\begin{prop}[\cf Barbasch-Vogan {\cite{BVUni}*{Proposition~5.28}}]
    \label{lem:Lcell}
    The
  map    \[
      \begin{array}{rcl}
        \bar{\mathrm A}(\ckcO) & \rightarrow & \LC_{\ckcO},\\
                       \bar \wp & \mapsto &  \tau_{\mathrm b} \otimes \tau_{\bar \wp}
      \end{array}
    \]
    is well-defined and bijective, and
        \[
      \tau_{\ckcO}:=\tau_{\mathrm b}\otimes \tau_{\emptyset}
    \] is the unique special representation in $\LC_{\ckcO}$.
 Moreover, 
 \begin{equation}\label{eq:dBV.W}
 \textrm{ $j_{W(\Lambda)}^{W}(\tau_{\ckcO})$ corresponds to $\dBV(\ckcO)$ under the Springer correspondence,}
 \end{equation}
    and
  \begin{equation}\label{eq:dBV.W2}
     \dBV(\ckcO)=   (\ckcOpb)^{t} \cupcol (\ckcOpb)^{t}\cupcol \dBV(\ckcOg)
   \end{equation}
as Young diagrams.
 \end{prop}

  \begin{proof}
    When  $\ckcO=\ckcO_{\mathrm g}$ and $\star\neq \wtC$, the proposition is proved in \cite{BVUni}*{Proposition~5.28}.  When $\star\neq \wtC$,
    \eqref{eq:dBV.W} and \eqref{eq:dBV.W2} are proved in \cite{BVUni}*{Proposition~A2}. In general, the proposition follows from
    %from an induction on number of columns of $\ckcO_{\mathrm g}$ using
    Lusztig's formula of $J$-induction in \cite{Lu}*{\S 4.4-4.6}, as well as the explicit form of the Springer correspondence \cite{Sho}. 
    \end{proof}

Recall that $\Wg' =\sfW'_{n_{\mathrm g}}$ when $\star \in \set{\wtC,D,D^{*}}$. Since the representation theory of $\sfW_{n}$ is more elementary than that of
$\sfW'_{n}$, we prefer to work with $\sfW_n$ instead of $\sfW_n'$ in some situations. For this reason, 
we also define for all cases
\begin{equation}\label{eq:ttauwp}
\wttau_{\wp} = (\imath_{\wp},\jmath_{\wp})\in \Irr(\sfW_{n_\mathrm g}), \qquad \wp\in {\mathrm A}(\ckcO).
\end{equation}

For later use, we record the following lemma, which follows immediately from our explicit descriptions of  $\tau_{\mathrm b}$ and $\tau_{\wp}$.

\begin{lem}\label{lem:WLcell}
Let $\wp\in {\mathrm A}(\ckcO)$. If $\star=\wtC$, then
  \[
    \Ind_{\sfW_{n_{\mathrm b}}\times \sfW'_{n_{\mathrm g}}}^{\sfW_{n_{\mathrm b}}\times \sfW_{n_{\mathrm g}}} \tau_{\mathrm b}\otimes\tau_{\wp}  \cong
    \begin{cases}
       \tau_{\mathrm b}\otimes \wttau_{\emptyset}, &\quad  \text{if } \CPPs(\ckcOg)=\emptyset ; \\%n_{\mathrm g}=0; \\
      (\tau_{\mathrm b}\otimes \wttau_{\wp}) \oplus ( \tau_{\mathrm b}\otimes \wttau_{\wp^{c}}),
      &\quad \text{otherwise}.
    \end{cases}
  \]
If $\star\in\set{D,D^{*}}$, then
  \[
    \Ind_{\sfW'_{n_{\mathrm b}}\times \sfW'_{n_{\mathrm g}}}^{\sfW_{n_{\mathrm b}}\times \sfW_{n_{\mathrm g}}} \tau_{\mathrm b}\otimes \tau_{\wp} \cong
    \begin{cases}
      \wttau_{\mathrm b}, & \quad \text{if } n_{\mathrm g}=0; \\
      ( \wttau_{\mathrm b}\otimes \wttau_{\wp})  \oplus ( \wttau_{\mathrm b}\otimes \wttau_{\wp}^{\epsilon}),
      &\quad \text{otherwise}.
    \end{cases}
  \]
  Here $\wttau_{\mathrm b} = \Ind_{\sfW'_{n_{\mathrm b}}}^{\sfW_{n_{\mathrm b}}}\tau_{\mathrm b}$ and $\wttau_{\wp}^{\epsilon}:= \wttau_{\wp}\otimes \varepsilon$ (recall the quadratic character $\varepsilon$ from \eqref{defep}).
\end{lem}

\subsection{From coherent continuation representation to counting}
\label{subsec:counting}

We have defined in \eqref{defpbp2222} the set $\PBPs(\ckcO)$ when $\ckcO$ has good parity. Similarly, we make the following definition in the bad parity case.
\begin{defn}
  Let $\PBP^*(\ckcOb)$ be the set of all triples
  $\uptau = (\imath,\cP)\times(\jmath,\cQ)\times \star $ where $(\imath,\cP)$ and
  $(\jmath,\cQ)$ are painted Young diagrams such that
  \begin{itemize}
    \item $(\imath,\jmath) = (\tau_{L,\mathrm b},\tau_{R,\mathrm b})$ (see \eqref{eq:taub1} and \eqref{eq:taub2});
    \item the symbols of $\cP$ are in
          \[
          \begin{cases}
            \set{\bullet, c,d},  & \text{if } \star\in \set{B,\wtC}; \\
            \set{\bullet, d},  & \text{if } \star\in \set{C,D};\\
            \set{\bullet},  & \text{if } \star\in \set{C^{*},D^{*}};\\
          \end{cases}
          \]
    \item the symbols of $\cQ$ are in 

          \[
          \begin{cases}
            \set{\bullet},  & \text{if } \star\in \set{B,\wtC, C^{*},D^{*}};\\
            \set{\bullet, c},  & \text{if } \star\in \set{C,D}.
            \\
          \end{cases}
          \]
  \end{itemize}
\end{defn}

 We introduce some additional notation. For each $\check \CO$ with good parity, and $\wp \subseteq \CPP(\ckcO)$, let
\begin{equation}\label{eq:PBPSP}
  \PBP_{\star}(\check \CO, \wp) := \Set{ \uptau \text{ is a painted bipartition }\mid  \star_{\uptau}=\star, (\imath_{\uptau},\jmath_{\uptau}) = (\imath_{\wp}, \jmath_{\wp})}
\end{equation}
and
\begin{equation}\label{eq:PBPGP}
  \PBP_{G}(\check \CO, \wp) := \Set{\uptau \text{ is a painted bipartition }\mid G_{\uptau}=G, (\imath_{\uptau},\jmath_{\uptau}) = (\imath_{\wp}, \jmath_{\wp})}.
  % \uptau=(\imath, \cP)\times (\jmath,\cP)\times \alpha|}
\end{equation}
Note that in the notation above, we have $\PBP_{\star}(\check \CO) =\PBP_{\star}(\check \CO, \emptyset)$ and
$\PBP_{G}(\check \CO) =\PBP_{G}(\check \CO, \emptyset)$.  

Put
\[
  \tPBP_{\star}(\ckcO) :=
  \bigsqcup_{\wp \subseteq \CPP(\ckcO)}\PBP_{\star}(\check \CO, \wp)
\]
and
\[
  \tPBP_{G}(\check \CO):=   \bigsqcup_{\wp \subseteq \CPP(\ckcO)}\PBP_{G}(\check \CO, \wp). 
 \]

Recall the notion of $\check \CO$ being $G$-relevant (\Cref{secrgp0}). Note that if  $\star =D^{*}$, then in view of the Remarks \ref{rm63} (c), $\check \CO$ is not $G$-relevant if and only if $\check \CO$ has bad parity and type II (see \eqref{assumi}). 

\begin{prop}\label{prop:countBCD}
 If  $\star =D^{*}$ and $\check \CO$ is not $G$-relevant, then
  \[
  \sum_{\bar \wp\in \bar{\mathrm A}(\ckcO)} [\tau_{\mathrm b}\otimes \tau_{\bar \wp}: \Coh_{\Lambda} (\CK'(G))]=0.
  \]
In all other cases,   
 \[
 \sum_{\bar \wp\in \bar{\mathrm A}(\ckcO)} [\tau_{\mathrm b}\otimes \tau_{\bar \wp}: \Coh_{\Lambda} (\CK'(G))]=
           \sharp (\PBP^{\star}(\ckcO_{\mathrm b}))\cdot \sharp(\tPBP_{G_{\mathrm g}}(\ckcO_{\mathrm g})).
           \]
 \end{prop}

 \begin{proof}
   We use the following
  formulas (\cite{Mc}*{p220 (6)})  to compute the multiplicities:
  \begin{equation}\label{eq:indSW}
    \begin{split}
      \Ind_{\sfH_{t}}^{\sfW_{2t}} \eta& \cong   \bigoplus_{\sigma\in \Irr(\sfS_{t})} (\sigma,\sigma),\\
      \Ind_{\sfH_{t}}^{\sfW'_{2t}} \eta &\cong  \bigoplus_{\sigma\in \Irr(\sfS_{t})} (\sigma,\sigma)_{I},\\
      \Ind_{\sfS_{t}}^{\sfW_{t}}\sgn &\cong \bigoplus_{s+r=t}\Ind_{\sfW_{s}\times \sfW_{\mathrm r}}^{\sfW_{t}} \bsgn\otimes \sgn \cong \bigoplus_{s+r=t} ([s]_{\mathrm{col}},[r]_{\mathrm{col}}),\\
      \Ind_{\sfS_{t}}^{\sfW_{t}} 1 &\cong \bigoplus_{c+d=t}\Ind_{\sfW_{c}\times \sfW_{\mathrm d}}^{\sfW_{t}} 1 \otimes \epsilon \cong \bigoplus_{c+d=t} ([c]_{\mathrm{row}},[d]_{\mathrm{row}}),
    \end{split}
\end{equation}
where $t\in \bN$.

We skip the details when $\star \in \set{B,\wtC, C,D,C^{*}}$, and present the computation for $\star = D^{*}$, which is the most complicated case (in certain aspect).
 Suppose that $\star=D^*$ so that $G=\SO^*(2n)$.  If $\check \CO$ has bad parity, then
  \begin{eqnarray*}
 && \sum_{\bar \wp\in \bar{\mathrm A}(\ckcO)} [\tau_{\mathrm b}\otimes \tau_{\bar \wp}: \Coh_{\Lambda} (\CK'(G))] \medskip \\
 &=&[\tau_{\mathrm b} : \Coh_{\Lambda} (\CK'(G))] \medskip \\
 &=& [\tau_\mathrm b:\Ind_{\sfH_{t}}^{\sfW'_{2t}} \eta] \medskip \\
& =& \begin{cases}
    1=\sharp (\PBP^{\star}(\ckcO_{\mathrm b}))\cdot \sharp(\tPBP_{G_{\mathrm g}}(\ckcO_{\mathrm g})),  \qquad &  \text{if $\check \CO$ is $G$-relevant;}\medskip \\
          0, \qquad & \text{if $\check \CO$ is not $G$-relevant.}
       \end{cases}
  \end{eqnarray*}

Now we assume that $\check \CO$ does not have bad parity so that $n_\mathrm g>0$.
 Put
 \[
   \sfW'_{\nbb, \ngg}:=\sfW'_n\cap (\sfW_{\nbb}\times \sfW_{\ngg})\supseteq \sfW'_{\nbb}\times \sfW'_{\ngg}= W_{\mathrm b}\times W_{\mathrm g}.
 \]

Recall that
  \[
    \wttau_{\mathrm b} := \Ind_{\sfW'_{\nbb}}^{\sfW_{\nbb}} \tau_{\mathrm b} = (\cOpb,\cOpb)\in \Irr(\sfW_{\nnb})
    \]
    and
    \[
    \wttau_{\wp}: = (\imath_{\wp},\jmath_{\wp})\in \Irr(\sfW_{\nng}) \quad \textrm{ for all } \wp \subseteq \CPP(\ckcOg).
  \]
  Note that $\imath_{\wp}\neq \jmath_{\wp}$ since
  $\bfcc_{1}(\imath_{\wp})> \bfcc_{1}(\jmath_{\wp})$, which implies that
  \begin{equation*}%\label{eq:W''}
    \Ind_{\sfW'_{\nbb}\times \sfW'_{\ngg}}^{\sfW_{n_{\mathrm b},n_{\mathrm g}}'} \tau_{\mathrm b}\otimes \tau_{\wp}
    \cong (\wttau_{\mathrm b}\otimes \wttau_{\wp})|_{\sfW_{n_{\mathrm b},n_{\mathrm g}}'}.
  \end{equation*}

  For ease of notation, write $\sfW'':=\sfW_{n_{\mathrm b},n_{\mathrm g}}'$. Put
  \[
    \cC_{\mathrm b}:=   \Ind_{\sfH_{\frac{n_\mathrm b}{2}}}^{\sfW_{n_\mathrm b}}\eta \quad\textrm{and}\quad  \cC_\mathrm g:=  \bigoplus_{\substack{2t+a=n_\mathrm g}} \Ind_{\sfH_{t} \times \sfS_{a}}^{\sfW'_{n_\mathrm g}}
         \eta \otimes \sgn,
           \]
           where $\cC_\mathrm b$ is viewed as a representation of $\sfW_{n_\mathrm b}'$ by restriction.
   For every finite group $E$ and any  two finite-dimensional  representations $V_1$ and $V_2$ of $E$, put
  \[
    [V_1, V_2]_E:=\dim \Hom_E(V_1, V_2).
  \]
   For each $\wp\in \mathrm A (\ckcOg)$, we have that
  \[
    \begin{split}
    [\tau_{\mathrm b}\otimes \tau_{\wp}: \Coh_{\Lambda} (\CK'(G))]
      = & [\tau_{\mathrm b}\otimes \tau_{\wp} :
      \Ind_{\sfW'_{n_{\mathrm b}}\times \sfW'_{n_{\mathrm g}}}^{\sfW''} \cC_{\mathrm b} \otimes \cC_{\mathrm g}]_{\sfW'_{n_{\mathrm b}}\times \sfW'_{n_{\mathrm g}}}\\
      = & [\Ind_{\sfW'_{n_{\mathrm b}}\times \sfW'_{n_{\mathrm g}}}^{\sfW''} \tau_{\mathrm b}\otimes \tau_{\wp} :
      \Ind_{\sfW'_{n_{\mathrm b}}\times \sfW'_{n_{\mathrm g}}}^{\sfW''} \cC_{\mathrm b} \otimes \cC_{\mathrm g}]_{\sfW''}\\
      = & [(\wttau_{\mathrm b}\otimes \wttau_{\wp})|_{\sfW''}:
      \Ind_{\sfW'_{n_{\mathrm b}}\times \sfW'_{n_{\mathrm g}}}^{\sfW''} \cC_{\mathrm b} \otimes \cC_{\mathrm g}]_{\sfW''}\\
      = & [\wttau_{\mathrm b}\otimes \wttau_{\wp}:
      \Ind_{\sfW'_{n_{\mathrm b}}\times \sfW'_{n_{\mathrm g}}}^{\sfW_{n_{\mathrm b}}\times \sfW_{n_{\mathrm g}}} \cC_{\mathrm b} \otimes \cC_{\mathrm g}]_{\sfW_{n_{\mathrm b}}\times \sfW_{n_{\mathrm g}}}\\
       = & [\wttau_{\mathrm b}:\Ind_{\sfW'_{n_{\mathrm b}}}^{\sfW_{n_{\mathrm b}}} \cC_{\mathrm b}]_{\sfW_{n_{\mathrm b}}}\cdot
           [\wttau_{\wp}:\Ind_{\sfW'_{n_{\mathrm g}}}^{\sfW_{n_{\mathrm g}}} \cC_{\mathrm g}]_{\sfW_{n_{\mathrm g}}}\\
      =& \sharp(\PBP^{\star}(\ckcOb)))\cdot \sharp(\PBP_{\Gg}(\check \CO_{\mathrm g}, \wp)). 
    \end{split}
  \]
  The last equality follows from induction in stages, the branching rules
  in \eqref{eq:indSW}, specifically $      \Ind_{\sfH_{t}}^{\sfW_{2t}} \eta \cong   \bigoplus_{\sigma\in \Irr(\sfS_{t})} (\sigma,\sigma)$ 
  and $      \Ind_{\sfS_{t}}^{\sfW_{t}}\sgn \cong \bigoplus_{s+r=t} ([s]_{\mathrm{col}},[r]_{\mathrm{col}})$, and
  the Pieri rule (\cite[Corollary 9.2.4]{GW}) (which is a special case of the Littlewood-Richard rule). For the case at hand, note that a painted bipartition 
  $\uptau=(\imath, \CP)\times (\jmath, \cQ)\times D^*$ in $\PBP_{\Gg}(\check \CO_{\mathrm g}, \wp)$ has symbols $\bullet, s$ in $\CP$ and symbols $\bullet, r$ in $\cQ$. In applying the branching rules, the  $\sigma $ in the summation $ \bigoplus_{\sigma\in \Irr(\sfS_{t})} (\sigma,\sigma)$ corresponds to filling an identical set of boxes in $\CP$ and $\cQ$ with the symbol $\bullet$, and the $[s]_{\mathrm{col}}$ (resp. $[r]_{\mathrm{col}}$) in the summation $\bigoplus_{s+r=t} ([s]_{\mathrm{col}},[r]_{\mathrm{col}})$ corresponds to filling a column consisting of symbols $s$ in $\CP$ (resp. a column consisting of symbols $r$ in $\cQ$), with each possible way of constructing a bipartition $\uptau $ in $\PBP_{\Gg}(\check \CO_{\mathrm g}, \wp)$ 
  contributing $1$ to the multiplicity
  $[\wttau_{\wp}:\Ind_{\sfW'_{n_{
  \mathrm g}}}^{\sfW_{n_{\mathrm g}}} \cC_{\mathrm g}]_{\sfW_{n_{\mathrm g}}}$. Likewise each possible way of constructing a bipartition in 
$\PBP^{\star}(\ckcOb)$ contributes $1$ to the multiplicity $[\wttau_{\mathrm b}:\Ind_{\sfW'_{n_{\mathrm b}}}^{\sfW_{n_{\mathrm b}}} \cC_{\mathrm b}]_{\sfW_{n_{\mathrm b}}}$.
  \end{proof}

Now the equality \eqref{boundc22}, 
\Cref{lem:Lcell} and \Cref{prop:countBCD} imply the following corollary.
\begin{cor}\label{prop:countBCD22}
 The following equality holds: 
  \begin{eqnarray*}
 && \sharp(\Unip_{\ckcO}(G))\\
  & = & \begin{cases}
    0,  \qquad &  \text{if $\star=D^*$ and $\check \CO$ is not $G$-relevant;}\medskip \\
          \sharp (\PBP^{\star}(\ckcO_{\mathrm b}))\cdot \sharp(\tPBP_{G_{\mathrm g}}(\ckcO_{\mathrm g})), \qquad & \text{otherwise.}
       \end{cases}
  \end{eqnarray*}
 \end{cor}

When $\ckcO $ has good parity, we will see from \Cref{prop:PBP1} and \Cref{prop:PBP2} that
  \[
    \sharp(\tPBP_{G}(\ckcO)) =
    \left\{
    \begin{array}{ll}
       \sharp (\PBP_{G}(\ckcO)),  & \hbox{if $\star\in \{C^*,D^*\}$}; \smallskip\\
       2^{\sharp(\CPPs(\check \CO))} \cdot \sharp (\PBP_{G}(\ckcO)),  &\hbox{if $\star\in \{B, C,D,\widetilde {C}\}$}.
    \end{array}
  \right.
  \]
 \Cref{prop:countBCD22} will thus imply Theorem \ref{countup}.

\section{Special unipotent representations in type $BCD$:  reduction to good parity}\label{sec:red}
\label{sec:proof.BCDred}

The goal of this section is to prove \Cref{reduction}.
\Cref{prop:countBCD22} implies that $\Unip_{\check \CO}(G)$ is empty if  $\check \CO$ is not $G$-relevant (this notion is defined in \Cref{secrgp0}). Thus we further assume that $\check \CO$ is $G$-relevant. With our earlier assumptions \eqref{nonemp0} and \eqref{nonemp00}, this is equivalent to saying that $\check \CO$ has type I when $\star=D^*$ and $\check \CO$ has bad parity.

If $\star=C^*$, or $\star=D^*$ and $\check \CO $ does not have bad parity, by possibly changing $\omega_{\check V_\mathrm b}$ and $\omega_{\check V_\mathrm g}$ defined in \eqref{omega12} to their negatives, we
assume without loss of generality that  ${\check \CO_\mathrm b}$ has type I.

\subsection{Separating bad parity and good parity for special unipotent representations}

By Proposition \ref{propKL33},  we have an injective linear map
\[
 \varphi_{\lambda_{\check \CO}} : \CK'_{\lambda_{\check \CO_\mathrm b}}(G_{\mathrm b})\otimes  \CK'_{\lambda_{\check \CO_\mathrm g}}(G_{\mathrm g})\rightarrow\CK'_{\lambda_{\check \CO}}(G)
\]
and an injective map
\be\label{injirr}
   \varphi_{\lambda_{\check \CO}} : \Irr'_{\lambda_{\check \CO_\mathrm b}}(G_\mathrm b)\times  \Irr'_{\lambda_{\check \CO_\mathrm g}}(G_\mathrm g)\rightarrow  \Irr'_{\lambda_{\check \CO}}(G).
\ee

\begin{prop}\label{propKL333}
The map  \eqref{injirr} restricts to a bijective map
\[
\varphi_{\lambda_{\check \CO}}:
\Unip_{\ckcO_{\mathrm b}}(G_{\mathrm b})\times \Unip_{\ckcO_{\mathrm g}}(G_{\mathrm g})
\rightarrow \Unip_{\ckcO}(G).
\]

\end{prop}
\begin{proof}
Let $\pi_\mathrm b\in  \Irr'_{\lambda_{\check \CO_\mathrm b}}(G_\mathrm b)$ and $\pi_\mathrm g\in  \Irr'_{\lambda_{\check \CO_\mathrm g}}(G_\mathrm g)$.  Pick a basal  element $\Psi_\mathrm b$ of  $\Coh_{ \Lambda_\mathrm b}(\CK'(G_\mathrm b))$ such that $\Psi_\mathrm b(\lambda_{\check \CO_\mathrm b})=\pi_\mathrm b$. Likewise pick a basal  element $\Psi_\mathrm g$ of $ \Coh_{ \Lambda_\mathrm b}(\CK'(G_\mathrm g))$ such that $\Psi_\mathrm g(\lambda_{\check \CO_\mathrm g})=\pi_\mathrm g$.
Write $\Psi:=\varphi(\Psi_\mathrm b\otimes \Psi_\mathrm g)$. Denote by $\CC$ the Harish-Chandra cell in $ \Coh_{\Lambda}(\CK'(G))$ containing $\Psi$, and similarly define the Harish-Chandra cells $\CC_\mathrm g\ni \Psi_{\mathrm g}$ and $\CC_\mathrm b\ni \Psi_{\mathrm b}$.

Put $\pi:= \varphi_{\lambda_{\check \CO}}(\pi_\mathrm b, \pi_\mathrm g)=\Psi(\lambda_{\check \CO})$. 
Recall the integral Weyl group $W(\Lambda)=W_\mathrm b\times W'_\mathrm g$. Then
 \begin{eqnarray*}
   && \pi\in  \Unip_{\ckcO}(G)\\
    &\Longleftrightarrow& \sigma_\CC\cong  \left( j_{W_{\lambda_\ckcO}}^{W(\Lambda)} \sgn\right)\otimes \sgn \quad (\textrm{by Proposition \ref{hcass222}})\\
   &\Longleftrightarrow& \sigma_{\CC_\mathrm b}\cong  \left( j_{W_{\mathrm b, \lambda_{\ckcO_\mathrm b}}}^{W_\mathrm b} \sgn\right)\otimes \sgn \ \ \textrm{and}\ \  \sigma_{\CC_\mathrm g}\cong  \left( j_{W'_{\mathrm g, \lambda_{\ckcO_\mathrm g}}}^{W'_\mathrm g} \otimes \sgn\right)\otimes \sgn\\
 &\Longleftrightarrow& \pi_\mathrm b\in  \Unip_{\ckcO_{\mathrm b}}(G_{\mathrm b}) \ \ \textrm{and}\ \  \pi_\mathrm g\in  \Unip_{\ckcO_{\mathrm g}}(G_{\mathrm g}) \quad (\textrm{by Proposition \ref{hcass222}}).
 \end{eqnarray*}
 Here $W_{\mathrm b, \lambda_{\ckcO_\mathrm b}}$ denotes the stabilizer of $\lambda_{\ckcO_\mathrm b}$ in $W_{\mathrm b}$, and likewise $W'_{\mathrm g, \lambda_{\ckcO_\mathrm g}}$ denotes the stabilizer of $\lambda_{\ckcO_\mathrm g}$ in $W'_{\mathrm g}$.

In particular, we have proved that 
\[
\varphi_{\lambda_{\check \CO}}(\Unip_{\ckcO_{\mathrm b}}(G_{\mathrm b})\times \Unip_{\ckcO_{\mathrm g}}(G_{\mathrm g}))\subseteq 
\Unip_{\ckcO}(G).
\]
On the other hand, \Cref{prop:countBCD22} implies that
\[
\sharp(\Unip_{\ckcO_{\mathrm b}}(G_{\mathrm b})\times \Unip_{\ckcO_{\mathrm g}}(G_{\mathrm g}))=\sharp
\Unip_{\ckcO}(G).
\]
This proves the proposition since the map  \eqref{injirr} is injective. 
\end{proof}

\subsection{The case of bad parity}

Recall from \eqref{Gpb} the group
\[
  G'_{\mathrm b} := \begin{cases}
    \GL_{n_{\mathrm b}}(\bR), & \text{if } \star \in \set{B,C,D}; \\
    \widetilde{\GL}_{n_{\mathrm b}}(\bR), & \text{if } \star = \wtC; \\
    \GL_{\frac{n_{\mathrm b}}{2}}(\bH), & \text{if } \star \in \set{C^{*},D^{*}}.\\
  \end{cases}
\]

\begin{lem}\label{prop:BP.PP} 
  The equalities
\[
    \begin{split}
      \sharp(\PBP^{\star}(\ckcO_{\mathrm b})) = \sharp(\PAP_{\star '}(\ckcO'_{\mathrm b})) = \sharp(\Unip_{\ckcO'_{\mathrm b}}(G'_{\mathrm b}))
    \end{split}
  \]
  hold, where
\begin{equation}\label{def:star'}
\star ':= \begin{cases}
    A^{\bR}, & \text{if } \star \in \set{B,C,\wtC,D}; \\
    A^{\bH}, & \text{if } \star \in \set{C^{*},D^{*}}.\\
  \end{cases}
  \end{equation}
\end{lem}

\begin{proof} Suppose that $\star \in \Set{C^{*},D^{*}}$. Then
\[
    \begin{split}
      \sharp(\PBP^{\star}(\ckcO_{\mathrm b})= \sharp(\PAP_{A^{\bH}}(\ckcO'_{\mathrm b}))= 1.
    \end{split}
\]
  Suppose that $\star \in \Set{B,C,\wtC,D}$.
     It is easy to see that we have a bijection
  \[
    \begin{array}{ccc}
      \PBP^{\star}(\ckcO_{\mathrm b}) &  \rightarrow & \PAP_{A^{\bR}}(\ckcO'_{\mathrm b}),\\
      (\tau_{L,\mathrm b},\cP)\times (\tau_{R,\mathrm b},\cQ)\times \star & \mapsto & \cP',
    \end{array}
  \]
  where $\cP'$ is defined by the condition that
  \[
    \cP(\bfcc_{j}(\tau_{L,\mathrm b}),j)=d \Longleftrightarrow \cP'(\bfrr_{j}(\ckcO'_{\mathrm b}),j)=d, \quad \textrm{ for all } j=1,2,\cdots, \bfcc_{1}(\ckcO'_{\mathrm b}).
  \]
  The last equality is in Theorem \ref{thm:mainR00} (which is proved in
  \Cref{sec:GLRH11}).  \end{proof}

Let $P_{\mathrm b}$ be a parabolic subgroup of $G_{\mathrm b}$ that is  $\check \CO_\mathrm b$-relevant as defined in Section \ref{secrgp0}. Then
  $G'_{\mathrm b}$ is naturally isomorphic to the Levi quotient of $P_\mathrm b$.

\def\fIb{\fI_{\mathrm b}}
\def\pib{\pi_{\mathrm b}}

\begin{prop}\label{lem:Unip.BP}
For every $\pi'\in \Unip_{\ckcOpb}(\Gpb)$, the normalized induced representation
$\Ind_{\Pb}^{\Gb}\pi'$ is irreducible and belongs to $\Unip_{\ckcOb}(\Gb)$.
Moreover, the map
  \be\label{badind}
    \begin{array}{rcl}
      \Unip_{\ckcO'_{\mathrm b}}(G'_{\mathrm b}) & \rightarrow & \Unip_{\ckcOb}(\Gb), \\
      \pi' & \mapsto & % \pib:=
                        \Ind_{\Pb}^{\Gb}\pi'
    \end{array}
  \ee
  is bijective.
\end{prop}

\begin{proof}
    By the construction of $\pi'$  (see \Cref{sec:GLRH}) and
  a result of Barbasch \cite[Corollary 5.0.10]{B.Orbit}, the
  wavefront cycle of $\Ind_{\Pb}^{\Gb}\pi'$ is a single orbit $\sO$ with
  multiplicity one and its complexification is $\dBV(\ckcOb)$.
  Note that every irreducible summand of  $\Ind_{\Pb}^{\Gb}\pi'$ belongs
  to $\Unip_{\ckcOb}(\Gb)$. Hence $\Ind_{\Pb}^{\Gb}\pi'$ has to be irreducible.

    We now suppose that $\star=D$ so that $\Gb=\SO(\nnb, \nnb)$. We will prove the injectivity of the map \eqref{badind}.
    Fix a split Cartan subgroup
    $$H=(\bR^\times)^{\nnb}\subseteq \Gpb\subseteq \Gb $$ and write $H = MA$ where
    $M = \set{\pm 1}^{\nnb}$ is the compact part of $H$ and
    $A = (\bR^\times_+)^{\nnb}$ is the split part of $H$. We identify a
    character of $A$ with a vector in $\bC^{\nnb}$ as usual. Let
    \[
    K:=\{(g_1, g_2)\in \rO(\nnb)\times \rO(\nnb)\mid \det(g_1)\cdot \det(g_2)=1\}
    \]
    be a maximal compact subgroup of $\Gb$, and let
    $B$ be a Borel subgroup of $\Gb$ containing $H$ and the unipotent radical of $\Pb$.

 For each integer $l$ such that $0\leq l \leq \nnb$, put
    \[
      \delta_{ l} = \underbrace{1_{\pm }\otimes \cdots \otimes 1_{\pm }}_{l} \otimes
      \underbrace{\sgn_{\pm } \otimes \cdots \otimes \sgn_{\pm }}_{\nnb-l}\in \Irr(M),
    \]
    where $1_{\pm }$ denotes the trivial character of $\{\pm 1\}$ and $\sgn_{\pm }$ denotes the non-trivial character of $\{\pm 1\}$.
It is a fine $M$-type (\cite{Vg}*{Definition~4.3.8}). Let $\lambda_{l}$ be
    the restriction to $K$ of the irreducible
    $\rO(\nnb)\times \rO(\nnb)$-representation
    $\wedge^{\nnb-l}\,\bC^{\nnb}\otimes \wedge^{0}\,\bC^{\nnb}$. Then $\lambda_{l}$ is irreducible and is a
    fine $K$-type, see \cite{BGG.M}*{\S 6}.

    \def\cusp{\fC}

    Vogan proves that there is a well-defined injective map
\[
      \begin{array}{lccc}
       \fX\colon &   W_H\backslash \Irr(H)&\rightarrow & \Irr(G_\mathrm b),\\
       & \textrm{the $W_H$-orbit of $\chi\in \Irr(H)$} & \mapsto & (\Ind_{B}^{\Gb}\chi)(\lambda_{l_\chi}),
      \end{array}
    \]
    where $W_H$ denotes the real Weyl group of $G_\mathrm b$ with respect to $H$,  $l_\chi\in\{0,1, \dots, n_\mathrm b\}$ is the integer such that $\chi|_M$ is conjugate to $\delta_{l_\chi}$ by $W_H$, and $(\Ind_{B}^{\Gb}\chi)(\lambda_{l_\chi})$ denotes the unique
    irreducible subquotient in the normalized induced representation
    $\Ind_{B}^{\Gb}\chi$ containing the $K$-type $\lambda_{l}$. See \cite{Vg}*{Theorem 4.4.8}.

As in \eqref{indgl}, every representation in  $\Unip_{\ckcOpb}(\Gpb)$ is a normalized smooth parabolic induction 
         \[
         \pi' = 1_{n_{1}}\times \cdots \times 1_{n_{r}}\times \sgn_{n_{r+1}}\times \cdots \times \sgn_{n_{k}},
         \]
       where   $k\geq r\geq 0$,   $n_1, n_2, \dots, n_k$ are positive integers with  $n_1+n_2+ \dots+ n_k=\nnb$, $1_{i}$ and $\sgn_{i}$ respectively denote  the trivial character and the sign character of $\GL_{i}(\R)$ ($i\in \bN^+$).

    Write
    \begin{eqnarray*}
    \delta_{\pi'}&:=&\delta_{n_{1}+n_{2}+\dots + n_{r}}\in \Irr(M),\\
        \nu_{\pi'} &:=& \big( {\tfrac{n_{1}-1}{2},\tfrac{n_{1}-3}{2},\dots ,\tfrac{1-n_{1}}{2}, \dots , \tfrac{n_{k}-1}{2},\tfrac{n_{k}-3}{2},\dots ,\tfrac{1-n_{k}}{2}} \big)\in \Irr(A).
   \end{eqnarray*}
   Clearly the map
   \[
   \begin{array}{rcl}
     \fP:   \Unip_{\ckcOpb}(\Gpb)&\rightarrow &  W_H\backslash \Irr(H),\\
      \pi' &\mapsto & \textrm{the $W_H$-orbit of $\delta_{\pi'}\otimes \nu_{\pi'}$}
      \end{array}
      \]
       is well-defined and
    injective.

   Note that the map \eqref{badind} equals the composition $\fX\circ \fP$, and hence it is also injective. This proves the injectivity of \eqref{badind} in the case when $\star=D$. The same proof works in the case when $\star\in \{B, C, \wtC\}$ and we omit the details. When $\star\in \{C^*, D^*\}$, the map \eqref{badind} has to be injective since its domain is a singleton. Thus \eqref{badind} is injective in all cases.

  The bijection of \eqref{badind} follows from the injectivity  and the counting
  inequalities below:
  \[
    \abs{\PAP_{\star'}(\ckcO'_{\mathrm b})}=\abs{\Unip_{\ckcO'_{\mathrm b}}(G'_{\mathrm b})}\leq \abs{\Unip_{\ckcO_\mathrm b}(G_\mathrm b)}
    \leq \abs{\PBP^{\star}(\ckcO_{\mathrm b})} = \abs{\PAP_{\star'}(\ckcO'_{\mathrm b})}.
  \]
  Here $\star' \in \set{A^\R,A^{\bH}}$ is as in \eqref{def:star'}, the first inequality follows from the injectivity of \eqref{badind}, the second inequality follows from  \Cref{prop:countBCD22}, and the last equality is in \Cref{prop:BP.PP}.
\end{proof}

\subsection{Coherent continuation representations and parabolic induction}

The normalized smooth parabolic induction from $P_{\mathrm b}$ to $G_\mathrm b$ yields a linear map
\[
 \Ind:  \CK'(G'_{\mathrm b})\rightarrow \CK'(G_\mathrm b).
\]
This induces a linear map
\[
\begin{array}{rcl}
 \Ind:  \Coh_{\Lambda_\mathrm b}(\CK'(G'_{\mathrm b}))&\rightarrow &  \Coh_{\Lambda_\mathrm b}( \CK'(G_\mathrm b)),\\
   \Psi'_\mathrm b&\mapsto & \left(\nu\mapsto \Ind(\Psi'_\mathrm b(\nu))\right).
   \end{array}
\]

Let $P$ be an $\check \CO$-relevant parabolic subgroup of $G$ as in \Cref{reduction}.
Then  $G'_\mathrm b\times G_\mathrm g$ (or its quotient by $\{\pm 1\}$ when $\star=\wtC$) is naturally isomorphic to the  Levi quotient  of $P$. The
 normalized smooth parabolic induction using $P$ yields a linear map
\[
 \Ind:  \CK'(G'_{\mathrm b})\otimes \CK'(G_{\mathrm g})\rightarrow \CK'(G),
\]
which further induces a linear map (see \cite{Vg}*{Corollary 7.2.10})
\[
\begin{array}{rcl}
 \Ind:  \Coh_{\Lambda_{\mathrm b}}(\CK'(G'_{\mathrm b}))\otimes  \Coh_{\Lambda_\mathrm g}(\CK'(G_{\mathrm g})) &\rightarrow &  \Coh_{\Lambda}( \CK'(G)),\\
  \Psi'_\mathrm b\otimes \Psi_\mathrm g&\mapsto & \left((\nu_\mathrm b,\nu_\mathrm g)\mapsto \Ind(\Psi'_\mathrm b(\nu_\mathrm b)\otimes \Psi_\mathrm g(\nu_\mathrm g))\right).
   \end{array}
\]

The following result follows from Langlands' construction of standard modules. See \cite{Vg}*{Lemma~6.6.12 and Theorem~6.6.15}.
\begin{lem}\label{lem:cohred000}
       The  diagram
    \be\label{cdt}
    \begin{tikzcd}[column sep={4cm,between origins}]
      & \Coh_{\Lambda_{\mathrm b}}(\CK'(G'_{\mathrm b}))\otimes  \Coh_{\Lambda_{\mathrm g}}(\CK'(G_{\mathrm g}))
      \ar[dl,"\Ind \otimes  \mathrm{Id }"'] \ar[dr,"\Ind "]&\\
   \Coh_{\Lambda_\mathrm b}(\CK'(G_{\mathrm b}))   \otimes     \Coh_{\Lambda_{\mathrm g}}(\CK'(G_{\mathrm g})) \ar[rr,"\varphi"]& & \Coh_{\Lambda}(\CK'(G))
    \end{tikzcd}
  \ee
   commutes.
 \end{lem}

  \begin{prop}\label{lem:cohred111}
    For all  $\pi'\in \Unip_{\ckcO'_{\mathrm b}}(G'_{\mathrm b})$ and  $\pi_{\mathrm g}\in  \Unip_{\ckcO_{\mathrm g}}(G_{\mathrm g})$, the  normalized smooth parabolic induction
    $\Ind_P^G (\pi_\mathrm g\widehat \otimes \pi')$ is irreducible and belongs to  $\Unip_{\ckcO}(G)$.
 \end{prop}
\begin{proof}
By using the evaluation maps, the commutative  diagram \eqref{cdt} descends to a commutative diagram
\be\label{tri}
\begin{tikzcd}[column sep={4cm,between origins}]
      & \CK'_{\lambda_{\check \CO_\mathrm b}}(G'_{\mathrm b})   \otimes \CK'_{\lambda_{\check \CO_\mathrm g}}(G_{\mathrm g})
      \ar[dl," \Ind \otimes \mathrm{Id}"']\ar[dr,"\Ind "]&\\
      \CK'_{\lambda_{\check \CO_\mathrm b}}(G_{\mathrm b})\otimes  \CK'_{\lambda_{\check \CO_\mathrm g}}(G_{\mathrm g})\ar[rr,"\varphi_{\lambda_{\check \CO}}"]& & \CK'_{\lambda_{\check \CO}}(G).
    \end{tikzcd}
\ee
Here $\varphi_{\lambda_{\check \CO}}$ is as in Proposition \ref{propKL33}. Thus  $\Ind_P^G (\pi' \widehat \otimes \pi_\mathrm g)$ is irreducible
    by Propositions \ref{lem:Unip.BP} and \ref{propKL33}. As before,   \cite[Corollary 5.0.10]{B.Orbit} implies that $\Ind_P^G (\pi' \widehat \otimes \pi_\mathrm g)\in \Unip_{\ckcO}(G)$.
\end{proof}

By \Cref{lem:cohred111}, we have a well-defined map
\be\label{ind000}
  \Ind: \Unip_{\ckcO'_{\mathrm b}}(G'_{\mathrm b})\times \Unip_{\ckcO_{\mathrm g}}(G_{\mathrm g})   \longrightarrow \Unip_{\ckcO}(G).
  \ee
In view of the commutative diagram \eqref{tri},   Propositions \ref{propKL33} and \ref{lem:Unip.BP} imply that the map \eqref{ind000} is injective.
On the other hand,
Propositions \ref{propKL333} and \ref{lem:Unip.BP} imply that the domain  
and codomain of the map  \eqref{ind000} have the same cardinality.
Thus the map  \eqref{ind000}  is bijective. This completes the proof of  \Cref{reduction}.

\def\fhhaso{(\fhh^a_1)^*}
\def\fhhast{(\fhh^a_2)^*}
\newcommand{\ff}{f}
\newcommand{\ffcoh}{\varphi}

\section{Combinatorics of painted bipartitions}
\label{sec:com}

Throughout this section, we assume that $\star \in \{B,C,\wtC,C^{*},D,D^{*}\}$, and $\ckcO = \ckcOg$, namely $\ckcO $ has good parity.

Recall the set  $\CPPs(\ckcO)$ of primitive $\star$-pairs in $\ckcO$ (\Cref{defn:PP}). For each subset $\wp$ of $\CPPs(\ckcO)$, we have defined a pair   $(\imath_{\wp},\jmath_{\wp})$ of Young diagrams in \Cref{sec:LCBCD}. Our main object of interest is the following set of disjoint union: 
\[
  \tPBP_{G}(\check \CO):=  \bigsqcup_{\wp \subseteq \CPP(\ckcO)}\PBPGOP ,
 \]
where $
\PBPGOP $ is defined in \eqref{eq:PBPGP}. 

\subsection{Main combinatorial results}\label{MainComb}

We say that the orbit $\check \CO$ is  quasi-distinguished 
if there is no $\star$-pair that is balanced in $\check \CO$ (see \Cref{defn:PP}). 

For the quaternionic groups, we have the following non-existence result on painted
bipartitions.

\begin{prop} \label{prop:PBP1} Suppose that $\star\in \set{C^{*}, D^{*}}$. 
If $\PBP_{G}(\check \CO, \wp)$ is nonempty for some $\wp\subset \CPPs(\ckcO)$, then    $\check \CO$ is quasi-distinguished and $\wp=\emptyset$. Consequently, 
     \[
     \sharp(\tPBP_{G}(\ckcO)) = \sharp(\PBP_{G}(\ckcO)).
  \]
\end{prop}

\begin{proof}
 Suppose that
  $\uptau = (\imath_{\wp}, \cP)\times (\jmath_{\wp},\cQ)\times \star\in \PBP_{G}(\check \CO, \wp)$.
Assume by contradiction that  $\wp \neq \emptyset$. 

   First assume that  $\star = C^{*}$.  Pick an element $(2i-1, 2i)\in \wp$. Then we have that
  \begin{equation}\label{eq:res.C*}
    \bfcc_{i}(\imath_{\wp}) = \half(\bfrr_{2i-1}(\ckcO)+1)>
    \half(\bfrr_{2i}(\ckcO)-1) = \bfcc_{i}(\jmath_{\wp}).
      \end{equation}
   By the requirements of a painted bipartition, we also have that
  \begin{eqnarray*}
    \bfcc_{i}(\imath_{\wp})& = &\sharp\set{j\in \BN^+ \mid (i,j)\in \BOX{\imath_\wp}, \, \cP(i,j)=\bullet} \\
    &= &\sharp\set{j\in \BN^+\mid  (i,j)\in \BOX{\jmath_\wp}, \, \cQ(i,j)=\bullet} \\
    &\leq & \bfcc_{i}(\jmath_{\wp}).
  \end{eqnarray*}
 This contradicts \eqref{eq:res.C*} and therefore $\wp = \emptyset$.

  Now assume that $\star = D^{*}$. Pick an element  $(2i, 2i+1)\in \wp$.
  Then we have that
  \begin{equation}\label{eq:res.D*}
    \bfcc_{i+1}(\imath_{\wp}) = \half(\bfrr_{2i}(\ckcO)+1)>
    \half(\bfrr_{2i+1}(\ckcO)-1) = \bfcc_{i}(\jmath_{\wp}).
  \end{equation}
      By the requirements of a painted bipartition, we also have that
  \begin{eqnarray*}
  \bfcc_{i+1}(\imath_{\wp})& \leq & \sharp\set{j\in \BN^+ \mid (i,j)\in \BOX{\imath_\wp}, \,  \cP(i,j)=\bullet} \\
  &=&\sharp\set{j\in \BN^+ \mid (i,j)\in \BOX{\jmath_\wp}, \,  \cQ(i,j)=\bullet} \\
  &\leq & \bfcc_{i}(\jmath_{\wp}).  \end{eqnarray*}
 This contradicts \eqref{eq:res.D*} and therefore  $\wp = \emptyset$.

It is also straightforward to show that $\check \CO$ is quasi-distinguished from the requirements of a painted bipartition. We omit the details. 
\end{proof}

In view of \Cref{prop:PBP1}, we  assume in the rest of this section that $\check \CO$ is quasi-distinguished when $\star\in \{C^*,D^*\}$. 

In the rest of the section, we will define a notion of descent for painted
bipartitions (in $\PBPGOP$) which will be used as an induction mechanism and we will deduce certain counting formulas by a generating function approach. The following
proposition will then be an immediate consequence of a certain independence property of the generating functions in \Cref{prop:gf.BD} and \Cref{prop:gf.C}.  

\begin{prop} \label{prop:PBP2}
Suppose that   $\star\in \set{B,C,\wtC,D}$. Then
  \[
    \sharp(\PBPop{G}{}{\ckcO}{\wp}) = \sharp(\PBPop{G}{}{\ckcO}{\emptyset}), \quad \textrm{for all } \wp \subset \CPPs(\ckcO).
  \]
 Consequently,
     \[
     \sharp(\tPBP_{G}(\ckcO)) = 2^{\sharp(\CPPs(\ckcO))}\cdot \sharp(\PBP_{G}(\ckcO)).
  \]
\end{prop}

\subsection{Type $C$ and $\wtC$: shape shifting}
 In this subsection, we assume that $\star \in \set{C,\wtC}$. We will define a shape shifting operation on painted bipartitions which  will be used in defining the notion of descent.

 Recall that $\ckcO$ is an orbit with good parity. In this subsection we assume that  $(1,2)\in \CPPs(\ckcO)$, namely the pair $(1,2)$ is primitive in $\check \CO$.
 Let $\wp$ be a subset of  $\CPPs(\ckcO)$ such that $(1,2)\notin \wp$. 
 Put  $\wp_{\uparrow}:= \wp\cup \set{(1,2)}\subset \CPPs(\ckcO)$. 
 
Note that 
  \[
      (\bfcc_{1}(\imath_{\wpu}), \bfcc_{1}(\jmath_{\wpu})) =
      \begin{cases}
        (\bfcc_{1}(\jmath_{\wpd})+1, \bfcc_{1}(\imath_{\wpd})-1), & \text{when
        } \star = C;\\
        (\bfcc_{1}(\jmath_{\wpd}), \bfcc_{1}(\imath_{\wpd})), & \text{when
        } \star = \wtC,\\
      \end{cases}
  \]
  and
  \[
      (\bfcc_{i}(\imath_{\wpu}),\bfcc_{i}(\jmath_{\wpu})) =(\bfcc_{i}(\imath_{\wpd}),\bfcc_{i}(\jmath_{\wpd})),\qquad \text{for $i=2, 3, 4, \dots$}.
  \]

 For
 $\uptau := (\imath_{\wp},\cP_{\uptau})\times (\jmath_{\wp},\cQ_{\uptau})\times \alpha \in \PBPop{G}{}{\ckcO}{\wp} $,
 we will define an element 
 \[
   \uptauu:= (\imath_{\wpu}, \cP_{\uptauu})\times (\jmath_{\wpu},\cQ_{\uptauu})\times \alpha\in \PBPop{G}{}{\ckcO}{\wpu}
 \]
 by the following recipe.

 \newcommand{\localtextbulletone}{\textcolor{black}{\raisebox{.45ex}{\rule{.6ex}{.6ex}}}}

\medskip
 {\bfseries The case when $\star = C$.}
 Note that $\bfcc_1(\imath_{\wpu})> \bfcc_1(\imath_\wp)\geq 1$ since $(1,2)\in \CPP_\star(\ckcO)$.
  For all $(i,j)\in \BOX{\imath_{\wpu}}$, we define $\cP_{\uptauu}(i,j)$ case by
  case as in what follows.
 \begin{enumerate}[label=(\alph*)]
   \item Suppose that
   $\cP_{\uptaud}(\bfcc_{1}(\imath_{\wpd}),1)\neq \bullet$.
   \begin{enumerate}[label={\localtextbulletone}]
     \item If $\bfcc_{1}(\imath_{\wpd})\geq 2$ and
     $\cP_{\uptau}(\bfcc_{1}(\imath_{\wpd})-1,1) = c$,
     we define
     \[
       \cP_{\uptauu}(i,j) := \begin{cases}
         r ,& \text{if $j=1$ and $\bfcc_{1}(\imath_{\wpd})-1
           \leq i \leq \bfcc_{1}(\imath_{\wpu})-2$};\\
         c ,& \text{if $(i,j)=(\bfcc_{1}(\imath_{\wpu})-1,1)$};\\
         d ,& \text{if $(i,j)=(\bfcc_{1}(\imath_{\wpu}),1)$};\\
         \cP_{\uptaud}(i,j) ,&\text{otherwise}.
       \end{cases}
     \]
     \item Otherwise, we define
     \[
       \cP_{\uptauu}(i,j) := \begin{cases}
         r ,& \text{if $j=1$ and $\bfcc_{1}(\imath_{\wp})
           \leq i \leq \bfcc_{1}(\imath_{\wpu})-1$};\\
         \cP_{\uptaud}(\bfcc_{1}(\imath_{\wpd}),1), &
         \text{if $(i,j)=(\bfcc_{1}(\imath_{\wpu}),1)$};\\
         \cP_{\uptaud}(i,j) ,&\text{otherwise},
       \end{cases}
     \]
   \end{enumerate}
   \item Suppose that $\cP_{\uptaud}(\bfcc_{1}(\imath_{\wpd}),1)=\bullet$.
   \begin{enumerate}[label={\localtextbulletone}]
     \item If $\bfcc_{2}(\imath_{\wpd}) = \bfcc_{1}(\imath_{\wpd})$
     and
     $\cP_{\uptaud}(\bfcc_{1}(\imath_{\wpd}),2) = r$,
     we define
     \[
       \cP_{\uptauu}(i,j) := \begin{cases}
         r ,& \text{if $j=1$ and $\bfcc_{1}(\imath_{\wpd})\leq i \leq \bfcc_{1}(\imath_{\wpu})-1$};\\
         c ,& \text{if $(i,j)=(\bfcc_{2}(\imath_{\wpd}),2)$};\\
         d ,& \text{if $(i,j)=(\bfcc_{1}(\imath_{\wpu}),1)$};\\
         \cP_{\uptaud}(i,j) ,&\text{otherwise}.
       \end{cases}
     \]
     \item Otherwise, we define
     \[
       \cP_{\uptauu}(i,j) := \begin{cases}
         r ,& \text{if $j=1$ and $\bfcc_{1}(\imath_{\wpd})\leq i \leq \bfcc_{1}(\imath_{\wpu})-2$};\\
         c ,& \text{if $(i,j)=(\bfcc_{1}(\imath_{\wpu})-1,1)$};\\
         d ,& \text{if $(i,j)=(\bfcc_{1}(\imath_{\wpu}),1)$};\\
         \cP_{\uptaud}(i,j) ,&\text{otherwise}.
       \end{cases}
     \]
   \end{enumerate}
 \end{enumerate}

Note that $\BOX{\jmath_{\wpu}}\subseteq \BOX{\jmath_{\wp}}$.   For all $(i,j)\in \BOX{\jmath_{\wpu}}$, we define
   \[
     \cQ_{\uptauu}(i,j) := \cQ_{\uptau}(i,j).
   \]

\medskip
{\bfseries The case when $\star = \wtC$.} 
  For all $(i,j)\in \BOX{\imath_{\wpu}}$, define 
  \begin{equation} \label{eq:T.M1}
     \cP_{\uptauu}(i,j) :=
     \begin{cases}
       s,& \text{if $j=1$ and  $\cQ_{\uptau}(i,j)=r$;}\\
       c,& \text{if $j=1$ and  $\cQ_{\uptau}(i,j)=d$;}\\
       \cP_{\uptau}(i,j), &\text{otherwise.}
     \end{cases}
   \end{equation}
  For all $(i,j)\in \BOX{\jmath_{\wpu}}$, define
  \begin{equation} \label{eq:T.M2}
    \cQ_{\uptauu}(i,j) :=
    \begin{cases}
      r,& \text{if $j=1$ and  $\cP_{\uptau}(i,j)=s$;}\\
      d,& \text{if $j=1$ and  $\cP_{\uptau}(i,j)=c$;}\\
      \cQ_{\uptau}(i,j), &\text{otherwise.}
    \end{cases}
  \end{equation}

In both cases, it is routine to check that $\uptauu$ is a painted 
bipartition and is an element of $\PBPop{G}{}{\ckcO}{\wpu}$. 

\begin{lem}\label{lem:sn}
 The map 
  \[
    T_{\wp,\wpu}\colon \PBPop{G}{}{\ckcO}{\wp} \rightarrow
    \PBPop{G}{}{\ckcO}{\wpu}, \quad \uptau \mapsto \uptauu
  \]
  defined above is bijective. 
\end{lem}
\begin{proof}
  Let $\uptau' = (\imath_{\wpu}, \cP_{\uptau'})\times (\jmath_{\wpu},\cQ_{\uptau'})\times \alpha \in\PBPop{G}{}{\ckcO}{\wpu}$. We will define a painted bipartition $\uptaud\in \PBPop{G}{}{\ckcO}{\wp}$ by the following recipe.

  \medskip

  {\bfseries The case when $\star = C$.}
  Note that $$\bfcc_1(\imath_{\wpu})= \bfcc_1(\jmath_{\wp})+1>  \bfcc_1(\imath_\wp)= \bfcc_1(\jmath_{\wpu})+1.$$
  We define $\cP_{\uptaud}$ and $\cQ_{\uptaud}$  case by case as in what follows. 
 \begin{enumerate}[label=(\alph*)]
   \item Suppose that
   $\cP_{\uptau'}(\bfcc_{1}(\imath_{\wpu})-1,1)=c$.
   \begin{enumerate}[label={\localtextbulletone}]
     \item
           If $\bfcc_{1}(\jmath_{\wpu})\geq 1$ and
           $\cP_{\uptau'}(\bfcc_{1}(\jmath_{\wpu}),1) = r$, 
         for all $(i,j)\in \BOX{\imath_{\wpd}}$  we define
          \[
       \cP_{\uptaud}(i,j) := \begin{cases}
         c ,& \text{if $(i,j)=(\bfcc_{1}(\imath_{\wpd})-1,1)$};\\
         d ,& \text{if $(i,j)=(\bfcc_{1}(\imath_{\wpd}),1)$};\\
         \cP_{\uptau'}(i,j) ,&\text{otherwise},
       \end{cases}
     \]
     and for all $(i,j)\in \BOX{\jmath_{\wpd}}$ we define
     \[
       \cQ_{\uptaud}(i,j) := \begin{cases}
         s ,& \text{if $j=1$ and $\bfcc_{1}(\imath_{\wpd})\leq i \leq \bfcc_{1}(\jmath_{\wpd})$};\\
         \cQ_{\uptau'}(i,j), &\text{otherwise}.
       \end{cases}
     \]
     
     \item Otherwise, for all $(i,j)\in \BOX{\imath_{\wpd}}$ we define
     \[
       \cP_{\uptaud}(i,j) := \begin{cases}
         \bullet ,& \text{if $(i,j)=(\bfcc_{1}(\imath_{\wpd}),1)$};\\
         \cP_{\uptau'}(i,j) ,&\text{otherwise},
       \end{cases}
     \]
     and for all $(i,j)\in \BOX{\jmath_{\wpd}}$ we define
     \[
       \cQ_{\uptaud}(i,j) := \begin{cases}
         \bullet ,& \text{if $(i,j)=(\bfcc_{1}(\imath_{\wpd}),1)$};\\
         s ,& \text{if $j=1$ and $\bfcc_{1}(\imath_{\wpd})+1\leq i \leq \bfcc_{1}(\jmath_{\wpd})$};\\
         \cQ_{\uptau'}(i,j) ,&\text{otherwise}. 
       \end{cases}
     \]
   \end{enumerate}
   \item Suppose that $\cP_{\uptau'}(\bfcc_{1}(\imath_{\wpu})-1,1)=r$.
   \begin{enumerate}[label={\localtextbulletone}]
     \item If $\bfcc_{2}(\imath_{\wpu}) = \bfcc_{1}(\jmath_{\wpu})+1$
     and
     $(\cP_{\uptau'}(\bfcc_{1}(\imath_{\wpu}),1), \cP_{\uptau'}(\bfcc_{2}(\imath_{\wpu}),2) ) = (d,c)$,
         for all $(i,j)\in \BOX{\imath_{\wpd}}$  we define
     \[
       \cP_{\uptaud}(i,j) := \begin{cases}
         \bullet ,& \text{if $(i,j)=(\bfcc_{1}(\imath_{\wpd}),1)$};\\
         r ,& \text{if $(i,j)=(\bfcc_{2}(\imath_{\wpd}),2)$};\\
         \cP_{\uptau'}(i,j) ,&\text{otherwise},
       \end{cases}
     \]
     and  for all $(i,j)\in \BOX{\jmath_{\wpd}}$  we define
     \[
       \cQ_{\uptaud}(i,j) := \begin{cases}
         \bullet ,& \text{if $(i,j)=(\bfcc_{1}(\imath_{\wpd}),1)$};\\
         s ,& \text{if $j=1$ and $\bfcc_{1}(\imath_{\wpd})+1\leq i \leq \bfcc_{1}(\jmath_{\wpd})$};\\
         \cQ_{\uptau'}(i,j), &\text{otherwise}.
       \end{cases}
     \]
     \item Otherwise, 
         for all $(i,j)\in \BOX{\imath_{\wpd}}$  we define
     \[
       \cP_{\uptaud}(i,j) := \begin{cases}
         \cP_{\uptau'}(\bfcc_{1}(\imath_{\wpu}),1), & \text{if $(i,j)=(\bfcc_{1}(\imath_{\wpd}),1)$};\\
         \cP_{\uptau'}(i,j) ,&\text{otherwise},
       \end{cases}
     \]
     and for all $(i,j)\in \BOX{\jmath_{\wpd}}$  we define
     \[
       \cQ_{\uptaud}(i,j) := \begin{cases}
         s ,& \text{if $j=1$ and $\bfcc_{1}(\imath_{\wpd})\leq i \leq \bfcc_{1}(\jmath_{\wpd})$};\\
         \cQ_{\uptau'}(i,j) ,&\text{otherwise}.
       \end{cases}
     \]
   \end{enumerate}
 \end{enumerate}

\medskip 

  {\bfseries The case when $\star = \wtC$.}
  The definition of $T_{\wpu,\wp}$ is given by the formulas \eqref{eq:T.M1} and
  \eqref{eq:T.M2} with the
  role of $\wp$ and $\wpu$ switched.

In both cases, it is routine to check that 
  \[
    T_{\wpu,\wp}\colon
    \PBPop{G}{}{\ckcO}{\wpu} \rightarrow \PBPop{G}{}{\ckcO}{\wp}  , \qquad \uptau'\mapsto \uptaud
  \]
 is well-defined and is the inverse of $T_{\wp,\wpu}$.
\end{proof}

\subsection{Naive descent of a painted bipartition}\label{sec:comb}

As before, let  $\star\in \{ B, C,  D, \widetilde{C},  C^*, D^*\}$ and let $\check \CO$ be a Young diagram that has good parity.

Define a label
\[
\star':=\widetilde{C}, \ D, \  C, \ B, \ D^*,\  \textrm{ or } \ C^*
\]
respectively if
\[
\star=B,\  C, \ D, \ \widetilde{C}, \ C^*, \ \textrm{ or }\  D^*.
\]
We call $\star'$
 the Howe dual of $\star$.

\def\bipartl{\mathrm{bi\cP_L}}
\def\bipartr{\mathrm{bi\cP_R}}
\def\dsdiagl{\mathrm{DS_L}}
\def\dsdiagr{\mathrm{DS_R}}
\def\DDl{\eDD_\mathrm{L}}
\def\DDr{\eDD_\mathrm{R}}

  The dual descent of the Young diagram $\check \CO$ is defined to be
  \[
   \check \CO':=\check \nabla(\check \CO):=\check \nabla_\star( \check \CO):=\begin{cases}
   \tytb{\   }\, , \quad& \textrm{if $\star\in \{D, D^*\}$ and $|\check \CO|=0$};\smallskip\\
   \check \nabla_{\mathrm{naive}}(\check \CO), \quad & \textrm{otherwise},
    \end{cases}
  \]
  where $\tytb{\   }$ denotes the Young diagram that has total size $1$,
  and $\ckDDn(\ckcO)$ denotes the Young diagram obtained from $\check \CO$ by removing the first row. We also view $\check \CO'$ as a nilpotent orbit in the Langlands dual (or metaplectic Langlands dual) of the complexified Lie algebra of a classical real Lie group of type $\star'$. 

For a Young diagram $\imath$, its naive descent, which is denoted by $\nabla_\mathrm{naive}(\imath)$, is defined to be the Young diagram obtained from $\imath$ by removing the first column. By convention, $\nabla_\mathrm{naive}(\emptyset)=\emptyset$.

We start with the notion of the naive descent of a painted bipartition. Let $\uptau=(\imath,\cP)\times (\jmath,\cQ)\times \gamma $ be a painted bipartition such that $\star_\uptau=\star$. Put
\delete{\begin{equation} \label{eq:def.alphap}
\alpha '=\begin{cases} B^+,
& \textrm{if $\gamma = \wtC$ and $\cP_\tau(l_{\star,\ckcO},1),1) \neq c$;}\\
B^-,
& \textrm{if $\gamma = \wtC$ and $\cP_\tau(l_{\star,\ckcO},1),1)  = c$;}\\
\star', & \textrm{if $\gamma \neq \widetilde C$}.
\end{cases}
\end{equation}
}
  \begin{equation} \label{eq:def.alphap}
    \gamma '=\begin{cases} B^+,
  & \textrm{if $\alpha=\widetilde{C}$ and $c$ does not occur in the first column of $(\imath,\cP)$}; \smallskip \\
  B^-,
  & \textrm{if $\alpha=\widetilde{C}$ and  $c$ occurs in the first column of $(\imath,\cP)$}; \smallskip \\
  \star', & \textrm{if $\alpha\neq \widetilde C$}.
  \end{cases}
  \end{equation}

\begin{lem}\label{lemDDn1}
  If $\star \in \set{B,C,C^*}$, then there is a unique painted bipartition of the form $\uptaupn= (\imath',\cPpn)\times (\jmath',\cQpn)\times \gamma '$ with the following properties:
  \begin{itemize}
        \item $
   (\imath',\jmath')= (\imath,\DD_\mathrm{naive}(\jmath)); \smallskip
   $
   \item for all $(i,j)\in \BOX{\imath'}$,
   \[
     \cPpn(i,j)=\begin{cases}
    \bullet \textrm{ or } s,&\textrm{ if  $\ \cP(i,j)\in \{\bullet, s\}$;} \smallskip \\
  \cP(i,j),& \textrm{ if $\ \cP(i,j)\notin \{\bullet, s\}$};\end{cases}
   \]
   \item for all $(i,j)\in \BOX{\jmath'}$,
   \[
     \cQpn(i,j)=\begin{cases}
    \bullet \textrm{ or } s,&\textrm{ if  $\ \cQ(i,j+1)\in \{\bullet, s\}$;} \smallskip \\
  \cQ(i,j+1), & \textrm{ if $\ \cQ(i,j+1)\notin \{\bullet, s\}$}.  \end{cases}
   \]
    \end{itemize}
    \end{lem}

   \begin{proof}
    First assume that the images of $\cP$ and $\cQ$ are both contained in $\{\bullet, s\}$. Then  the image of $\cP$  is in fact contained in $\{\bullet\}$, and $(\imath, \jmath)$ is  right interlaced in the sense that
 \[
 \mathbf{c}_1(\jmath)\geq \mathbf{c}_1(\imath)\geq \mathbf{c}_2(\jmath)\geq \mathbf{c}_2(\imath)\geq \mathbf{c}_3(\jmath)\geq \mathbf{c}_3(\imath) \geq \cdots.
 \]
 Hence $ (\imath',\jmath'):= (\imath,\DD(\jmath))$ is left interlaced in the sense that
 \[
 \mathbf{c}_1(\imath')\geq \mathbf{c}_1(\jmath')\geq \mathbf{c}_2(\imath')\geq \mathbf{c}_2(\jmath')\geq \mathbf{c}_3(\imath')\geq \mathbf{c}_3(\jmath') \geq \cdots.
 \]
 Then it is clear that there is a unique painted bipartition of the form  $\uptaupn=(\imath',\cPpn)\times (\jmath',\cQpn)\times \gamma '$ such that symbols of $\cPpn$ and $\cQpn$ are both in $\{\bullet, s\}$. This proves the lemma in the special case when the images of $\cP$ and $\cQ$ are both contained in $\{\bullet, s\}$.

 The proof of the lemma in the general case is easily reduced to this special case.
   \end{proof}
    \begin{lem}\label{lemDDn2}
    If $\star \in \set{ \wtC, D,D^*}$, then there is a unique painted bipartition of the form $\uptaupn= (\imath',\cPpn)\times (\jmath',\cQpn)\times \gamma '$ with the following properties:
  \begin{itemize}
        \item $
   (\imath',\jmath')= (\DD_\mathrm{naive}(\imath),\jmath); \smallskip
   $
   \item for all $(i,j)\in \BOX{\imath'}$,
   \[
     \cPpn(i,j)=\begin{cases}
    \bullet \textrm{ or } s,&\textrm{ if  $\ \cP(i,j+1)\in \{\bullet, s\}$;} \smallskip \\
  \cP(i,j+1),& \textrm{ if $\ \cP(i,j+1)\notin \{\bullet, s\}$};\end{cases}
   \]
   \item for all $(i,j)\in \BOX
          {\jmath'}$,
   \[
     \cQpn(i,j)=\begin{cases}
    \bullet \textrm{ or } s,&\textrm{ if  $\ \cQ(i,j)\in \{\bullet, s\}$;} \smallskip \\
  \cQ(i,j), & \textrm{ if $\ \cQ(i,j)\notin \{\bullet, s\}$}.  \end{cases}
   \]

    \end{itemize}
\end{lem}
\begin{proof}
  The proof is similar to that of \Cref{lemDDn1}.
\end{proof}
 In the notation of \Cref{lemDDn1} and \Cref{lemDDn2}, we call $\uptaupn$ the naive descent of $\uptau$, to be denoted by $\DDn(\uptau)$.

\begin{eg}  
 Suppose that the nonzero row lengths of $\check \CO$ are $8,6,6,6,4,4,2$. Then 
    \[
     \uptau = \ytb{\bullet\bullet\bullet {c},\bullet {s} {c},{s} {c},{c}}
    \times \ytb{\bullet\bullet\bullet ,\bullet {r} {d},{d}{d}, \none}
    \times \widetilde C \in \PBP_{\widetilde{\Sp}_{36}(\bR)}(\check \CO).     
    \]
  We have that  
   \[
    \nabla_{\mathrm{naive}}(\uptau) =\ytb{\bullet\bullet{c} ,\bullet{c},{c} }
    \times  \ytb{\bullet\bullet {s} ,\bullet {r} {d},{d}{d}}\times B^- \in \PBP_{\SO(14,15)}(\check \CO').
    \]
\end{eg}

\subsection{Descent of a painted bipartition}

Given a subset $\wp\subset \mathrm{PP}_{\star}(\check \CO)$, the dual descent of $\wp$ is defined to be the subset
\begin{equation*}\label{eq:DD.wp}
  \wp':=\ckDD(\wp):=\{(i,i+1)\mid i\in \bN^+, \, (i+1, i+2)\in \wp \}\subset \mathrm{PP}_{\star'}(\check \CO').
\end{equation*}
Note that 
\[
 (\imath_{\wp'}, \jmath_{\wp'})=\begin{cases}
    (\imath_\wp, \nabla_{\mathrm{naive}}(\jmath_{\wp})), & \textrm{if $\star=B$, or $\star\in\{C, C^*\}$ and $(1,2)\notin \wp$;} \\
    ( \imath_{\wpm},  \nabla_{\mathrm{naive}}(\jmath_{\wpm})), & \textrm{if $\star\in \{C, C^*\}$ and $(1,2)\in \wp$}; \\
        ( \nabla_{\mathrm{naive}}(\imath_{\wp}),\jmath_{\wp}), & \textrm{if $\star\in \{D, D^*\}$, or $\star=\wtC$ and $(1,2)\notin \wp$}; \\
          ( \nabla_{\mathrm{naive}}(\imath_{\wpm}),\jmath_{\wpm}), & \textrm{if  $\star=\wtC$ and $(1,2)\in \wp$}, \\
 \end{cases}
\]
where $\imath_{\wp'}:=\imath_{\star'}(\check \CO', \wp')$, $\jmath_{\wp'}:=\jmath_{\star'}(\check \CO', \wp')$, and $\wpm=\wp\setminus \{(1,2)\}$.

Recall from \eqref{eq:PBPSP} the set \[
\PBPOP := \set{\uptau \text{ is a painted bipartition}|  \star_{\uptau}=\star, (\imath_{\uptau},\jmath_{\uptau}) = (\imath_{\wp},\jmath_{\wp})}. 
\]
Suppose that
\[
\uptau=(\imath_{\wp},\cP)\times(\jmath_{\wp},\cQ)\times \gamma  \in  \PBPOP. 
\]
We will use the notion   of \Cref{lemDDn1} and \Cref{lemDDn2}, and define an element 
\[
  \uptau' := (\imathwpp, \cP')\times (\jmathwpp, \cQ')\times \gamma ' \in \PBP_{\star'}(\ckcOp,\wp'),
\]
to be called the descent of $\uptau$. 

\medskip

{\bfseries The case when $\star = B$. }
\begin{enumerate}[label=(\alph*),wide=0pt]
  \item If  
    \[
    \begin{cases}
      \gamma = B^+; & \\
      (2,3) \notin \wp ;& \\
      \bfrr_2(\ckcO)>0; & \\
      \cQ(\bfcc_1(\imath_{\wp}),1)\in \set{r,d},
    \end{cases}
    \]
    we define %
        \[
        \cP'(i,j) := \begin{cases}
          s, & \ \text{ if $(i,j) = (\bfcc_1(\imathwpp),1)$;}\\
          \cPpn(i,j), & \ \text{ otherwise}
        \end{cases}
        \] for all $(i,j)\in \BOX{\imathwpp}$, and
        $\cQ' := \cQpn $.
  \item If  
    \[
    \begin{cases}
      \gamma = B^+; & \\
      (2,3)\in \wp ;& \\
      \cQ(\bfcc_2(\jmath_{\wp}),1)\in \set{r,d},
    \end{cases}
    \]
    we define
    $\cP' := \cPpn $, and
    \[
    \cQ'(i,j) := \begin{cases}
      r, & \ \text{ if $(i,j) = (\bfcc_1(\jmathwpp),1)$;}\\
      \cQpn(i,j), & \ \text{ otherwise}
    \end{cases}
    \] for all $(i,j)\in \BOX{\jmathwpp}$.
    \trivial[h]{
      Note that $\bfcc_{1}(\jmathwpp) = \bfcc_{2}(\jmath_{\wp})$
    }
        \item Otherwise, we define $\cP' := \cPpn$ and $\cQ':= \cQpn$.
\end{enumerate}

\medskip 
{\bfseries The case when $\star = D$. }
\begin{enumerate}[label=(\alph*),wide=0pt]
  \item If  
  \[
    \begin{cases}
      \bfrr_2(\ckcO)=\bfrr_{3}(\ckcO)>0; & \\
      \cP(\bfcc_{2}(\imath_{\wp}),2) = c;  &\\
      \cP(i,1)\in \set{r,d}, & \text{for all
        $\bfcc_{2}(\imath_{\wp})\leq i\leq \bfcc_{1}(\imath_{\wp})$},
    \end{cases}
  \]
  we define %
  \[
    \cP'(i,j) := \begin{cases}
      r, & \ \text{ if $(i,j) = (\bfcc_1(\imathwpp),1)$;}\\
      \cPpn(i,j), & \ \text{ otherwise}
    \end{cases}
    % \quad \text{for all $(i,j)\in \BOX{\imathwpp}$,}
  \] for all $(i,j)\in \BOX{\imathwpp}$, and
  $\cQ' := \cQpn $.
  \item If 
  \[
    \begin{cases}
      (2,3)\in \wp ;& \\
      \cP(\bfcc_2(\imath_{\wp})-1,1)\in \set{r,c},
    \end{cases}
  \]
  we define
  \[
    \cP'(i,j) := \begin{cases}
      r, & \ \text{ if $(i,j) = (\bfcc_1(\imathwpp)-1,1)$;}\\
      \cP(\bfcc_2(\imath_{\wp})-1,1), & \ \text{ if $(i,j) = (\bfcc_1(\imathwpp),1)$;}\\
      \cPpn(i,j), & \ \text{ otherwise}
    \end{cases}
  \] for all $(i,j)\in \BOX{\imathwpp}$,  and 
  $\cQ' := \cQpn $. 
  \item Otherwise, we define $\cP' := \cPpn$ and $\cQ':= \cQpn$.
\end{enumerate}

\medskip

{\bfseries The case when $\star \in \set{C,\wtC,C^{*},D^{*}}$}. 
\begin{enumerate}[label=(\alph*),wide=0pt]
  \item If $(1,2)\notin \wp$, we define
  \[
    \uptau' := \DDn(\uptau).
  \]
  \item If $(1,2)\in \wp$, then $\star \in \set{C,\wtC}$ (by \Cref{prop:PBP1}) and we define
  \[
  \uptau' := \DDn(\uptau_{\wpm}),
  \]
  where $\wpm := \wp \setminus \set{(1,2)}$,  
  $\uptau_{\wpm}:= T_{\wpm,\wp}^{-1}(\uptau)$,
  and $T_{\wpm,\wp}$ is as in \Cref{lem:sn}.
 \end{enumerate}

\medskip 
In all cases, it is routine to check that $\uptau'$ thus constructed is an element in $\PBPOPp$, which is called the descent of $\uptau$. 
To summarize, we have a map 
  \begin{equation}
    \label{eq:DD.CC}
    \DD \colon \PBPOP \rightarrow \PBPOPp.
  \end{equation}
We call $\DD$ the descent map of painted bipartitions. 

\begin{eg} Suppose that the nonzero row lengths of $\check \CO$ are $7,7,7,3$. 
Then
 \[
 \uptau= \ytb{\bullet\bullet, \bullet s, \bullet s, r c} \times \ytb{\bullet\bullet,\bullet,\bullet, \none }\times
  D \in \PBP_{\SO(11,13)}(\ckcO).
 \]
We have that  
\[
 \DDn(\uptau)=\ytb{\bullet, \bullet , \bullet ,  c} \times \ytb{\bullet s,\bullet,\bullet, \none } \times
  C\qquad\textrm{and}\qquad \DD(\uptau)=\ytb{\bullet, \bullet , \bullet ,  r} \times \ytb{\bullet s,\bullet,\bullet, \none } \times
  C,
 \]
both in $\PBP_{\Sp_{16}(\bR)}(\check \CO')$.
\end{eg}

\subsection{Properties of the descent map}

We first define the notion of the tail of a painted bipartition, when $\star\in\set{B, D, C^*}$. 

\def\startt{\star_{\mathrm t}}
Put
\[
  \startt:= \begin{cases}
  D, & \ \text{ if $\star\in \set{B,D}$}; \\
  C^*, &\  \text{ if $\star=C^*$},
\end{cases}
\quad
\textrm{and}\quad 
k := \begin{cases}
  \frac{\bfrr_{1}(\ckcO)-\bfrr_{2}(\ckcO)}{2} + 1, &\ 
    \text{if $\star\in \{B,D\}$}; \\
\lfloor\frac{\bfrr_{1}(\ckcO)-\bfrr_{2}(\ckcO)-1}{2}\rfloor,  & \  \text{if $\star=C^*$}.
  \end{cases}
\]
Note that $k$ is a non-negative integer, and is positive when $\star \in \set{B,D}$. 

From the pair $(\star, \ckcO)$, we define a Young diagram $\ckcO_{\bftt}$ as follows.
\begin{itemize}
  \item If $\star \in \set{B,D}$, then $\ckcO_{\bftt}$ consists of two rows with
        lengths $2k-1$ and $1$.
  \item If $\star =C^*$, then $\ckcO_{\bftt}$ consists of one row with length
        $2k+1$.
\end{itemize}
Note that in all the three cases
 $\check \CO_{\mathbf t}$ has good parity (with respect to $\star_{\mathbf t}$), $\mathrm{PP}_{\star_{\mathbf t}}(\check \CO_{\mathbf t})=\emptyset$,  and every element in $\PBP_{\star_\bftt}(\ckcO_\bftt)$ has the form
 \begin{equation}
 \label{tail0}
  \ytb{{x_1} , {x_2} , {\enon\vdots},{\enon{\vdots}},{x_k}  } \times \emptyset \times
  D,\qquad \qquad  \ytb{{x_1} , {x_2} , {\enon\vdots},{\enon{\vdots}},{x_k}  } \times \emptyset \times
  D,\qquad\textrm{or}\qquad \emptyset \times  \ytb{{x_1} , {x_2} , {\enon\vdots},{\enon{\vdots}},{x_k}  } \times
 C^*,
\end{equation}
according to $\star=B, D$, or $C^*$, respectively. When $k=0$ (and thus $\star=C^*$), the element in $\PBP_{\star_\bftt}(\ckcO_\bftt)$ is understood to be 
$\emptyset \times  \emptyset \times C^*$. 

\medskip

 When $\star\in\set{B, D, C^*}$, the tail $\uptau_\bftt\in \PBP_{\star_\bftt}(\ckcO_\bftt)$ of an element $ \uptau=(\imath,\cP)\times(\jmath,\cQ)\times \gamma \in
\PBP_\star(\check \CO, \wp) $ will be the painted bipartition in \eqref{tail0}, specified by the multiset $\{x_1, x_2, \cdots, x_k\}$ case by case as follows. 

\medskip

{\bfseries The case when $\star = B$.} In this case,  $\mathbf c_1(\imath)\leq \mathbf c_1(\jmath)$. The multiset $\{x_1, x_2, \cdots, x_k\}$ is the
union of the multiset
\[
  \set{\cQ(j,1)| \bfcc_{1}(\imath)+1 \leq j \leq  \bfcc_{1}(\jmath) }
\]
with the set
\[
  \begin{cases}
 \set{c}, &
 \qquad
  \text{if $\alpha = B^+$, and either $\bfcc_{1}(\imath)=0$ or $\cQ(\bfcc_{1}(\imath),1)\in \set{\bullet,s}$};  \\
 \set{s},&
  \qquad \text{if $\alpha = B^-$, and either $\bfcc_{1}(\imath)=0$ or $\cQ(\bfcc_{1}(\imath),1)\in \set{\bullet,s}$}; \\
\set{\cQ(\bfcc_{1}(\imath),1)},&
\qquad
\text{otherwise.}
\end{cases}
\]

{\bfseries The case when $\star = D$.} In this case,  $\mathbf c_1(\imath)>\mathbf c_1(\jmath)$ when $|\check \CO|>0$. 
The multiset $\{x_1, x_2, \cdots, x_k\}$
is the union of the multiset
\[
\set{\cP(j,1)| \bfcc_{1}(\jmath)+2 \leq j \leq \bfcc_{1}(\imath)}
\]
with the set
\[
\begin{cases}
 \set{d}, &\textrm{ if $|\check \CO|$=0};\\
    \set{c},                          &
    \ \text{if $\bfrr_2(\ckcO)=\bfrr_3(\ckcO)>0$, $\ \cP(\bfcc_1(\imath),1)\in \set{r,d}$,   and }\\   &\qquad ( \cP(\bfcc_{1}(\jmath)+1,1), \cP(\bfcc_{1}(\jmath)+1,2)) = (r,c);                                                                     \\
    \set{\cP(\bfcc_{1}(\jmath)+1,1)}, &
    \    \text{otherwise.}
  \end{cases}
\]

{\bfseries The case $\star = C^*$.}  In this case,  $\check \CO$ is quasi-distinguished,
$\wp=\emptyset$,
and $\mathbf c_1(\imath)\leq \mathbf c_1(\jmath)$.
The multiset $\{x_1, x_2, \cdots, x_k\}$ 
equals the multiset 
\[
  \set{\cQ(j,1)| \bfcc_{1}(\imath)+1 \leq j \leq  \bfcc_{1}(\jmath) }.
\]

We introduce one final notation when $\star \in \set{B,D}$: 
\[
x_\uptau := \cP_{\uptau_\bftt}(k,1), 
\]
which is the symbol in the last box of the tail $\uptau_\bftt$.

The following two propositions summarize key properties of the descent map. Both readily follow from the detailed description of the descent algorithm.

\begin{prop} \label{lem:PBPd.C}
  Suppose that  $\star\in \set{C,\wtC, D^*}$. 

 \noindent (a) If $\star = D^{*}$ or $\bfrr_1(\ckcO)>\bfrr_2(\ckcO)$,
    then the descent map \eqref{eq:DD.CC} is bijective. %Case $D^{*}$.

\noindent (b) If $\star\in \{C,\widetilde C\}$ and  $\bfrr_1(\ckcO)=\bfrr_2(\ckcO)$, then
    the descent map \eqref{eq:DD.CC} is injective with image
    \[
      \Set{\uptau'\in \PBPOPp |  x_{\uptau'}\neq s}.
    \]
\end{prop}

Recall that $\cO = \dBV(\ckcO)$ is the Barbasch-Vogan dual of $\ckcO$. 
Write 
\[\ckcOpp := \ckDD^2(\ckcO):=\ckDD(\ckcO') \AND \wp'':=\ckDD^2(\wp):= \ckDD(\wp').\]
Then we have the double descent  map
\[
\DD^2:=\DD\circ \DD\, \colon\,  \PBPOP \longrightarrow
    \PBPOPpp.
\]

\begin{prop}
\label{lem:delta}
Suppose that $\star \in \set{B,D,C^*}$ and $(\star, \abs{\check \CO})\neq (D, 0)$. 
Consider the map
\begin{equation}\label{eq:delta}
  \PBPOP \longrightarrow
    \PBPOPpp \times \PBP_{\star_\bftt}(\ckcO_\bftt),
    \qquad \uptau \mapsto (\DD^2(\uptau),\uptau_\bftt).
\end{equation}
\begin{enumerate}[label=(\alph*),wide=2pt]
\item 
Suppose that
$\star = C^*$ or $\bfrr_2(\ckcO)>\bfrr_3(\ckcO)$. Then the map \eqref{eq:delta} is bijective, and for every $\uptau\in  \PBP_\star(\ckcO) $,
\begin{equation*}\label{eq:sign.D}
\Sign(\uptau)
=(\bfcc_2(\cO),\bfcc_2(\cO))+\Sign(\DD^2(\uptau))+\Sign(\uptau_\bftt).
\end{equation*}

\item 
Suppose that  $\star \in \set{B,D}$ and $\bfrr_2(\ckcO)=\bfrr_3(\ckcO)>0$.
Then the map \eqref{eq:delta} is injective with image %equals
\begin{equation*}\label{eq:delta.I}
    \Set{ (\uptau'',\uptau_0)  \in \PBPOPpp \times \PBP_D(\ckcO_\bftt)  |
    \begin{array}{l}
        \text{either
    $x_{\uptau''} = d$, or} \\
    \text{$x_{\uptau''}\in \set{r,c}$  and
    $\cP_{\uptau_0}^{-1}(\set{s,c})\neq \emptyset$}
    \end{array}
},
\end{equation*}
and for every $\uptau\in  \PBPOP$,
\begin{equation*}\label{eq:sign.GD}
\Sign(\uptau)
=(\bfcc_2(\cO)-1,\bfcc_2(\cO)-1)+\Sign(\DD^2(\uptau))+\Sign(\uptau_\bftt).
\end{equation*}
\end{enumerate}
\end{prop}

The following corollary will be used in \cite{BMSZ2}. 
\begin{cor}\label{prop:DD.BDinj}
  Suppose that $\star \in \set{B, D,C^*}$.
 Denote 
  \[
    \varepsilon_{\uptau} := \begin{cases}
      0, & \text{if $\star\in \set{B,D}$ and $x_{\uptau}=d$;}\\
      1, & \text{otherwise.}
    \end{cases}
  \]
  Then the map
\begin{equation*}
  \begin{array}{rcl}
   \PBPOP&\rightarrow&
   \PBPOPp \times \BN\times \bN\times \Z/2\Z, \smallskip\\
   \uptau& \mapsto & (\DD(\uptau), p_\uptau, q_\uptau, \varepsilon_\uptau)
   \end{array}
\end{equation*}
is injective.
\end{cor}
\begin{proof}
This is easy to verify in the special  case when $\bfrr_{3}(\ckcO)=0$ and
    \[
      (\bfrr_{1}(\ckcO), \bfrr_{2}(\ckcO))
      =
      \begin{cases}
        (2k-2,0), & \text{if } \star=B;\\
        (2k-1,1), & \text{if } \star=D;\\
        (2k-1,0), & \text{if } \star=C^{*}.\\
      \end{cases}
    \]
    In general, it follows from the injectivity of $\DD$ in \Cref{lem:PBPd.C}
    and the signature formula in \Cref{lem:delta}.
\end{proof}

\subsection{Generating functions}\label{sec: generating}

We introduce some generating functions to count painted bipartitions. 
When $\star \in  \set{B,D, C^{*}}$, we define
\[
   f_{\star,\ckcO,\wp}(\bfpp,\bfqq) := \sum_{\uptau \in \PBPOP} \bfpp^{p_{\uptau}} \bfqq^{q_{\uptau}},
\]
where $(p_{\uptau}, q_{\uptau})$ is as in \eqref{ptqt}.  
 This is an element in the polynomial ring $\bZ[\bfpp,\bfqq]$, where $\bfpp$ and $\bfqq$ are indeterminants. The coefficient of $\bfpp^{p}\bfqq^{q}$ in $f_{\star,\ckcO,\wp}(\bfpp,\bfqq)$  equals the cardinality of
 $\PBPop{G}{}{\ckcO}{\wp}$ (note that $G=\SO(p,q)$ or $\Sp(\frac{p}{2},\frac{q}{2})$).

When $\star\in \{B, D\}$, for each subset $S\subset \set{c,d,r,s}$, we also define
\[
  \PBPOP[S] := \set{\uptau\in \PBPOP|x_{\uptau}\in S}
\]
and the corresponding generating function
\[
   f_{\star,\ckcO,\wp}^{S}(\bfpp,\bfqq) := \sum_{\uptau \in \PBPOP[S]} \bfpp^{p_{\uptau}} \bfqq^{q_{\uptau}}. 
\]
Note that \[
  f_{D,\yrow{0}, \emptyset}^{\set{d}} := 1, \quad
  f_{D,\yrow{0}, \emptyset}^{\set{c,r}} := 0,
  \AND
  f_{D,\yrow{0}, \emptyset}^{\set{s}} := 0.
\]
In the following we will compute
$f_{\star,\ckcO,\wp}^{\set{d}}$, $f_{\star,\ckcO,\wp}^{\set{c,r}}$, and
$f_{\star,\ckcO,\wp}^{\set{s}}$ whose sum is the desired function
$f_{\star, \ckcO,\wp}$.

\def\CSk#1#2{h_{#1}^{#2}}
\def\TSk#1#2{g_{#1}^{#2}}
\def\RS{\nu}

For every integer $k$, we define 
\[
\RS_{k} := \begin{cases}
    \sum_{i=0}^{k} \bfpp^{2i}\bfqq^{2(k-i)}, \quad &\textrm{if $k\geq 0$};\\
    0,\quad & \textrm{if $k<0$.}
\end{cases}
\]

It is straightforward to check the following identities: for every $k\in \BN$,
\[
  f_{B,\yrow{2k},\emptyset}^{S}  = \begin{cases}
      (\bfpp^{2}\bfqq+ \bfpp\bfqq^{2})\RS_{k-1}, \quad &\textrm{if $S=\set{d}$};\\
      \bfpp\RS_{k} + \bfpp^{2}\bfqq\RS_{k-1}, \quad &\textrm{if $S=\set{c,r}$};\\
       \bfqq^{2k+1}, \quad &\textrm{if $S=\set{s}$}, 
       \end{cases}
\]
and for every $k\in \BN^+$,  
\[
 f_{D,\yrow{2k-1,1},\emptyset}^{S} = \begin{cases}
      \bfpp\bfqq\RS_{k-1} + \bfpp^{2}\bfqq^{2}\RS_{k-2}, \quad &\textrm{if $S=\set{d}$};\\
      (\bfpp^{2}+\bfpp\bfqq)\RS_{k-1}, \quad &\textrm{if $S=\set{c,r}$};\\
      \bfqq^{2k}, \quad &\textrm{if $S=\set{s}$}. 
      \end{cases}
\]

For every $k\in \BN$, we also define
\[
  \begin{split}
    \CSk{k}{\set{d}} &:= (\bfpp^{2}\bfqq^{2} + \bfpp\bfqq^{3})\RS_{k-2}, \\
  \CSk{k}{\set{c,r}} &:= \bfpp\bfqq\RS_{k-1}+\bfpp^{2}\bfqq^{2}\RS_{k-2},\\
  \CSk{k}{\set{s}} &:= \bfqq^{2k}.\\
  \end{split}
\]

\begin{prop}\label{prop:gf.BD}
  Suppose that $\star\in\set{B,D}$ and 
  $\bfrr_{2}(\ckcO)>0$. Let $k:=\frac{\bfrr_{1}(\ckcO)-\bfrr_{2}(\ckcO)}{2}+1$ and suppose that $S = \set{d}, \set{c,r}$, or $ \set{s}$. 

\begin{enumerate}[label=(\alph*),wide=0pt]
\item If $(2,3)\in \CPPs(\ckcO)$, then
    \[
      f_{\star,\ckcO, \wp}^{S} = (\bfpp\bfqq)^{\bfcc_{1}(\ckcO)} 
      f_{D,\yrow{2k-1,1},\emptyset}^{S}
        % \TSk{k}{S} 
      \cdot  f_{\star,\ckDD^{2}(\ckcO),\ckDD^{2}(\wp)}.
    \]
  
\item  If $(2,3)\notin \CPPs(\ckcO)$,
    then
    \[
      f_{\star,\ckcO, \wp}^{S} := (\bfpp\bfqq)^{\bfcc_{1}(\ckcO)-1} (
      f_{D,\yrow{2k-1,1},\emptyset}^{S}
      %\TSk{k}{S} 
      \cdot f_{\star,\ckDD^{2}(\ckcO),\ckDD^{2}(\wp)}^{\set{d}}
      + \CSk{k}{S} \cdot f_{\star,\ckDD^{2}(\ckcO),\ckDD^{2}(\wp)}^{\set{c,r}}).
    \]
\end{enumerate}
\end{prop}
\begin{proof}
  This follows from \Cref{lem:delta}, after a routine calculation. 
\end{proof}

\begin{prop}\label{prop:gf.C}
  Suppose that $\star\in\set{C,\wtC}$ and 
  $\bfrr_{1}(\ckcO)>0$.
  \begin{enumerate}[label=(\alph*),wide=0pt]
    \item If $(1,2)\in \CPPs(\ckcO)$, then
    \[
      \# \PBPOP =  f_{\star',\ckDD(\ckcO), \ckDD(\wp)}(1,1).
    \]
    \item If $(1,2)\notin \CPPs(\ckcO)$,
    then
    \[
      \# \PBPOP = f_{\star',\ckDD(\ckcO), \ckDD(\wp)}^{\set{c,r}}(1,1) + f_{\star',\ckDD(\ckcO), \ckDD(\wp)}^{\set{d}}(1,1).
    \]
  \end{enumerate}
 \end{prop}
 
\begin{proof}
  This follows from   
  \Cref{lem:PBPd.C}, after a routine calculation.  
\end{proof}

As a consequence of \Cref{prop:gf.BD} and \Cref{prop:gf.C},
when $\star \in \set{B,C,\wtC, D}$, the generating functions $f_{\star,\ckcO,\wp}$ are independent of $\wp$ and so
  \[
    f_{\star,\ckcO,\wp} = f_{\star,\ckcO,\emptyset}, \quad \text{for all $\wp\subseteq \CPPs(\ckcO)$.}
  \]
Consequently, 
\[
    \# \PBPOP = \# \PBPop{\star}{}{\ckcO}{\emptyset}, \quad \text{for all $\wp\subseteq  \CPPs(\ckcO)$.}
\]
This proves \Cref{prop:PBP2}.

The following two results are easy to check using \Cref{lem:PBPd.C} and \Cref{lem:delta}.

\begin{prop}\label{lem:gf.C*}
  Suppose that $\star = C^{*}$. Then
  \[
      f_{\star, \ckcO,\emptyset} = 
     \begin{cases}
     \nu_{\frac{\bfrr_1(\ckcO)-1}{2}}, &\quad \text{if }\bfrr_2(\ckcO) = 0;  \\
     (\bfpp\bfqq)^{\bfrr_{2}(\ckcO)+1} \RS_{k}\cdot f_{\star,\ckDD^{2}(\ckcO),\emptyset},
     & \quad \text{if } \bfrr_2(\ckcO)>0. \\
     \end{cases} 
  \]
  where $k := \frac{\bfrr_{1}(\ckcO)-\bfrr_{2}(\ckcO)}{2}-1$.
\end{prop}

\begin{prop}\label{lem:gf.D*}
  Suppose that $\star = D^{*}$. Then 
  \[
    \# \PBP_G(\ckcO,\emptyset) = f_{\star',\ckDD(\ckcO),\emptyset}(1,1).
  \]
  \end{prop}

The following result is clear by \Cref{lem:gf.C*} and \Cref{lem:gf.D*}. 

\begin{prop}\label{prop:count.CD*}
  Suppose that $\star\in \set{C^*, D^*}$. % and $\ckcO$ is quasi-distinguished.
  Then the cardinality of $\PBP_{G}(\ckcO)$ equals the number of $G$-orbits in 
  $(\sqrt{-1}\g_0^*)\cap \CO$, where $\g_0$ denotes the Lie algebra of $G$ and $\CO:=\dBV(\ckcO)$.
\end{prop}

\section*{Acknowledgements}
The authors have learned much about  coherent continuation representations  from David Vogan, and are grateful to him for many helpful  discussions/clarifications. 
The authors also thank Devra Garfinkle Johnson for fruitful discussions  concerning \Cref{thm:LC.ABC} and \Cref{rem:LC.ABC}.  
The authors are further grateful to the referees for their helpful comments and detailed suggestions.

D. Barbasch is supported by NSF grant, Award Number 2000254. 
J.-J. Ma is supported by the National Natural Science Foundation of China (Grant No. 11701364 and Grant No. 11971305), the Fundamental Research Funds for the Central Universities (Grant No. 20720230022) and  Xiamen University
Malaysia Research Fund (Grant No. XMUMRF/2022-C9/IMAT/0019).
B. Sun is supported by  National Key R \& D Program of China (No. 2022YFA1005300 and 2020YFA0712600) and New Cornerstone Investigator Program.  C.-B. Zhu is supported by MOE AcRF Tier 1 grants A-0004280-00-00 and A-8002490-00-00, and
Provost’s Chair grant E-146-00-0018-01 in NUS. 

C.-B. Zhu is grateful to Max Planck Institute for Mathematics in Bonn, for its warm hospitality and conducive work environment, where he spent the academic year 2022/2023 as a visiting scientist.

%\appendix


\begin{bibdiv}
  \begin{biblist}

\bib{AAM}{article}{
       author = {Adams, J.},
       author = {Arancibia Robert, N.},
       author = {Mezo, P.},
        title = {Equivalent definitions of Arthur packets for real classical groups},
        journal = {Mem. Amer. Math. Soc.},
        volume = {300},
    number = {1503},
    year = {2024},
    pages ={vi+110 pp.},
    }

\bib{ABV}{book}{
  title={The Langlands classification and irreducible characters for real reductive groups},
  author={Adams, J.},
  author={Barbasch, D.},
  author={Vogan, D. A.},
  series={Progr. Math.},
  volume={104},
  year={1991},
  publisher={Birkh\"auser}
}


\bib{AC}{article}{
  title={Algorithms for representation theory of
    real reductive groups},
  volume={8},
  %DOI={10.1017/S1474748008000352},
  number={2},
  journal={J. Inst. Math. Jussieu},
  publisher={Cambridge University Press},
  author={Adams, J.},
  author={du Cloux, F.},
  year={2009},
  pages={209--259}
}

\bib{AM}{article}{
       author = {Arancibia Robert, N.},
       author = {Mezo, P.},
        title = {Equivalent definitions of Arthur packets for real unitary groups},
        eprint = {arXiv:2204.19715},
}


\bib{ArPro}{article}{
author = {Arthur, J.},
title = {On some problems suggested by the trace formula},
bookTitle = {Lie group representations, II (College Park, Md.), Lecture Notes in Math., vol. 1041},
pages = {1--49},
year = {1984}
}


\bib{ArUni}{article}{
author = {Arthur, J.},
title = {Unipotent automorphic representations: conjectures},
journal = {Orbites unipotentes et repr\'esentations, II, Ast\'erisque},
pages = {13--71},
volume = {171-172},
year = {1989}
}

\bib{ArEnd}{article}{
    author = {Arthur, J.},
    title = {The Endoscopic Classification of Representations: Orthogonal and Symplectic Groups},
    journal = {Amer. Math. Soc. Colloq. Publ., vol. 61, Amer. Math. Soc., Providence, RI},
    year = {2013},
    }



\bib{B89}{article}{
  author = {Barbasch, D.},
  title = {The unitary dual for complex classical Lie groups},
  journal = {Invent. Math.},
  volume = {96},
  number = {1},
 pages = {103--176},
 year = {1989}
}

\bib{B.Orbit}{article}{
  author = {Barbasch, D.},
  title = {Orbital integrals of nilpotent orbits},
    journal = {The mathematical legacy of {H}arish-{C}handra, Proc. Sympos. Pure Math.},
    volume = {68},
     pages = {97--110},
 publisher = {Amer. Math. Soc., Providence, RI},
      year = {2000},
}


\bib{B17}{article}{
  author = {Barbasch, D.},
  title = {Unipotent representations and the dual pair correspondence},
  bookTitle = {J. Cogdell et al. (eds.), Representation Theory, Number Theory, and Invariant Theory, In Honor of Roger Howe. Progr. Math., vol. 323, Birkh\"auser},
  %series ={Progress in Math.},
  pages = {47--85},
  year = {2017},
}





\bib{BMSZ2}{article}{
      title={Special unipotent representations of real classical groups: construction and unitarity},
      author={Barbasch, D.},
      author = {Ma, J.-J.},
      author = {Sun, B.},
      author = {Zhu, C.-B.},
      year={2021},
      eprint={arXiv:1712.05552},
}

\bib{BMSZ0}{article}{
      title={On the notion of metaplectic Barbasch-Vogan duality},
      %year={2020},
      author={Barbasch, D.},
      author = {Ma, J.-J.},
      author = {Sun, B.},
      author = {Zhu, C.-B.},
      %volume = {2023},
      number = {20},
    pages = {17822-17852},
    year = {2023},
      month = {05},
      journal={Int. Math. Res. Not. IMRN},
}


\bib{BV1}{article}{
   author={Barbasch, D.},
   author={Vogan, D. A.},
   title={Primitive ideals and orbital integrals in complex classical
   groups},
   journal={Math. Ann.},
   volume={259},
   date={1982},
   number={2},
   pages={153--199},
   issn={0025-5831},
 %  review={\MR{656661}},
   %doi={10.1007/BF01457308},
}

\bib{BV2}{article}{
   author={Barbasch, D.},
   author={Vogan, D. A.},
   title={Primitive ideals and orbital integrals in complex exceptional
   groups},
   journal={J. Algebra},
   volume={80},
   date={1983},
   number={2},
   pages={350--382},
   issn={0021-8693},
 %  review={\MR{691809}},
   %doi={10.1016/0021-8693(83)90006-6},
}

\bib{BV.W}{article}{
  author={Barbasch, D.},
  author={Vogan, D. A.},
  editor={Trombi, P. C.},
  title={Weyl Group Representations and Nilpotent Orbits},
  bookTitle={Representation Theory of Reductive Groups:
    Proc. Univ. Utah Conference (1982)},
  year={1983},
  publisher={Birkh{\"a}user Boston},
  address={Boston, MA},
  pages={21--33},
  %doi={10.1007/978-1-4684-6730-7_2},
}


\bib{BVUni}{article}{
 author = {Barbasch, D.},
 author = {Vogan, D. A.},
 journal = {Ann. of Math.},
 number = {1},
 pages = {41--110},
 title = {Unipotent representations of complex semisimple groups},
 volume = {121},
 year = {1985}
}

\bib{BGG.M}{article}{
   author={Bernstein, I. N.},
   author={Gel'fand, I. M.},
   author={Gel'fand, S. I.},
   title={Models of representations of compact Lie groups},
   language={Russian},
   journal={Funkcional. Anal. i Prilo\v{z}en.},
   volume={9},
   date={1975},
   number={4},
   pages={61--62},
   %issn={0374-1990},
   %review={\MR{0414792}},
}

\delete{
\bib{Bor}{article}{
 author = {Borho, W.},
 journal = {S\'eminaire Bourbaki, Exp. No. 489},
 pages = {1--18},
 title = {Recent advances in enveloping algebras of semisimple Lie-algebras},
 year = {1976/77}
}}

\bib{Borel}{article}{
  title={Automorphic $L$-functions},
  author={Borel, A.},
  book = {
    title={Automorphic Forms, Representations and $L$-functions},
    series={Proc. Sympos. Pure Math.},
    volume={33, Part 2},
  },
  pages={27--61},
  year={1979},
  publisher={Amer. Math. Soc.},
}

\bib{Bor}{article}{
 author = {Borho, W.},
 journal = {S\'eminaire Bourbaki, Exp. No. 489},
 pages = {1--18},
 title = {Recent advances in enveloping algebras of semisimple Lie-algebras},
 year = {1976/77}
}


\bib{BK}{article}{
author={Borho, W.},
author={Kraft, H.},
title={\"{U}ber die Gelfand-Kirillov-Dimension},
journal={Math. Ann.},
volume={220},
date={1976},
number={1},
pages={1--24},
issn={0025-5831},
% review={\MR{412240}},
%doi={10.1007/BF01354525},
}


\bib{BM}{article}{
   author={Borho, W.},
   author={MacPherson, R.},
   title={Repr\'{e}sentations des groupes de Weyl et homologie d'intersection
   pour les vari\'{e}t\'{e}s nilpotentes},
   language={French, with English summary},
   journal={C. R. Acad. Sci. Paris S\'{e}r. I Math.},
   volume={292},
   date={1981},
   number={15},
   pages={707--710},
   issn={0249-6291},
  % review={\MR{618892}},
}


\bib{Bo}{article}{
   author={Bo\v{z}i\v{c}evi\'{c}, M.},
   title={Double cells for unitary groups},
   journal={J. Algebra},
   volume={254},
   date={2002},
   number={1},
   pages={115--124},
   issn={0021-8693},
 %  review={\MR{1927434}},
   %doi={10.1016/S0021-8693(02)00070-4},
}


\bib{Carter}{book}{
   author={Carter, R. W.},
   title={Finite groups of Lie type},
   series={Wiley Classics Library},
   %note={Conjugacy classes and complex characters;
   %Reprint of the 1985 original;
   %A Wiley-Interscience Publication},
   publisher={John Wiley \& Sons, Ltd., Chichester},
   date={1993},
   pages={xii+544},
   isbn={0-471-94109-3},
   %review={\MR{1266626}},
}

\bib{Cas}{article}{
   author={Casian, L. G.},
   title={Primitive ideals and representations},
   journal={J. Algebra},
   volume={101},
   date={1986},
   number={2},
   pages={497--515},
   issn={0021-8693},
  % review={\MR{847174}},
   %doi={10.1016/0021-8693(86)90208-5},
}



\bib{CM}{book}{
  title = {Nilpotent orbits in semisimple Lie algebra: an introduction},
  author = {Collingwood, D. H.},
  author = {McGovern, W. M.},
  year = {1993},
  publisher = {Van Nostrand Reinhold Co.},
}

\bib{Dix}{book}{
  title={Enveloping algebras},
  author={Dixmier, J.},
  year={1996},
  publisher={Grad. Stud. Math., vol. 11, Amer. Math. Soc.},
}


\bib{Du77}{article}{
  author = {Duflo, M.},
  journal = {Ann. of Math.},
  number = {1},
  pages = {107-120},
  title = {Sur la Classification des Ideaux Primitifs Dans
    L'algebre Enveloppante d'une Algebre de Lie Semi-Simple},
  volume = {105},
  year = {1977}
}


\bib{FJMN}{article}{
   author={Folz{-Donahue}, T.},
   author={Jackson, S. G.},
   author={Milev, T.},
   author={No\"{e}l, A. G.},
   title={Sign signatures and characters of Weyl Groups},
   journal={Adv. in Appl. Math.},
   volume={130},
   date={2021},
   pages={Paper No. 102225, 73},
   issn={0196-8858},
  % review={\MR{4272991}},
   %doi={10.1016/j.aam.2021.102225},
}

\bib{GI}{article}{
   author={Gan, W. T.},
   author={Ichino, A.},
   title={On the irreducibility of some induced representations of real
   reductive Lie groups},
   journal={Tunis. J. Math.},
   volume={1},
   date={2019},
   number={1},
   pages={73--107},
   issn={2576-7658},
  % review={\MR{3907735}},
  % doi={10.2140/tunis.2019.1.73},
}

\bib{G3}{article}{
   author={Garfinkle, D.},
   title={On the classification of primitive ideals for complex classical
   Lie algebras. III},
   journal={Compos. Math.},
   volume={88},
   date={1993},
   number={2},
   pages={187--234},
   issn={0010-437X},
   % review={\MR{1237920}},
}

\bib{GV}{article}{
   author={Garfinkle, D.},
   author={Vogan, D. A.},
   title={On the structure of Kazhdan-Lusztig cells for branched Dynkin
   diagrams},
   journal={J. Algebra},
   volume={153},
   date={1992},
   number={1},
   pages={91--120},
   issn={0021-8693},
   %review={\MR{1195408}},
   %doi={10.1016/0021-8693(92)90150-K},
}

\bib{Gar.D4}{article}{
    author = {Garfinkle Johnson, D.},
    title = {Edge Transport from Parabolic Subgroups of Type $D_4$},
    year = {2019},
    eprint={arXiv:1907.09717},
}




\bib{GW}{book}{
   author={Goodman, R.},
   author={Wallach, N. R.},
   title={Symmetry, representations, and invariants},
   series={Grad. Texts in Math.},
   volume={255},
   publisher={Springer, Dordrecht},
   date={2009},
   pages={xx+716},
   isbn={978-0-387-79851-6},
   %review={\MR{2522486}},
   %doi={10.1007/978-0-387-79852-3},
}


\bib{Ho}{article}{
author={Hotta, R.},
title={On Joseph's construction of Weyl group representations},
journal={Tohoku Math. J.},
volume={36},
%number = {3},
pages={49--74 },
year={1984},
}


\bib{Howe79}{article}{
  title={$\theta$-series and invariant theory},
  author={Howe, R.},
  book = {
    title={Automorphic Forms, Representations and $L$-functions},
    series={Proc. Sympos. Pure Math.},
    volume={33, Part 1},
  },
  pages={275--285},
  year={1979}, 
  publisher={Amer. Math. Soc.},
}



\bib{Howe89}{article}{
author={Howe, R.},
title={Transcending classical invariant theory},
journal={J. Amer. Math. Soc.},
volume={2},
pages={535--552},
year={1989},
}


 \bib{H}{book}{
   author={Humphreys, J. E.},
    title={Representations of semisimple Lie algebras in the BGG category
    $\scr{O}$},
    series={Grad. Stud. Math.},
    volume={94},
   publisher={Ameri. Math. Soc., Providence, RI},
   date={2008},
   pages={xvi+289},
  isbn={978-0-8218-4678-0},
  %  review={\MR{2428237}},
    doi={10.1090/gsm/094},
}


\bib{Jan}{book}{
   author={Jantzen, J. C.},
   title={Moduln mit einem h\"{o}chsten Gewicht},
   series={Lecture Notes in Math.},
   volume={750},
   publisher={Springer-Verlag, Berlin/Heidelberg/New York},
   date={1979},
  % pages={xvi+289},
   % isbn={978-0-8218-4678-0},
   % review={\MR{2428237}},
   % doi={10.1090/gsm/094},
}



\delete{
\bib{J79W}{article}{
   author={Joseph, A.},
   title={$W$-module structure in the primitive spectrum of the enveloping
   algebra of a semisimple Lie algebra},
   conference={
      title={Noncommutative harmonic analysis},
      address={Proc. Third Colloq., Marseille-Luminy},
      date={1978},
   },
   book={
      series={Lecture Notes in Math.},
      volume={728},
      publisher={Springer, Berlin},
   },
   date={1979},
   pages={116--135},
  % review={\MR{548328}},
}
}

 \bib{J79D}{article}{
   author={Joseph, A.},
   title={Dixmier's problem for Verma and principal series submodules},
   journal={J. Lond. Math. Soc. (2)},
   volume={20},
   date={1979},
   number={2},
   pages={193--204},
   issn={0024-6107},
  % review={\MR{551445}},
  % doi={10.1112/jlms/s2-20.2.193},
}
	

\bib{J1}{article}{
   author={Joseph, A.},
   title={Goldie rank in the enveloping algebra of a semisimple Lie algebra. I},
   journal={J. Algebra},
   volume={65},
   date={1980},
   number={2},
   pages={269--283},
   issn={0021-8693},
  % review={\MR{585721}},
  % doi={10.1016/0021-8693(80)90217-3},
}

\bib{J2}{article}{
   author={Joseph, A.},
   title={Goldie rank in the enveloping algebra of a semisimple Lie algebra. II},
   journal={J. Algebra},
   volume={65},
   date={1980},
   number={2},
   pages={284--306},
   issn={0021-8693},
  % review={\MR{585721}},
  % doi={10.1016/0021-8693(80)90217-3},
}

\bib{J.hw}{article}{
   author={Joseph, A.},
   title={On the variety of a highest weight module},
   journal={J. Algebra},
   volume={88},
   date={1984},
   number={1},
   pages={238--278},
   issn={0021-8693},
  % review={\MR{741942}},
  % doi={10.1016/0021-8693(84)90100-5},
}


\bib{J.av}{article}{
   author={Joseph, A.},
   title={On the associated variety of a primitive ideal},
   journal={J. Algebra},
   volume={93},
   date={1985},
   number={2},
   pages={509--523},
   issn={0021-8693},
  % review={\MR{786766}},
 % doi={10.1016/0021-8693(85)90172-3},
}

 \bib{KL}{article}{
    author={Kazhdan, D.},
    author={Lusztig, G.},
    title={Representations of Coxeter groups and Hecke algebras},
    journal={Invent. Math.},
    volume={53},
    date={1979},
    number={2},
    pages={165--184},
    issn={0020-9910},
  %  review={\MR{560412}},
   % doi={10.1007/BF01390031},
 }

 
\bib{Ki62}{article}{
author={Kirillov, A. A.},
title={Unitary representations of nilpotent Lie groups},
journal={Uspehi Mat. Nauk},
volume={17},
issue ={4},
pages={57--110},
year={1962},
}


\bib{Kn}{book}{
title={Lie groups beyond an introduction},
author={Knapp, A. W.},
year={2002},
publisher={Progr. Math., vol. 140, Birkhäuser Boston},
}



\bib{Ko70}{article}{
author={Kostant, B.},
title={Quantization and unitary representations},
bookTitle ={Lectures in Modern Analysis and Applications III, Lecture Notes in Math., vol. 170},
pages={87--208},
year={1970},
}


\bib{Lsp}{article}{
   author={Lusztig, G.},
   title={A class of irreducible representations of a Weyl group},
   journal={Nederl. Akad. Wetensch. Indag. Math.},
   volume={41},
   date={1979},
   number={3},
   pages={323--335},
   issn={0019-3577},
  % review={\MR{546372}},
}
	

\bib{Lcell}{article}{
   author={Lusztig, G.},
   title={A class of irreducible representations of a Weyl group. II},
   journal={Nederl. Akad. Wetensch. Indag. Math.},
   volume={44},
   date={1982},
   number={2},
   pages={219--226},
   issn={0019-3577},
  % review={\MR{662657}},
}


\bib{Lu}{book}{
   author={Lusztig, G.},
   title={Characters of reductive groups over a finite field},
   series={Ann. of Math. Stud.},
   volume={107},
   publisher={Princeton University Press, Princeton, NJ},
   date={1984},
   pages={xxi+384},
   isbn={0-691-08350-9},
   isbn={0-691-08351-7},
  % review={\MR{742472}},
   doi={10.1515/9781400881772},
}



\bib{Mat96}{article}{
   author={Matumoto, H.},
   title={On the representations of ${\rm U}(m,n)$ unitarily induced from
   derived functor modules},
   journal={Compos. Math.},
   volume={100},
   date={1996},
   number={1},
   pages={1--39},
   issn={0010-437X},
  % review={\MR{1377407}},
}

\bib{Mat04}{article}{
    title={On the representations of $\mathrm{Sp}(p, q)$ and $\mathrm{SO}^*(2n)$ unitarily induced from derived functor modules},
    volume={140},
    %DOI={10.1112/S0010437X03000629}, 
    number={4},
    journal={Compos. Math.}, 
    author={Matumoto, H.},
    year={2004},
    pages={1059–1096}
    }


\bib{Mc}{article}{
   author = {McGovern, W. M.},
    title = {Cells of Harish-Chandra modules for real classical groups},
    journal = {Amer. J. Math.},
    volume = {120},
    issue = {01},
    year = {1998},
    pages ={211--228},
}


\bib{GJMP}{article}{
      title={On the classification of primitive ideals for complex classical Lie algebras, IV}, 
      %author={Garfinkle Johnson, D.},
      author={McGovern, W. M.},
      author = {Pietraho, T.},
      year={2023},
      eprint={arXiv:2309.02363},
      archivePrefix={arXiv},
}

\bib{Mo11}{article}{
   author={M{\oe}glin, C.},
   title={Multiplicit\'{e} 1 dans les paquets d'Arthur aux places $p$-adiques},
   language={French, with English summary},
   bookTitle={On certain $L$-functions, Clay Math. Proc., vol. 13, Amer. Math. Soc., Providence, RI},
   year={2011},
   pages={333--374},
}
	

\bib{Mo17}{article}{
  author={M{\oe}glin, C.},
  title = {Paquets d'Arthur Sp\'eciaux Unipotents aux Places Archim\'ediennes et Correspondance de Howe},
  bookTitle = {J. Cogdell et al. (eds.), Representation Theory, Number Theory, and Invariant Theory, In Honor of Roger Howe. Progr. Math., vol. 323, Birkh\"auser},
  %series ={Progress in Math.},
  pages = {469--502},
  year = {2017},
}

\bib{MR.C}{article}{
   author={M{\oe}glin, C.},
   author={Renard, D.},
   title={Paquets d'Arthur des groupes classiques complexes},
   language={French, with English and French summaries},
   conference={
      title={Around Langlands correspondences},
   },
   book={
      series={Contemp. Math.},
      volume={691},
      publisher={Amer. Math. Soc., Providence, RI},
   },
   date={2017},
   pages={203--256},
  % review={\MR{3666056}},
  % doi={10.1090/conm/691/13899},
}

\bib{MR.U}{article}{
   author={M{\oe}glin, C.},
   author={Renard, D.},
   title={Sur les paquets d'Arthur des groupes unitaires et quelques
   cons\'{e}quences pour les groupes classiques},
   language={French, with English and French summaries},
   journal={Pacific J. Math.},
   volume={299},
   date={2019},
   number={1},
   pages={53--88},
   issn={0030-8730},
  % review={\MR{3947270}},
  % doi={10.2140/pjm.2019.299.53},
}

\delete{
\bib{MR}{article}{
  author={M{\oe}glin, C.},
    author={Renard, D.},
  title = {Sur les paquets d'Arthur des groupes classiques r\'eels},
  journal = {J. Eur. Math. Soc. },
  volume = {22},
    issue = {6},
    year = {2020},
    pages ={1827--1892}
    }
}

\bib{Mok}{article}{
author={Mok, C. P.},
  title = {Endoscopic classification of representations of quasi-split unitary groups},
  journal = {Mem. Amer. Math. Soc.},
  volume = {235},
    number = {1108},
    year = {2015},
   pages ={vi+248 pp.},
    }


\bib{RT1}{article}{
   author={Renard, D.},
   author={Trapa, P.},
   title={Irreducible genuine characters of the metaplectic group:
   Kazhdan-Lusztig algorithm and Vogan duality},
   journal={Represent. Theory},
   volume={4},
   date={2000},
   pages={245--295},
  % review={\MR{1795754}},
  % doi={10.1090/S1088-4165-00-00105-9},
}


\bib{RT2}{article}{
   author={Renard, D.},
   author={Trapa, P.},
   title={Irreducible characters of the metaplectic group. II.
   Functoriality},
   journal={J. Reine Angew. Math.},
   volume={557},
   date={2003},
   pages={121--158},
   issn={0075-4102},
  % review={\MR{1978405}},
  % doi={10.1515/crll.2003.028},
}



\bib{Sch}{article}{
  author={Schmid, W.},
  title = {Two character identities for semisimple Lie groups},
  bookTitle = {Non-commutative harmonic analysis (Actes Colloq., Marseille-Luminy, 1976).
Lecture Notes in Math., Vol. 587, Springer, Berlin},
  pages = {196--225},
  year = {1977},
}

\bib{SV}{article}{
author = {Schmid, W.},
author = {Vilonen, K.},
journal = {Ann. of Math.},
number = {3},
pages = {1071--1118},
%publisher = {Princeton University, Mathematics Department, Princeton, NJ; Mathematical Sciences Publishers, Berkeley},
title = {Characteristic cycles and wave front cycles of representations of reductive Lie groups},
volume = {151},
year = {2000},
}

\bib{Sho}{article}{
author = {Shoji, T.},
title = {On the Springer representations of the Weyl groups of classical algebraic groups},
journal = {Comm. Algebra},
volume = {7},
pages = {1713--1745},
year = {1979},
}


\bib{Soergel}{article}{
   author={Soergel, W.},
   title={Kategorie $\scr O$, perverse Garben und Moduln \"{u}ber den
   Koinvarianten zur Weylgruppe},
   language={German, with English summary},
   journal={J. Amer. Math. Soc.},
   volume={3},
   date={1990},
   number={2},
   pages={421--445},
   issn={0894-0347},
  % review={\MR{1029692}},
  % doi={10.2307/1990960},
}


\bib{So}{article}{
author = {Sommers, E.},
title = {Lusztig's canonical quotient and generalized duality},
journal = {J. Algebra},
volume = {243},
number = {2},
pages = {790--812},
year = {2001},
}


\bib{Spa}{book}{
   author={Spaltenstein, N.},
   title={Classes unipotentes et sous-groupes de Borel},
   series={Lecture Notes in Math.},
   volume={946},
   publisher={Springer-Verlag, Berlin-New York},
   date={1982},
}


\bib{SpVo}{article}{
author = {Speh, B.},
author = {Vogan, D. A.},
title = {Reducibility of generalized principal series representations},
journal = {Acta Math.},
volume = {145},
number = {3--4},
pages = {227--299},
year = {1980},
}


\bib{Spr}{article}{
author = {Springer, T. A.},
title = {A construction of representations of Weyl groups},
journal = {Invent. Math.},
volume = {44},
pages = {279--293},
year = {1978},
}


\bib{Tr.U}{article}{
   author={Trapa, P.},
   title={Annihilators and associated varieties of $A_{\mathfrak q}(\lambda)$
   modules for $\mathrm U(p,q)$},
   journal={Compos. Math.},
   volume={129},
   date={2001},
   number={1},
   pages={1--45},
   issn={0010-437X},
  % review={\MR{1856021}},
   %doi={10.1023/A:1013115223377},
}

\bib{Tr.H}{article}{
  title={Special unipotent representations and the Howe correspondence},
  author={Trapa, P.},
  year = {2004},
  journal={Univ. Aarhus Publ. Series},
  volume = {47},
  pages= {210--230}
}

\bib{VGK}{article}{
   author={Vogan, D. A.},
   title={Gel\cprime fand-Kirillov dimension for Harish-Chandra modules},
   journal={Invent. Math.},
   volume={48},
   date={1978},
   number={1},
   pages={75--98},
   issn={0020-9910},
  % review={\MR{506503}},
  % doi={10.1007/BF01390063},
}

\bib{Vtau}{article}{
   author={Vogan, D. A.},
   title={A generalized $\tau $-invariant for the primitive spectrum of a
   semisimple Lie algebra},
   journal={Math. Ann.},
   volume={242},
   date={1979},
   number={3},
   pages={209--224},
   issn={0025-5831},
  % review={\MR{545215}},
  % doi={10.1007/BF01420727},
}
		

\bib{V1}{article}{
   author={Vogan, D. A.},
   title={Irreducible characters of semisimple Lie groups. I},
   journal={Duke Math. J.},
   volume={46},
   date={1979},
   number={1},
   pages={61--108},
   issn={0012-7094},
  % review={\MR{523602}},
}

\bib{Vg}{book}{
   author={Vogan, D. A.},
   title={Representations of real reductive Lie groups},
   series={Progr. Math.},
   volume={15},
   publisher={Birkh\"{a}user, Boston, Mass.},
   date={1981},
   pages={xvii+754},
   isbn={3-7643-3037-6},
  % review={\MR{632407}},
}


\bib{V4}{article}{
   author={Vogan, D. A.},
   title={Irreducible characters of semisimple Lie groups. IV.
   Character-multiplicity duality},
   journal={Duke Math. J.},
   volume={49},
   date={1982},
   number={4},
   pages={943--1073},
   issn={0012-7094},
  % review={\MR{683010}},
}

\bib{V3}{article}{
   author={Vogan, D. A.},
   title={Irreducible characters of semisimple Lie groups. III. Proof of Kazhdan-Lusztig conjecture in the integral case},
   journal={Invent. Math.},
   volume={71},
   date={1983},
   number={2},
   pages={381--417},
}


\bib{V.GL}{article}{
   author={Vogan, D. A.},
   title={The unitary dual of ${\rm GL}(n)$ over an Archimedean field},
   journal={Invent. Math.},
   volume={83},
   date={1986},
   number={3},
   pages={449--505},
   issn={0020-9910},
  % review={\MR{827363}},
  % doi={10.1007/BF01394418},
}

\bib{VoICM}{article}{
author={Vogan, D. A.},
   title={Representations of reductive Lie groups},
   bookTitle={Proc. Intern. Congr. Math. (Berkeley, Calif., 1986)},
   publisher = {Amer. Math. Soc.},
  year = {1987},
pages={246--266},
}

\bib{VoBook}{book}{
author = {Vogan, D. A.},
  title={Unitary representations of reductive Lie groups},
  year={1987},
  series = {Ann. of Math. Stud.},
 volume={118},
  publisher={Princeton University Press}
}


\bib{Vo89}{article}{
  author = {Vogan, D. A.},
  title = {Associated varieties and unipotent representations},
 %booktitle ={Harmonic analysis on reductive groups, Proc. Conf., Brunswick/ME (USA) 1989,},
  bookTitle = {Harmonic analysis on reductive groups, Proc. Conf., Brunswick/ME
    (USA) 1989, Progr. Math., vol. 101},
  publisher = {Birkh\"{a}user, Boston-Basel-Berlin},
  year = {1991},
pages={315--388},
  editor = {W. Barker and P. Sally},
}



\bib{Wa2}{book}{
title={Real reductive groups II},
author={Wallach, N. R.},
year={1992},
publisher={Academic Press Inc. }
}

\bib{Weis}{article}{
 author={Weissman, M. H.},  
 title={L-groups and parameters for covering groups},
 journal={Astérisque},
 volume={398},
 pages={33--186},
 year={2018},
 }


\bib{Zu}{article}{
  title={Tensor products of finite and infinite dimensional representations of semisimple Lie groups},
  author={Zuckerman, G.},
  year = {1977},
  journal={Ann. of Math.},
  volume = {106},
  pages= {295--308}
}



\end{biblist}
\end{bibdiv}
\end{document}